\documentclass[10pt]{amsart}
\usepackage{amsmath,amssymb}
\usepackage{amsthm}
\usepackage{amscd}
\usepackage{graphicx}
\usepackage[pagebackref]{hyperref}

\newcommand{\cF}{{\mathcal F}}
\newcommand{\cE}{{\mathcal E}}
\newcommand{\cK}{{\mathcal K}}

\newcommand{\bbR}{{\mathbb R}}
\newcommand{\bbC}{{\mathbb C}}
\newcommand{\C}[1]{\ensuremath{\mathbb C^{\,#1}}} 
\newcommand{\bbZ}{{\mathbb Z}}
\newcommand{\bbP}{{\mathbb P}}

\newcommand{\euso}{\operatorname{\mathfrak{so}}}

\newcommand{\euu}{\operatorname{\mathfrak u}}
\newcommand{\eugl}{\operatorname{\mathfrak{gl}}}
\newcommand{\eug}{\operatorname{\mathfrak g}}
\newcommand{\euh}{\operatorname{\mathfrak h}}
\newcommand{\eua}{\operatorname{\mathfrak a}}
\newcommand{\euz}{\operatorname{\mathfrak z}}

\newcommand{\Un}{\operatorname{U}}
\newcommand{\GL}{\operatorname{GL}}
\newcommand{\Gr}{\operatorname{Gr}}
\newcommand{\Or}{\operatorname{O}}
\newcommand{\spn}{\operatorname{span}}
\newcommand{\Cof}{\operatorname{Cof}}

\newcommand{\euX}{\operatorname{\frak X}}

\newcommand{\Lie}{\operatorname{\frak L}}
\newcommand{\Hom}{\operatorname{Hom}}
\newcommand{\End}{\operatorname{End}}
\newcommand{\Aut}{\operatorname{Aut}}
\newcommand{\Ric}{\operatorname{Ric}}
\DeclareMathOperator{\tr}{tr}
\newcommand{\I}{\operatorname{I}}

\newcommand{\ts}{\textstyle }
\newcommand{\la}{\langle}
\newcommand{\ra}{\rangle}
\newcommand{\lhk}{\mathbin{\hbox{\vrule height1.4pt width4pt depth-1pt 
                  \vrule height4pt width0.4pt depth-1pt}}}
\newcommand{\w}{{\mathchoice{\,{\scriptstyle\wedge}\,}{{\scriptstyle\wedge}}
      {{\scriptscriptstyle\wedge}}{{\scriptscriptstyle\wedge}}}}
\newcommand{\smcirc}{{\mathchoice{\,{\scriptstyle\circ}\,}
{{\scriptstyle\circ}}{{\scriptscriptstyle\circ}}{{\scriptscriptstyle\circ}}}}

\renewcommand{\Re}{\operatorname{Re}}
\renewcommand{\Im}{\operatorname{Im}}

\numberwithin{equation}{section}

\newtheorem{theorem}{Theorem}

\newtheorem{proposition}{Proposition}
\newtheorem{corollary}{Corollary}
\newtheorem{apptheorem}{Theorem}[section]

\theoremstyle{remark}

\newtheorem{remark}{Remark}
\newtheorem{example}{Example}

\begin{document}

\author[R. Bryant]{Robert L. Bryant}
\address{Duke University Mathematics Department\\
         P.O. Box 90320\\
         Durham, NC 27708-0320}
\email{\href{mailto:bryant@math.duke.edu}{bryant@math.duke.edu}}
\urladdr{\href{http://www.math.duke.edu/~bryant}%
         {http://www.math.duke.edu/\lower3pt\hbox{\symbol{'176}}bryant}}

\title[Bochner-Kahler metrics]{Bochner-K\"ahler metrics}

\date{March 25, 2000}

\begin{abstract}
A K\"ahler metric is said to be \emph{Bochner-K\"ahler} if its Bochner
curvature vanishes.  This is a nontrivial condition when the
complex dimension of the underlying manifold is at least~$2$. 
In this article it will be shown that, in a certain well-defined
sense, the space of Bochner-K\"ahler metrics in complex dimension~$n$
has real dimension $n{+}1$ and a recipe for an explicit formula
for any Bochner-K\"ahler metric will be given.  

It is shown that any Bochner-K\"ahler metric in
complex dimension~$n$ has local (real) cohomogeneity
at most~$n$.  The Bochner-K\"ahler metrics that can be
`analytically continued' to a complete metric, free of singularities,
are identified.  In particular, it is shown that the only 
compact Bochner-K\"ahler manifolds are the discrete quotients of 
the known symmetric examples.  However, there are compact Bochner-K\"ahler
orbifolds that are not locally symmetric.  In fact, every weighted
projective space carries a Bochner-K\"ahler metric.

The fundamental technique is to construct a canonical infinitesimal 
torus action on a Bochner-K\"ahler metric whose associated momentum
mapping has the orbits of its symmetry pseudo-groupoid as fibers.
\end{abstract}



\thanks{ The research for this article was made possible by support 
from the National Science Foundation through grant DMS-9870164
and from Duke University.  The most current version
of this article can be found at {\tt arXiv:math.DG/0003099}.}

\maketitle

\tableofcontents

\listoffigures

\section[Introduction]{Introduction}\label{sec:intro}

In Riemannian geometry, the decomposition of the curvature tensor
into its irreducible summands under the orthogonal group is regarded 
as fundamental.  There are three such summands, the scalar curvature,
the traceless Ricci curvature, and the Weyl curvature.%
\footnote{The Weyl curvature exists as a nontrivial summand only when
the dimension~$n$ of the underlying manifold is 4 or more.  When $n=4$, 
the Weyl curvature is further reducible under the special orthogonal
group, but not the full orthogonal group.} 
The metrics for which one or more of these irreducible tensors vanishes 
have been the subject of much research and a great deal is now known 
about restrictions on the topology of the complete or compact examples.
For example, consult~\cite{Be}, where the bulk of the work is 
devoted to studying the metrics for which the traceless Ricci
curvature vanishes, i.e., the Einstein metrics.  The metrics 
in dimensions 4 or higher for which the Weyl curvature vanishes
are the conformally flat metrics.  While such metrics are trivial
to describe locally, their global geometry is rather delicate, so 
that classifying the complete or compact examples remains a challenge.

In K\"ahler geometry, the corresponding decomposition of the 
curvature tensor into its irreducible summands under the unitary 
group is not quite as familiar, although it has been known
since the 1949 work of Bochner~\cite{Bo}. (For a more recent
treatment, see~\cite[2.63]{Be}.) The K\"ahler decomposition
bears some resemblance to the Riemannian one, there being three
irreducible summands, the scalar curvature, the traceless Ricci curvature, 
and what has become known as the \emph{Bochner} curvature.%
\footnote{N.B.:  The Bochner curvature is one component of the 
Weyl curvature, but not the only component. For example, in complex 
dimension~$2$ the Bochner curvature is the anti-self-dual part
of the Weyl curvature.  See~\S\ref{sssec:  curv decomp}.}
Bochner's interest in this latter tensor was due to its appearance 
in certain Weitzenbock-type formulae.  In~\cite{Bo}, he proved some
cohomological vanishing theorems for compact K\"ahler manifolds
with vanishing Bochner tensor or, more generally, for manifolds
for which the pointwise norm of the Bochner tensor was sufficiently
small relative to the smallest eigenvalue of the Ricci tensor.

While K\"ahler metrics with vanishing scalar curvature or vanishing
traceless Ricci curvature (i.e., K\"ahler-Einstein metrics) have
been much studied, those with vanishing Bochner tensor, now known
as \emph{Bochner-K\"ahler} metrics, have received considerably 
less attention.  For surveys of what has been known up to now
about these metrics, the reader might consult~\cite{DSV}, \cite{KK}, 
\cite{LPV}, \cite{TW}, or~\cite{VV} 
in addition to~\S\ref{sec: BK structure equations}
of the present article.
One will be struck by the paucity of examples.  For example, 
up until now, every known complete Bochner-K\"ahler metric was also 
locally symmetric.  (The symmetric examples are the products
of the form~$M^p_c\times M^{n-p}_{-c}$ where~$M^p_c$ denotes the
$p$-dimensional complex space form of constant holomorphic 
sectional curvature~$c$.)

At first glance, one might expect the theory of Bochner-K\"ahler 
manifolds to parallel the theory of conformally flat manifolds.  
However, this expectation is quickly abandoned.  Unlike the 
local description of conformally flat metrics, a local description
of Bochner-K\"ahler metrics is far from trivial.  In fact,
no such description was known until now.  

Theorem~\ref{thm: existence}  and Corollary~\ref{cor: BK germ moduli}
show that the space of isometry classes of germs of
$C^5$ Bochner-K\"ahler metrics in complex dimension~$n$ can be naturally 
regarded as a closed semi-algebraic subset~$F_n\subset\bbR^{2n+1}$ 
(with a nonempty interior).  More precisely, if~$M$ is 
a complex $n$-manifold endowed with a $C^5$ Bochner-K\"ahler metric~$g$, 
there is a mapping~$f:M\to F_n\subset\bbR^{2n+1}$ (which is a polynomial 
function of the curvature tensor of~$g$ and its first two 
covariant derivatives) with the property that~$f(x)=f(y)$ for~$x,y\in M$ 
if and only if the germ of~$g$ at~$x$ is holomorphically isometric 
to the germ of~$g$ at~$y$.  Moreover, I show that for every~$v\in F_n$, 
there is a Bochner-K\"ahler metric~$g$ on a neighborhood~$U$ of~$0\in\C{n}$ 
so that the associated classifying map~$f:U\to F_n$ satisfies~$f(0)=v$.
(This existence theorem relies on some old results of \'Elie Cartan
that are not readily available in the current literature, so I have
included an appendix that exposes these results in a form convenient
for the applications in this article.)
A by-product of this analysis is that any $C^5$ Bochner-K\"ahler
metric is necessarily real-analytic.%
\footnote{Presumably, any $C^2$ Bochner-K\"ahler
metric is real-analytic, but I have not shown this.}
Accordingly, for the rest of the article, I assume that the 
Bochner-K\"ahler metrics under consideration are real-analytic.

Theorem~\ref{thm: existence}  suggests that a notion of `analytic continuation' 
of Bochner-K\"ahler metrics might be useful.
Elements~$v_1,v_2\in F_n$ are said to be \emph{analytically
connected} if there is a connected Bochner-K\"ahler manifold~$(M^n,g)$
for which~$f(M)$ contains both~$v_1$ and~$v_2$.  This is an equivalence
relation, so denote the analytically connected equivalence 
class of~$v\in F_n$ by~$[v]\subset F_n$.  In Theorem~\ref{thm:  BK constants}, 
I construct a polynomial submersion~$C:\bbR^{2n+1}\to\bbR^{n+1}$ 
and show that it is constant on each~$[v]$.  Eventually, 
Theorem~\ref{thm: analytically connected classes} will show that each
fiber~$C^{-1}(c)\cap F_n$ consists of a finite number of analytically
connected equivalence classes and explicitly identify each one as
a (not necessarily closed) semi-algebraic set of (real) dimension at most~$n$.  
Thus, the components of~$C$ furnish a set of `coarse moduli' for 
Bochner-K\"ahler metrics.  The image~$C(F_n)\subset\bbR^{n+1}$ (which will 
be explicitly identified below) has nonempty interior, so it makes sense to 
say that, roughly speaking, the moduli space of Bochner-K\"ahler metrics in 
complex dimension~$n$ has real dimension~$n{+}1$. 

Since each equivalence class~$[v]\subset F_n$ has real dimension~$m\le n$
at its smooth points, this suggests that a connected Bochner-K\"ahler manifold
of complex dimension~$n$ must always have a non-trivial local isometry
`group', acting with some cohomogeneity~$m\le n$.  
In Theorem~\ref{thm:  sym alg of dim n } and
Proposition~\ref{prop: dim sym alg} , I show that when~$M^n$ is 
simply-connected, the Lie
algebra~$\eug$ of Killing fields for a Bochner-K\"ahler structure 
on~$M$ does indeed have dimension at least~$n$ and I compute its precise
dimension for each analytically connected equivalence class~$[v]\subset F_n$.
Moreover, for each~$v\in F_n$, I compute the dimension of the orbit
of the local isometry pseudogroup through an~$x\in M$ with~$f(x)=v$.  
In particular, I show in \S\ref{sssec: orbit dim n slices} 
how to compute the cohomogeneity~$m$ 
for each~$v\in F_n$.  (Interestingly enough, it turns out that~$m$ 
cannot, in general, be computed from the coarse moduli~$C(v)$ alone.  
This is a reflection of the fact that not all of the equivalence 
classes~$[v]$ are closed sets in~$F_n$.)  The ultimate conclusion is 
that a Bochner-K\"ahler metric always possesses a rather high degree of 
infinitesimal symmetry.   

Perhaps the greatest surprise and what, ultimately, turns out to be
the key to understanding the geometry of Bochner-K\"ahler metrics
is that the Lie algebra~$\eug$ contains a canonical central 
subalgebra~$\euz$ whose dimension~$m$ is the same as that 
of~$[f(x)]\subset F_n$ for some (and hence any)~$x\in M$.  This
infinitesimal torus action can be described explicitly as follows:
Let~$\Omega$ be the K\"ahler form and let~$\rho = \Ric(\Omega)$
be its associated Ricci form~\cite[2.44]{Be}.  Define a `renormalized'
Ricci 2-form~$\eta$ by
\begin{equation*}
\eta = \frac{1}{2(n{+}1)(n{+}2)}\,(\tr_\Omega\rho)\,\,\,\Omega 
        - \frac1{2(n{+}2)}\,\rho\,
\end{equation*}
and define~$p_h(t)$ by the 
formula~$\bigl(t\Omega-\eta)^n=p_h(t)\,\Omega^n$.  Thus,
\begin{equation*}
p_h(t) = t^n - h_1\,t^{n-1} + \cdots + (-1)^n h_n
\end{equation*}
where~$h_j:M\to\bbR$ is a certain symmetric
polynomial of degree~$j$ in the eigenvalues of the Ricci tensor.  
Then Theorem~\ref{thm: generators of center} asserts that the 
$\Omega$-Hamiltonian vector 
fields~$X_j$ defined by~$X_j\lhk\Omega = -dh_j$ for~$1\le j\le n$
are Killing fields for the metric~$g$ and that they Lie commute, i.e.,
span a torus~$\euz\subset\eug$.  Of course, this infinitesimal
action is Poisson since~$h = (h_1,\ldots, h_n):M\to\bbR^n$ is a
momentum mapping by definition.

As is shown in~\S\ref{sssec: mod map ii}, 
the map~$h:M\to \bbR^n$ can be written
as~$\Psi\circ f$  where~$\Psi$ is a weighted 
homogeneous polynomial mapping from~$\bbR^{2n+1}$ to~$\bbR^n$. 
When~$M$ is connected, the maps~$h$ and~$f$ have the
same fibers.  The image of~$h$ is $m$-dimensional and lies in an affine 
subspace~$\eua\subset\bbR^n$ of dimension~$m$ (the same~$m\le n$ as
defined above).  This number~$m$ is defined to be the \emph{cohomogeneity}
of the Bochner-K\"ahler structure.  Theorem~\ref{thm: cohomogeneity m}
 shows that,
in fact, $p_h(t)$ has a polynomial factor~$p_{h''}(t)$ with 
constant coefficients and of degree~$n{-}m$.  Thus, $p_h(t)
= p_{h''}(t)\,p_{h'}(t)$ where
\begin{equation*}
p_{h'}(t) = t^m - h'_1\,t^{m-1} + \cdots + (-1)^m h'_m
\end{equation*}
and the functions~$h'_j:M\to\bbR^m$ for~$1\le j\le m$ are 
smooth.  Theorem~\ref{thm: cohomogeneity m} also shows that, 
outside a (possibly singular) 
complex submanifold~$N\subset M$ (called the \emph{exceptional locus}), 
the \emph{reduced momentum mapping}~$h'=(h'_1,\ldots, h'_m):M\to \bbR^m$ 
is a submersion.  This singular locus is the union of a number of
totally geodesic complex submanifolds of~$M$. Let~$M^\circ = M\setminus N$ 
be its complement, the \emph{regular locus}.

Theorem~\ref{thm: analytically connected classes} 
yields a polynomial embedding~$\iota_v:[v]\to\bbR^m$ of 
each $m$-dimensional analytically connected equivalence class~$[v]$ 
into~$\bbR^m$ as a convex polytope, i.e., an intersection of half-spaces
(which can be open or closed).  The embedding~$\iota_v$ satisfies
$h' = \iota_v\circ f$ when~$f(M)$ lies in~$[v]$.  Moreover,~$h'$ 
maps $M^\circ$ into the interior of the polytope.  
Theorem~\ref{thm: cell metric from pD} 
shows that the interior of $\iota_v\bigl([v]\bigr)$
carries a canonical Riemannian metric so that~
$h':M^\circ \to \iota_v\bigl([v]\bigr)^\circ$ 
is a Riemannian submersion.  In fact, this metric on~
$\iota_v\bigl([v]\bigr)^\circ$ has rational polynomial coefficients 
when expressed in terms of linear coordinates on~$\bbR^m$.
These metrics are related to certain metrics considered
by Guillemin in his study of K\"ahler structures on toric
varieties~\cite{Gu}, as will be explained.

Since the metric on the polytope is very explicitly computed, 
this allows conclusions to be drawn about the existence of
complete Bochner-K\"ahler metrics based on the geometry of the
polytopes.  In Proposition~\ref{prop: metric complete implies cell bounded}, 
I show that if there is a complete
Bochner-K\"ahler metric whose moduli image lies in~$[v]$, then~$[v]$
must be bounded (which turns out to be the same as saying that its
corresponding polytope is bounded).  Essentially, it turns out that
when~$[v]$ is unbounded, any attempt to `analytically continue' the
metric to a maximal domain will run into curvature blow-up at finite
distance.  Since there are very few~$[v]$
that are bounded, this considerably narrows the search 
for complete examples. 

On the other hand, Proposition~\ref{prop: distinct roots is not complete} 
shows that if~$[v]$ is compact but is not a single point, 
then there is no complete Bochner-K\"ahler manifold whose moduli
image lies in~$[v]$.  In this case, the problem is not curvature
blow-up but is, instead, the presence of essential orbifold
singularities in any attempted completion.

A corollary of Proposition~\ref{prop: distinct roots is not complete} 
is that the only compact 
Bochner-K\"ahler manifolds are the compact quotients of the
known symmetric ones.  This result renders vacuous or trivial many 
of the results in the literature about Bochner-K\"ahler metrics.
For example, the only K\"ahler $n$-manifold satisfying the conditions 
of~\cite[Theorems~8.25 and~8.26]{BY} is $\bbC\bbP^n$ 
endowed with a constant multiple of the Fubini-Study metric.  
The conclusions of these theorems (which concern the vanishing
of various cohomology groups) are trivial for these manifolds.

Theorem~\ref{thm: leaf metric} 
provides explicit models for the Bochner-K\"ahler
metrics in dimension~$n$ that are of cohomogeneity~$n$ (i.e., 
the \emph{least} symmetric ones) on the regular locus.  
It constructs, for each $n$-dimensional class~$[v]\subset F_n$,
a Bochner-K\"ahler metric on~$\iota_v\bigl([v]\bigr)^\circ\times\bbR^n$
with the following universal embedding property:  
If~$(M,g)$ is a Bochner-K\"ahler $n$-manifold with~$f(M)\subset[v]$, 
then the universal cover~$\widetilde{M^\circ}$ can be isometrically
immersed into~$\iota_v\bigl([v]\bigr)^\circ\times\bbR^n$, lifting
the momentum submersion~$h:M^\circ\to\iota_v\bigl([v]\bigr)^\circ$. 
Completeness issues can then be addressed by studying the model metric 
on~$\iota_v\bigl([v]\bigr)^\circ\times\bbR^n$.  These metrics
are closely related to the metrics studied in~\cite{Gu} and~\cite{Ab}.
In particular, Abreu's results in~\cite{Ab} can be generalized
to show that the above metrics are actually extremal in the sense of
Calabi.  

Theorem~\ref{thm: complete metrics} provides 
a contractible $n$-parameter family of 
complete Bochner-K\"ahler metrics on~$\C{n}$ and proves that every
simply-connected, complete Bochner-K\"ahler manifold that is 
not homogeneous is isometric to a unique member of this family.

Thus, the set of complete Bochner-K\"ahler manifolds is very restricted.
However, if one is willing to consider orbifolds, it turns out that
there are many nontrivial complete Bochner-K\"ahler metrics on 
orbifolds.  I include some discussion of these at the end of the article.
In fact, by Theorem~\ref{thm: wtd proj space}, 
every weighted projective space carries
a Bochner-K\"ahler metric,%
\footnote{A natural guess would be that this metric is the one that 
comes by symplectic reduction from the standard metric on~$\C{n+1}$
via the weighted $S^1$-action that defines the weighted projective
space.  However, this `reduced' metric is never Bochner-K\"ahler except
in the case of equal weights.}
presumably unique up to constant multiples, though I have not shown this.
For example, the Fubini-Study metric is, up to isometry and
constant multiples, the unique Bochner-K\"ahler metric on~$\bbC\bbP^n$.
For more detail on this, see~\S\ref{sssec: explicit forms} 
and~\S\ref{sssec: wtd proj spaces}.

Finally, in~\S\ref{sec: final remarks}, 
I collect some miscellaneous and incidental remarks
about generalizations and related problems.  In particular, I comment
on how this work in the dimension~2 case is related to the recent
work of Apostolov and Gauduchon~\cite{AG} that classifies the self-dual 
Hermitian Einstein metrics in (real) dimension~4 and use the normal forms 
constructed in this article to produce the first known complete examples 
of such metrics that are of cohomogeneity~2 (the maximum possible,
as it turns out).

\section[Bochner-K\"ahler structure equations]
{The Structure Equations of Bochner-K\"ahler Metrics}
\label{sec: BK structure equations}

First, some standard notation. Let~$\C{n}$ (thought of as columns of 
height~$n$ whose entries are complex numbers) be endowed with its usual 
Hermitian inner product, in which~$\la z, w\ra = {}^t\bar z w$ 
for all~$w,z\in\C{n}$.   Let~$\Un(n)\subset M_n(\bbC)$ denote
the group of unitary matrices and let~$\euu(n)\subset M_n(\bbC)$ denote
its Lie algebra, i.e., the space of skew-Hermitian $n$-by-$n$ matrices.  
As is customary, the conjugate transpose operation will be denoted 
by a superscript asterisk.  Thus, $\la z,w\ra = z^*w$, and~$a\in M_n(\bbC)$ 
lies in~$\euu(n)$ if and only if~$a^* = -a$.

\subsection{The unitary coframe bundle}\label{ssec: Un coframe bundle}

Let~$(M,g,\Omega)$ be a K\"ahler manifold, i.e., $M$ is an 
$n$-dimensional complex manifold and $g$ is an Hermitian metric on~$M$ 
whose associated K\"ahler 2-form~$\Omega$ is closed.  As is customary,
let~$J:TM\to TM$ be the associated almost complex structure endomorphism.

For $x\in M$, 
let~$P_x$ be the set of unitary isomorphisms~$u:T_xM\to\C{n}$.  
Then~$P=\cup_{x\in M} P_x$ is a principal right~$\Un(n)$-bundle over~$M$, 
with the basepoint projection~$\pi:P\to M$ given by~$\pi\bigl(P_x\bigr)=x$
and~$\Un(n)$-action given by~$u\cdot a = a^{-1}\smcirc u$ for~$a\in\Un(n)$.

\subsubsection{The first and second structure equations}
\label{sssec: 1st and 2nd str eqs}

Let~$\omega$ be the $\C{n}$-valued $1$-form on~$P$ defined
by the rule~$\omega(v) = u\bigl(\pi'(v)\bigr)$ for all~$v\in T_uP$.
Then~$\pi^*\Omega = -\frac{i}2\,\omega^*\w\,\omega$.

Because the structure~$(M,g,\Omega)$ is K\"ahlerian, there exists
a unique $\euu(n)$-valued 1-form~$\phi$ on~$P$ satisfying the
\emph{first structure equation} of \'E. Cartan,
\begin{equation}\label{eq: 1st str eq}
d\omega = - \phi\w\omega.
\end{equation}
The \emph{second structure equation} of \'E. Cartan takes the form
\begin{equation}\label{eq: 2nd str eq}
d\phi = -\phi\w\phi + {\ts\frac12}R(\omega\w\omega^*),
\end{equation}
where~$R:P\to\Hom\bigl(\euu(n),\euu(n)\bigr)$ 
is the \emph{K\"ahler curvature function}.  The adjoint representation 
of~$\Un(n)$ on~$\euu(n)$ induces a representation~$\rho$ of~$\Un(n)$
on~$\Hom\bigl(\euu(n),\euu(n)\bigr)$. The curvature function~$R$ is 
equivariant with respect to this action,
i.e., $R(u\cdot a) = \rho(a^{-1})\bigl(R(u)\bigr)$ for~$a\in\Un(n)$.

The first Bianchi identity is $0 = d(d\omega)=-R(\omega\w\omega^*)\w\omega$.  
Thus,~$R$ takes values in the subspace~$\cK\bigl(\euu(n)\bigr)$ 
consisting of those elements~$r\in\Hom\bigl(\euu(n),\euu(n)\bigr)$ 
that satisfy~
\begin{equation*}
r(xy^*{-}yx^*)z+r(yz^*{-}zy^*)x+r(zx^*{-}xz^*)y=0,
\quad\forall\,x,y,z\in\C{n}.
\end{equation*}

\subsubsection{Tensors, vector fields, and symmetries}
\label{sssec: tensors vfs and symms}

The reader will recall that any (real or complex)
representation~$\chi:\Un(n)\to \Aut(V)$ defines a (tensor) vector bundle
$P_\chi = P\times_\chi V$ over~$M$. A section~$\sigma$ of~$P_\chi$ is 
then uniquely defined by a function~$s:P\to V$ that satisfies
the equivariance condition~$s(u{\cdot} a) = \chi(a^{-1})\bigl(s(u)\bigr)$
for all~$u\in P$ and~$a\in\Un(n)$ and~$\sigma(x) = \bigl[u,s(u)\bigr]_\chi$
for some (and hence any)~$u\in P_x$.  The function~$s$ is said
to \emph{represent}~$\sigma$.  For example,~$R$ represents the K\"ahler
curvature tensor.

For notational simplicity, I will use~$\chi$ also to denote the 
induced map on Lie algebras; thus, $\chi:\euu(n)\to\End(V)$.  The
$\Un(n)$-equivariance of a representative function~$s:P\to V$ implies that
the $1$-form~$ds{+}\chi(\phi)\,s$  is~$\pi$-semibasic.   Thus, 
there exists a linear mapping~$Ds:P\to\Hom_\bbR(\C{n},V)$ satisfying
\begin{equation*}
ds + \chi(\phi)\,s = Ds(\omega).
\end{equation*}
Naturally, $Ds$ represents the covariant derivative of the section~ 
$\sigma$ represented by~$s$.

For example, the standard inclusion~$\iota:\Un(n)\hookrightarrow\Aut(\C{n})$
yields~$P_\iota\simeq TM$.  A vector field~$Z$ on~$M$
is represented by the function~$z:P\to\C{n}$ defined by~$z(u) 
= u\bigl(Z_{\pi(u)}\bigr)$.
Now, $\Hom_\bbR(\C{n},\C{n})=\Hom_\bbC(\C{n},\C{n})
\oplus \Hom_\bbC(\C{n},\C{n})\,C$ where~$C:\C{n}\to\C{n}$
is conjugation. Thus, since~$\Hom_\bbC(\C{n},\C{n}) = M_n(\bbC)$,
there are functions~$z'$ and $z''$ on~$P$ with values in~$M_n(\bbC)$ 
so that
\begin{equation*}
dz + \phi\,z = z'\,\omega + z''\,\bar\omega\,.
\end{equation*}
These functions have the $\Un(n)$-equivariance
\begin{equation*}
z'(u\cdot a) = a^{-1}z'(u)a,\qquad
z''(u\cdot a) = a^{-1} z''(u) \bar a
\end{equation*}
and thus represent tensors on~$M$.  In fact,~$z'$ 
represents~$\nabla^{1,0}(Z{-}iJZ)$ while~$z''$ 
represents~$\nabla^{0,1}(Z{-}iJZ) = \bar\partial(Z{-}iJZ)$. 

In particular,~$Z$ is the real part of a holomorphic vector field,
namely~$Z{-}iJZ$, if and only if~$z''=0$.  
Moreover, computation shows that
\begin{equation*}
\pi^*\bigl(Z\lhk\Omega\bigr) = -{\ts\frac i2}(z^*\,\omega - \omega^*\,z),
\end{equation*}
implying, in particular, that
\begin{equation*}
\pi^*\bigl(\Lie_Z\Omega\bigr) 
= -{\ts\frac i2}\,\omega^*\bigl(z'+(z')^*\bigr)\,\omega 
-{\ts\frac i2}\,\omega^*\,z''\,\bar\omega
-{\ts\frac i2}\,\bar\omega^*\,(z'')^*\,\omega.
\end{equation*}

Thus, the flow of~$Z$ is both holomorphic
and symplectic (and hence an infinitesimal symmetry of the K\"ahler
structure) if and only if~$z'' = 0$ and $z'+(z')^* = 0$.  In such a
case,~$Z=\pi'(Z')$ where~$Z'$ is the vector field on~$P$ that satisfies
\begin{equation*}
\omega(Z') = z,\qquad \phi(Z') = z'.
\end{equation*}
The flow of~$Z'$ preserves both~$\omega$ and~$\phi$.  In fact,
\begin{equation*}
\Lie_{Z'}\omega 
= d\bigl(\omega(Z')\bigr) + Z'\lhk\bigl(-\phi\w\omega)
= dz  + \phi\,z - z'\,\omega = 0,
\end{equation*}
so the flow of~$Z'$ does indeed preserve~$\omega$.  Moreover,
since~$\phi$ is the unique $\euu(n)$-valued $1$-form that
satisfies~$d\omega=-\phi\w\omega$, the flow of~$Z'$
must preserve~$\phi$ as well.

Conversely, any vector field on~$P$ whose flow preserves both
$\omega$ and~$\phi$ is of the form~$Z'$ where~$Z$ is a symmetry vector
field of the K\"ahler structure.

If~$Z$ is a symmetry vector field of the K\"ahler structure and~$Z$ 
vanishes at~$x\in M$, then $\nabla Z(x)\in T_xM\otimes T^*_xM$ is
both skew-symmetric and commutes with the complex structure~$J_x$.
Moreover, the flow~$\Phi_Z$ of~$Z$ is complete on the open geodesic 
ball~$B_\delta(x)$ for all sufficiently small~$\delta>0$ and is
isometric there.  Let~$z:P\to\C{n}$ represent~$Z$.  Then~$z(u)=0$ for
all~$u\in P_x$ and, by the above discussion,~$z'(u)$ 
belongs to~$\euu(n)$.  In particular, the linear 
transformation~$a = u^{-1}\circ z'(u)\circ u:T_xM\to T_xM$ is a 
well-defined skew-Hermitian transformation of~$T_xM$.

Then, for all~$v\in T_xM$ with~$|v|<\delta$,
\begin{equation*}
\Phi_Z\bigl(t,\exp_x(v)\bigr) 
     = \exp_x\bigl(e^{-a\,t}v\bigr),
\end{equation*}
i.e., the map~$\exp_x:B_\delta(0_x)\to B_\delta(x)$ intertwines
the linear 1-parameter subgroup action on~$T_xM$ generated by
exponentiating~$-a$ with the flow of~$Z$.  

This has two consequences that will be needed in this article 
(see~\S\ref{sssec: ness conds for complete}). 
First, exponentiating the kernel of~$a$ gives the component of
the fixed locus of the flow of~$Z$ that passes through~$x$, which
is therefore a totally geodesic complex submanifold of~$M$.
Second, when~$M$ is connected, the flow of~$Z$ will be periodic of 
period~$T$ if and only if the eigenvalues of~$z'(u)$ generate the discrete
subgroup of~$i\bbR\subset\bbC$ that consists of the integral multiples 
of~$2\pi i/T$. 

A `micro-local' version of symmetry will be useful.
Two coframes~$u,v\in P$ are said to be \emph{equivalent}
if there is a connected $u$-neighborhood~$U$, 
a connected~$v$-neigh\-borhood~$V$ and a diffeomorphism~$\psi:U\to V$ 
that satisfies~$p(u)=v$ and $p^*(\omega_V)=\omega_U$.  
(It follows, as a consequence, that
$p^*(\phi_V)=\phi_U$.)   Such a $p$, when it exists, is 
unique once~$U$ is specified and is locally of the form
$p(w) = w\circ({\bar p}')^{-1}$ for some local 
isomorphism~$\bar p:\pi(U)\to\pi(V)$ of the K\"ahler structure on~$M$.

If~$u$ and~$v$ are equivalent, then~$R(u)=R(v)$; 
in fact, $D^kR(u)=D^kR(v)$ for all~$k\ge0$.%
\footnote{The converse is not generally true, 
though it is when the K\"ahler structure is real-analytic.}
Let~$\Gamma\subset P\times P$ consist of the equivalent pairs.
Then the set~$\bar\Gamma = \Gamma/\Un(n)$ (where the 
$\Un(n)$-action is the diagonal one on $P\times P$) can be
identified with the set of pointed local isomorphisms of
the K\"ahler structure on~$M$.  For want of a better name,
I will refer to~$\bar\Gamma$ as the \emph{symmetry pseudo-groupoid}
of the K\"ahler structure.   

For any~$x$, the set~$\bar\Gamma\cdot x$
is defined to consist of the points~$\pi(v)$ where~$\pi(u)=x$ and~$(u,v)$
lies in~$\Gamma$.  Thus, $\bar\Gamma\cdot x\subset M$ consists of
the points~$y\in M$ about which the K\"ahler structure is locally
isomorphic to the K\"ahler structure about~$x$.  Even though~$\bar\Gamma$
is not a group, I will, by an extension of the usual language, refer
to~$\bar\Gamma\cdot x$ as the \emph{$x$-orbit} of the symmetry
pseudo-groupoid of the K\"ahler structure.  For any~$x\in M$,
the $x$-orbit is a smooth (but not necessarily closed) 
submanifold of~$M$.

The subset~$\bar\Gamma_x = \bigl(\Gamma\cap (P_x{\times} P_x)\bigr)/\Un(n)$
actually is a group in a natural way, canonically represented 
as a closed subgroup of~$\Un(T_xM)$ as the (local) rotations about~$x$
that preserve the metric and complex structure.  This group will be known
as the \emph{stabilizer} of~$x$.

\subsubsection{Curvature decomposition}
\label{sssec: curv decomp}

Now, the curvature representation $\cK\bigl(\euu(n)\bigr)$ 
is a $\Un(n)$-invariant subspace 
of~$\Hom\bigl(\euu(n),\euu(n)\bigr)$.
It is known~\cite{KN} that~$\cK\bigl(\euu(n)\bigr)$ is isomorphic as
a~$\Un(n)$-module to~$S^{2,2}_\bbR(\C{n}) 
 = \bigl(S^{2,0}(\C{n})\otimes_\bbC S^{0,2}(\C{n})\bigr)_\bbR$, 
the real-valued quartic functions on~$\C{n}$ 
that are complex quadratic and complex conjugate quadratic. 

Now, for each~$p>0$, the $\Un(n)$-invariant Hermitian inner product
on~$\C{n}$ induces a surjective $\Un(n)$-equivariant `trace'
(also called a `contraction' or `Laplacian')
\begin{equation*}
\tr:S^{p,p}_\bbR(\C{n})\to S^{p-1,p-1}_\bbR(\C{n}).
\end{equation*}  
Its kernel~$S^{p,p}_{\bbR,0}(\C{n})\subset
S^{p,p}_\bbR(\C{n})$ is an irreducible~$\Un(n)$-module~\cite{FH}.  

It follows that there is an isomorphism of~$\Un(n)$-modules
\begin{equation}\label{eq: Un decomp of K}
\cK\bigl(\euu(n)\bigr)\simeq S^{2,2}_{\bbR}(\C{n}) 
       \simeq \bbR\oplus S^{1,1}_{\bbR,0}(\C{n})
         \oplus S^{2,2}_{\bbR,0}(\C{n}),
\end{equation}
where the $\Un(n)$-irreducible modules on the right hand side have (real)
dimensions~$1$, $n^2{-}1$, and $\frac14 n^2(n{-}1)(n{+}3)$, respectively.
Thus, there are unique $\Un(n)$-invariant 
subspaces~$\cK_i \subset \cK\bigl(\euu(n)\bigr)$ 
satisfying~$\cK_i\simeq S^{i,i}_\bbR(\C{n})$
for $i=0$, $1$, and $2$. 
 
The K\"ahler curvature function~$R$ can therefore be written as a sum
\begin{equation*}
R = R_0 + R_1 + R_2
\end{equation*}
where~$R_i$ takes values in~$\cK_i$ and represents a 
section of the bundle~$S^{i,i}_\bbR(TM)$, i.e., a tensor 
associated to the K\"ahler structure~$\Omega$. 

The function~$R_0$ represents the scalar curvature,~$R_1$
represents the traceless Ricci tensor, and~$R_2$ represents
the \emph{Bochner tensor}, identified in~1949 by S. Bochner~\cite{Bo}.
When~$n=1$, both~$R_1$ and~$R_2$ are zero by definition, but
when~$n\ge 2$, all three tensors are nonzero for the 
generic K\"ahler metric.  

The K\"ahler structures for which~$R_0$ vanishes are the scalar-flat
K\"ahler structures.  When~$n\ge 2$, those for which~$R_1$ vanishes are
the K\"ahler-Einstein structures and those for which~$R_2$ vanishes are 
known as \emph{Bochner-K\"ahler} structures.

\begin{remark}[The Riemannian analogy]  
Bochner's decomposition of the K\"ahler 
curvature bears a resemblance to the more familiar decomposition of the 
Riemann curvature tensor of a Riemannian metric into the scalar curvature, 
the traceless Ricci tensor, and the Weyl curvature tensor.  However, this 
resemblance is somewhat misleading.  

While the scalar curvature and the Ricci curvature in the two cases do 
correspond, the Weyl curvature tensor of a K\"ahler metric is not simply
the Bochner curvature tensor. For example, when~$n=2$, so that the underlying
manifold has dimension~$4$ and is canonically oriented, the Bochner
tensor turns out to be~$W^-$, the anti-self-dual part of the
Weyl curvature.  Thus, in complex dimension~$2$, the Bochner-K\"ahler
metrics are the same as the self-dual K\"ahler metrics.%
\footnote{These metrics have been studied from this point of view.  
For example, see~\cite{De} and the forthcoming~\cite{AG}.
For further comments on this relationship, 
see~\S\ref{ssec: self-dual K metrics}.}

Bochner observed~\cite{Bo} that the Weyl curvature of a K\"ahler metric
breaks up into two or three irreducible components under the action 
of~$\Un(n)\subset\Or(2n)$,  one of which is the Bochner 
curvature tensor.  One of the other components is equivalent to the 
scalar curvature while, when~$n>2$, 
another is equivalent to the traceless Ricci curvature.  
Thus, when~$n>2$, the vanishing of the Weyl curvature of a 
K\"ahler metric implies that the metric is flat.  In particular, when
$n>2$, a conformally flat K\"ahler metric is flat.  When $n=2$, the
conformal flatness of a K\"ahler metric implies only 
that the structure is Bochner-K\"ahler, with vanishing scalar curvature.%
\footnote{However, as will be seen in Example~\ref{ex: loc sym} 
below, when~$n=2$ there are essentially only two conformally flat 
K\"ahler structures up to local isomorphism and homothety.}
\end{remark}

\subsection{Explicit Bochner-K\"ahler structures}
\label{ssec: explicit BKs}

Few explicit examples of Bochner-K\"ahler structures have been found
up to now.  The main strategy for constructing examples so far has been to
look for examples that satisfy conditions sufficiently stringent
to reduce the construction to an ODE problem.

\begin{example}[Locally symmetric]\label{ex: loc sym} 
The simplest Bochner-K\"ahler metric is the complex $n$-dimensional 
space~$M^n_c$ of constant holomorphic sectional curvature~$c\in\bbR$.  
(In fact, $R_1=R_2=0$ characterizes these metrics.)  

Tachibana and Liu~\cite[\S2]{TL} showed that the 
products~$M^p_c\times M^{n-p}_{-c}$ are Bochner-K\"ahler for 
any~$n$, $p$, and~$c$.  Moreover, they showed that any Bochner-K\"ahler
structure that is a product in a nontrivial way is locally isomorphic 
to~$M^p_c\times M^{n-p}_{-c}$. 

Matsumoto~\cite[Theorem~2]{Ma} proved that a Bochner-K\"ahler structure
with constant scalar curvature is locally symmetric.
Matsumoto and Tanno~\cite{MT} then proved that any locally symmetric 
Bochner-K\"ahler structure is locally isomorphic to one of the above
examples. (For a simple proof, see Proposition~\ref{prop: BK loc symm} below.)  

Note that their results, combined with the preceding remark, imply
the well-known result that the only conformally flat K\"ahler structures 
in dimension~$n=2$ are those that are locally isometric 
to~$M^1_c\times M^{1}_{-c}$ for some~$c\ge0$.
\end{example}

\begin{example}[Rotationally symmetric]\label{ex: rot sym} 
The first
examples with nonconstant scalar curvature appear to be due to 
Tachibana and Liu~\cite{TL}, who considered K\"ahler structures of the form
\begin{equation}\label{eq: rot inv BK}
\Omega = {\ts\frac{i}2}\,\partial\bar\partial f\bigl(|z|^2\bigr)
=-{\ts\frac{i}2}\,\,dz^*\w\left[ f'\bigl(|z|^2\bigr)\,\I_n
                         +f''\bigl(|z|^2\bigr)\,z\,z^*\right]\,dz
\end{equation}
where $f$ is smooth and real-valued on some interval~$I\subset\bbR$.  
The $(1,1)$-form~$\Omega$ is positive 
on~$D = \{\,z\in\C{n}\,\,\mid\,\,|z|^2\in I\,\}$
if and only if~$f'(t)+tf''(t)>0$ and~$f'(t)>0$ (when~$n>1$) 
for all $t\in I\cap [0,\infty)$.  

For~$n\ge2$, they showed that~$\Omega$ is Bochner-K\"ahler on~$D$
if and only if~$f'$ satisfies 
\begin{equation}\label{eq: rot sym ode}
f''(t) = \bigl(a\,tf'(t) + k\bigr){f'(t)^2}
\end{equation}
for some constants~$a$ and~$k$.%
\footnote{ While equation~\eqref{eq: rot sym ode} 
makes sense even when $n=1$,
`Bochner-K\"ahler' has not yet been defined for~$n=1$.  
This will be remedied in~\S\ref{sssec: dim 1} in such a way
that the present discussion extends without change 
to the case $n=1$.} 

For such an~$\Omega$, the eigenvalues of~$\Ric(\Omega)$ 
with respect to~$\Omega$ are
\begin{equation*}
\begin{split}
\rho_1 &= -2(n{+}1)k-2(n{+}2)a\,|z|^2f'\bigl(|z|^2\bigr)\,,\\
\rho_2 &= -2(n{+}1)k-4(n{+}2)a\,|z|^2f'\bigl(|z|^2\bigr)\,,\\
\end{split}
\end{equation*}
with $\rho_1$ having multiplicity~$n{-}1$,
representing the $(n{-}1)$-plane orthogonal to the radial direction, 
and~$\rho_2$ having multiplicity~$1$,
representing the radial direction.  
Thus, the solutions of~\eqref{eq: rot sym ode} for 
which~$a\not=0$ yield Bochner-K\"ahler structures that are not homogeneous.

Tachibana and Liu integrated the above equation when~$k=0$, 
thereby giving explicit examples of Bochner-K\"ahler structures that are not 
homogeneous.  They do not discuss completeness issues, but it is evident 
from their formulae that none of their explicit examples are complete.

\subsubsection{Further analysis}
\label{ssec: further analysis}
Now, \eqref{eq: rot sym ode} can be integrated even when~$k\not=0$.  
Set~$x(t) = tf'(t)$, so that~\eqref{eq: rot sym ode} becomes
\begin{equation}\label{eq: change of vars}
t\,x'(t) = x(t)\bigl(1+k\,x(t)+a\,x(t)^2\bigr).
\end{equation}
Admissible solutions must satisfy~$x>0$ when~$t>0$ 
and~$x'=f'+tf''>0$. 
Now,~\eqref{eq: change of vars} can be integrated by separation of variables
\begin{equation}\label{eq: sep of vars}
\frac{dx}{x\,(1{+}k\,x{+}a\,x^2)} = \frac{dt}t\,.
\end{equation}

\emph{Scaling equivalences.}  
Relation~\eqref{eq: sep of vars} is invariant under scaling $t$, 
which corresponds geometrically to homothety in~$\C{n}$.  Thus, solutions of~
\eqref{eq: sep of vars} that differ by constant scaling 
in~$t$ represent isomorphic K\"ahler 
structures and can be regarded as equivalent.  Similarly, multiplying~$x$ 
by a positive constant corresponds to multiplying the K\"ahler form~$\Omega$ 
by that constant, so solutions for a given pair of constants~$(k,a)$ can be 
regarded as equivalent to the solutions 
for any other pair~$(\lambda k, \lambda^2 a)$ with~$\lambda\in\bbR^+$.

\emph{The two types of solutions.}
For any fixed~$(k,a)\in\bbR^2$, let~$J_{k,a}\subset\bbR$ be the
maximal $x$-interval containing~$0$ on which~$(1{+}k\,x{+}a\,x^2)$
is positive.  Define a positive function~$F$ on~$J_{k,a}$ by the formula
\begin{equation*}
\log F(x) = -\int_0^x \frac{k{+}a\,\xi}{(1{+}k\,\xi{+}a\,\xi^2)}\,d\xi\,.
\end{equation*}
One solution to~\eqref{eq: sep of vars} can then be written implicitly in 
the form
\begin{equation*}
x F(x) = t.
\end{equation*}
The expression on the left hand side of this equation defines a function 
on~$J_{k,a}$ that has positive derivative and that vanishes at~$x=0$. 
Let~$I_{k,a}\subset\bbR$ denote the range of this function.  (More will
be said about this range below.)  The above relation can then be solved 
for~$x$, yielding a real-analytic solution~$x:I_{k,a}\to J_{k,a}$ 
to~\eqref{eq: change of vars}.  

This solution satisfies~$x'(0)=1$.  Any other solution 
to~\eqref{eq: change of vars} whose range 
lies in~$J_{k,a}$ differs from this one by scaling in~$t$.  These solutions
will be said to be of \emph{type one}.

When~$a>0$ and $k\le -2\sqrt a$, there is a second, geometrically distinct, 
admissible solution to~\eqref{eq: change of vars}. 
Under these assumptions, let~$J^*_{k,a}$ be the 
interval~$(p,\infty)$ where~
\begin{equation*}
p = \frac{-k+\sqrt{k^2-4a}}{2a}>0
\end{equation*} 
is the larger root of~$(1{+}k\,p{+}a\,p^2)=0$.
(When the two roots are equal, $p$ is simply the root.)
Define a function~$F^*$ on~$J^*_{k,a}$ by the formula
\begin{equation*}
\log F^*(x) = -\int_x^\infty \frac1{\xi\,(1{+}k\,\xi{+}a\,\xi^2)}\,d\xi.
\end{equation*}
Then~\eqref{eq: sep of vars} can be integrated in the form~$F^*(x) = t$.  
Since the
integral diverges to infinity as $x$ approaches $p$ from above,
the function~$F^*$  maps~$(p,\infty)$ diffeomorphically onto~$(0,1)$. 
Thus, the equation~$F^*(x)=t$ can be solved for~$x$, yielding a real 
analytic solution~$x^*:(0,1)\to(p,\infty)$ to~\eqref{eq: change of vars}.  
Any other solution 
to~\eqref{eq: change of vars} whose range lies in~$J^*_{k,a}$ 
differs from this one by scaling 
in~$t$.  These solutions will be said to be of \emph{type two}.

\emph{Completeness.}  For any admissible solution~$x:I\to J$ 
to~\eqref{eq: change of vars}, 
consider the Bochner-K\"ahler structure~$\Omega$ got by 
setting~$f'(t)=x(t)/t$.  The differential of arc length~$\sigma$ 
along a radial curve~$\{s\,v\mid s^2\in I\}$
for any fixed~$v\in\C{n}$ with~$|v|=1$ can be calculated to be
\begin{equation}\label{eq: radial ds}
d\sigma 
= \frac{d\bigl(x(s^2)\bigr)}
    {2\,\sqrt{x(s^2)\bigl(1+k\,x(s^2)+a\,x(s^2)^2\bigr)}}\,.
\end{equation}
This formula permits an analysis of the completeness properties 
of~$\Omega$ without having to write down an explicit formula for~$x$.

When~$a=0$, the solution of type one is~$x(t)=t/(1{-}kt)$
and~$\Omega$ has constant holomorphic sectional curvature.  
This metric is complete on~$\C{n}$ when~$k=0$.  When~$k>0$,
it is complete on the ball~$|z|^2<k^{-1}$.
When~$k<0$ it is not complete, since the radial arc length
\begin{equation*}
\int_0^{-k^{-1}}\frac{d\xi}{2\,\sqrt{\xi(1+k\,\xi)}}
\end{equation*}
is finite.  However, in this case, the metric on~$\C{n}$ 
extends smoothly to (a multiple of) the Fubini-Study metric 
on~$\bbC\bbP^n$.  

When~$a>0$ and~$1+k\,p+a\,p^2=0$ has no positive root~$p$, the interval~
$J_{k,a}$ contains some interval of the form~$(\alpha,\infty)$
for~$\alpha<0$.  Because the integral
\begin{equation*}
\int_1^\infty\frac{d\xi}{\xi(1+k\,\xi+a\,\xi^2)}
\end{equation*} 
converges, the type one solution to~\eqref{eq: change of vars} 
is defined on an 
interval~$I_{k,a} =(-\delta,R^2)$ for~$\delta$ and~$R^2$ positive, 
with $x(t)$ tending to infinity as~$t$ approaches~$R^2$.  
Thus,~$\Omega$ is defined and nondegenerate on a ball~$|z|<R$.  Since
\begin{equation*}
\int_0^\infty\frac{d\xi}{2\,\sqrt{\xi(1+k\,\xi+a\,\xi^2)}} <\infty
\end{equation*}
the metric is not complete.  Yet, $\Omega$ cannot be extended 
beyond~$|z|<R$ because the two curvatures~$\rho_1$
and~$\rho_2$ tend to $-\infty$ as~$|z|$ approaches~$R$.

Suppose now that~$1+k\,p+a\,p^2=0$ does have at least one
positive root.  By the $x$-scaling argument, it can be assumed 
that~$p=1$ is a root and that there is no root in the interval~$(0,1)$.
Thus,~$k = -(1{+}a)$, so that~$(1+k\,p+a\,p^2) = (1{-}p)(1{-}a\,p)$, 
and~$a\le1$.

Suppose first that~$a=1$ (the extreme value), 
so that $(1{+}k\,p{+}a\,p^2) = (1{-}p)^2$.  Since the integral 
\begin{equation*}
\int_0^1\frac{d\xi}{\xi(1-\xi)^2}
\end{equation*} 
diverges at both endpoints, the type one solution~$h$ 
to~\eqref{eq: change of vars} 
is defined on all of~$\bbR$ and maps~$[0,\infty)$ to~$[0,1)$.  
Thus, $\Omega$ is defined and nondegenerate on all of~$\C{n}$.  
Since
\begin{equation*}
\int_0^1\frac{d\xi}{2\,\sqrt{\xi(1-\xi)^2}} =\infty,
\end{equation*} 
this metric is complete on~$\C{n}$.  As~$|z|^2$ goes to infinity,  
the curvatures~$\rho_1$ and $\rho_2$ approach~$2n$ and~$-4$, 
respectively. 

Still assuming~$(1{+}k\,p{+}a\,p^2)=(1{-}p)^2$, 
consider the type two solution to~\eqref{eq: change of vars}.  
The form~$\Omega$ is 
defined and nondegenerate on the punctured ball~$0<|z|<1$.  
The arc length integral shows that this metric 
is complete on a neighborhood of the puncture but not complete near the 
boundary~$|z|=1$. Since~$x$ goes to infinity as~$|z|$ tends to~$1$, the 
curvatures~$\rho_1$ and~$\rho_2$ tend to $-\infty$ near this 
boundary. Thus, $\Omega$ cannot be extended beyond the punctured 
ball~$0<|z|<1$.

Now suppose~$a<1$. The integral
\begin{equation*}
\int_0^1\frac{d\xi}{\xi(1-\xi)(1-a\xi)}
\end{equation*} 
still diverges at both endpoints, so the type one solution 
to~\eqref{eq: change of vars}
is defined on an open interval in~$\bbR$ that contains~$[0,\infty)$.  
Moreover,~$x$ maps~$[0,\infty)$ to~$[0,1)$.  
Again, $\Omega$ is defined and positive definite on all of~$\C{n}$.  
However, now, the elliptic integral
\begin{equation*}
\int_0^1\frac{d\xi}{2\,\sqrt{\xi(1-\xi)(1-a\xi)}}
\end{equation*}
is finite, so the metric is not complete.  The curvatures~$\rho_1$ 
and~$\rho_2$ approach the limits~$2(n{+}1)-2a$ and~$2(n{+}1)(1{-}a)$, 
respectively, as~$|z|^2$ goes to infinity.   It can be shown that
this  Bochner-K\"ahler structure extends to an `orbifold' 
Bochner-K\"ahler structure on~$\bbC\bbP^n$ even when~$a\not=0$.  
I will not discuss this extension here since its nature will be more 
clear after the considerations to be taken up in the next section.  
Unless~$a=0$ (the Fubini-Study case), this is not a homogeneous metric.

Finally, when~$0<a<1$, consider the type two solution 
to~\eqref{eq: change of vars},
whose range is~$(a^{-1},\infty)$. Since the integral
\begin{equation*}
\int_{a^{-1}{+}1}^\infty\frac{d\xi}{\xi(1-\xi)(1-a\xi)}
\end{equation*} 
converges, the domain of this solution is~$(0,1)$.  Then $\Omega$ is 
defined and nondegenerate on the punctured unit ball~$0<|z|<1$.
When~$0<a<1$, the elliptic integral
\begin{equation*}
\int_{a^{-1}}^\infty\frac{d\xi}{2\,\sqrt{\xi(1-\xi)(1-a\xi)}}
\end{equation*}
is finite, so the metric is not complete at either the puncture or 
the boundary.  The curvatures~$\rho_1$ and~$\rho_2$ approach~$-\infty$ 
as~$|z|$ approaches~$1$.  However, these curvatures remain bounded and 
approach a limit when~$z$ approaches~$0$.  The nature of the 
singularity at~$|z|=0$ and whether or not it can be removed 
will be discussed in~\S\ref{sec: glo geo and sym}.

\emph{Conclusion.} 
Up to constant multiples and scaling, the Ansatz of Tachibana
and Liu provides exactly one example of a complete Bochner-K\"ahler
metric (on~$\C{n}$) that is not locally symmetric.
\end{example}

\begin{example}[Ejiri metrics]\label{ex: ejiri metrics} 
Ejiri~\cite{Ej} considered a somewhat more general Ansatz, seeking
Bochner-K\"ahler metrics for which the Ricci tensor has at most two
distinct eigenvalues, an evident property of the Tachibana-Liu examples
and the locally symmetric examples.  He showed that when~$n\ge3$, such 
examples that are not locally symmetric have cohomogeneity one and that the
isometry stabilizer of the general point is $\Un(n{-}1)\subset\Un(n)$.  
Thus, the problem of describing these examples reduces to an ODE problem, 
which Ejiri integrated up to a Weierstra\ss{}-type equation, 
thereby producing the desired examples.  

In [Ej,\S4], Ejiri remarked that none of his examples (aside from the locally
symmetric ones) were known to be complete.  However, since the Tachibana-Liu
examples are special cases of his examples, at least one of his examples 
is complete.  In fact, Ejiri's example in~\cite[\S4]{Ej} of a complete, $C^2$
Bochner-K\"ahler metric on~$\bbR^{2n}$ turns out to be the complete example
of Tachibana and Liu on~$\C{n}$, but presented in unusual coordinates 
in which it is not fully regular at the origin.  This will become
apparent in~\S\ref{ssec: geod fol}, when all of the complete
Bochner-K\"ahler metrics in dimension~$n$ will be classified.
\end{example}

\subsection{The differential analysis}
\label{ssec: diff analysis}

Now suppose that $M$ is a complex manifold of complex dimension~$n\ge2$
endowed with a Bochner-K\"ahler structure~$\Omega$.  As before,
let~$\pi:P\to M$ be the unitary coframe bundle of~$\Omega$ and
denote its canonical forms by~$\omega$, with values in~$\C{n}$
and~$\phi$, with values in~$\euu(n)$.  Let~$R:P\to\cK(\euu(n))$
be the K\"ahler curvature function.  

\subsubsection{Simplification of the curvature}
\label{sssec: simp of curv}

By definition,~$\Omega$ is Bochner-K\"ahler if and only if~$R_2$
vanishes identically.  The curvature decomposition
of~\S\ref{sssec:  curv decomp}
shows that the remaining part of~$R$ takes values in a representation
isomorphic to~$S^{1,1}_\bbR(\C{n})$, the Hermitian symmetric
quadratic forms.  Now, for any function~
$S = S^*: P\to i\euu(n)\subset M_n(\bbC)$, the $2$-form
\begin{equation*}
\Phi = S\,\omega^*\w\omega - S\,\omega\w\omega^* 
                         - \omega\w\omega^*S + \omega^*\w S\omega\,\I_n
\end{equation*}
takes values in $\euu(n)$ and satisfies~$\Phi\w\omega = 0$ (which is
the first Bianchi identity).  Moreover,
$\Phi$ vanishes if and only if~$S$ vanishes.

It follows that the assumption that~$\Omega$ be Bochner-K\"ahler is 
equivalent to the existence of a function~$S:P \to i\euu(n)\subset M_n(\bbC)$
for which
\begin{equation}\label{eq: simplified curvature}
d\phi + \phi\w\phi
                      = S\,\omega^*\w\omega - S\,\omega\w\omega^* 
                         - \omega\w\omega^*S + \omega^*\w S\omega\,\I_n\,.
\end{equation}

Now,~$S$ does not represent the Ricci tensor \emph{per se}.  
However, the identity
\begin{equation*}
\pi^*\bigl(\Ric(\Omega)\bigr) = i\,\tr\bigl(d\phi{+}\phi\w\phi) 
= i\,\bigl(\tr(S) \omega^*\w\omega + (n{+}2)\,\omega^*\w S\omega\bigr)
\end{equation*}
shows how~$S$ is related to the Ricci form.  In particular, the
scalar curvature of the underlying metric 
is~$2\tr_\Omega\bigl(\Ric(\Omega)\bigr) = -8(n{+}1)\tr S$.

\subsubsection{Higher Bianchi identities}
\label{sssec: higher Bianchi}

Now, consider the consequences of 
differentiating~\eqref{eq: simplified curvature}.
Setting~$\sigma = dS +\phi\,S - S\,\phi$ and taking the exterior derivative
of~\eqref{eq: simplified curvature} leads to the identity
\begin{equation*}
\sigma\w\omega^*\w\omega - \sigma\w\omega\w\omega^* 
    - \omega\w\omega^*\w\sigma - \omega^*\w\sigma\w\omega\,\I_n = 0.
\end{equation*}
This, coupled with the evident identity~$\sigma=\sigma^*$ implies, by
a straightforward variant of Cartan's Lemma, that there must exist a 
function~$T:P\to\C{n}$ so that
\begin{equation}\label{eq: second Bianchi}
dS +\phi\,S - S\,\phi = \sigma = 
T\,\omega^* + \omega\,T^* + {\ts\frac12}(T^*\omega+\omega^*T)\,\I_n\,.
\end{equation}
(Equation~\eqref{eq: second Bianchi} is the second Bianchi identity 
for Bochner-K\"ahler structures.)

Setting~$\tau = dT + \phi\,T - S^2\,\omega$ and computing the exterior
derivative of~\eqref{eq: second Bianchi} yields
\begin{equation*}
\tau\w\omega^* - \omega\w\tau^* 
+ {\ts\frac12}(\tau^*\w\omega-\omega^*\w\tau)\,\I_n = 0.
\end{equation*}
By another variant of Cartan's Lemma, there is a function~$U:P\to\bbR$ so that
\begin{equation}\label{eq: third Bianchi}
dT + \phi\,T - S^2\,\omega = \tau = U\omega.
\end{equation}
(This might be thought of as a sort of \emph{third} Bianchi identity.)

Finally, setting~$\upsilon = dU - (T^*S\omega+\omega^*ST)$ and 
differentiating~\eqref{eq: third Bianchi} yields~$\upsilon\w\omega=0$, 
implying that~$\upsilon = 0$, i.e., that
\begin{equation}\label{eq: fourth Bianchi}
dU  = T^*S\omega+\omega^*ST.
\end{equation}
(This is a \emph{fourth} Bianchi identity.)  The exterior derivative
of~\eqref{eq: fourth Bianchi} is an identity.

The collection of formulae
\begin{equation}\label{eq: structure equations i}
\begin{split}
d\omega &= -\phi\w\omega,\\
d\phi &= -\phi\w\phi + S\,\,\omega^*\w\omega - S\,\omega\w\omega^* 
                         - \omega\w\omega^*S
         + {\ts\omega^*{\w}S\omega}\,\I_n\,,\\
\noalign{\vskip3pt}
dS &= -\phi\,S + S\,\phi + T\,\omega^* + \omega\,T^* 
                + {\ts\frac12}(T^*\omega+\omega^*T)\,\I_n\,,\\
dT &= -\phi\,T +(U\,\I_n +  S^2)\,\omega,\\
dU &= T^*S\omega+\omega^*ST.
\end{split}
\end{equation}
will be referred to as the \emph{structure equations}
of a Bochner-K\"ahler structure.

\subsubsection{First consequences}
\label{sssec: 1st conseqs}

The equations~\eqref{eq: structure equations i} 
allow simple proofs of some known results
about Bochner-K\"ahler structures.

The first part of the following result is due to Matsumoto~\cite{Ma}
and the second part is due to Matsumoto and Tanno~\cite{MT}.

\begin{proposition}\label{prop: BK loc symm}
If a Bochner-K\"ahler structure has constant scalar curvature,
then it is a locally symmetric space.  Any locally symmetric
Bochner-K\"ahler structure is locally isometric to~$M^p_c\times M^{n-p}_{-c}$
for some~$n$, $p$, and $c$.
\end{proposition}

\begin{proof}
Since the pullback of the scalar curvature to~$P$ is~$-8(n{+}1)\tr S$,
the hypo\-thesis of constant scalar curvature is equivalent 
to~$d(\tr S)=0$. Now, by the structure equations
\begin{equation*}
d(\tr S) = {\ts\frac12}(n{+}2)\bigl(T^*\omega+\omega^*T\bigr),
\end{equation*}
so $d(\tr S)=0$ implies that~$T$ vanishes identically.  However,
if~$T$ vanishes identically, then~$dS = -\phi S + S\phi$, so 
that the curvature tensor is parallel.  Thus, the structure is 
locally symmetric.  In particular, the eigenvalues of~$S$ are
all constant. 

Also, $T=0$ implies that ~$S^2 = -U \I_n$.
This, combined with the constancy of the eigenvalues of~$S$ implies
that~$U$ is constant and equal to~$-s^2$ for some real number~$s\ge0$.
This, in turn, implies that~$(S-s\I_n)(S+s\,\I_n)=0$.  Consequently,
$S$ has at most two distinct eigenvalues.  It follows that either
$S=\pm s \I_n$, in which case the structure has constant holomorphic 
sectional curvature~$\mp4s$, or else that there is a symmetric frame 
reduction of~$P$ to a $\bigl(\Un(p){\times}\Un(n{-}p)\bigr)$-subbundle~
$P'\subset P$ on which
\begin{equation*}
S = \begin{pmatrix} -s\,\I_p & 0\\ 0 & s\,\I_{n-p}\end{pmatrix}.
\end{equation*}
Thus, the structure is a locally isomorphic to~$M^p_c\times M^{n-p}_{-c}$
where $c=4s$. 
\end{proof}

The structure equations also yield a simple proof of the
following result of Tachibana and Liu.

\begin{proposition}\label{prop: BK loc prod}
If a Bochner-K\"ahler structure is locally a nontrivial product, then
it is locally isometric to~$M^p_c\times M^{n-p}_{-c}$
for some~$n$, $p$, and $c$.
\end{proposition}

\begin{proof} Assume that the Bochner-K\"ahler structure is locally 
a nontrivial product.  Then for some~$1\le p\le n/2$,
there is a $\bigl(\Un(p){\times}\Un(n{-}p)\bigr)$-subbundle~$P'\subset P$ 
on which~$\phi$ is blocked in the form
\begin{equation*}
\phi = \begin{pmatrix} \phi_1 & 0\\ 0 & \phi_2 \end{pmatrix},
\end{equation*}
where~$\phi_1$ takes values in~$\euu(p)$ and~$\phi_2$ takes values
in~$\euu(n{-}p)$.  This forces~$S$ to be blocked in the corresponding form
\begin{equation*}
S = \begin{pmatrix} S_1 & 0\\ 0 & S_2\end{pmatrix}.
\end{equation*}
The vanishing of the off-diagonal blocks of the structure equation 
for~$dS$ then shows that~$T$ must be zero, thus implying that the structure 
is locally symmetric. Now apply Proposition~\ref{prop: BK loc symm}.
\end{proof}

\subsubsection{The structure function}
\label{sssec: str fctn}

It turns out%
\footnote{This was only noticed in hindsight, after the momentum
mapping construction of~\S\ref{ssec: central syms}.} 
to be more convenient to work with~$H = S - \frac{1}{n+2}(\tr S)\I_n$
than to work with~$S$ directly.  Thus,~$S = H + \frac{1}{2}(\tr H)\I_n$,
and the structure equations~\eqref{eq: structure equations i}
 assume the form
\begin{equation}\label{eq: structure equations ii}
\begin{split}
d\omega &= -\phi\w\omega,\\
d\phi &= -\phi\w\phi + H\,\,\omega^*\w\omega - H\,\omega\w\omega^* 
                         - \omega\w\omega^*H
         + {\ts\omega^*{\w}H\omega}\,\I_n\\
   &\qquad\qquad\qquad\qquad
         +(\tr H)\bigl(\omega^*\w\omega\,\I_n-\omega\w\omega^*\bigr),\\
\noalign{\vskip3pt}
dH &= -\phi\,H + H\,\phi + T\,\omega^* + \omega\,T^* \,,\\
dT &= -\phi\,T +\bigl(H^2+(\tr H)\,H + V\,\I_n\bigr)\,\omega,\\
dV &= (\tr H)\bigl(T^*\omega+\omega^*T\bigr) 
           +\bigl(T^*H\omega+\omega^*HT\bigr).
\end{split}
\end{equation}
where I have also set~$V = U + \frac14\,(\tr H)^2$.
The map~$(H,T,V):P\to i\euu(n)\oplus\C{n}\oplus\bbR$ will be known
as the \emph{structure function}.  

While several of these equations seem more complicated than their 
counterparts in~\eqref{eq: structure equations i}, 
the decisive simplification is the formula for~$dH$ 
versus the formula for~$dS$, as will be seen.  For later use, I record the
identity
\begin{equation}\label{eq: Ricci as H}
\pi^*\bigl(\Ric(\Omega)\bigr) = (n{+}2)\,i\,
\omega^*\w \bigl(H + (\tr H)\,\I_n\bigr)\omega
\end{equation}
which follows from the earlier formula for the Ricci form in terms
of~$S$.

\subsubsection{Scaling weights}
\label{sssec: scal wts}

If~$\Omega$ is a Bochner-K\"ahler structure on a complex manifold~$M$,
then so is~$c\,\Omega$ for any constant~$c>0$.  The unitary coframe
bundle of this scaled structure is
\begin{equation*}
\sqrt{c}\,P = \left\{\,\sqrt{c}\,u \mid  u\in P\,\right\}.
\end{equation*}
The structure functions on the two bundles~$P$ and~$\sqrt{c}\,P$ then satisfy
\begin{equation}\label{eq: scal wts}
H\bigl(\sqrt{c}\,u\bigr) = c^{-1  }\,H(u),\quad
T\bigl(\sqrt{c}\,u\bigr) = c^{-3/2}\,T(u),\quad
V\bigl(\sqrt{c}\,u\bigr) = c^{-2  }\,V(u).\quad
\end{equation}
This motivates assigning `scaling weights' to the components of the 
structure function as follows:  $H$ has scaling weight~1, $T$ has 
scaling weight~$\frac32$, and $V$ has scaling weight~$2$.
(Taking positive, rather than negative, scaling weights is
a simplifying convention.)

\subsubsection{Dimension~$1$}
\label{sssec: dim 1}

Equations~\eqref{eq: structure equations ii} still make sense when~$n=1$, 
i.e., when~$M$ is 
a complex curve endowed with a positive 2-form~$\Omega$ 
and~$\pi:P\to M$ is its $\Un(1)$-coframe bundle.  In this case,
$H$ is an $\bbR$-valued function while~$T$ is~$\bbC$-valued.
The equations~\eqref{eq: structure equations ii} 
then simplify to the scalar equations
\begin{equation}\label{eq: dim 1 str eqs}
\begin{split}
d\omega &= -\phi\w\omega\,,\\
  d\phi &=  - 6H\,\omega\w\bar\omega\,, \\
\noalign{\vskip3pt}
dH &=  \bar T\,\omega + T\,\bar\omega\,,\\
dT &= -\phi\,T +\bigl(2H^2+V\bigr)\,\omega\,,\\
dV &= 2H\,\bigl(\bar T\,\omega + T\,\bar\omega\bigr) = 2H\,dH \\
\end{split}
\end{equation}
Accordingly, when~$n=1$, the satisfaction of these structure
equations can be taken to be the \emph{definition} of the 
Bochner-K\"ahler property.  Throughout this article, this will
be done.  It is not difficult to check that the rotationally
symmetric analysis of Example~\ref{ex: rot sym}, extends to the case~$n=1$ 
when one takes this as the definition of Bochner-K\"ahler.

The Gaussian curvature of the associated metric~$g$ is~$K=-12H$.
In fact, the geometric interpretation of the 
equations~\eqref{eq: dim 1 str eqs} is
just that the $\Omega$-Hamiltonian flow associated to~$K$
should be $g$-isometric. (Compare~\S\ref{sssec:  tensors vfs and symms}.)  
Thus, any constant 
curvature metric in (complex) dimension~$1$ is Bochner-K\"ahler.
Moreover, any nonconstant curvature metric in dimension~1 
that is Bochner-K\"ahler has a canonically defined nontrivial 
Killing field.

Assume that~$M$ is connected, which implies that~$P$ is also connected.
The last structure equation of~\eqref{eq: dim 1 str eqs} implies that~$V-H^2$
is a constant~$C_2$ (the index denotes the scaling weight), and 
the next-to-last equation of~\eqref{eq: dim 1 str eqs} 
then implies that there is
a constant~$C_3$ so that $|T|^2 = H^3 + C_2\,H + C_3$, or equivalently,
that~$|T|^2-VH = C_3$.

These two `constants of the structure' will be generalized 
considerably in higher dimensions, as will the existence of nontrivial
symmetry vector fields.

\section{Existence and Moduli}\label{sec: exist and mod}

\subsection{Existence}\label{ssec: existence}

In~\cite{Ca}, \'Elie Cartan proved a powerful existence and
uniqueness theorem that generalizes Lie's Third Fundamental 
Theorem from the case of a transitive group action to the
case of an intransitive group action.  

For the convenience of the reader and because Cartan's
rather sketchy treatment needs amplification on some 
minor points, a discussion of his theorem is included 
in the Appendix.

Cartan's conditions for the existence of a (local) coframing 
and system of functional invariants satisfying a given set
of structure equations are satisfied by the 
system~\eqref{eq: structure equations ii}.
The following result is then an immediate consequence of his 
general theorem.

\begin{theorem}\label{thm: existence} 
For any~$(H_0,T_0,V_0)\in i\euu(n)\oplus\C{n}\oplus\bbR$, there
exists a Bochner-K\"ahler structure~$\Omega$ on a
neighborhood~$V$ of~$0\in\C{n}$ whose unitary coframe bundle~$\pi:P\to V$ 
contains a~$u_0\in P_0=\pi^{-1}(0)$ for which~$H(u_0) = H_0$, 
$T(u_0)=T_0$, and $V(u_0)=V_0$. Any two real-analytic Bochner-K\"ahler 
structures with this property are isomorphic on a neighborhood 
of~$0\in\C{n}$. Finally, any Bochner-K\"ahler structure that
is~$C^5$ is real-analytic.
\end{theorem}

\begin{proof}
Since the exterior derivatives of the 
equations~\eqref{eq: structure equations ii} are
identities, Cartan's conditions (i.e., his generalization
of the Jacobi conditions) are satisfied for these equations as 
structure equations of a coframing.   

Thus, by Theorem~\ref{thm: Cartan existence} (see the Appendix),  
for any~$(H_0,T_0,V_0)\in i\euu(n){\oplus}\C{n}{\oplus}\bbR$,
there exists a real-analytic manifold~$N$ of dimension~$n^2{+}2n$ 
on which there are two real-analytic 1-forms~$\omega$ and~$\phi$, 
taking values in~$\C{n}$ and~$\euu(n)$, respectively, and a real-analytic 
function~$(H,T,V):P\to i\euu(n)\oplus\C{n}\oplus\bbR$
with the properties that~$(\omega,\phi)$ is a $\C{n}{\oplus}\euu(n)$-valued 
coframing on~$N$, that the 
equations~{eq: structure equations ii} are satisfied on~$N$, and
that there exists a~$u_0\in N$ for which~$\bigl(H(u_0),T(u_0),V(u_0)\bigr)
= (H_0,T_0,V_0)$.  

Since~$d\omega = -\phi\w\omega$, the equation~$\omega=0$ defines
an integrable plane field of codimension~$2n$ on~$N$.  After shrinking~$N$ 
to an open neighborhood of~$u_0$ if necessary, an application of the
complex Frobenius theorem shows that there is a submersion~$z:N\to\C{n}$
with~$z(u_0) = 0$ so that the leaves of this integrable plane field 
are the fibers of~$\pi$ and, moreover, that~$dz = p\,\omega$
for some function~$p:N\to\GL(n,\bbC)$ that satisfies~$p(u_0) = \I_n$.

Since~$\phi=-\phi^*$, the 2-form
\begin{equation*}
\Omega = -{\ts\frac i2}\,\omega^*\w\omega 
   =-{\ts\frac i2}\,dz^*\w (pp^*)^{-1}\,dz  
\end{equation*}
is closed.  Since $\Omega$ is $z$-semibasic and since, by
definition, the fibers of~$z$ are connected, it follows that~$\Omega$
is actually the pullback to~$N$ of a closed, positive (1,1)-form 
on the open set~$V = z(N)\subset\C{n}$, i.e., a K\"ahler structure on~$V$.

Let~$\pi:P\to V$ be the unitary coframe bundle of this K\"ahler structure.
Define a mapping~$\tau:N\to P$ as follows:  If~$z(u)=x\in V$, then
$dz_u:T_uN\to T_x\C{n}\simeq\C{n}$ is surjective and, by construction, 
has the same kernel as~$\omega_u:T_uN\to\C{n}$.  Thus, there
is a unique linear isomorphism~$\tau(u):T_x\C{n}\to\C{n}$ 
so that~$\omega_u=\tau(u)\circ dz_u$.  In fact,~$\tau(u)$ is complex
linear; using the standard identification~$T_x\C{n}\simeq\C{n}$,
one sees that~$\tau(u)$ becomes~$p(u)^{-1}\in M_n(\bbC)$.

The equation~$\Omega = -{\ts\frac i2}\,\omega^*\w\omega$ implies that
$\tau(u)$ is a unitary coframe for all~$u\in N$.  Since~$(\omega,\phi)$
is a coframing, it follows that~$\tau:N\to P$ is an open immersion of~$N$
into~$P$.  Shrinking~$N$ again if necessary, it can be assumed that~$\tau$
embeds~$N$ as an open subset of~$P$.  Thus, nothing is lost by
identifying~$N$ with this open subset of~$P$.  

The structure equations~\eqref{eq: structure equations ii} 
now become identified with the structure
equations of the unitary coframe bundle~$P$, implying that the underlying
K\"ahler structure on~$V$ is, in fact Bochner-K\"ahler, and that the
structure function~$(H,T,V)$ takes on the value~$(H_0,T_0,V_0)$ at~$u_0\in P$,
as desired.  Further details are left to the reader.  This completes
the existence proof.

Uniqueness in the real-analytic category now follows directly from 
Theorem~\ref{thm: Cartan existence}. 
Now, while Cartan states the uniqueness part of 
Theorem~\ref{thm: Cartan existence} only in
the real-analytic category, uniqueness can actually be proved 
using only ordinary differential equations (i.e., the Frobenius theorem);
the Cauchy-Kowalewski or Cartan-K\"ahler Theorems are not needed.  Thus, 
his uniqueness result is valid as long as the form~$\Omega$ is 
sufficiently differentiable for~$P$ to exist as a differentiable bundle 
and for~$H$, $T$, and $V$ to be defined and differentiable.  For this to
be true, it certainly suffices for~$\Omega$ to be~$C^5$. 

Since Cartan's existence proof produces a real-analytic example, 
uniqueness then implies that any $C^5$
Bochner-K\"ahler structure is real-analytic. 
\end{proof}

\begin{remark}[Minimal Regularity]
With some work, one can show that if~$H$ 
and~$T$ are differentiable, then~$V$ (which, 
by~\eqref{eq: third Bianchi}, must exist) 
must be differentiable as well, thus reducing the regularity needed to
apply Cartan's Theorem to~$C^4$.  However, this is almost certainly
not optimal since, presumably, when~$n\ge2$, any $C^2$ K\"ahler structure 
that is Bochner-K\"ahler is real-analytic.  However, the above proof 
does not show this.  

\emph{From now on, I will assume that the Bochner-K\"ahler 
structures under consideration are real-analytic.}
\end{remark}

\subsection{Local moduli}\label{ssec: local mods}

The group~$\Un(n)$ acts on the space~$i\euu(n)\oplus\C{n}\oplus\bbR$
in the usual way:
\begin{equation}\label{eq: Un action}
a\cdot(h,t,v) = (aha^*,\,at,\,v)
\end{equation}
for~$a\in\Un(n)$.  This action makes the structure 
function of a Bochner-K\"ahler structure~
$(H,T,V):P\to i\euu(n)\oplus\C{n}\oplus\bbR$ equivariant
with respect to the right bundle action, i.e.,
\begin{equation}\label{eq: Un equivariance}
\bigl(H(u{\cdot}a),T(u{\cdot}a),V(u{\cdot}a)\bigr)
= a^{-1}\cdot \bigl(H(u),T(u),V(u)\bigr).
\end{equation}
Consequently, it will be useful to have an understanding of the 
orbits of~$\Un(n)$ acting on this space.  

\subsubsection{Orbits}
\label{sssec: orbits}

Let~$W\subset i\euu(n)\oplus\C{n}\oplus\bbR$ be the linear
subspace consisting of the triples~$(h,t,v)$ for which~$h$ is diagonal
and~$t$ is real.  Then~$W$ is a linear subspace of (real) 
dimension~$2n{+}1$.  Let~$C\subset W$ be the `chamber' defined
by the inequalities~$h_{1\bar1}\ge h_{2\bar2}\ge \cdots\ge h_{n\bar n}$
augmented by the conditions that~$t_j\ge0$, with equality 
if~$h_{j\bar\jmath}=h_{i\bar\imath}$ for any~$i<j$.  
N.B.: The set~$C$ has nonempty interior in~$W$. Note,
however, that~$C$ is not closed when~$n\ge2$.

\begin{proposition}\label{prop: Un orbits}
Each $\Un(n)$-orbit in $i\euu(n)\oplus\C{n}\oplus\bbR$ meets~$C$
in exactly one point.
\end{proposition}

\begin{proof}
Consider any~$(h,t,v)\in i\euu(n)\oplus\C{n}\oplus\bbR$. 
Act by an element~$a\in\Un(n)$ so as to reduce to the case
where~$h$ is diagonal and its (real) eigenvalues are arranged in 
decreasing order down the diagonal.  If there are integers~$i\le j$ 
so that~$h_{j\bar\jmath}=h_{i\bar\imath}$, suppose that $i,i{+}1,\ldots, j$ 
is a maximal unbroken string with this property.
Then the stabilizer of~$h$ in~$\Un(n)$ will contain a subgroup
isomorphic to $\Un(j{-}i{+}1)$ that will act as unitary rotations on 
the subvector~$(t_i,\ldots, t_j)$.  Acting by an element of
the stabilizer of~$h$, one can then reduce to the case where~$t_i$ is real
and nonnegative while~$t_{i+1}=\cdots=t_j=0$.  By definition, the 
resulting new~$(h,t,v)$ is an element of~$C$.  
It is clear from the construction that this
element is unique. 
\end{proof}

\begin{corollary}\label{cor: BK germ moduli}
The set of isomorphism classes of germs of Bochner-K\"ahler 
structures in dimension~$n$ is in one-to-one correspondence
with the elements of~$C$. 
\end{corollary}

\subsubsection{Invariant polynomials}
\label{sssec: inv polys}

It is not difficult to exhibit enough $\Un(n)$-invariant 
polynomials on~$i\euu(n)\oplus\C{n}\oplus\bbR$ to separate the 
$\Un(n)$-orbits.  For~$k\ge0$, define the $\Un(n)$-invariant polynomials
\begin{equation}\label{eq: define a and b}
a_k(h,t,v) = \tr(h^k), \quad\qquad b_{k+3}(h,t,v) = t^* h^k\,t\,,
\end{equation}
and set~$b_2(h,t,v) = v$. (The indexing is chosen so as to indicate the
scaling weight as defined in~\S\ref{sssec: scal wts}.  The anomalous definition 
of~$b_2$ will be explained below.)
Then an easy argument using Proposition~\ref{prop: Un orbits} shows
that the collection of~$2n{+}1$ functions
\begin{equation}\label{eq: varphi components}
\varphi  = (a_1,\ldots,a_n,b_2,b_3,\ldots,b_{n+2})
\end{equation}
separates the $\Un(n)$-orbits in~$i\euu(n){\oplus}\C{n}{\oplus}\bbR$.%
\footnote{In fact, by~\cite[Theorem~12.1]{Pr}, the components of~$\varphi$ 
generate the ring of $\Un(n)$-invariant polynomials 
on~$i\euu(n){\oplus}\C{n}{\oplus}\bbR$.}  

When~$n=1$, the function~$a_1^2 + {b_2}^2 + b_3\ge0$ is evidently a
proper function on~$i\euu(1){\oplus}\bbC{\oplus}\bbR$ while, 
for~$n\ge2$, the function~$a_2 + {b_2}^2 + b_3\ge0$ is a proper function 
on~$i\euu(n){\oplus}\C{n}{\oplus}\bbR$.  

It follows that~$\varphi$ is a proper mapping, implying that~$F_n 
= \varphi\bigl(i\euu(n){\oplus}\C{n}{\oplus}\bbR\bigr)$ is closed 
in~$\bbR^{2n+1}$.  The set~$F_n\subset\bbR^{2n+1}$ is thus 
the proper moduli space of orbits.  

If~$(h_0,t_0,v_0)\in i\euu(n){\oplus}\C{n}{\oplus}\bbR$ is such that~$h_0$ 
has~$n$ distinct eigenvalues and~$t_0$ is not orthogonal to any of the 
eigenvectors of~$h_0$, then an elementary computation shows that~
$\varphi'(h_0,t_0,v_0):i\euu(n){\oplus}\C{n}{\oplus}\bbR\to\bbR^{2n+1}$
is surjective.  It follows from this that~$F_n$ is the closure 
of its interior. Of course,~$\varphi:C\to F_n$ is a bijection. 

\subsubsection{The moduli mapping}
\label{sssec: moduli map}

This description of~$F_n$ can be interpreted as saying that the germs 
of Bochner-K\"ahler structures in dimension~$n$ form a singular space of 
real dimension~$2n{+}1$.  It is~$F_n$ that is the natural
moduli space for germs of Bochner-K\"ahler structures in the following 
sense: For any Bochner-K\"ahler structure~$(M,g,\Omega)$, there is
a commutative diagram
\begin{equation}\label{eq: comm di of maps}
\begin{CD}
P @>(H,T,V)>> i\euu(n){\oplus}\C{n}{\oplus}\bbR\\
@V\pi VV @VV{\varphi}V\\
M @>f>> {F_n\rlap{${}\subset\bbR^{2n+1}$}}\\
\end{CD}
\end{equation}
where~$f:M\to F_n$ is a real-analytic map each of whose fibers 
is an orbit of the symmetry pseudo-groupoid of the Bochner-K\"ahler
structure on~$M$.  This function will be known as the \emph{moduli
mapping} of the Bochner-K\"ahler structure.

\subsubsection{Analytic connectedness}
\label{sssec: anal con}

However, this description does not really say `how many' Bochner-K\"ahler
structures there are locally since, for a given Bochner-K\"ahler structure,
the map~$f:M\to F_n$ might have rather large image in~$F_n$.  A priori, the
image could even have dimension as large as~$2n$, in which case one would
be tempted to say that the `generic' Bochner-K\"ahler structures
depend on only one parameter, the parameter that distinguishes the
`hypersurfaces' in~$F_n$ that are the images of generic Bochner-K\"ahler
structure maps.  However, as will be shown in the next subsection, 
this is not the case.  Instead, the dimension of the image~$f(M)$ 
turns out to be no more than~$n$ for any Bochner-K\"ahler structure.

One of the difficulties that arises in discussing this `how many' question
is that it turns out that not every connected Bochner-K\"ahler structure
can be regarded as an open subset of a unique `maximal' connected
Bochner-K\"ahler structure (cf. the discussion of the 
dimension~$n=1$ at the end of~\S\ref{sssec: explicit dim 1}).  
Even when one restricts
attention to the simply-connected, connected Bochner-K\"ahler structures,
this difficulty persists.  Compare this situation with that of locally
symmetric spaces:  Every simply-connected, connected locally symmetric
space has an isometric open immersion (sometimes called a developing
map) into a unique (complete) simply-connected symmetric space 
and this immersion is unique up to ambient isometry.  The discussion
carried out in Example~\ref{ex: rot sym} and 
in~\S\ref{sssec: explicit dim 1} below shows that no 
such result could hold for Bochner-K\"ahler structures.

Two elements~$v_1,v_2\in F_n$ will be said to be 
\emph{analytically connected}
if there exists a connected Bochner-K\"ahler manifold~$(M,\Omega)$
so that both $v_1$ and~$v_2$ lie in~$f(M)$.  An elementary argument
shows that this is an equivalence relation.  
One of the tasks of this article is to describe these equivalence 
classes explicitly.

\subsubsection{Dimension~$1$}
\label{sssec: explicit dim 1}

Now,~$i\euu(1)=\bbR$ and the map~$\varphi:\bbR\oplus\bbC\oplus\bbR\to\bbR^3$
takes the form
\begin{equation*}
\varphi(h,t,v) = \bigl(h,v,|t|^2\bigr).
\end{equation*}
Thus~$F_1\subset\bbR^3$ is the closed upper half-space.  
The fiber~$\varphi^{-1}(x,y,0)$ is a single point for each~$(x,y,0)$ on
the boundary of~$F_1$ while the fiber~$\varphi^{-1}(x,y,z)$ is a circle
when~$(x,y,z)$ lies in the interior~$F^\circ_1$, i.e., when~$z>0$.

By Theorem~\ref{thm: existence}, every point of~$F_1$ lies in the $f$-image 
of some Bochner-K\"ahler structure in dimension~1.

Let~$(M,g,\Omega)$ be a connected Bochner-K\"ahler manifold of dimension~$1$,
so that $\varphi\circ(H,T,V) = (H,V,|T|^2)$.  As was pointed out 
in~\S\ref{sssec: dim 1}, there are constants~$C_2$ and~$C_3$ so that 
\begin{equation*}
V-H^2 = C_2\,\qquad\text{and}\qquad |T|^2 - HV = C_3\,.
\end{equation*}
In other words, $\varphi\circ(H,T,V) = (H,\,H^2{+}C_2,\,H^3{+}C_2H{+}C_3)$,
implying that the map~$f:M\to \bbR^3$ has its image either a point (if
$H$ is constant) or a curve.  

For any~$C = (C_2,C_3)$, let~$p_C(t) = t^3{+}C_2\,t{+}C_3$ and set
\begin{equation*}
\Gamma_C = \bigl\{(t,\,t^2{+}C_2,\,t^3{+}C_2\,t{+}C_3)\ 
            \vrule\ p_C(t)\ge0\,\bigr\}.
\end{equation*}
Since~$dH=\bar T\,\omega+T\,\bar\omega$,
it follows that~$df$ vanishes only at those~$x\in M$ where~$|T|^2 = 0$.  
In other words, if
$M^\circ = f^{-1}(F_1^\circ)$ is the locus where~$|T|^2$ is nonzero, then
$f:M^\circ\to F_1^\circ$ is a submersion onto an open subset of 
$\Gamma^\circ_C = \Gamma_C\cap F_1^\circ$. 

Since~$4|T|^2$ is the squared norm of the gradient of~$dH$ and
the $\Omega$-Hamiltonian of~$H$ is a Killing field on the surface, it
follows that either~$|T|^2$ vanishes identically or else it vanishes
only at isolated points in~$M$ and then only to second order.  

In the former case, $H$ is constant on~$M$. 
By the structure equations~\eqref{eq: dim 1 str eqs}, 
since~$T$ vanishes identically
it follows that~$V\equiv -2H^2$.  Thus, 
each of the points~$v = (r,-2r^2,0)\in F_1$ constitutes a single 
analytically connected equivalence class that is the $f$-image 
of any surface endowed with a metric of constant curvature~$K=-12r$.
Note that, in this case, the constants~$C_2$ and~$C_3$ assume the values
$C_2 = -3r^2$ and $C_3 = 2r^3$, so that $p_C(t) = (t-r)^2(t+2r)$
has either a double or triple root (if~$r=0$).
Let~$\Pi = \bigl\{(r,-2r^2,0)\ \vrule\ r\in\bbR\,\bigr\}$ be the 
parabola of `isolated' classes.  These are the only points that can be
the value of a constant~$f$.

Now suppose that~$|T|^2$ is not identically zero, so that~$M^\circ$
is simply~$M$ minus a set of isolated points.

When~$p_C(t)$ has only one real simple root, say~$r_0$, then~$\Gamma_C$
is connected and homeomorphic to a closed half-line.  Call this Case 1.
In this case~$\Gamma_C\cap \Pi=\emptyset$, so that it is not possible
for $M$ to satisfy~$f(M)\subset\Gamma_C$ and have~$f(M)$ be a point.
Since~$f:M^\circ\to\Gamma^\circ_C$ is a submersion, it follows that
if~$f(M)$ lies in~$\Gamma_C$, then~$f(M)$
is a open subset of~$\Gamma_C$.  Since~$\Gamma_C$ is connected, it
follows that~$\Gamma_C$ must constitute a single analytically
connected equivalence class.

When~$p_C(t)$ has only one real root, but this root is multiple, the
only possibility is that this root is~$t=0$ and, in fact, $p_C(t) = t^3$.
Call this Case~2. In this case, $\Gamma_C = \{(0,0,0)\}\cup\Gamma^\circ_C$
where~$\Gamma^\circ_C = \bigl\{(t,\,t^2,\,t^3)\,\vrule\,t>0\,\bigr\}$.
In this case, $f(M)$ can lie in~$\Gamma_C$ only
if either~$f(M)=\{(0,0,0)\}$ or~$f(M)$ is an open subset of~$\Gamma^\circ_C$.
Since~$\Gamma^\circ_C$ is connected, it follows that~$\Gamma^\circ_C$
constitutes a single analytically connected equivalence class.

When~$p_C(t)$ has two real distinct roots, say~$r_1>r_2$, one must be
double, so there are two possibilities.  
Case 3-$i$ will be that in which~$r_i$ is the double root.

In Case 3-1, $p_C(t) = (t-r)^2(t+2r)$ where~$r>0$.  
Since~$\Gamma_C \cap \Pi = \{(r,-2r^2,0)\}$, define
\begin{equation*}
\begin{split}
\Gamma^a_C &= \bigl\{\bigl(t,\,t^2{-}3r^2,\,(t{-}r)^2(t{+}2r)\bigr)\ 
                   \vrule\ t>r\,\bigr\},\\
\Gamma^b_C &= \bigl\{\bigl(t,\,t^2{-}3r^2,\,(t{-}r)^2(t{+}2r)\bigr)\ 
                   \vrule\ -2r\le t < r\,\bigr\}.\\
\end{split}
\end{equation*}
Then~$\Gamma_C = \Gamma^b_C\cup\{(r,-2r^2,0)\}\cup\Gamma^a_C$, and 
each of~$\Gamma^a_C$, $\{(r,-2r^2,0)\}$, and $\Gamma^b_C$ is evidently
a single analytically connected equivalence class.

In Case 3-2, $p_C(t) = (t-r)^2(t+2r)$ where~$r<0$.  
Still,~$\Gamma_C \cap \Pi = \{(r,-2r^2,0)\}$, but now
$\Gamma_C = \{(r,-2r^2,0)\}\cup \Gamma^a_C$ where
\begin{equation*}
\Gamma^a_C = \bigl\{\bigl(t,\,t^2{-}3r^2,\,(t{-}r)^2(t{+}2r)\bigr)\ 
                   \vrule\ -2r\le t\,\bigr\},
\end{equation*}
and each of~$\{(r,-2r^2,0)\}$ and $\Gamma^a_C$ is evidently
a single analytically connected equivalence class.

When~$p_C(t)$ has three distinct real roots, say~$r_0>r_1>r_2$,
then~$r_0+r_1+r_2=0$ and again~$\Gamma_C\cap \Pi=\emptyset$. 
Call this Case 4.  In this case,~$\Gamma_C$ has two components
\begin{equation*}
\begin{split}
\Gamma^0_C &= 
  \bigl\{\bigl(t,\,t^2+(r_0r_1{+}r_0r_2{+}r_1r_2),
          \,(t{-}r_0)(t{-}r_1)(t{-}r_2)\bigr)\ 
                   \vrule\ r_0 \le t\,\bigr\},\\
\Gamma^1_C &=   
  \bigl\{\bigl(t,\,t^2+(r_0r_1{+}r_0r_2{+}r_1r_2),
          \,(t{-}r_0)(t{-}r_1)(t{-}r_2)\bigr)\ 
                   \vrule\ r_2 \le t\le r_1\,\bigr\},\\
\end{split}
\end{equation*}
each of which is a single analytically connected equivalence class.

Now, in all these cases, the metric~$g = \omega\circ\bar\omega$
can be expressed directly in terms of the invariants.  Restrict
attention to~$M^\circ\subset M$ and note that, by the structure equations,
the complex-valued $1$-form~$\omega/T$ is closed and therefore
a nowhere vanishing holomorphic $1$-form on~$M^\circ$.  
Since~$|T|^2$ vanishes only to second order at each of its zeroes, 
$\omega/T$ extends to all of~$M$ as a meromorphic $1$-form 
with simple poles at the places where~$\omega/T$ vanishes.

Also, since~$|T|^2 = H^3+C_2\,H + C_3$, it follows that
\begin{equation*}
\frac{dH}{H^3+C_2\,H + C_3} = \frac{dH}{|T|^2} 
= \frac{\omega}{T}+\frac{\bar\omega}{\bar T}.
\end{equation*}
Thus,
\begin{equation*}
\frac{\omega}{T} = \frac{dH}{2(H^3+C_2\,H + C_3)} + 2i\,d\theta
\end{equation*}
where~$\theta$ is locally well-defined on~$M^\circ$ 
up to a (real) additive constant (the factor of~$2$ in front of
the~$d\theta$ term provides for consistency with later notation).  
Consequently, one has the formula
\begin{equation*}
g = \omega\circ\bar\omega 
= \frac{dH^2}{4(H^3+C_2\,H + C_3)} + 4(H^3+C_2\,H + C_3)\,d\theta^2\,.
\end{equation*}
More precisely, the simply-connected cover~$\widetilde{M^\circ}$ admits
a developing map~$(H,\theta):\widetilde{M^\circ}\to\bbR^2$ that
isometrically embeds~$\widetilde{M^\circ}$ into the region~$R_C$ 
in the~$H\theta$-plane defined by the inequality~$H^3+C_2\,H + C_3>0$,
endowed with the above metric.

Using this representation, one can determine which of the Cases above
can allow complete Bochner-K\"ahler metrics in dimension~$1$.  

For example, let~$R$ be the largest real root of~$p_C(t)$. 
Then because the integral
\begin{equation*}
\int_{R+1}^\infty \frac{dH}{\,2\sqrt{H^3+C_2\,H + C_3}\,}
\end{equation*}
converges, the metric~$g$ defined above is not complete at the `edge' 
$H=\infty$ of the half-plane~$H>R$.  This implies that if~$(M,g,\Omega)$ 
is a Bochner-K\"ahler metric with characteristic polynomial $p_C$, 
satisfying~$H\ge R$, and having $H$ non-constant, then the length of
the gradient lines of~$H$ would be finite in the increasing direction
and so could not be complete.

Consequently, a complete Bochner-K\"ahler metric must have its image lie in a 
bounded region of~$F_1$.  In particular, the analytically connected component
that contains~$f(M)$ must be bounded.  The only bounded analytically 
connected equivalence classes are 
\begin{enumerate}
\item Case~3 with~$f(M)=\{(r,-2r^2,0)\}$;
\item Case~3-1 with~$f(M) = \Gamma^b_C$; and 
\item Case~$4$ with~$f(M)=\Gamma^1_C$.
\end{enumerate}

The case of a single point has already been discussed:  There is a
unique connected and simply-connected complete example for each~$r$.

In Case~3-1, with~$f(M)=\Gamma^b_C$ with~$r>0$, the metric~$g$ 
on the region~$-2r<H<r$ in the $H\theta$-plane is of the form
\begin{equation*}
g  = \frac{dH^2}{4(H{-}r)^2(H{+}2r)} + 4(H{-}r)^2(H{+}2r)\,d\theta^2\,.
\end{equation*}
Because
\begin{equation*}
\int_0^r \frac{dH}{\,2\sqrt{(H{-}r)^2(H{+}2r)}\,} = \infty,
\end{equation*}
this metric is complete near the `edge' $H=r$.  However, since
\begin{equation*}
\int_{-2r}^0 \frac{dH}{\,2\sqrt{(H{-}r)^2(H{+}2r)}\,} <\infty,
\end{equation*}
the metric is not complete near the `edge' $H=-2r$.  In fact, making
the substitution~$H+2r = 3r\rho^2$, the metric takes the form
\begin{equation*}
g  
= \frac{d\rho^2 + \rho^2(1-\rho^2)^4\,(18r^2\,d\theta)^2}{(1-\rho^2)^2}\,,
\end{equation*}
and one recognizes that $g$ will extend to a smooth metric at~$\rho=0$
in polar coordinates~$(\rho,\theta)$ on the disk~$\rho<1$ 
if and only if $\theta$ is taken to be periodic with period~$\pi/(9r^2)$.
This disk endowed with this complete metric is conformally equivalent 
to~$\bbC$.  The Gaussian curvature decreases monotonically from~$24r$ 
at~$\rho=0$ to a limiting value of~$-12r$ as~$\rho$ 
approaches~$1$.

Finally, consider Case 4 with image in~$\Gamma^1_C$.  
Let~$r_0>r_1>r_2$ be the three roots
satisfying~$r_0 = -(r_1+r_2)$, so that~$H^3+C_2\,H+C_3 
= (H{-}r_0)(H{-}r_1)(H{-}r_2)$. Consider the metric on the strip~
$r_2<H<r_1$ in the $H\theta$-plane given by
\begin{equation*}
g 
= \frac{dH^2}{4(H-r_0)(H-r_1)(H-r_2)} + 4(H-r_0)(H-r_1)(H-r_2)\,d\theta^2\,.
\end{equation*}
Since
\begin{equation*}
\int_{r_2}^{r_1} \frac{dH}{2\sqrt{(H-r_0)(H-r_1)(H-r_2)}} <\infty,
\end{equation*}
this metric is not complete at either edge~$H=r_i$ for~$i=1,2$.

Letting~$H=r_2 + v^2$ and computing as above, one finds that the metric
will extend to a smooth metric on a disk about~$v=0$ in $(v,\theta)$
polar coordinates if and only if~$\theta$ is taken to be periodic with
period
\begin{equation*}
\tau_2 = \frac{\pi}{3{r_2}^2+C_2}>0\,.
\end{equation*}
Similarly, setting~$H = r_1 - w^2$, and computing as above,
one finds that the metric will extend to a smooth metric on a disk 
about~$w=0$ in $(w,\theta)$ polar coordinates if and only if~$\theta$ 
is taken to be periodic with period
\begin{equation*}
\tau_1 = \frac{-\pi}{3{r_1}^2+C_2}>0\,.
\end{equation*}
Now, computation shows that~$\tau_1 = \tau_2$ has no
solutions with~$r_1>r_2$.  

\emph{Consequently, there is no complete Bochner-K\"ahler metric 
on a surface whose moduli image is~$\Gamma^1_C$.}

However, complete Bochner-K\"ahler metrics on orbifolds do exist:  
Taking~$r_1 = r(q{-}2p)$ and $r_2 = r(p{-}2q)$ where~$0<p<q$ 
are relatively prime integers and~$r$ is a positive real number, 
one can choose a period for~$\theta$ so that the resulting quotient 
completes to an orbifold metric on~$S^2$ with one conical point of 
order $1/q$ and the other of order~$1/p$.  

This orbifold is the weighted projective 
line~$\bbC\bbP^{(p,q)}$, i.e., $\C{2}$ minus the origin modulo
the $\bbC^*$-action~$\lambda\cdot(z,w)=\bigl(\lambda^pz,\lambda^qw\bigr)$.
This compact Riemannian orbifold could reasonably be regarded as the 
natural complete model for this case.  Note that the Gaussian curvature 
of this metric will be strictly positive if and only if~$q<2p$.

\subsection{Infinitesimal symmetries}
\label{ssec: inftess symms}

It turns out that any Bochner-K\"ahler structure 
has a nontrivial symmetry pseudo-groupoid~$\bar\Gamma$. 
In this subsection, some useful information about the `dimension'
and orbits of~$\bar\Gamma$ will be collected.

For~$(h,t,v)\in i\euu(n){\oplus}\C{n}{\oplus}\bbR$, 
let~$G^0_{(h,t,v)}\subset \Un(n)$ be the stabilizer of~$(h,t,v)$ 
under the action defined in~\S\ref{ssec: local mods}.  
Since~$a\in~\Un(n)$ lies
in~$G^0_{(h,t,v)}$ if and only if~$aha^* = h$ and $at = t$, it 
follows that $G^0_{(h,t,v)}$ is a closed, connected subgroup of~$\Un(n)$.

In fact, $G^0_{(h,t,v)}$ is a product of unitary groups and can be described 
as follows:  Let~$h_1>h_2>\cdots> h_\delta$ be the distinct eigenvalues 
of~$h$ and, for~$1\le\alpha\le n$, 
let~$L_\alpha\subseteq\C{n}$ be the $h_\alpha$-eigenspace 
of~$h$.  Since~$h$ is Hermitian symmetric, there is an orthogonal direct
sum decomposition
\begin{equation*}
\C{n} = L_1\oplus L_2\oplus \cdots \oplus L_\delta
\end{equation*}
with $\dim L_\alpha = n_\alpha\ge 1$. 
Write~$t = t_1 + \cdots + t_\delta$ where $t_\alpha$ lies in~$L_\alpha$
and let~$t_\alpha^\perp\subseteq L_\alpha$ be the subspace of~$L_\alpha$
that is perpendicular to~$t_\alpha$.  Then, using obvious notation,
\begin{equation*}
G^0_{(h,t,v)} 
= \Un(t^\perp_1)\times\Un(t^\perp_2)\times\cdots\times\Un(t^\perp_\delta).
\end{equation*}
The uniqueness part of Cartan's Theorem~\ref{thm: Cartan existence}  
then has the following useful corollary.

\begin{corollary}\label{cor: sym stab alg}  
Let~$P\to M$ be a Bochner-K\"ahler structure.
Then for any~$u\in P_x$, the unitary isomorphism $u:T_xM\to\C{n}$
induces an isomorphism
\begin{equation*}
\bar\Gamma_x \simeq G^0_{(H(u),T(u),V(u))} .
\end{equation*}
Thus, $\bar\Gamma_x$ is isomorphic to a product of unitary
groups and, in particular, is connected. 
\end{corollary}

\subsubsection{Existence and lower bounds}
\label{sssec: exist and lower bds}

Roughly speaking, a Bochner-K\"ahler structure has at least an
$n$-dimensional `infinitesimal symmetry group'.  As will be 
seen below, this lower bound is reached for the `generic' 
Bochner-K\"ahler structure.

\begin{theorem}\label{thm: sym alg of dim n }
Let~$M$ be a simply-connected complex $n$-manifold endowed with
a Bochner-K\"ahler structure~$\Omega$.  Let~$\eug\subset\euX(M)$ denote
the Lie algebra of vector fields on~$M$ whose flows preserve the complex
structure and~$\Omega$.  Then $\dim_\bbR\eug\ge n$.
\end{theorem}

\begin{proof}
Let~$(M,\Omega)$ satisfy the assumptions of the theorem, let~$\pi:P\to M$ 
be the unitary coframe bundle, with canonical forms~$\omega$
and~$\phi$, and let~$(H,T,V):P\to i\euu(n){\oplus}\C{n}{\oplus}\bbR$ be 
the structure function.

Because~$M$ is simply-connected and the Bochner-K\"ahler structure~$\Omega$
is real-analytic, any symmetry vector field of the structure defined
on a connected open subset of~$M$ can be uniquely analytically
continued to a symmetry vector field on all of~$M$.  
Moreover if~$Z\in\euX(M)$ is such a symmetry 
vector field, then, as discussed in~\S\ref{ssec: Un coframe bundle}, 
there is a unique vector field~$Z'$ on~$P$ satisfying~$\pi'(Z')=Z$ 
and $\Lie_{Z'}\omega = \Lie_{Z'}\phi=0$.  Conversely, if~$Y$
is a vector field on~$P$ satisfying $\Lie_Y\omega = \Lie_Y\phi=0$, 
then~$Y = Z'$ where $Z = \pi'(Y)$ is a symmetry vector field on~$M$.  

In other words, the mapping~$Z\mapsto Z'$ defines an 
embedding~$\eug\hookrightarrow\euX(P)$ that realizes~$\eug$ as the
Lie algebra of vector fields on~$P$ whose flows preserve the 
coframing~$\eta = (\omega,\phi)$.  
By the structure equations~\eqref{eq: structure equations ii}
the flow of such a vector field must necessarily preserve the structure
function~$(H,T,V): P\to i\euu(n){\oplus}\C{n}{\oplus}\bbR$, which 
is a submersion onto its (connected) image.
 
Applying Cartan's Theorem~\ref{thm: symmetry group} 
(see the Appendix), for any~$u\in P$ 
the evaluation map $e_u:\eug\to T_uP$ defined by~$e_u(Z) = Z'(u)\in T_uP$ 
is a vector space isomorphism between~$\eug$ and the kernel 
of~$(H,T,V)'(u):T_uP\to i\euu(n){\oplus}\C{n}{\oplus}\bbR$.  
Let~$K_u\subset T_uP$ denote this kernel.  Then, 
by~\eqref{eq: structure equations ii}, the
image~$(\omega,\phi)(K_u)\subset \C{n}{\oplus}\,i\euu(n)$ consists of
the pairs~$(w,f)\in \C{n}{\oplus}\,i\euu(n)$ that satisfy
\begin{equation*}
\begin{split}
0 &= H(u)\,f - f\,H(u) + T(u)\,w^* + w\,T(u)^*\,,\\
0 &= -f\,T(u) + \bigl({H(u)}^2+(\tr H(u))\,H(u) + V(u)\,\I_n\bigr)\,w\,,\\
0 &= (\tr H(u))\bigl(T(u)^*\,w+w^*T(u)\bigr)+T(u)^*H(u)\,w+w^* H(u) T(u).\\
\end{split}
\end{equation*}
By the first of these equations, 
\begin{equation*}
T(u)^*\,w + w^*T(u) = \tr\bigl([f,H(u)]\bigr) = 0
\end{equation*}
and
\begin{equation*}
2(T(u)^*H(u)\,w + w^* H(u) T(u)) 
= 2\tr\bigl([f,H(u)]H(u)\bigr) =  \tr\bigl([f,{H(u)}^2]\bigr) = 0.
\end{equation*}
Thus, the third equation is a consequence of the first and so
can be ignored for the rest of this discussion.

Let~$(H,T,V)(u_0) = (H_0,T_0,V_0)\in i\euu(n){\oplus}\C{n}{\oplus}\bbR$.
By $\Un(n)$-equivariance, it is enough to 
show that the dimension of $K_{u_0}$ is at least~$n$
at any point~$u_0$ where~$H_0$ and $T_0$ are both real, so assume this
for the rest of the argument.

By the structure equations~\eqref{eq: structure equations ii}, 
the dimension of $K_{u_0}$ 
is equal to the dimension of the space of solutions
of the linear equations
\begin{equation}\label{eq: the kernel K}
\begin{split}
0 &= H_0\,f - f\,H_0 + T_0\,w^* + w\,T_0^*\\
0 &= -f\,T_0 + \bigl({H_0}^2+(\tr H_0)\,H_0 + V_0\,\I_n\bigr)\,w\\
\end{split}
\end{equation}
for~$w\in\C{n}$ and~$f\in\euu(n)$.  

Consider the solutions of~\eqref{eq: the kernel K} for which
$f$ and~$w$ are purely imaginary, i.e., where~$f = i s$ and $w = i y$
for some symmetric (real) matrix~$s$ and some~$y\in\bbR^n$.  Then the
equations in~\eqref{eq: the kernel K} reduce to
\begin{equation}\label{eq: special K}
\begin{split}
0 &= H_0\,s - s\,H_0 - T_0\,{}^ty + y\,{}^tT_0\,,\\
0 &= -s\,T_0 + \bigl({H_0}^2+(\tr H_0)\,H_0 + V_0\,\I_n\bigr)\,y\,.\\
\end{split}
\end{equation}
The right hand side of the first equation of~\eqref{eq: special K} 
takes values in~$\euso(n)$
and the right hand side of the second equation of~\eqref{eq: special K} 
takes values in~$\bbR^n$.
Thus, this is $\frac12n(n{-}1)+n$ equations for the $\frac12n(n{+}1)+n$
components of~$s={}^ts$ and~$y$.  Consequently, the space of solutions is
at least of dimension~$n$. 
\end{proof}

\subsubsection{The symmetry algebra}
\label{sssec: sym alg}

For~$(h,t,v)\in i\euu(n){\oplus}\C{n}{\oplus}\bbR$, 
let~$\eug_{(h,t,v)}\subset \C{n}{\oplus}\euu(n)$ 
be the space of solutions~$(w,f)\in \C{n}{\oplus}\euu(n)$
of the linear equations
\begin{equation}\label{eq: define eug}
\begin{split}
0 &= h\,f - f\,h + t\,w^* + w\,t^*\\
0 &= -f\,t + \bigl({h}^2+(\tr h)\,h + v\,\I_n\bigr)\,w\\
\end{split}
\end{equation}
As was established in the course of the above proof,
$\eug_{(h,t,v)}$ is isomorphic as a vector space to the 
symmetry algebra~$\eug$ of any simply-connected 
Bochner-K\"ahler manifold whose structure function assumes the
value~$(h,t,v)$.  In fact, the structure 
equations~\eqref{eq: structure equations ii}
show that, if~$\bigl(H(u),T(u),V(u)\bigr)=(h,t,v)$,
then the vector space isomorphism~$\eug\to \eug_{(h,t,v)}$ defined
by~$X\mapsto\bigl(\omega_u(X'),\phi_u(X')\bigr)$ induces a Lie algebra
structure on~$\eug_{(h,t,v)}$ that is given by the formula
\begin{equation*}
\bigl[(x,x'),(y,y')\bigr] = \bigl(-x'y+y'x,\,-[x',y']+\{x,y\}_h\bigr)
\end{equation*}
where, for~$x,y\in\C{n}$, the element~$\{x,y\}_h$ in~$\euu(n)$ is 
defined by
\begin{equation*}
\begin{split}
\{x,y\}_h &= {} - h\,(xy^* - yx^*) - (xy^* - yx^*)\,h + (x^*y-y^*x)\,h\\
    &\qquad         + (x^*hy-y^*hx)\,\I_n
           +(\tr h)\bigl((x^*y-y^*x)\,\I_n-xy^*+yx^*\bigr).\\
\end{split}
\end{equation*}

For~$x\in M$, let~$\eug_x\subset \eug$ denote the subalgebra that
consists of the vector fields in~$\eug$ that vanish at~$x$.  Under
the vector space isomorphism $\eug\to \eug_{(h,t,v)}$ defined above,
$\eug_x$ maps into the subalgebra~$\eug^0_{(h,t,v)}\subset\eug_{(h,t,v)}$
defined by~$w=0$.

Since information about~$\eug^0_{(h,t,v)}$ and~$\eug_{(h,t,v)}$ 
will be needed later, these spaces will now be described more fully.

Fix~$(h,t,v)\in i\euu(n){\oplus}\C{n}{\oplus}\bbR$.
Suppose that~$h_1>h_2>\cdots>h_\delta$ are the distinct eigenvalues
of~$h$, that~$h$ has~$L_\alpha\subset\C{n}$ as 
its~$h_\alpha$-eigenspace, and that~$n_\alpha\ge1$ is the
(complex) dimension of~$L_\alpha$.  Write
\begin{equation*}
t = t_1 + \cdots+t_\delta
\end{equation*}
where~$t_\alpha$ lies in~$L_\alpha$.  Define the quantities
\begin{equation*}
v_\alpha = {h_\alpha}^2 + (\tr h)\,h_\alpha + v 
         +\sum_{\beta\not=\alpha}\frac{|t_\beta|^2}{(h_\alpha - h_\beta)}\,.
\end{equation*}
Now define
\begin{equation*}
\tau_\alpha
= \begin{cases}
  1 & \\
  0 & \\
  2n_\alpha& \\
  \end{cases}
\quad \text{and}\qquad
\rho_\alpha
= \begin{cases}
  (n_\alpha{-}1)^2 &\qquad \text{if $t_\alpha\not=0$;}\\
  (n_\alpha)^2 &\qquad \text{if $t_\alpha = 0$ but~$v_\alpha\not=0$;}\\
  (n_\alpha)^2 &\qquad \text{if $t_\alpha = 0$ and~$v_\alpha=0$.}\\
  \end{cases}
\end{equation*}

\begin{proposition}\label{prop: dim sym alg} 
For any~$(h,t,v)\in i\euu(n){\oplus}\C{n}{\oplus}\bbR$,
\begin{equation*}
\dim \eug^0_{(h,t,v)} = \rho_1 + \cdots + \rho_\delta
\end{equation*}
and
\begin{equation*}
\dim\eug_{(h,t,v)} = \dim\eug^0_{(h,t,v)}\ +\ \tau_1 + \cdots + \tau_\delta\,.
\end{equation*}
\end{proposition}

\begin{proof}  Because all the integers involved are invariant
under the action of~$\Un(n)$, it suffices to prove this formula
in the case that~$(h,t,v)$ lies in~$C$.  Maintaining the
notation introduced above, this means that
\begin{equation*}
h = \begin{pmatrix} 
        h_1\,\I_{n_1} & 0 & \cdots & 0\\
        0 & h_2\,\I_{n_2} & \cdots & 0\\
        \vdots & \vdots & \ddots & \vdots \\
        0 & 0 & \cdots & h_\delta\,\I_{n_\delta} \\
     \end{pmatrix}\quad\text{and}\quad
t = \begin{pmatrix} t_1\\ t_2\\ \vdots\\ t_\delta\end{pmatrix},
\end{equation*}
where~$t_\alpha$ takes values in~$\bbR^{n_\alpha}$ for~$1\le\alpha\le\delta$ 
and has all of its entries equal to zero except possibly the top one, 
which is nonnegative. 

For~$f\in\euu(n)$, write~$f$ in `block' form as~$f = (f_{\alpha\bar\beta})$ 
where~$f_{\alpha\bar\beta} = -{f_{\beta\bar\alpha}}^*$ takes values in
$n_\alpha$-by-$n_\beta$ complex matrices for~$1\le \alpha,\beta\le\delta$.
Correspondingly, write~$w\in\C{n}$ in `block' form as~$w = (w_\alpha)$ 
where~$w_\alpha$ takes values in~$\C{n_\alpha}$.
Then the first equation of~\eqref{eq: define eug} breaks into blocks as
\begin{equation*}
0 = (h_\alpha{-}h_\beta) f_{\alpha\bar\beta} 
         + t_\alpha\,{w_\beta}\!^* + w_\alpha\,{t_\beta}\!^*.
\end{equation*}
When~$\alpha=\beta$, this forces~
$t_\alpha\,{w_\alpha}\!^* + w_\alpha\,{t_\alpha}\!^* = 0$, so that either
$t_\alpha=0$, in which case this places no restriction on~$w_\alpha$
or else~$t_\alpha\not=0$, in which case $w_\alpha$ must be a purely
imaginary multiple of~$t_\alpha$, say $w_\alpha = i\,r_\alpha\,t_\alpha$
for some~$r_\alpha\in\bbR$.  In either case,~${w_\alpha}\!^*\,t_\alpha$ 
is purely imaginary.

When~$\alpha\not=\beta$, the above equation can be written as
\begin{equation*}
f_{\alpha\bar\beta}  
= \frac{t_\alpha\,{w_\beta}^* + w_\alpha\,{t_\beta}^*}{(h_\beta-h_\alpha)},
\qquad \alpha\not=\beta.
\end{equation*}
Substituting this equation into the second equation 
of~\eqref{eq: define eug} yields
\begin{equation*}
0 = -f_{\alpha\bar\alpha}\,t_\alpha
-\sum_{\beta\not=\alpha}
\frac{\ t_\alpha\,{w_\beta}^*+w_\alpha\,{t_\beta}^*\ }
                             {(h_\beta-h_\alpha)}\,t_\beta 
   +({h_\alpha}^2 + (\tr h)\,h_\alpha + v )\,w_\alpha\,,
\end{equation*}
which, by the definition of~$v_\alpha$ and the purely
imaginary nature of~${w_\beta}\!^*\,t_\beta$, can be written in the
form
\begin{equation*}
0 = -f_{\alpha\bar\alpha}\,t_\alpha
+\left(\sum_{\beta\not=\alpha}
\frac{{t_\beta}\!^*\,w_\beta}{(h_\beta{-}h_\alpha)} \right)\,t_\alpha
   +v_\alpha\,w_\alpha\,,
\end{equation*}

Now, if~$t_\alpha\not=0$, then this equation can be written in the
form
\begin{equation*}
0 = \left(-f_{\alpha\bar\alpha} 
+\left(iv_\alpha\,r_\alpha + \sum_{\beta\not=\alpha}
\frac{{t_\beta}\!^*\,w_\beta}{(h_\beta{-}h_\alpha)} \right) 
\,\frac{t_\alpha\,{t_\alpha}\!^*}{|t_\alpha|^2}
\right)\,t_\alpha\,,
\end{equation*}
which is $2n_\alpha{-}1$ real equations for 
the~${n_\alpha}^2$ entries of~$f_{\alpha\bar\alpha}$.  In fact, 
the solutions of this equation can be written in the form
\begin{equation*}
f_{\alpha\bar\alpha} = f^\perp_{\alpha\bar\alpha} + 
\left(iv_\alpha\,r_\alpha + \sum_{\beta\not=\alpha}
\frac{{t_\beta}\!^*\,w_\beta}{(h_\beta{-}h_\alpha)} \right) 
\,\frac{t_\alpha\,{t_\alpha}\!^*}{|t_\alpha|^2}
\end{equation*}
where~$f^\perp_{\alpha\bar\alpha}\in\euu(n_\alpha)$ is any solution 
to~$f^\perp_{\alpha\bar\alpha}\,t_\alpha = 0$, an equation that defines
the stabilizer subalgebra of~$t_\alpha$ in~$\euu(n_\alpha)$ and so 
has a solution space of dimension~$(n_\alpha{-}1)^2$.

If~$t_\alpha=0$, then the equation above simplifies to~$v_\alpha\,w_\alpha=0$.
If~$v_\alpha\not=0$, then this implies that~$w_\alpha=0$ while 
if $v_\alpha=0$, the equation degenerates to an identity.

In particular, it follows that the equations~\eqref{eq: define eug} 
impose no interrelations
among the~$w_\alpha$, just the condition~$w_\alpha = ir_\alpha\,t_\alpha$
if~$t_\alpha\not=0$, the condition~$w_\alpha=0$ if~$t_\alpha=0$ 
but~$v_\alpha\not=0$, and no condition on~$w_\alpha$ if~$t_\alpha=v_\alpha=0$.

Moreover, once the $w_\alpha$ have been chosen subject to these conditions,
the~$f_{\alpha\bar\beta}$ for~$\alpha\not=\beta$ are completely determined
while the~$f_{\alpha\bar\alpha}\in\euu(n_\alpha)$ are determined up to 
a choice of~$f^\perp_{\alpha\bar\alpha}$ if~$t_\alpha\not=0$ or are freely
specifiable if~$t_\alpha=0$. 

The desired dimension formulae follow immediately. 
\end{proof}

\subsubsection{Orbit dimension and slices}
\label{sssec: orbit dim n slices}

The proof of Proposition~\ref{prop: dim sym alg} shows how to compute 
the dimension of the $x$-orbit for any~$x\in M$ for which there
is a coframe~$u\in P_x$ with~$\bigl(H(u),T(u),V(u)\bigr) = (h,t,v)$.
Maintain the notation introduced above for the invariants of~$(h,t,v)$.  

Let~$O_x\subset T_xM$ be the tangent to the orbit through~$x$.
Then~$u(O_x)\subset\C{n}$ is the direct 
sum of the lines~$\bbR{\cdot}t_\alpha\subset L_\alpha$ for those~$\alpha$ 
with~$t_\alpha\not=0$ and the subspaces~$L_\alpha$ for those~$\alpha$ with
$t_\alpha=0$ and~$v_\alpha=0$. Thus, the dimension this of the $x$-orbit 
is equal to~$\tau_1+\cdots+\tau_\delta$.

A more interesting result is the calculation of a near-slice to the orbits
near~$x$.  For each~$\alpha$, let~$t^\perp_\alpha\subset L_\alpha$ be
the (complex) subspace~$\{\ w\in L_\alpha\,\mid\, t_\alpha^*w=0\,\}$.

If~$O^\perp_x\subset T_xM$ is the perpendicular to~$O_x$,
then~$u(O^\perp_x)\subset\C{n}$ is the direct sum 
of the subspaces~$i\bbR{\cdot}t_\alpha{\oplus}t^\perp_\alpha$ 
for those~$\alpha$ with~$t_\alpha\not=0$ together 
with the subspaces~$L_\alpha$ for those~$\alpha$ with~$t_\alpha=0$ 
and $v_\alpha\not=0$.

Now, from the description of~$\eug^0_{(h,t,v)}$, it follows that the
flows of the vector fields in~$\eug_x$ generate a group of rotations
about~$x$ that, via the unitary identification~$u:T_xM\to\C{n}$,
is carried isomorphically into the product of the unitary 
groups~$\Un(t^\perp_\alpha)$ for all~$\alpha$.  This is a closed 
subgroup of~$\Un(n)$ that evidently preserves the subspace~$u(O^\perp_x)$.

A near-slice to this action can be constructed as follows:  
For each~$\alpha$ with~$t_\alpha\not=0$ and~$t^\perp_\alpha=0$
(i.e., $n_\alpha = 1$), let~$S_\alpha
= i\bbR{\cdot}t_\alpha$.  If~$t_\alpha\not=0$ and~$t^\perp_\alpha\not=0$
(i.e., $n_\alpha>1$), choose a unit vector~$s_\alpha\in t^\perp_\alpha$
and let~$S_\alpha = i\bbR{\cdot}t_\alpha\oplus \bbR{\cdot}s_\alpha$.
If~$t_\alpha=0$ and~$v_\alpha\not=0$, choose a unit vector~$s_\alpha\in 
t^\perp_\alpha = L_\alpha$ and let~$S_\alpha=\bbR{\cdot}s_\alpha$.
Finally, if~$t_\alpha = 0$ and~$v_\alpha=0$, then set~$S_\alpha=0$.

Then the direct sum~$S_1\oplus\cdots\oplus S_\delta\subset u(O^\perp_x)$ 
is of the form~$u(S_x)$ for a subspace~$S_x\subset O^\perp_x$ that is
a near-slice to the action of the isometric rotations about~$x$
Consequently, the submanifold~$\exp_x(S_x)$ near~$x$ meets each orbit 
in a finite number of points and meets the generic orbit transversely.  
Let~$m_\alpha = \dim S_\alpha$, so that
\begin{equation*}
m_\alpha
= \begin{cases}
  2 &\quad \text{if $t_\alpha\not=0$ and $n_\alpha>1$;}\\
  1 &\quad \text{if $t_\alpha\not=0$ and $n_\alpha=1$;}\\
  1 &\quad \text{if $t_\alpha = 0$ but~$v_\alpha\not=0$;}\\
  0 &\quad \text{if $t_\alpha = 0$ and~$v_\alpha=0$.}\\
  \end{cases}
\end{equation*}
Then the `generic' orbit in~$M$ has codimension~$m = m_1+\cdots+m_\delta$.
Since~$m_\alpha\le n_\alpha$ for all~$\alpha$, it follows that~$m\le n$.

\subsubsection{Minimal symmetry}
\label{sssec: min sym}

By Proposition~\ref{prop: dim sym alg}, if~$(h,t,v)\in C$ 
satisfies~$h_{i\bar\imath}>h_{j\bar\jmath}$ 
for~$i<j$ and $t_i>0$ for all~$i$ then~$\dim \eug_{(h,t,v)}=n$.  
Thus, any simply-connected
Bochner-K\"ahler manifold whose structure function assumes such 
a~$(h,t,v)$ must have its symmetry algebra~$\eug$ be of 
dimension~$n$ exactly.  Moreover, from the above discussion, it
follows that the generic orbits of such a Bochner-K\"ahler structure
have codimension~$n$, the maximum possible. 

Thus, a `generic' Bochner-K\"ahler structure has its infinitesimal 
symmetry algebra of dimension~$n$ as well as cohomogeneity equal to~$n$.

\subsection{Constants of the structure}
\label{ssec: consts of str}

In the previous subsection, it was shown that the structure 
function~$(H,T,V):P\to i\euu(n){\oplus}\C{n}{\oplus}\bbR$ has rank
at most~$n^2{+}n$.  Since the image of~$P$ under this map is
$\Un(n)$-invariant, it is natural to look for a set of 
$\Un(n)$-invariant polynomials whose simultaneous level sets
will contain the images of structure functions.  In this section, 
I will exhibit $n{+}1$ such polynomials and show that they are
independent.

\subsubsection{Conserved polynomials}
\label{sssec: cons polys}

Let~$\Omega$ be any Bochner-K\"ahler structure on a connected
complex manifold~$M$ and let~$\pi:P\to M$ be the unitary coframe bundle, 
with canonical forms~$\omega$ and~$\phi$ and structure 
functions~$S$, $T$, and $V$ as above.  

By~\eqref{eq: structure equations ii}, there are identities
\begin{equation*}
\begin{split}
d(\tr H) &= (T^*\,\omega+\omega^*\,T),\\
d(\tr H^2) &= 2\,(T^*\,H\,\omega+\omega^*\,H\,T),\\
\end{split}
\end{equation*}
Thus, by the last equation of~\eqref{eq: structure equations ii}
\begin{equation*}
dV = (\tr H)\,d(\tr H) + {\ts\frac12}\,d(\tr H^2).
\end{equation*}
Since~$P$ is connected, there is a constant~$C_2$ for which
\begin{equation*}
V -  {\ts\frac12}\,\tr(H^2) - {\ts\frac12}\,(\tr H)^2 = C_2.
\end{equation*}

I am now going to show that this example can be generalized 
by constructing $n$ additional polynomials 
on~$i\euu(n){\oplus}\C{n}{\oplus}\bbR$ that have this constancy
property.  

Define~$A_k$ and~$B_k$ for~$k\ge0$ by the formulae
\begin{equation}\label{eq: define A and B}
\begin{split}
A_k &= \tr(H^k), \\
B_0 &= 1,\quad B_1 = \tr H,\quad B_2 = V, \\
B_k &= T^* H^{k-3} T,\quad k\ge3.\\
\end{split}
\end{equation}
Because these functions are constant on the fibers of~$\pi:P\to M$, 
they can be regarded as the pullbacks to~$P$ of well-defined smooth
functions on~$M$.  In what follows, I will usually treat them as 
functions on~$M$.  For convenience, define~$A_k=B_k=0$ when~$k<0$.

Also, for~$0\le k\le n$, let~$h_k$ denote the $k$-th elementary symmetric
function of the eigenvalues of~$H$.  These functions can
be expressed as polynomials in the~$A_k$ and hence are
smooth functions on~$M$.  For example,~$h_0 = 1$, $h_1 = A_1$,
$h_2 = \frac12({A_1}^2-A_2)$, etc.   For convenience, set~$h_k=0$
for~$k<0$ or~$k>n$.

\begin{theorem}\label{thm: BK constants}
For any connected Bochner-K\"ahler $n$-manifold~$(M,\Omega)$, 
the functions
\begin{equation}\label{eq: define C}
C_k =  B_k - h_1\,B_{k-1} +  h_2\,B_{k-2}
         - \cdots +(-1)^{k-1} h_{k-1}\,B_1 + (-1)^{k} h_{k}\,B_0\,.
\end{equation}
are locally constant for~$2\le k\le n+2$.
\end{theorem}

\begin{example}[Lowest constants]\label{ex: low consts}
For example, in addition to the evident constancy of the function
\begin{equation*}
C_2 = B_2 - h_1\,B_1+h_2\,B_0 
    = B_2 -  {\ts\frac12}\,A_2 - {\ts\frac12}\,{A_1}^2,
\end{equation*}
one has the constancy of
\begin{equation*}
C_3 = B_3 - h_1\,B_2+h_2\,B_1-h_3\,B_0 
    = B_3 - A_1\,B_2-{\ts\frac13}\,\bigl(A_3-{A_1}^3\bigr).  
\end{equation*}

The reader may notice that the above formula for~$C_k$ makes sense
for~$k=1$ and for~$k>n{+}2$.  Now, the expression~$C_1$
is just~$B_1 - h_1 = B_1 - A_1$, which vanishes by 
definition.  When~$k\ge n+3$, applying the Cayley-Hamilton theorem
to the definition of~$C_k$ yields
\begin{equation*}
C_{k} 
= T^*H^{k-n-3}\,(H^{n}-h_1\,H^{n-1}+h_2\,H^{n-2}-\cdots+(-1)^nh_n\,\I_n)T 
= 0.
\end{equation*}
However, when $2\le k\le n{+}2$, the expression~$C_k$ is a nontrivial
polynomial of weighted degree~$k$ in the variables~$A_j$ and~$B_j$. 
In fact, the above expressions for~$(C_2,\ldots, C_{n+2})$ 
can obviously be solved for~$(B_2,\ldots, B_{n+2})$.
\end{example}

\begin{proof} Define 1-forms~$\alpha_0 = 0$ and 
\begin{equation}\label{eq: define alpha}
\alpha_{k+1}\ = T^*H^k\omega + \omega^* H^k T,\qquad\text{for~$k\ge0$}.
\end{equation}
(The indexing is determined by `scaling weight' considerations.)
The~$\alpha_k$ are visibly $\pi$-semibasic, but they are also invariant
under the~$\Un(n)$-action on~$P$. Thus, they are the $\pi$-pullbacks 
of well-defined 1-forms on~$M$.  Consequently, they will, by
abuse of language, be treated as 1-forms on~$M$.

The first step will be to prove the following identities for all~$k\ge 0$:
\begin{equation}\label{eq: DEs for A and B}
\begin{split}
dA_{k} &= k\,\alpha_k\,,\qquad d\alpha_{k}  = 0\,,\\
dB_k & = B_0\,\alpha_k + B_1\,\alpha_{k-1} + B_2\,\alpha_{k-2}
                  + \cdots + B_{k-1}\,\alpha_1\,.\\
\end{split}
\end{equation}

Now, the first set of identities is just a calculation.  The case~$k=0$
is obvious, so assume~$k>0$.
Using~$\tr(PQ) = \tr(QP)$ and the structure equation, one has
\begin{equation*}
\begin{split}
dA_k &= d\bigl(\tr (H^{k})\bigr) = k \tr(H^{k-1}\,dH)\\
 &= k\tr\bigl(H^{k-1} (T\,\omega^* + \omega\, T^*)\bigr)
 = k\,\alpha_{k}\,. \\
\end{split}
\end{equation*}
(The terms in~$dH$ involving~$\phi$ cancel since~$A_{k}$ is 
constant on the fibers of~$\pi$.) Taking the exterior derivative of 
this relation and dividing by~$k$ yields
\begin{equation*}
0 = d\alpha_k\,.
\end{equation*}

The second set is a little more complicated, but still just a calculation.
The case~$k=0$ is trivial, and the case~$k=1$ follows from the
fact that~$B_1 = A_1$, so~$dB_1 = dA_1 = \alpha_1 = B_0\,\alpha_1$. 

Because $B_2 = V$, the second identity for~$k=2$ is just the 
structure equation for~$dV$.  Also, 
\begin{equation*}
\begin{split}
dB_3 &= dT^*\,T + T^* dT = \omega^*\bigl(H^2{+}(\tr H)\,H{+}V\,\I_n\bigr)T 
      + T^*\bigl(H^2{+}(\tr H)\,H{+}V\,\I_n\bigr)\,\omega\\
    &= B_0\,\alpha_3 +B_1\,\alpha_2 + B_2\,\alpha_1\,,
\end{split}
\end{equation*}
verifying the formula when~$k=3$. Thus, suppose from now
on that~$k>3$ and compute
(again ignoring terms involving~$\phi$, which must cancel)
\begin{equation*}
\begin{split}
dB_k &= dT^* \,H^{k-3}\,T +  T^*\,H^{k-3}\,dT
     + \sum_{l=0}^{k-4} T^*\,H^l\,dH\,H^{k-l-4}\,T
               \\
&= \omega^*\bigl(H^2{+}(\tr H)\,H{+}V\,\I_n\bigr) H^{k-3}\,T 
+  T^*\,H^{k-3} \bigl(H^2{+}(\tr H)\,H{+}V\,\I_n\bigr)\omega\\
&\qquad\qquad +\sum_{l=0}^{k-4} T^*\,H^l\,(T\,\omega^* + \omega\,T^*
                              )\,H^{k-l-4}\,T\\
&= \alpha_k + (\tr H)\,\alpha_{k-1}+ V\,\alpha_{k-2} 
       +\sum_{l=0}^{k-4} B_{l+3}\,\omega^*\,H^{k-l-4}\,T
                              + T^*\,H^l\,\omega\,B_{k-l-1}\\
&= B_0\,\alpha_k + B_1\,\alpha_{k-1}+ B_2\,\alpha_{k-2}
       +\sum_{l=0}^{k-4} B_{l+3}\,\alpha_{k-l-3}\\
&=B_0\,\alpha_{k+2} + B_1\,\alpha_{k-1} + B_2\,\alpha_{k} + B_3\,\alpha_{k-3}
                  + \cdots + B_{k-1}\,\alpha_1\,.\\
\end{split}
\end{equation*}
Thus, the formulae~\eqref{eq: DEs for A and B} are established.

Now, I claim that the functions~$h_k$ satisfy the differential equations
\begin{equation}\label{eq: DEs for h}
dh_k = h_{k-1}\,\alpha_1-h_{k-2}\,\alpha_2+\cdots+(-1)^{k+1}h_0\,\alpha_k\,.
\end{equation}
Granting~\eqref{eq: DEs for h} for the moment, computation gives
\begin{equation}\label{eq: dC vanishes}
\begin{split}
dC_k &= d\left(\sum_{j=0}^k (-1)^jh_jB_{k-j}\right)
     = \sum_{j=0}^k (-1)^j\bigl(dh_j\,B_{k-j}+h_j\,dB_{k-j}\bigr)\\
     &= \sum_{j=0}^k (-1)^j\left(
        \sum_{l=0}^j\,(-1)^{j-l+1}h_l\alpha_{j-l}\,B_{k-j}
 +\sum_{l=0}^{k-j}\,h_j\,B_{l}\alpha_{k-j-l}\right)\\
 &= \sum_{j=0}^k \sum_{l=0}^j\,(-1)^{-l+1}h_lB_{k-j}\,\alpha_{j-l}
 +\sum_{j=0}^k \sum_{l=0}^{k-j}\,(-1)^jh_jB_{l}\,\alpha_{k-j-l}\\
 &= \sum_{l=0}^k \sum_{j=l}^k\,(-1)^{-l+1}h_lB_{k-j}\,\alpha_{j-l}
 +\sum_{l=0}^k \sum_{j=0}^{k-l}\,(-1)^lh_lB_{j}\,\alpha_{k-j-l}\\
 &= \sum_{l=0}^k \sum_{j=l}^k\,(-1)^{-l+1}h_lB_{k-j}\,\alpha_{j-l}
 +\sum_{l=0}^k \sum_{j=l}^{k}\,(-1)^lh_lB_{k-j}\,\alpha_{j-l}\\
&=0.
\end{split}
\end{equation}

It remains to verify~\eqref{eq: DEs for h}.  
This is a classical identity:  Let
$\lambda_1,\ldots,\lambda_n$ be free variables, let~$s_k$ be
the $k$-th elementary function of the~$\lambda_i$, and let~$p_k$ be 
the $k$-th power function of the~$\lambda_i$, i.e., 
$p_k = {\lambda_1}^k+\cdots+{\lambda_n}^k$.  For any constant~$t$,
one has
\begin{equation*}
(1+s_1\,t + s_2\,t^2 + \cdots + s_n\,t^n) 
= (1+\lambda_1\,t)\cdots(1+\lambda_n\,t).
\end{equation*}
Taking the logarithm and then computing the differential of both sides
yields
\begin{equation*}
\frac{(ds_1\,t + ds_2\,t^2 + \cdots + ds_n\,t^n) }
{(1+s_1\,t + s_2\,t^2 + \cdots + s_n\,t^n) }
= \sum_{i=0}^n \frac{t\,d\lambda_i}{1+\lambda_i\,t}.
\end{equation*}
Expanding the right hand side out as a formal geometric power series in~$t$ 
and collecting like powers of~$t$ yields
\begin{equation*}
\frac{(ds_1\,t + ds_2\,t^2 + \cdots + ds_n\,t^n) }
{(1+s_1\,t + s_2\,t^2 + \cdots + s_n\,t^n) }
= \sum_{k=1}^\infty \left( \frac{(-1)^{k+1}}{k} dp_k\right) t^k\,.
\end{equation*}
It follows that
\begin{equation*}
\frac{(dh_1\,t + dh_2\,t^2 + \cdots + dh_n\,t^n) }
{(1+h_1\,t + h_2\,t^2 + \cdots + h_n\,t^n) }
= \sum_{k=1}^\infty \left( \frac{(-1)^{k+1}}{k} dA_k\right) t^k
= \sum_{k=1}^\infty (-1)^{k+1}\alpha_k\,t^k\,.
\end{equation*}
Thus,
\begin{equation}\label{eq: differential of ph}
dh_1\,t  + \cdots + dh_n\,t^n
= (1 + h_1\,t  + \cdots + h_n\,t^n)
 ( \alpha_1\,t -\alpha_2\,t^2 + \alpha_3\,t^3 - \cdots ),
\end{equation}
which, after equating coefficients of like powers of~$t$ on each side,
is~\eqref{eq: DEs for h}. 
\end{proof}

\subsubsection{The moduli map}
\label{sssec: mod map ii}

The map~$f:M\to F_n\subset\bbR^{2n+1}$ defined in~\S\ref{ssec: local mods} 
satisfies
\begin{equation*}
f = (A_1,\ldots,A_n,B_2,B_3,\ldots, B_{n+2}).
\end{equation*}

As is well known,%
\footnote{
In fact, this mapping can be computed by comparing like powers of~$t$
in the formal series expansions of the identity
\begin{equation*}
1+h_1\,t+\cdots+h_n\,t^n 
= \exp(A_1\,t -{\ts\frac12}A_2\,t^2 + {\ts\frac13}A_3\,t^3 - \cdots ).
\end{equation*} 
The mapping in the other direction can be computed by taking the
logarithm of both sides, expanding the left side as a series in~$t$,
and then comparing like powers.} 
there is a unique, invertible weighted-homogeneous
polynomial mapping~$\Sigma:\bbR^n\to\bbR^n$ that satisfies
\begin{equation*}
\Sigma(A_1,\ldots,A_n) = (h_1,\ldots,h_n).
\end{equation*}
 From the definition of the~$C_k$ as polynomials in~$B_j$ and~$h_j$,
it follows easily that~$\Sigma$
can be extended to an invertible, weighted-homogeneous 
polynomial mapping~$\Delta:\bbR^{2n+1}\to\bbR^{2n+1}$ so that
\begin{equation*}
\Delta{\circ}f = (h_1,\ldots,h_n,C_2,C_3,\ldots,C_{n+2}).
\end{equation*}
Thus, the fibers and the rank of the map~$\Delta{\circ}f:M\to\bbR^{2n+1}$ 
are the same as for the map~$f:M\to\bbR^{2n+1}$. In particular, the
fibers of~$\Delta{\circ}f$ are the orbits of the symmetry pseudo-groupoid
of the Bochner-K\"ahler structure. 

By Theorem~\ref{thm: BK constants}, when~$M$ is connected, the functions~$C_k$ 
are constants, so the fibers and rank of~$f$ are the same as the fibers and 
rank of the map~$h:M\to\bbR^n$ defined by
\begin{equation*}
h = (h_1,\ldots, h_n).
\end{equation*}

\subsection{Central symmetries}\label{ssec: central syms}

The symmetry algebra~$\eug$ of a Bochner-K\"ahler structure~$\Omega$ on a
connected~$M$ turns out to contain a canonical central subalgebra~$\euz$, 
whose dimension is equal to the infinitesimal cohomogeneity of the structure.

\subsubsection{An isometry vector field}\label{sssec: isom vf}

Consider the vector field~$Z'_2$ on~$P$ that satisfies
\begin{equation}\label{eq: define Z-2}
\omega(Z'_2) = 2i\,T,\qquad\qquad \phi(Z'_2) = 2i(H^2+\tr(H)\,H+V\,\I_n).
\end{equation}
(The indexing is dictated by scaling weight considerations.)
Since~$\omega{\oplus}\phi:T_uP\to\C{n}{\oplus}\euu(n)$ is an
isomorphism for all~$u\in P$ and since~$H$ is Hermitian symmetric
while~$V$ is real, this does indeed define a unique vector field~$Z'_2$.
Because of the $\Un(n)$-equivariance of the structure functions,
the vector field~$Z'_2$ is invariant under the right $\Un(n)$-action.  
Consequently, there is a unique vector field~$Z_2$ on~$M$ 
that is~$\pi$-related to~$Z'_2$, i.e., 
that satisfies~$\pi'\bigl(Z'_2(u)\bigr)=Z_2\bigl(\pi(u)\bigr)$. 
The structure equation~$dT+\phi\,T = (H^2{+}\tr(H)\,H+V\,\I_n)\,\omega$ 
coupled with the discussion in~\S\ref{sssec:  tensors vfs and symms} shows
that~$Z_2$ is the real part of a holomorphic vector field on~$M$.

Since
\begin{equation}\label{eq: Z-2 is Ham}
\pi^*(Z_2\lhk\Omega) = Z'_2\lhk\bigl(-{\ts\frac{i}2}\,\omega^*\w\omega\bigr)
= -\bigl(T^*\omega + \omega^* T\bigr) 
= -\, d\bigl(\tr(H)\bigr),
\end{equation}
the flow of~$Z_2$ is the $\Omega$-Hamiltonian flow associated to~$h_1$.

Alternatively, one can see directly that the (local) flow of~$Z_2$ preserves
the Bochner-K\"ahler structure.  It has already been observed that~$Z_2$
is $\pi$-related to the~$\Un(n)$-invariant vector field~$Z'_2$ on~$P$.  
The defining formulae for~$Z'_2$ yields
\begin{equation*}
\begin{split}
\Lie_{Z'_2}\omega &= d\bigl(\omega(Z'_2)\bigr)+ Z'_2\lhk d\omega
 = d\bigl(\omega(Z'_2)\bigr) + Z'_2\lhk (-\phi\w\omega) \\
 &= d(2i\,T) -\phi(Z'_2)\,\omega + \phi\,\omega(Z'_2)\\
 &= 2i\left(dT - (H^2{+}\tr(H)\,H{+}V\,\I_n)\,\omega +\phi\,T\right) = 0.\\
\end{split}
\end{equation*}
Thus, the flow of~$Z'_2$ preserves~$\omega$.  In turn, this implies
that the flow of~$Z'_2$ preserves~$\phi$ (since~$\phi$ is the unique
$\euu(n)$-valued 1-form that satisfies~$d\omega = -\phi\w\omega$).  Thus,
the vector field~$Z'_2$ is $\pi$-related to the symmetry vector field~$Z_2$
of the~$\Un(n)$-structure that~$P$ defines, i.e., the original 
Bochner-K\"ahler structure.

\begin{remark}[Matsumoto's observation]
To my knowledge it was Matsumoto~\cite[Theorem~2]{Ma} who first observed
that one could construct a holomorphic vector field on~$M$ 
by~$\Omega$-dualizing the exterior derivative of the scalar curvature,
at least when~$M$ is compact.  The vector field that he constructs
is, up to a constant complex multiple, the same as the one whose real part
is~$Z_2$.  The above argument shows that compactness actually 
plays no role; the holomorphicity of~$Z_2 - iJZ_2$ is a purely local fact.
Apparently, Matsumoto did not realize that some complex multiple of
his vector field had a real part whose flow of~$Z_2$ was not only 
holomorphic but isometric as well.
\end{remark}

\subsubsection{The central algebra}\label{sssec: central alg}

I am now going to show that~$Z_2$ is the first of a sequence 
of real parts of holomorphic vector fields on~$M$ whose representative
functions can be written down explicitly in terms of~$H$ and $T$.  For
example, the next term in the sequence will be seen to be 
the vector field~$Z_3$ whose representative function 
is~$z_3 = -2i\left(H {-}\,\tr(H)\,\I_n\right) T$.

\begin{theorem}\label{thm: generators of center}
For every~$k$ in the range~$0\le k\le n{-}1$, the function
\begin{equation*}
z_{k+2} = 2i(-1)^k\bigl(H^k -h_1\,H^{k-1}+ h_2\,H^{k-2}
    + \cdots + (-1)^k\,h_k\,\I_n\bigr)\,T
\end{equation*}
is the representative function of a vector field~$Z_{k+2}\in\eug$.
Moreover, the span~$\euz$ of~$Z_2,\ldots,Z_{n+1}$ lies in the
center of~$\eug$.  

For~$x\in M$, the subspace~$\euz_x =
\spn\{Z_2(x),\ldots,Z_{n+1}(x)\}$ is the $\Omega$-complement to~$\ker df_x$,
where~$f:M\to\bbR^{2n+1}$ is the moduli mapping 
of~\S{\upshape\ref{sssec: moduli map}}.

If~$M$ is connected, then~$\dim\euz\le n$ is the maximum over~$x\in M$
of~$\dim\euz_x$.  If~$\dim\euz=n$, then~$\eug=\euz$.
\end{theorem}

\begin{remark}
While the formula for~$z_{k+2}$ makes sense for all~$k\ge0$,
this expression vanishes identically when~$k\ge n$, due to the
Cayley-Hamilton Theorem.
\end{remark}

\begin{proof}
A computation like that done in the proof 
of~\eqref{eq: DEs for A and B} shows that, 
for~$k\ge0$,
\begin{equation}\label{eq: DEs for HT}
d(H^kT) + \phi\,H^kT 
\equiv H^kT\,\alpha_0 + H^{k-1}T\,\alpha_1 +\cdots + H^0T\,\alpha_k 
\mod \omega.
\end{equation}
Using this identity, a calculation analogous 
to~\eqref{eq: dC vanishes} yields
\begin{equation*}
d(z_{k+2}) + \phi\,z_{k+2} \equiv 0 \mod \omega.
\end{equation*}
Thus, according to~\S\ref{sssec:  tensors vfs and symms}, 
each~$z_{k+2}$ is the representative function 
of a vector field~$Z_{k+2}$ on~$M$ whose local flow is holomorphic.
 
Letting~$\tilde Z_{k+2}$ be any vector field on~$P$ so 
that~$\omega(\tilde Z_{k+2})= z_{k+2}$,
\begin{equation}\label{eq: Z-k is Ham}
\begin{split}
\pi^*\bigl(Z_{k+2}\lhk\Omega \bigr)
&= \tilde Z_{k+2}\lhk\bigl(-{\ts\frac{i}2}\,\omega^*\w\omega\bigr)
= -{\ts\frac{i}2}\,\bigl(z_{k+2}^*\w\omega-\omega^*\,z_{k+2}\bigr)\\
&= -(-1)^k\bigl( T^*(H^k -h_1\,H^{k-1}
    + \cdots + (-1)^k\,h_k\,\I_n)\,\omega \\
&\qquad\qquad {} +
 \omega^*(H^k -h_1\,H^{k-1}+\cdots+(-1)^k\,h_k\,\I_n)\,T\bigr)\\
&= (-1)^{k+1}\bigl(\alpha_{k+1} - h_1\,\alpha_k + h_2\,\alpha_{k-1}
      + \cdots+ (-1)^kh_k\,\alpha_1\bigr)\\
&=  -dh_{k+1} = \pi^*\bigl(-dh_{k+1}\bigr)\,,\\
\end{split}
\end{equation}
by virtue of~\eqref{eq: DEs for h}.  Thus, $Z_{k+2}\lhk\Omega = - dh_{k+1}$,  
so that~$Z_{k+2}$ is the $\Omega$-Hamiltonian vector field 
associated to~$h_{k+1}$.

Since the flow of~$Z_{k+2}$ is both holomorphic and symplectic, 
$Z_{k+2}$ belongs to~$\eug$, as claimed.

Since the representative function~$z_{k+2}$ is constructed as a 
polynomial in~$H$ and~$T$, which are invariant under the $Y'$-flow on~$P$
for any vector field~$Y\in\eug$, it follows that~$Z_{k+2}$ is
invariant under the flow of any~$Y\in\eug$, i.e.,~$[Y,Z_{k+2}] = 0$
for any~$Y\in\eug$.  Thus, the~$Z_{k+2}$ for~$0\le k\le n{-}1$ 
span a central subalgebra~$\euz\subseteq\eug$.

For any~$x\in M$, the nondegeneracy of~$\Omega$ implies 
that~$\euz_x\subset T_xM$ is 
the $\Omega$-dual of
\begin{equation*}
\spn\{dh_1(x),\ldots,dh_n(x)\} = \spn\{dA_1(x),\ldots,dA_n(x)\}
\end{equation*}
(see \S\ref{sssec: mod map ii}).
The map~$f:M\to\bbR^{2n+1}$ has components given by~$A_i$ and~$B_i$. 
By Theorem~\ref{thm: BK constants} and~\S\ref{sssec: mod map ii}, 
each~$B_i$ can be written on each connected 
component of~$M$ as a weighted polynomial in~$A_1,\ldots,A_i$ with constant 
coefficients. Thus, the kernel of~$df_x$ is 
the same the kernel of~$(dA)_x$ where~$A = (A_1,\ldots, A_n)$,
which establishes the stated $\Omega$-complementarity. 

Now suppose that~$M$ is connected.  For each~$k\ge0$, each~$x\in M$,
and each~$u\in P_x$,
\begin{equation*}
u\bigl(\euz_x\bigr) 
= \spn\bigl\{\ iH(u)^k\,T(u)\mid\ 0\le k\le n{-}1\ \bigr\}.
\end{equation*}
Since~$H(u)$ is Hermitian symmetric, it follows that the dimension 
of~$\euz_x$ over~$\bbR$ is the largest integer~$m_x\le n$ so that the 
vectors~$T(u), H(u)T(u),\ldots, H(u)^{m_x-1} T(u)$
are linearly independent (over either~$\bbC$ or~$\bbR$) 
and, moreover, that $\euz_x\,\cap\,J\euz_x = (0)_x$.

Let~$M^\circ\subset M$ be the (nonempty) open set consisting of those~$x\in M$
for which~$m_x$ achieves the maximum value~$m\le n$.  If~$m=n$, 
then~$\dim\euz_x=n$ for all~$x\in M^\circ$ and, 
since~$\dim\euz_x\le\dim\euz\le n$ for all~$x\in M$, 
it follows that~$\dim\euz = n$.  If~$m=0$, then~$T$ vanishes identically,
implying that~$\euz = (0)$.  If~$0<m<n$, let~$k$ satisfy~$m\le k<n$.
Then, on~$M^\circ$, the vector field~$Z_{k+2}$
is a linear combination of the independent vector 
fields~$Z_2,\ldots,Z_{m+1}$.  Thus, there exist
smooth real-valued functions~$w_2,\ldots w_{m+1}$ on~$M^\circ$ so that
\begin{equation*}
Z_{k+2}= w_2\,Z_{2} + \cdots + w_{m+1}\,Z_{m+1}.
\end{equation*}
Consequently,
\begin{equation*}
Z_{k+2}{-}iJZ_{k+2} 
= w_2\,(Z_{2}{-}iJZ_{2}) + \cdots + w_{m+1}\,(Z_{m+1}{-}iJZ_{m+1}).
\end{equation*}
However, the left hand side of this equation is a holomorphic 
vector field while the holomorphic vector 
fields~$(Z_{2}{-}iJZ_{2}), \ldots,(Z_{m+1}{-}iJZ_{m+1})$ are 
linearly independent (over~$\bbC)$ at each point of~$M^\circ$.  It follows
that the functions~$w_2,\ldots,w_{m+1}$ are real-valued holomorphic
functions on~$M^\circ$ and hence must be constants.  Since~$M$ is 
connected, the identity
\begin{equation*}
Z_{k+2}= w_2\,Z_{2} + \cdots + w_{m+1}\,Z_{m+1}.
\end{equation*}
must hold on all of~$M$. In other words, the vector 
fields~$Z_2,\ldots,Z_{m+1}$ are a basis of~$\euz$, as desired.

If~$\dim\euz=n$, then at any point~$x\in M^\circ$, 
the vectors~$Z_2(x),\ldots,Z_{n+1}(x)$
are linearly independent. If~$\pi(u)=x$, 
then~$\{\,H^k(u) T(u)\ \mid\ 0\le k\le n{-}1\ \}$ are linearly
independent, implying (since~$H(u)$ is diagonalizable) that~$T(u)$ does not 
lie any sum of fewer that~$n$ distinct eigenspaces of~$H(u)$.  
By~\S\ref{sssec: min sym}, this implies that the differential of the 
mapping~$(H,T,V):P\to i\euu(n){\oplus}\C{n}{\oplus}\bbR$ has kernel 
of dimension equal to~$n$, so $\dim\eug = n$, as desired. 
\end{proof}

\subsubsection{The momentum mapping}\label{sssec: mom map}

The proof of Theorem~\ref{thm: generators of center} 
shows that the map~$h:M\to\bbR^n$ defined by
\begin{equation}\label{eq: define h-tot}
h = \bigl(h_1,h_2,h_3,\ldots,h_n\bigr)
\end{equation}
is a momentum mapping for the infinitesimal torus
action generated by~$\euz$.%
\footnote{More precisely, the infinitesimal $n$-torus action on~$M$ is
given by the Lie algebra homomorphism~$\bbR^n\to\euz\subset\euX(M)$
defined by the explicit generators~$Z_2,\ldots,Z_{n+1}$.}  

As already remarked, there is an invertible 
weighted-homogeneous polynomial mapping~$\Delta:\bbR^{2n+1}\to\bbR^{2n+1}$ 
that satisfies
\begin{equation*}
\Delta\bigl(A_1,\ldots,A_n,B_2,\ldots,B_{n+2}\bigl) 
= \bigl(h_1,\ldots,h_n,C_2,\ldots, C_{n+2}\bigr).
\end{equation*}
Thus, the fibers of~$\Delta\circ f:M\to\bbR^{2n+1}$
are the orbits of the symmetry pseudo-groupoid of the underlying
Bochner-K\"ahler structure. By Theorem~\ref{thm: BK constants}, 
if~$M$ is connected, 
then the functions~$C_k$ are constant.  Thus, for connected~$M$, 
the fibers of~$h$ are the orbits of this symmetry pseudo-groupoid.

\subsection{Cohomogeneity and the momentum polynomial}
\label{ssec: coho mo po}

Assume that~$M$ is connected and that~$\dim\euz = m$.

\subsubsection{Cohomogeneity}\label{sssec: cohomo}

The proof of Theorem~\ref{thm: generators of center} shows that $dh_k$ 
for~$k > m$ is
a constant linear combination of the differentials~$dh_1,\ldots, dh_m$
and that these latter $1$-forms are linearly independent on an open
subset of~$M$.  Thus, $h(M)$ lies in an $m$-dimensional 
affine subspace~$\eua\subset\bbR^n$ and, moreover, 
contains an open subset of~$\eua$.  
Since the fibers of~$h$ are the orbits of the symmetry pseudo-groupoid,
it is reasonable to call the number~$m$ the \emph{cohomogeneity} 
of the Bochner-K\"ahler structure.

\subsubsection{Cohomogeneity}\label{sssec: mo po}

Let~$t$ be a parameter and define the \emph{momentum polynomial}~$p_h(t)$ 
of~$M$ by the formula
\begin{equation*}
p_h(t) = t^n - h_1\,t^{n-1} + \cdots + (-1)^n h_n\,.
\end{equation*}
Of course, $p_h(t)= \det(t\,\I_{n}-H)$ is the characteristic
polynomial of the Hermitian symmetric matrix~$H$, so all
of its roots are real. 

\begin{theorem}\label{thm: cohomogeneity m}
If~$(M,\Omega)$ is a connected Bochner-K\"ahler manifold of
cohomogeneity~$m$, then $n{-}m$ of the roots of~$p_h(t)$
are constant and, outside a closed, proper, complex analytic 
subvariety~$N\subset M$, the remaining~$m$ roots are
distinct, real-analytic, and functionally independent.
\end{theorem}

\begin{proof}
If~$m=0$, then, in particular, $h_1$ is constant, so~$(M,\Omega)$
has constant scalar curvature.  By Proposition~\ref{prop: BK loc symm}, 
$(M,\Omega)$ is
locally homogeneous, so all of the eigenvalues of~$H$ are constant.
(By Proposition~\ref{prop: BK loc symm}, there are at most two distinct 
eigenvalues.)  In this case, $N$ can be taken to be empty.

Suppose from now on that~$m>0$. Technically, I should treat the cases~$m=n$
and~$m<n$ separately, but the argument for~$m=n$ differs from that for
$m<n$ by trivial notational changes, so I will not explicitly 
assume~$m<n$ but, rather, let the reader make the necessary modifications 
for the case~$m=n$.

By Theorem~\ref{thm: generators of center}, the 
differentials~$dh_1,\ldots,dh_m$ are linearly
independent exactly where the vector fields~$Z_2,\ldots,Z_{m+1}$
are linearly independent.  In particular, the locus~$N\subset M$ where
$dh_1\w\cdots\w\,dh_m$ vanishes is also where the holomorphic
$m$-vector
\begin{equation*}
(Z_2-\imath JZ_2)\w (Z_3-\imath JZ_3)\w \cdots\w (Z_{m+1}-\imath JZ_{m+1})
\end{equation*} 
vanishes.  Thus, $N$
is a closed, proper, complex analytic subvariety of~$M$ and so has real
codimension at least~$2$.  Consequently, its complement~$M^\circ\subset M$ 
is a connected, open, dense subset of~$M$.  

Since~$\dim\euz_x = m$ for all~$x\in M^\circ$, a  
subbundle~$P_0\subset P$ can be defined over~$M^\circ$ by saying 
that~$u\in\pi^{-1}(M^\circ)$ lies in~$P_0$ if and only 
if~$u(\euz_x)= i\bbR^m\subset\C{m}\subset\C{n}$.  
Then~$\pi:P_0\to M^\circ$
is a smooth principal $\Or(m){\times}\Un(n{-}m)$-bundle over~$M^\circ$.

Pull the forms~$\omega$ and~$\phi$ and the functions~$H$, $T$,
and~$V$ back to~$P_0$. By definition, for every~$u\in P_0$, 
the vectors~$T(u),H(u)T(u),\ldots, H(u)^{m-1}T(u)$ 
span~$\bbR^m\subset\C{n}$ and, moreover, $H(u)\cdot\bbR^m\subset\bbR^m$.
Thus, there exists a
function~$T':P_0\to\bbR^m$, a function~$H'$ on~$P_0$
with values in the open set of symmetric~$m$-by-$m$ (real) 
matrices with $m$ distinct eigenvalues, and a
function~$H''$ on~$P_0$ with values in 
Hermitian symmetric~$(n{-}m)$-by-$(n{-}m)$ matrices so that
\begin{equation}\label{eq: reduced TH}
T = \begin{pmatrix} T'\\0\\\end{pmatrix},\qquad
H = \begin{pmatrix} H'&0\\0&H''\end{pmatrix}.
\end{equation}
Write $\phi=-\phi^*$ in $(m,n{-}m)$-block form as
\begin{equation}\label{eq: reduced phi}
\phi = \begin{pmatrix} \phi' & \tau^* \\ -\tau & \phi'' \\ \end{pmatrix}
\end{equation}
where, of course, ~$\phi'$ and ~$\phi''$ take values in skew-Hermitian
matrices of dimensions~$m$ and~$n{-}m$, respectively.  The lower right-hand
$(n{-}m)$-by-$(n{-}m)$ block of the $dH$ equation 
in~\eqref{eq: structure equations ii} then becomes
\begin{equation}\label{eq: H'' DE}
dH'' = -\phi''\,H'' + H''\,\phi''.
\end{equation}
Consequently, the eigenvalues of~$H''$ are constant on~$P_0$.  Let
\begin{equation*}
\det(t\,\I_{n-m}-H'') 
= p_{h''}(t) = t^{n-m} - h''_1\,t^{n-m} + \cdots + (-1)^{n-m}\,h''_{n-m}
\end{equation*}
be the characteristic polynomial of~$H''$, where the~$h''_i$ are constants.
Then, on~$M^\circ$ at least, $p_{h''}(t)$ divides~$p_h(t)$.  
Using the Euclidean algorithm, write
\begin{equation*}
p_h(t) = p_{h''}(t)\,q(t) + r(t)
\end{equation*}
where~$q$ and~$t$ are polynomials in~$t$ and where the degree of~$r$ is 
at most $n{-}m{-}1$.  The coefficients of~$q$ and~$r$ are constant
linear combinations of the coefficients in~$p_h$ and so are continuous.
Since the coefficients of~$r$ vanish on~$M^\circ$, which is dense in~$M$, it
follows that~$r$ vanishes identically on~$M$.  Thus, 
$p_{h''}(t)$ divides~$p_h(t)$ on all of~$M$.

Defining real-analytic functions~$h'_1,\ldots, h'_m$ on~$M$ by
\begin{equation*}
q(t) =  t^m - h'_1\,t^{m-1} +\cdots+ (-1)^m\,h'_m\,,
\end{equation*}
one sees that~$q(t) = p_{h'}(t) = \det(t\,\I_{m}-H')$ on~$M^\circ$, i.e.,
that~$p_h(t) = p_{h'}(t)\,p_{h''}(t)$.

Of course the roots of~$p_{h'}(t)$ on~$M^\circ$ 
equal the eigenvalues of~$H'$ on~$P_0$ and so are distinct and 
therefore real-analytic on~$M^\circ$. Since the~$h_i$ are constant
coefficient linear combinations of the~$h'_j$, it follows that
there is a constant~$a$ so that
\begin{equation*}
dh_1\w\cdots\w dh_m = a\,dh'_1\w\cdots\w dh'_m\,.
\end{equation*}
Obviously, $a$ is nonzero and $dh'_1\w\cdots\w dh'_m$ is nonvanishing
on~$M^\circ$.  Since the roots of~$p_{h'}(t)$ are distinct on~$M^\circ$, 
and the $h'_i$ are the elementary symmetric functions of these roots, it
follows that these roots must be functionally independent on~$M^\circ$,
as claimed.  
\end{proof}

\begin{remark}
Theorem~\ref{thm: cohomogeneity m} 
accounts for the $n{-}m$ constant coefficient 
linear relations among the momenta~$h_1,\ldots h_n$ 
implicit in the initial discussion.  They are just the $n{-}m$ coefficients
of the remainder polynomial~$r(t)$.
\end{remark}

\subsubsection{Reduced momentum}\label{sssec: red mo}

The mapping~$h' = (h'_1,\ldots,h'_m):M\to \bbR^m$ will be known as
the \emph{reduced momentum mapping} of~$M$.  
The proof of Theorem~\ref{thm: cohomogeneity m} shows 
that~$h':M^\circ\to \bbR^m$ is a submersion onto its image.

The polynomial~$p_{h'}(t)$ will be referred to as the \emph{reduced
momentum polynomial} of~$M$.  The roots of~$p_{h'}(t)$ are real 
on~$M^\circ$, which is dense in~$M$, so the roots of~$p_{h'}(t)$
are real at every point of~$M$. For each~$x\in M$, let
\begin{equation*}
\lambda_1(x)\ge \lambda_2(x) \ge \cdots \ge \lambda_m(x)
\end{equation*}
be the roots of~$p_{h'}(t)$, counted with multiplicity.   
By a standard argument based on the Stone-Weierstra\ss{} theorem,
the functions~$\lambda_i:M\to\bbR$ are continuous.%
\footnote{Note that $\lambda_i$ will be real-analytic even at~$x\in N$
as long as it is a simple root of~$p_{h'}(t)$ at~$x$.}  
Thus, the reduced momentum polynomial
factors continuously as
\begin{equation*}
p_{h'}(t) = (t-\lambda_1)(t-\lambda_2)\cdots(t-\lambda_m).
\end{equation*}

\begin{example}[Low cohomogeneity]\label{ex: low cohom} 
Suppose~$M$ is locally isometric 
to~$M^p_c\times M^{n-p}_{-c}$. 
Then, looking back at the proof of
Proposition~\ref{prop: BK loc symm} and the definition of~$H$, 
one computes that
\begin{equation*}
p_h(t) = \left(t + \frac{c(n-p+1)}{2(n+2)}\right)^{p}
         \left(t - \frac{c(  p+1)}{2(n+2)}\right)^{n-p}.
\end{equation*}
Since this example is locally homogeneous, i.e., $m=0$, it 
follows that~$p_h(t) = p_{h''}(t)$.

On the other hand, for Example~\ref{ex: rot sym}
(i.e., rotationally symmetric), 
\begin{equation*}
p_h(t) = \left(t - \frac{k}{(n+2)}\right)^{n-1}
         \left(t - \frac{k}{(n+2)}- a\,|z|^2f'\bigl(|z|^2\bigr)\right).
\end{equation*} 
As long as $a\not=0$, these examples have cohomogeneity~$m=1$, with
\begin{equation*}
p_{h''}(t) = \left(t - \frac{k}{(n+2)}\right)^{n-1}
\quad\text{and}\qquad
p_{h'}(t) = \left(t - \frac{k}{(n+2)}- a\,|z|^2f'\bigl(|z|^2\bigr)\right).
\end{equation*}
\end{example}

\section{Global Geometry and Symmetries}\label{sec: glo geo and sym}

Throughout this section it will be assumed that~$M$ is a connected
complex $n$-manifold endowed with a Bochner-K\"ahler structure~$\Omega$.
All the notation introduced earlier will be retained.

\subsection{The characteristic polynomials}\label{ssec: char polys}

In this section, two constant coefficient polynomials
will be introduced that are invariants of the analytically
connected equivalence class of the Bochner-K\"ahler structure.
Also a formula (Theorem~\ref{thm: point data to polys}) will be developed 
to compute them 
from the value of the structure function at a single point.

\subsubsection{The characteristic polynomial}\label{sssec: char poly i}

Let~$C_k$ for~$k = 2,\ldots,n{+}2$ be the constants introduced
in Theorem~\ref{thm: BK constants}. For the sake of convenience, 
set $C_0 = 1$ and $C_k=0$
for~$k=1$ and~$k>n{+}2$.  Let~$t$ be a real parameter.  
Then by Theorem~\ref{thm: BK constants} (plus the remark following it),
\begin{equation*}
\begin{split}
\sum_{k=0}^\infty C_k\,t^k 
  &= \sum_{k=0}^\infty \sum_{l=0}^\infty (-1)^l h_l\,B_{k-l}\,t^k \\
  &=\left( \sum_{k=0}^\infty (-1)^l h_l\,t^l\right)
     \left(\sum_{j=0}^\infty B_{j}\,t^j\right) \\
  &= \det(\I_n{-}t\,H) \left(1 + h_1\,t + V\,t^2 
          +  t^3 \sum_{k=3}^\infty T^*(tH)^{k-3}T\right)\\
  &= \det(\I_n{-}t\,H) \left(1 + h_1\,t + V\,t^2 
          +  t^3\,T^*(\I_n{-}tH)^{-1}T\right)\\
  &= \det(\I_n{-}t\,H)\,(1 + h_1\,t + V\,t^2) 
          +  t^3\,T^*\,\Cof(\I_n{-}tH)\,T\\
\end{split}
\end{equation*}
where $\Cof(\I_n{-}tH)$ is the signed cofactor matrix%
\footnote{The signed cofactor matrix of any $n$-by-$n$ matrix~$R$
is the (unique) homogeneous polynomial matrix of degree~$(n{-}1)$ 
that satisfies the identity~$R\,\Cof(R) = \det(R) \I_n$.}
of~$\I_n{-}tH$.

The cautious reader may object that the second factors on the second
and third lines need not converge for all~$t$.  However, every~$u\in P$ 
has an open neighborhood on which $T^*T$ and~$\tr(H^*H)$ are bounded,
so that the series is bounded by a geometric series and hence
converges for~$|t|$ sufficiently small. The upshot of this
is that the two series converge absolutely and uniformly 
on compact subsets of a certain open neighborhood 
of~$P\times{0}$ in~$P\times\bbR$, so equality of the first
and last terms holds on that open subset.  The left hand side is 
evidently a polynomial in~$t$ of degree at most~$n{+}2$ and the
final form of the right hand side is also a polynomial in~$t$, so
it follows that these first and last expressions are equal for all~$t$.

Replacing $t$ by~$t^{-1}$ and multiplying through by~$t^{n+2}$ gives
the form of the identity that will be most useful, namely:
\begin{equation}\label{eq: define pC}
\sum_{k=0}^{n+2} C_k\,t^{n+2-k} 
  = \det(t\,\I_n-H)\,(t^2 + h_1\,t + V) 
          +  T^*\,\Cof(t\,\I_n-H)\,T.
\end{equation}
The polynomial~$p_C(t) = t^{n+2}+C_2\,t^n+C_3\,t^{n-2}+\cdots+C_{n+2}$
will be said to be the \emph{characteristic polynomial} of the
Bochner-K\"ahler structure.

\begin{example}[Low cohomogeneity]\label{ex: pC in low cohom} 
Suppose~$M$ is locally isometric 
to~$M^p_c\times M^{n-p}_{-c}$. 
Then, looking back at the proof of
Proposition~\ref{prop: BK loc symm} and the definition of~$H$, 
one computes that
\begin{equation*}
p_C(t) = \left(t + (n{-}p{+}1)\,r\right)^{p+1}
         \left(t - (p{+}1)\,r\right)^{n-p+1},
\quad\text{where}\qquad
r = \frac{c}{2(n+2)}\,.
\end{equation*}
For Example~\ref{ex: rot sym}
(i.e., rotationally symmetric), the formula is
\begin{equation*}
p_C(t) = \left(t - 2r\right)^n
         \left[\left(t + nr\right)^2 
                  - {\ts\frac14}k^2 + a \right],
\quad\text{where}\qquad
r = \frac{k}{2(n+2)}\,.
\end{equation*} 
\end{example}

\subsubsection{The reduced characteristic polynomial}
\label{sssec: red char po}

Let~$P_1$ be the set of those~$u\in P_0$ that satisfy the condition 
that~$H'(u)$ be diagonal, with eigenvalues arranged in descending order, 
and that each of the entries~$T_i(u)$ of~$T'(u)\in\bbR^m$ 
be positive.  (See \S\ref{ssec: coho mo po} for definitions.)  
Then $P_1$ is a 
$\{\I_m\}{\times}\Un(n{-}m)$-bundle over~$M^\circ$.  Using the identities
derived in \S\ref{ssec: coho mo po}, equation~\eqref{eq: define pC} 
can be written as
\begin{equation}\label{eq: pC to ph'' ratio}
\frac{p_C(t)}{p_{h''}(t)} 
= \bigl(t^2 + h_1\,t + V\bigr)\,\prod_{j=1}^m(t-\lambda_j)
         \ +\  \sum_{i=1}^m T_i^2\,\prod_{{j\not=i}}(t-\lambda_j).
\end{equation}
In particular,~$p_{h''}(t)$ divides~$p_C(t)$.  Denote the quotient
by~$p_D(t)$.  It is a monic polynomial with constant coefficients
of degree~$m{+}2$ and will be called the \emph{reduced characteristic 
polynomial}.   

\begin{example}[Low cohomogeneity]\label{ex: pD in low cohom} 
Suppose~$M$ is locally isometric 
to~$M^p_c\times M^{n-p}_{-c}$. 
Then, $m=0$ and
\begin{equation*}
p_D(t) = \left(t + (n{-}p{+}1)\,r\right)
         \left(t - (p{+}1)\,r\right),
\quad\text{where}\qquad
r = \frac{c}{2(n+2)}\,.
\end{equation*}
For Example~\ref{ex: rot sym}, where~$m=1$, 
the formula is
\begin{equation*}
p_D(t) = \left(t - 2r\right)
         \left[\left(t + nr\right)^2 
                  - {\ts\frac14}k^2 + a \right],
\quad\text{where}\qquad
r = \frac{k}{2(n+2)}\,.
\end{equation*} 
\end{example}

\begin{proposition}\label{prop: ph'' roots are pD roots}
Every root of~$p_{h''}(t)$ is also a root of~$p_D(t)$.
\end{proposition}

\begin{proof}
Let~$p_{h''}(t)$ have roots $\lambda_{m+1}\ge\lambda_{m+2}\ge\cdots
\ge \lambda_n$, counting multiplicity, and set
\begin{equation}\label{eq: define Lambda}
\Lambda = \begin{pmatrix} 
             \lambda_{m+1} & 0 & \cdots & 0\\
              0 & \lambda_{m+2} & \cdots & 0\\
              \vdots & \vdots & \ddots & \vdots\\
              0 & 0 & \cdots & \lambda_n\\
          \end{pmatrix} .
\end{equation} 
Let $P_2\subset P_1$ consist of the coframes~$u\in P_1$ for which 
$H''(u) = \Lambda$.  
This~$P_2$ is a bundle over~$M^\circ$ with structure 
group~$\{\text{I}_m\}\times G_\Lambda$, where~$G_\Lambda\subset\Un(n{-}m)$
is the group of unitary matrices commuting with~$\Lambda$.
All calculations will now take place on~$P_2$. 

Adopt the index range convention $1\le i,j,k\le m < a,b,c < n$.
Thus, for example~$T_i>0$ but~$T_a = 0$.  Also, $H_{a\bar\imath}=0$.
The $(a,i)$-entry of the structure 
equation~\eqref{eq: structure equations ii} for~$dH$ becomes
(no sum over~$i$)
\begin{equation}\label{eq: off block vanishing}
(\lambda_a - \lambda_i)\phi_{a\bar\imath} + T_i\,\omega_a = 0.
\end{equation}
Meanwhile, the structure equation for~$dT_a$ becomes
\begin{equation}\label{eq: DE for T}
({\lambda_a}^2+h_1\,\lambda_a + V)\,\omega_a 
-\sum_{i=1}^m T_i\,\phi_{a\bar\imath} = 0.
\end{equation}
Combining these equations yields
\begin{equation*}
\begin{split}
\left(({\lambda_a}^2+h_1\,\lambda_a + V)\,
    \prod_{i=1}^m(\lambda_a-\lambda_i)\right)\,
\omega_a 
&= \left(\prod_{i=1}^m(\lambda_a-\lambda_i)\right)
    \,\sum_{j=1}^m T_j\,\phi_{a\bar\jmath}\\
&= \sum_{j=1}^m \left(\prod_{i\not=j}(\lambda_a-\lambda_i)\right)
            \,T_j\,(\lambda_a{-}\lambda_j)\phi_{a\bar\jmath}\\
&= -\left(\sum_{j=1}^m {T_j}^2\,\prod_{i\not=j}(\lambda_a-\lambda_i)\right)
            \,\omega_a\,.\\
\end{split}
\end{equation*}

Since~$\omega_a$ is nonzero on~$P_2$, it follows that
\begin{equation*}
p_D(\lambda_a) 
= \bigl({\lambda_a}^2 + h_1\,\lambda_a + V\bigr)
\,\prod_{j=1}^m(\lambda_a-\lambda_j)
         \ +\  \sum_{i=1}^m T_i^2\,\prod_{{j\not=i}}(\lambda_a-\lambda_j) 
= 0,
\end{equation*}
as desired. 
\end{proof}

\subsubsection{Point data}\label{sssec: point data}

In this subsubsection, a formula will be developed for the characteristic 
polynomials of a connected Bochner-K\"ahler structure in terms of
a single value of its structure function.  The formula for~$p_C$ is,
of course, already given by~\eqref{eq: define pC}.  
However, the formula for~$p_D$ is
somewhat more subtle.

First, recall the concepts introduced in~\S\ref{ssec: inftess symms}, 
suitably modified for the present section.  
For any~$(H_0,T_0,V_0)\in i\euu(n){\oplus}\C{n}{\oplus}\bbR$, 
let~$H_1,\,H_2,\,\cdots,\,H_\delta$ be the distinct eigenvalues
of~$H_0$. Let~$L_\alpha\subset\C{n}$ be the eigenspace of
$H_0$ belonging to the eigenvalue~$H_\alpha$ and let~$n_\alpha\ge1$ be 
the (complex) dimension of~$L_\alpha$.  Write
\begin{equation}\label{eq: decomp T}
T_0 = T_1 + \cdots + T_\delta
\end{equation}
where~$T_\alpha$ lies in~$L_\alpha$ for~$1\le\alpha\le\delta$.  
Define the quantities
\begin{equation}\label{eq: define V-alpha}
V_\alpha = {H_\alpha}^2 + (\tr H_0)\,H_\alpha + V_0 
         +\sum_{\beta\not=\alpha}\frac{|T_\beta|^2}{(H_\alpha - H_\beta)}
\end{equation}
and
\begin{equation}\label{eq: define m-alpha}
m_\alpha
= \begin{cases}
  2 &\quad \text{if $T_\alpha\not=0$ and $n_\alpha>1$;}\\
  1 &\quad \text{if $T_\alpha\not=0$ and $n_\alpha=1$;}\\
  1 &\quad \text{if $T_\alpha = 0$ and $V_\alpha\not=0$;}\\
  0 &\quad \text{if $T_\alpha = 0$ and~$V_\alpha=0$.}\\
  \end{cases}
\end{equation}

\begin{theorem}\label{thm: point data to polys}  
If $(M^n,g,\Omega)$ is a connected Bochner-K\"ahler
manifold whose structure function~$(H,T,V)$ assumes the value~$(H_0,T_0,V_0)$,
then
\begin{equation}\label{eq: pC from points}
p_C(t) = \prod_{\alpha=1}^\delta (t-H_\alpha)^{n_\alpha}
          \left[t^2 + (\tr H_0)\,t + V_0 
            + \sum_{\alpha=1}^\delta \frac{|T_\alpha|^2}{(t-H_\alpha)}\right]
\end{equation}
and
\begin{equation}\label{eq: pD from points}
p_D(t) 
= \prod_{\alpha=1}^\delta (t-H_\alpha)^{m_\alpha}
   \left[t^2 + (\tr H_0)\,t + V_0
     + \sum_{\alpha=1}^\delta \frac{|T_\alpha|^2}{(t-H_\alpha)}\right]\rlap{.}
\end{equation}
\end{theorem}

\begin{proof}  The formula for~$p_C$ follows directly 
from~\eqref{eq: define pC}, so
the formula for~$p_D$ will follow from the equivalent statement
\begin{equation}\label{eq: ph'' product}
p_{h''}(t) = \prod_{\alpha=1}^\delta (t-H_\alpha)^{n_\alpha-m_\alpha},
\end{equation}
and this is what will be proved. 

By \S\ref{sssec: orbit dim n slices}, the generic orbit of the symmetry
pseudo-groupoid has codimension equal to $m_1+\cdots+m_\delta$.  
By~\S\ref{sssec: mom map} and Theorem~\ref{thm: cohomogeneity m}, the orbits of
the symmetry pseudo-groupoid in~$M^\circ$ (which is open and dense in~$M$)
have codimension~$m$.  Consequently, $m=m_1+\cdots+m_\delta$, so 
that~$n-m = (n_1{-}m_1) + \cdots + (n_\delta{-}m_\delta)$.  Moreover,
the inequality~$n_\alpha\ge m_\alpha$ follows immediately from the 
definitions.

Thus, by the very definition of~$p_{h''}(t)$, it will suffice to show
that, for each~$\alpha$, the polynomial~$p_h(t)$ has 
a constant root~$H_\alpha$ of multiplicity at least~$n_\alpha{-}m_\alpha$.
By Theorem~\ref{thm: cohomogeneity m}, it suffices to to show this constancy 
in an open neighborhood of the point~$x\in M$ for which there exists 
a~$u\in P_x$ satisfying~$\bigl(H(u),T(u),V(u)\bigr) = (H_0,T_0,V_0)$.

Thus, let~$w\in\C{n}$ be a nonzero vector and 
let~$c:(-\varepsilon,\varepsilon)\to M$ be the constant speed
geodesic satisfying~$u\bigl(\dot c(0)\bigr) = w$.  Then~$c$
can be lifted uniquely to a curve~$\gamma:(-\varepsilon,\varepsilon)\to P$
that satisfies~$\gamma(0) = u$ and~$\gamma^*\phi=0$ (i.e., the 
coframe field~$\gamma$ is parallel along~$c$).  Because $c$ is
a constant speed geodesic, $\gamma$ also 
satisfies~$\gamma^*(\omega) = w\,ds$, where~$s$ is the parameter 
on~$(-\varepsilon,\varepsilon)$.

Because the polynomial~$p_h(t)$ is
invariant under the action of the symmetry pseudo-groupoid, it suffices
to consider only geodesics with initial velocities orthogonal to the
subspace~$O_x\subset T_xM$ that is the tangent to the orbit through~$x$.
Thus, I will assume that if~$T_\alpha\not=0$, then~$T_\alpha^*w$ is real
and that if~$T_\alpha=V_\alpha=0$, then $w$ is orthogonal to $L_\alpha$.

For simplicity, 
set~$H(s) = H\bigl(\gamma(s)\bigr)$, $T(s) = T\bigl(\gamma(s)\bigr)$, 
and~$V(s) = V\bigl(\gamma(s)\bigr)$.
Then these functions on~$(-\varepsilon,\varepsilon)$
satisfy the initial conditions~$\bigl(H(0),T(0),V(0)\bigr) = (H_0,T_0,V_0)$
and the system of ordinary differential equations
\begin{equation}\label{eq: odes for HTV}
\begin{split}
\dot H &=  T\,w^* + w\,T^* \,,\\
\dot T &=  \bigl(H^2+(\tr H)\,H + V\,\I_n\bigr)\,w,\\
\dot V &= (\tr H)\bigl(T^*w+w^*T\bigr) 
           +\bigl(T^*Hw+w^*HT\bigr).
\end{split}
\end{equation}

Let~$L\in i\euu(n)$ be any fixed element that
satisfies~$Lw = LT_0=0$ and~$[L,H_0] = 0$.  Because of the latter
equation,~$L$ preserves each of the eigenspaces of~$H_0$, i.e., 
the subspaces~$L_\alpha$.  Consequently,~$LT_\alpha=0$
and~$Lw_\alpha=0$ for all~$\alpha$.  Conversely, if~$L\in i\euu(n)$ 
preserves the eigenspaces of~$H_0$ and annihilates~$T_\alpha$
and~$w_\alpha$ for all~$\alpha$, then it satisfies
satisfies~$Lw = LT_0=0$ and~$[L,H_0] = 0$.

Now, the above differential equations imply the differential equations
\begin{equation}\label{eq: odes for LT and LH}
\begin{split}
[L,\dot H] &=  LT\,w^* - w\,(LT)^* \,,\\
L\dot T &=  \bigl([L,H]H+H[L,H]+(\tr H)\,[L,H]\bigr)\,w,\\
\end{split}
\end{equation}
so that the quantities~$\bigl([L,H(s)],LT(s)\bigr)$ satisfy a 
linear system of ordinary differential equations with vanishing 
initial condition at~$s=0$.  Consequently~$[L,H(s)]$ and $LT(s)$ 
vanish identically for all~$s$, as does~$Lw$ (for trivial reasons).
By the above characterization of those~$L\in i\euu(n)$ that
satisfy~$[L,H_0] = LT_0 = Lw = 0$, this implies that the 
subspace~$K_\alpha\subset L_\alpha$ that is perpendicular 
to~$T_\alpha$ and~$w_\alpha$ is necessarily an eigenspace 
of~$H(s)$ that is perpendicular to~$T_\alpha(s)$ and~$w$ 
for all~$s$.  

If~$K_\alpha\not=0$, then there is a well-defined
eigenvalue~$H_\alpha(s)$ of~$H(s)$ associated to~$K_\alpha$.  
In particular,
\begin{equation}\label{eq: H-alpha is const}
\dot H_\alpha(s) y = \dot H(s) y 
    = \bigl(T(s)\,w^* + w\,T(s)^*\bigr)y = 0
\end{equation}
for all~$y\in K_\alpha$.  Of course, this implies 
that~$\dot H_\alpha(s) = 0$, i.e., 
that~$H_\alpha(s) = H_\alpha(0) = H_\alpha$ for all~$s$.

There are now four cases to consider:

If $\alpha$ is such that~$T_\alpha\not=0$ and~$n_\alpha>1$,
then $\dim K_\alpha$ is either~$n_\alpha{-}1$ or $n_\alpha{-}2$,
depending on whether~$w_\alpha$ is zero or not.  In either case,
$\dim K_\alpha\ge n_\alpha-2 = n_\alpha-m_\alpha$, so $H_\alpha$ 
is a root of~$p_h(t)$ of multiplicity at least~$n_\alpha{-}m_\alpha$,
as desired.

If $\alpha$ is such that~$T_\alpha\not=0$ and~$n_\alpha=1$, then
$m_\alpha=1$ and~$K_\alpha=0$.  In this case, of course, 
$n_\alpha{-}m_\alpha = 0$, so $H_\alpha$ is trivially a root of
of~$p_h(t)$ of multiplicity at least~$n_\alpha{-}m_\alpha$,
as desired.

If~$\alpha$ is such that~$T_\alpha=0$ but~$V_\alpha\not=0$, then
$\dim K_\alpha$ is either~$n_\alpha$ or~$n_\alpha{-}1$, 
depending on whether~$w_\alpha$ is zero or not.  In either case,
$\dim K_\alpha\ge n_\alpha-1 = n_\alpha-m_\alpha$, so $H_\alpha$ 
is a root of~$p_h(t)$ of multiplicity at least~$n_\alpha{-}m_\alpha$,
as desired.

Finally, if $T_\alpha=0$ and~$V_\alpha=0$, then~$w_\alpha=0$ by
the above condition on~$c$.  Thus~$K_\alpha = L_\alpha$, so that
$\dim K_\alpha = n_\alpha{-}0 = n_\alpha-m_\alpha$ and, again,
$H_\alpha$ 
is a root of~$p_h(t)$ of multiplicity at least~$n_\alpha{-}m_\alpha$,
as desired.  
\end{proof}

The following result will be needed in the next subsection.  Its
proof follows by inspection of the formula for~$p_D(t)$ and
the definition of the~$m_\alpha$ and so will be omitted.

\begin{corollary}\label{cor: no mult roots} 
No root of~$p_{h'}(t)$ is a multiple 
root of~$p_D(t)$. 
\end{corollary}

\subsection{Momentum cells}\label{ssec: mom cells}

In general, the two characteristic polynomials do not completely determine 
the analytically connected equivalence class of a Bochner-K\"ahler structure.  
However, as will be seen in this subsection, they do 
determine it up to a finite number (at most~$m{+}1$) of 
possibilities (Theorem~\ref{thm: analytically connected classes}).

\subsubsection{The roots of~$p_D$}\label{sssec: roots of pD}

It turns out that the reality and multiplicity properties of the
roots of~$p_D$ are severely constrained.

\begin{proposition}\label{prop: pD root cases}
One of the following cases holds:
\begin{enumerate}
\item $p_D$ has $\phantom{+}m\phantom{1}$ real, distinct roots, 
       all of order~$1$;
\item $p_D$ has $\phantom{+}m\phantom{1}$ real, distinct roots, 
       one of order~$3$ and the rest of order~$1$;
\item $p_D$ has $m{+}1$ real, distinct roots, 
       one of order~$2$ and the rest of order~$1$;
\item $p_D$ has $m{+}2$ real, distinct roots, 
       all of order~$1$.
\end{enumerate}
\end{proposition}

\begin{proof}
Substituting~$t=\lambda_i$ into~\eqref{eq: pC to ph'' ratio} yields
\begin{equation}\label{eq: pD at lambda-i}
p_D(\lambda_i) = T_i^2\,\prod_{{j\not=i}}(\lambda_i-\lambda_j).
\end{equation} 
Since~$\lambda_1>\lambda_2>\cdots>\lambda_m$ and~$T_i>0$ on~$M^\circ$, 
it follows that~$(-1)^{i-1}p_D(\lambda_i)>0$ for $1\le i\le m$ holds
on~$M^\circ$.

Equivalently, for every~$x\in M^\circ$, the polynomial $p_D(t)$ 
has an even number of real roots (counted with multiplicity) 
greater than~$\lambda_1(x)$ and an odd number of real roots 
(counted with multiplicity) 
in each open interval~$\bigl(\lambda_i(x),\lambda_{i+1}(x)\bigr)$ 
for~$1\le i<m$. Moreover, since~$(-1)^mp_D(t)$ is positive for all~$t$ 
sufficiently negative, $p_D(t)$ has an odd number of real roots 
(counted with multiplicity) less than~$\lambda_m(x)$.  
These considerations imply that $p_D(t)$ has at least~$m$ 
distinct real roots.

If~$p_D$ has exactly $m$ real roots, then the above parity conditions
show that they must all have odd order.  Since~$p_D(t)$ has degree~$m{+}2$,
either Case~1 or Case~2 must hold.

If $p_D$ has exactly $m{+}1$ real, distinct roots, 
then Case~3 must hold.

If $p_D$ has exactly $m{+}2$ real, distinct roots, 
then Case 4 must hold. 
\end{proof}

\begin{remark}[Global inequalities]
By continuity, the 
inequality $(-1)^{i-1}p_D(\lambda_i)\ge0$  holds on~$M$.
\end{remark}

\begin{remark}[Root labeling]
I will use the following convention to label the real roots of~$p_D$:  
When~$p_D$ has~$m$ distinct real roots, denote them by~$r_1>r_2>\cdots>r_m$;
when $p_D$ has $m{+}1$ distinct real roots, denote them 
by~$r_1>\cdots>r_{m+1}$; 
and when $p_D$ has $m{+}2$ distinct real roots, 
denote them by~$r_0>r_1>\cdots>r_{m+1}$. In all cases,
the list of real roots of~$p_D$, in descending order, will be denoted by~$r$.

If $r_i$ is any real root of~$p_D(t)$, the number of
roots~$\{\lambda_1(x),\ldots,\lambda_m(x)\}$
that are strictly greater than~$r_i$ is independent of~$x\in M^\circ$, 
so I will denote this common value by~$\mu_i$.  The function~$\mu$ is
constrained as follows:

Cases 1 and 2: Necessarily, $\mu_i = i$.

Case 3: Let~$r_i$ be the double root.  
Then $\mu_i$ is either~$i$ (SubCase~(3-$i$,$a$); impossible when~$i=m{+}1$) 
or~$i{-}1$ (Subcase (3-$i$,$b$)).
Moreover, $\mu_j=j$ for $j<i$, while~$\mu_j = j{-}1$ for~$j>i$.  

Case 4: There is an integer~$i\le m$ so that~$\mu_i=i$ (Subcase~4-$i$). 
Then~$\mu_j = j{+}1$ for~$j<i$ while~$\mu_j = j{-}1$ for~$j>i$.  
\end{remark}

\subsubsection{Momentum cells}\label{sssec: mom cells ii}

Since
\begin{equation*}
p_{h'}(r_i) = \prod_{j=1}^m(r_i-\lambda_j),
\end{equation*}
it follows that $(-1)^{\mu_i}p_{h'}(r_i)>0$ on~$M^\circ$.  By continuity,
the inequality
\begin{equation*}
(-1)^{\mu_i}
\bigl({r_i}^m-h'_1\,{r_i}^{m-1}
             +h'_2\,{r_i}^{m-2}-\cdots+(-1)^m\,h'_m\bigr)\ge0
\end{equation*}
holds on~$M$, with strict inequality on~$M^\circ$.  Moreover, by 
Corollary~\ref{cor: no mult roots},
if~$r_i$ is a multiple root of~$p_D$, then it is not a root of
$p_{h'}(t)$ at any point of~$M$, so that the above inequality is strict
on all of~$M$.

The image~$h'(M)\subset\bbR^m$  therefore lies in the intersection 
of the closed half-spaces~$\overline{H(r_i,\mu_i)}$ 
defined by the inequalities
\begin{equation}\label{eq: simple root inequality}
(-1)^{\mu_i} 
\bigl({r_i}^m-{r_i}^{m-1}\,x_1+{r_i}^{m-2}\,x_2-\cdots+(-1)^m\,x_m\bigr)\ge0
\end{equation}
as~$r_i$ ranges over the simple real roots of~$p_D$
and the open half-space~$H(r_i,\mu_i)$ defined by
\begin{equation}\label{eq: multi root inequality}
(-1)^{\mu_i}
\bigl({r_i}^m-{r_i}^{m-1}\,x_1+{r_i}^{m-2}\,x_2-\cdots+(-1)^m\,x_m\bigr)>0
\end{equation}
if $r_i$ is a multiple root of~$p_D$.  This intersection will
be referred to as the \emph{momentum cell}~$C(p_D,\mu)\subset\bbR^m$.
Note that~$h'(M^\circ)$ lies in~$C(p_D,\mu)^\circ$, 
the interior of~$C(p_D,\mu)$, and that~$h':M^\circ\to C(p_D,\mu)^\circ$
is a submersion onto its image.

\subsubsection{Possible momentum cells}\label{sssec: pos mo cells}

More generally, if $p_D(t)$ is any monic polynomial of degree~$m{+}2$ 
with real coefficients that falls into one of the Cases 1 through 4 of
Proposition~\ref{prop: pD root cases}, define the \emph{possible momentum 
cells of~$p_D$ }
as follows:

\begin{figure}
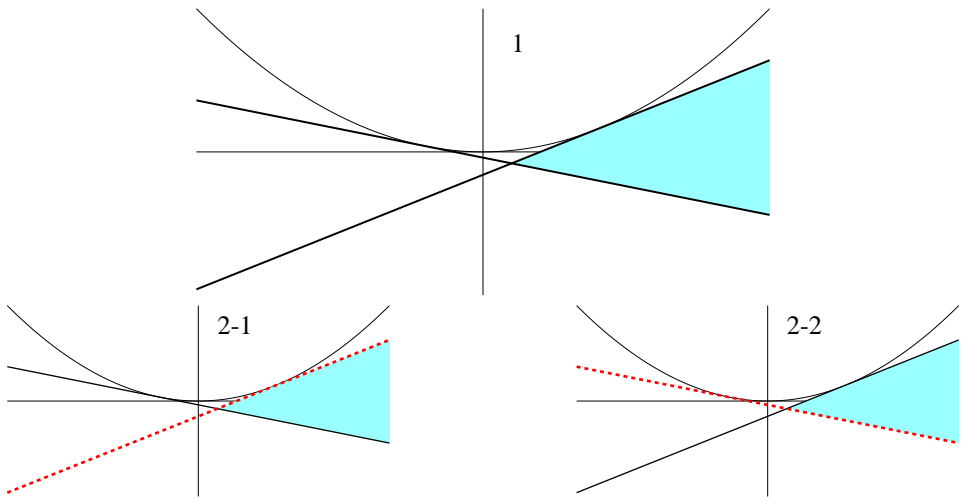

\includegraphics[width=3.0in]{Figure1a.eps}
\vskip3pt
\includegraphics[width=2.0in]{Figure1b.eps}
\hfill
\includegraphics[width=2.0in]{Figure1c.eps}
\caption[Two distinct real roots, $m=2$]{\label{fig: two roots}
Possible momentum cells with $m=2$ and two distinct real roots:
Case 1 (both roots simple),  Case~2-1 ($r_1$ triple),
and Case~2-2 ($r_2$ triple).}
\end{figure}

	If $p_D$ falls into Case~1 or Case~2, 
define~$\mu_i=i$ for~$1\le i\le m$ and let~$C(p_D,\mu)$ 
be defined by the inequalities~(16) (strict or not depending on 
the multiplicity of the roots).  This cell $C(p_D,\mu)$ is
a closed, unbounded, convex polytope in Case~1, but is not
closed in Case 2, since one of the faces is missing.

\begin{figure}
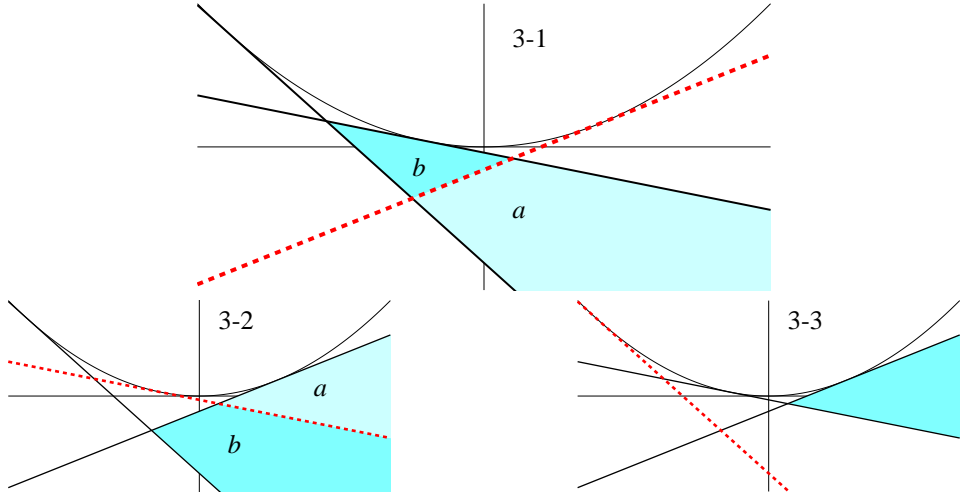

\includegraphics[width=3.0in]{Figure2a.eps}
\vskip3pt
\includegraphics[width=2.0in]{Figure2b.eps}
\hfill
\includegraphics[width=2.0in]{Figure2c.eps}
\caption[Three distinct real roots, $m=2$]{\label{fig: three roots}
Possible momentum cells with $m=2$ and three distinct real roots:
Case~3-1 ($r_1$ double; two cells, $a$ unbounded, $b$ bounded), 
Case~3-2 ($r_2$ double; two cells, both unbounded),
and Case~3-3 ($r_3$ double; one cell, unbounded).
}
\end{figure}

If $p_D$ falls into Case~3, with $r_i$ being the double root, 
define $\mu_j=j$ for $j<i$ and $\mu_j = j{-}1$ for~$j>i$, while~$\mu_i$ 
is allowed to be one of~$i$ (type~$a$) or~$i{-}1$ (type~$b$).  
Let~$C(p_D,\mu)$ be defined by the inequalities~(16) 
(strict or not depending on the multiplicity of the roots).  
Thus, there are two possible momentum cells except in the case
that~$r_{m+1}$ is the double root, in which case, there is only 
one possible cell. When there are two cells, neither is closed
and their closures share the missing face.  The only subcase 
with a bounded cell is Subcase~(3-1,$b$).

\begin{figure}
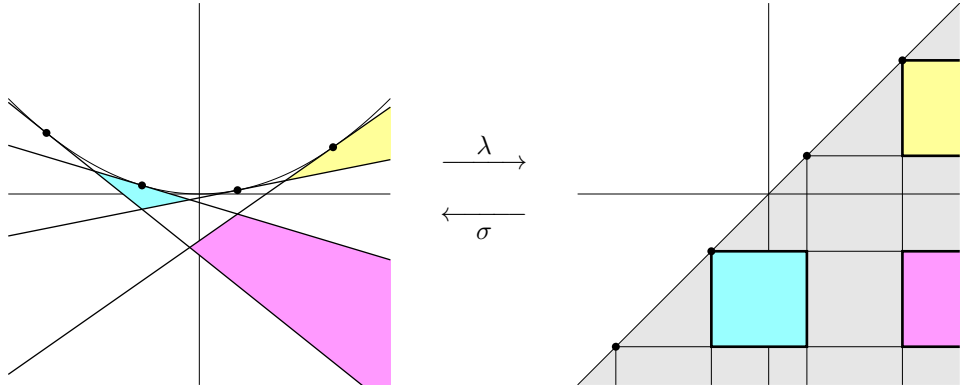

\includegraphics[width=2.0in]{Figure3a.eps}
\hfill
\raisebox{1.0in}{$\begin{CD}@>\quad{\ts\lambda}\quad>>\\
           @<<\quad{\ts\sigma}\quad<\end{CD}$}
\hfill
\includegraphics[width=2.0in]{Figure3b.eps}
\caption[Four distinct real roots, $m=2$]{\label{fig: four roots}
Possible momentum cells in~Case~4 with $m=2$ and four distinct real roots
and their corresponding spectral products. 
The `highest' cell (SubCase 4-2) has only two faces.  
The `lowest' cell (SubCase 4-0) is the only bounded one. 
}
\end{figure}

When~$p_D$ falls into Case~4, choose an integer~$i$ in the 
range~$0\le i\le m$ and define~$\mu$ so that~$\mu_j = j{+}1$
for~$j<i$, while $\mu_i=i$ and~$\mu_j = j{-}1$ for~$j>i$. 
Let~$C(p_D,\mu)$ be defined by the 
inequalities~\eqref{eq: simple root inequality}. 
Thus, there are $m{+}1$ possible momentum cells, one for each 
possible choice of~$i$.  Each of these cells is a closed polytope 
and they are mutually disjoint.  When~$\mu_{m}=m$
(i.e., $\lambda_m>r_{m-1}$ on~$M^\circ$), this `highest' cell has $m$ faces.
Each of the other cells has $m{+}1$ faces.  The only bounded cell
falls in Subcase~4-0, i.e.,~$\mu_0=0$ (implying 
that~$\mu_i = i{-}1$ for all~$1\le i\le m{+}1$).

Figures~\ref{fig: two roots}, \ref{fig: three roots}, 
and \ref{fig: four roots} show the possible momentum cells when~$m=2$
in the four cases.  The drawn axes are~$u_1$ ($= h_1'$) and $u_2$ ($=h'_2$).
Of course, all of these cells lie below the discriminant parabola~
$u_2 = \frac14\,{u_1}^2$ and are bounded by its tangent lines.

\subsubsection{The product representation}\label{sssec: prod rep}

It will be useful to have another description of the
possible momentum cells.

Let~$\sigma:\bbR^m\to\bbR^m$ be the standard symmetrizing map, so
that~$\sigma = (\sigma_1,\ldots,\sigma_m)$ where~$\sigma_k:\bbR^m\to\bbR$
is the $k$-th elementary symmetric function of its arguments.  
The map~$\sigma$ is one-to-one on the closed set
\begin{equation*}
\bbR^m_\ge=\{\ (x_1,\ldots,x_m)\in\bbR^m
                \ \mid\ x_1\ge x_2\ge\cdots\ge x_m\ \}.
\end{equation*}
Moreover, $\sigma:\bbR^m_\ge\to\sigma(\bbR^m_\ge)\subset\bbR^m$ is 
a homeomorphism onto its image and a real-analytic diffeomorphism 
on its interior, which will be denoted~$\bbR^m_>\subset\bbR^m_\ge$.
Denote the inverse of~$\sigma$ 
by~$\lambda:\sigma(\bbR^m_\ge)\to \bbR^m_\ge$.

For any momentum cell~$C(p_D,\mu)$, the 
subset~$\lambda\bigl(C(p_D,\mu)\bigr)$ is a product of the form
\begin{equation}\label{eq: spectral bands}
\lambda\bigl(C(p_D,\mu)\bigr) = I_1\times I_2 \times\cdots\times I_m
\end{equation}
where~$I_1,I_2,\cdots,I_m$ are (non-empty) intervals in~$\bbR$ 
with non-overlapping interiors.  The endpoints of the 
closure~$\overline{I_i}$ are roots of~$p_D$ and~$I_i$ 
contains such an endpoint if and only if that endpoint 
is a simple root of~$p_D$.  

In Case~1, where the distinct roots~$r_1>\cdots>r_m$ are all
simple, $I_1 = [r_1,\infty)$ and $I_k = [r_k,r_{k-1}]$
for~$1<k\le m$.  

In Case~2 in which, say, $r_1$ 
is the triple root, $I_1 = (r_1,\infty)$, $I_2 = [r_2,r_1)$, 
and (assuming~$m>2$)~$I_k = [r_k,r_{k-1}]$ for~$3\le k\le m$.

The intervals~$(I_1,\ldots,I_m)$ will be referred to as
the \emph{spectral bands} associated to~$C(p_D,\mu)$ 
and~$I_1\times I_2 \times\cdots\times I_m$ will be 
known as the \emph{spectral product}.  Since
\begin{equation}\label{eq: inv spectral bands}
\sigma(I_1\times I_2 \times\cdots\times I_m) = C(p_D,\mu),
\end{equation}
specifying the spectral bands is equivalent to specifying~$C(p_D,\mu)$. 
Also, note that $\sigma$ 
maps~$I^\circ_1\times I^\circ_2 \times\cdots\times I^\circ_m$ 
diffeomorphically onto~$C(p_D,\mu)^\circ$, the interior of~$C(p_D,\mu)$.

\begin{proposition}\label{prop: values of ph'}
Suppose that $p_D(t)$ is a polynomial of degree~$m{+}2$ that 
falls into one of the Cases of 
Proposition~{\upshape\ref{prop: pD root cases}}.  
Suppose that 
there exists a monic polynomial~$p_C(t)$ of degree~$n{+}2$ 
with the properties that~$p_C(t)/p_D(t)$ is a polynomial 
all of whose roots are real roots of~$p_D(t)$ 
and that its~$t^{n+1}$-coefficient vanishes.  

Then for every~$k'$ in a possible momentum 
cell~$C(p_D,\mu)\subset\bbR^m$, there exists a Bochner-K\"ahler 
$n$-manifold~$(M,g,\Omega)$ whose characteristic polynomials
are~$p_C(t)$ and~$p_D(t)$ and whose reduced momentum 
mapping~$h':M\to\bbR^m$ assumes the value~$k'$.
\end{proposition}

\begin{proof}
This will be a matter of checking cases.

First, some generalities.  Given polynomials~$p_D(t)$ and~$p_C(t)$
satisfying the hypotheses of the proposition and a~$k'\in\bbR^m$
lying in a possible momentum cell~$C(p_D,\mu)$ for~$p_D$, define
$p_{h''}(t) = p_C(t)/p_D(t)$.  By hypothesis, all the roots 
of~$p_{h''}(t)$ are real and are roots of~$p_D(t)$ as well.

Define
\begin{equation*}
p_{k'}(t) = t^m-k'_1\,t^{m-1}+\cdots+(-1)^mk'_m.
\end{equation*}
Let~$\lambda\bigl(C(p_D,\mu)\bigr) = I_1\times I_2 \times\cdots\times I_m$,
and let~$\lambda(k') = (s_1,s_2,\ldots,s_m)$, so that
there exists a real factorization of the form
\begin{equation*}
p_{k'}(t) = (t-s_1)(t-s_2)\cdots(t-s_m)
\end{equation*}
with~$s_k\in I_k$ for $1\le k\le m$.

First, assume that~$k'$ lies in~$C(p_D,\mu)^\circ$, 
the interior of~$C(p_D,\mu)$.   Then~$s_k$ lies in~$I^\circ_k$ 
for~$1\le k\le m$ and, in particular, the~$s_k$ are all distinct and 
are not roots of~$p_D(t)$.  It follows that the rational function~
$p_D(t)/p_{k'}(t)$ has a simple pole at~$t=s_k$ for~$1\le k\le m$.
Since~$p_D(t)$ has degree~$m{+}2$ and is monic, there is 
a partial fractions expansion of the form
\begin{equation*}
\frac{p_D(t)}{p_{k'}(t)}
= t^2 + b_1\,t + b_2 + \sum_{k=1}^{m} \frac{q_k}{(t-s_k)}.
\end{equation*}
Because of the way that the possible momentum cells were defined, 
the inequality $(-1)^{k-1}p_D(s_k)>0$ holds for~$s_k\in I^\circ_k$,
so it follows that~$q_k>0$ for~$1\le k\le m$.

If~$n>m$, define~$s_{m+1}\ge s_{m+2}\ge \cdots \ge s_n$ so that
\begin{equation*}
p_{h''}(t) = (t-s_{m+1})(t-s_{m+2})\cdots(t-s_m).
\end{equation*}
By hypothesis, each root of~$p_{h''}(t)$ is a real root of~$p_D$ and so
is not equal to any of the roots of~$p_{k'}(t)$.

Consider the element~$(s,t,v)\in i\euu(n)\oplus\C{n}\oplus\bbR$
defined by letting~$s$ be the diagonal matrix with entries~$s_{i\bar\imath}
= s_i$ for~$1\le i\le n$; letting~$t_i = \sqrt{q_i}$ for~$1\le i\le m$
and~$t_i = 0$ for~$m<i\le n$ (if~$n>m$); and letting~$v = b_2$. 
The hypothesis that $p_C$ have no~$t^{n+1}$-term is then seen to
be equivalent to the condition that~$b_1 = \tr s$, while the condition
that each~$s_a$ for~$a>m$ be a root of~$p_D$ is then equivalent to the
condition that
\begin{equation*}
v_a = {s_a}^2 + b_1\,s_a + v + \sum_{i=1}^m \frac{{t_i}^2}{(s_a-s_i)} = 0.
\end{equation*}
By Theorem~\ref{thm: point data to polys}, the Bochner-K\"ahler structure on
a neighborhood~$M$ of~$0\in\C{n}$ that has a unitary 
coframe~$u_0:T_0M\to\C{n}$ with~$\bigl(H(u_0),T(u_0),V(u_0)\bigr)=(s,t,v)$
has~$p_C(t)$ and $p_D(t)$ as its characteristic polynomials and satisfies~
$h'(0) = k'$.  This establishes existence for the 
interior points of~$C(p_D,\mu)$.

It remains to treat the boundary cases, i.e., cases in which one or
more of the~$s_i$ are actually roots of~$p_D(t)$.

Now, if~$s_j = s_{j+1}$ for
any~$j$, then~$\{s_j\} = I_j\cap I_{j+1}$, so that~$s_j$ is a 
simple root of~$p_D$.  In such a case, necessarily, $s_{j-1}>s_j$
(if~$j>1$) since~$I_{j-1}\cap I_{j+1} = \emptyset$
and~$s_{j+1}>s_{j+2}$ (if ~$j<m{-}1$) since $I_{j}\cap I_{j+2} = \emptyset$.  
Consequently, each $s_j$ is at most a double root of $p_{k'}(t)$ 
and, if so, it must also be a root of~$p_D(t)$.  It follows that the
rational function~$p_D(t)/p_{k'}(t)$ has a simple pole at~$t=s_j$ 
for~$1\le j\le m$.  Since~$p_D(t)$ has degree~$m{+}2$ and is monic, 
there is a partial fractions expansion
\begin{equation*}
\frac{p_D(t)}{p_{k'}(t)}
= t^2 + b_1\,t + b_2 + \sum_{j=1}^{m} \frac{q_j}{(t-s_j)},
\end{equation*}
where, in order to make the~$q_j$ unique, it is now necessary to add
the condition that~$q_{j+1} = 0$ if~$s_{j+1} = s_j$.
If~$j$ is such that~$s_j$ is not a root of~$p_D(t)$, then
the inequality $(-1)^{j-1}p_D(s_j)>0$ holds so that~$q_j>0$.
If~$j$ is such that~$s_j$ is a simple root of both~$p_{k'}(t)$ and~$p_D(t)$,
then~$q_j=0$.  If~$j$ is such that~$s_j = s_{j+1}$, then 
$(-1)^{j-1}p_D(t)$ and ~$(-1)^{j-1}p_{k'}(t)$ are both positive 
on~$I^\circ_j$. Since $s_j$ is a double root of~$p_{k'}(t)$ and
a simple root of~$p_D(t)$, it follows that
\begin{equation*}
\lim_{t\to s_j^+}\frac{p_D(t)}{p_{k'}(t)} = +\infty,
\end{equation*}
which can only hold if~$q_j>0$.  In particular, $q_j\ge0$
has been defined for~$1\le j\le m$ so that the above partial 
fractions expansion is valid.

If~$n>m$, again define~$s_{m+1}\ge s_{m+2}\ge \cdots \ge s_n$ so that
\begin{equation*}
p_{h''}(t) = (t-s_{m+1})(t-s_{m+2})\cdots(t-s_m).
\end{equation*}
Again, each root of~$p_{h''}(t)$ is a real root of~$p_D$ but now
it may also be a root of~$p_{k'}(t)$.

Define an element~$(s,t,v)\in i\euu(n)\oplus\C{n}\oplus\bbR$
by letting~$s$ be the diagonal matrix with entries~$s_{i\bar\imath}
= s_i$ for~$1\le i\le n$; letting~$t_i = \sqrt{q_i}$ for~$1\le i\le m$
and~$t_i = 0$ for~$m<i\le n$ (if~$n>m$); and letting~$v = b_2$. 
It must now be verified that the element~$(s,t,v)$ does indeed have
$p_C(t)$ and $p_D(t)$ as its characteristic polynomials.

Now, the hypothesis that $p_C$ have no~$t^{n+1}$-term is again seen to 
be equivalent to the condition that~$b_1 = \tr s$, and the condition
that~$t_i = 0$ for~$i>m$ or when~$s_i = s_{i-1}$ implies that
\begin{equation*}
\begin{split}
t^2 + (\tr s)\,t + v + \sum_{j=1}^n\frac{{t_j}^2}{(t-s_j)}
&= t^2 + b_1\,t + b_2 + \sum_{j=1}^{m} \frac{q_j}{(t-s_j)} \\
&= \frac{p_D(t)}{p_{k'}(t)} =  \frac{p_C(t)}{p_{h''}(t)p_{k'}(t)},\\
\end{split}
\end{equation*}
so that
\begin{equation*}
p_C(t) = {\prod_{i=1}^n(t-s_i)}\,
    \left[t^2+(\tr s)\,t+v+\sum_{j=1}^{n}\frac{{t_j}^2}{(t-s_j)}\right],
\end{equation*}
as desired.  It remains to verify that~$p_D(t)$ is the reduced
characteristic polynomial associated to~$(s,t,v)$, i.e., to compute the
numbers~$n_i$ and~$m_i$ for each eigenvalue~$s_i$ of~$s$ according to
the recipe of~\S\ref{ssec: inftess symms} 
and show that $s_i$ is a root of~$p_{h''}(t)$ 
of multiplicity exactly equal to~$n_i-m_i$.  This will be done by 
breaking it down into a number of cases.

If~$s_a$ is not~$s_i$ for any~$1\le i\le m$, 
then $p_D(s_a) =0$ is equivalent to
\begin{equation*}
v_a = {s_a}^2 + (\tr s)\,s_a + v + \sum_{i=1}^m \frac{{t_i}^2}{(s_a-s_i)} = 0,
\end{equation*}
and this implies that~$s_a$ is an eigenvalue of~$s$ of some multiplicity
$n_a\ge 1$ that satisfies~$t_a = v_a = 0$, so that~$m_a =0$.  Thus,~$s_a$
is a root of~$p_{h''}(t)$ and has multiplicity~$n_a{-}m_a = n_a$, 
as desired.

If~$s_i$ is not a root of~$p_D(t)$, then, by construction, it is
a simple eigenvalue of~$s$ and also satisfies~$t_i = \sqrt{q_i}>0$,
so~$m_i = 1$ and $n_i-m_i = 0$, so that~$(t-s_i)$ is not a factor
of~$p_{h''}(t)$, again, as desired.

If~$s_i$ is a simple root of~$p_D(t)$ and a simple root of~$p_{k'}(t)$,
then, by construction,~$t_i=0$.  The quantity~$v_i$ is
calculated to be
\begin{equation*}
\begin{split}
v_i &= {s_i}^2 + (\tr s)\,s_i + v 
            + \sum_{j\not=i} \frac{{t_j}^2}{(s_i-s_j)}\\
&=\lim_{t\to s_i}\left({t}^2 + (\tr s)\,t + v 
            + \sum_{j=1}^m \frac{{t_j}^2}{(t-s_j)}\right)\\
&=\lim_{t\to s_i} \frac{p_D(t)}{p_{k'}(t)}\not=0.\\
\end{split}
\end{equation*}
Thus, the recipe gives~$m_i = 1$, again as desired.

If~$s_i$ is a simple root of~$p_D(t)$ and a double root of~$p_{k'}(t)$,
then it can be assumed that~$s_i = s_{i+1}$, so that~$t_i>0$
(and~$t_{i+1}=0$).  Since~$n_i\ge 2$, the recipe gives~$m_i=2$,
again, as desired. 

Finally, if~$s_i$ is a multiple root of~$p_D(t)$, then it cannot
be a root of~$p_{k'}(t)$ at all, by the definition of the momentum
cell~$C(p_D,\mu)$.  Consequently,~$t_i=0$ by definition and calculation
shows that~$v_i=0$ as well.  Thus~$m_i=0$, as desired.

By Theorem~\ref{thm: point data to polys}, the Bochner-K\"ahler structure on
a neighborhood~$M$ of~$0\in\C{n}$ that has a unitary 
coframe~$u_0:T_0M\to\C{n}$ with~$\bigl(H(u_0),T(u_0),V(u_0)\bigr)=(s,t,v)$
has~$p_C(t)$ and $p_D(t)$ as its characteristic polynomials and satisfies~
$h'(0) = k'$.  This establishes existence for the 
boundary points of~$C(p_D,\mu)$. 
\end{proof}

The way is now paved for the following result, which, together
with the previous proposition, classifies the analytically
connected equivalence classes of Bochner-K\"ahler structures.

\begin{theorem}\label{thm: analytically connected classes}
The analytically connected class of a Bochner-K\"ahler
structure is determined by~$p_C$, $p_D$, and the momentum cell~$C(p_D,\mu)$
that contains the reduced momentum image.  Moreover, for any 
Bochner-K\"ahler structure with data~$\bigl(p_C,p_D,\mu\bigr)$, the
union of the reduced momentum images of the Bochner-K\"ahler
structures that are analytically connected to it is the entire
momentum cell~$C(p_D,\mu)$.
\end{theorem}

\begin{proof}
It has been established that~$p_C$ and~$p_D$ and the momentum 
cell~$C(p_D,\mu)$ are invariants of the analytically connected equivalence 
class.  Moreover, by Proposition~\ref{prop: values of ph'}, every point 
of~$C(p_D,\mu)$ lies in the image of the reduced momentum mapping of 
some Bochner-K\"ahler structure.

To prove Theorem~\ref{thm: analytically connected classes}, 
it will thus suffice to show that any two
Bochner-K\"ahler structures with the same data~$(p_C,p_D,\mu)$
are analytically connected.

Now, if~$(M,g,\Omega)$ and~$(\tilde M,\tilde g,\tilde \Omega)$ are
connected Bochner-K\"ahler manifolds with the same data~$(p_C,p_D,\mu)$
and their reduced momentum images~$h'(M)$ and $h'(\tilde M)$ 
have nontrivial intersection, then they contain points~$x\in M$
and~$\tilde x\in\tilde M$ so that~$f(x) = \tilde f(\tilde x)$ 
where~$f:M\to\bbR^{2n+1}$ and~$\tilde f:\tilde M\to \bbR^{2n+1}$ are
the corresponding moduli maps.  By Theorem~\ref{thm: existence} 
 and Corollary~\ref{cor: BK germ moduli},
the germs of Bochner-K\"ahler structures around~$x\in M$
and~$\tilde x\in \tilde M$ are isomorphic.  Since~$M$ and $\tilde M$
are connected, the germ of the Bochner-K\"ahler structure
around any~$y\in M$ is analytically connected to the germ of the
Bochner-K\"ahler structure around any~$\tilde y\in M$.  

Now, from Theorem~\ref{thm: cohomogeneity m}, it follows that~$h'(M^\circ)$ 
lies in the
interior of~$C(p_D,\mu)$ and that~$h':M^\circ\to C(p_D,\mu)^\circ$ is
a submersion onto its image, which is therefore open.

The union of the open sets~$h'(\tilde M^\circ)$ 
as~$(\tilde M,\tilde g,\tilde\Omega)$ 
ranges over the Bochner-K\"ahler structures that are analytically
connected to any given~$(M,g,\Omega)$ is a connected component 
of~$C(p_D,\mu)^\circ$.  Since~$C(p_D,\mu)^\circ$ is convex and hence
connected, this union must be all of~$C(p_D,\mu)^\circ$.

By Proposition~\ref{prop: values of ph'}, the union of all the sets~$h'(M)$ 
as~$(M,g,\Omega)$ 
ranges over the Bochner-K\"ahler structures with data~$(p_C,p_D,\mu)$ 
is equal to the entire cell~$C(p_D,\mu)$. Since $h'(M^\circ)$
is a nonempty subset of~$C(p_D,\mu)^\circ$ for any 
such~$(M,g,\Omega)$, it follows that all of these are analytically 
connected, as desired. 
\end{proof}

\begin{remark}[Coarse moduli and polytope embeddings]
By Theorem~\ref{thm: analytically connected classes}, 
the analytically connected equivalence classes 
in~$F_n$ correspond to the data~$(p_C,p_D,\mu)$ that satisfy the 
conditions of Proposition~\ref{prop: values of ph'}.  
Note that, for any given~$p_C(t)$, 
there are at most a finite number of choices of~$(p_D,\mu)$ that
will satisfy these constraints.  Thus, each value of~
$C = (C_2,\ldots,C_{n+2})$ corresponds to only a finite number of
equivalence classes.  It is in this sense that the 
functions~$C_i:F_n\to\bbR$ furnish the complete set of `coarse moduli'
for Bochner-K\"ahler structures in dimension~$n$.

Moreover, the mapping~$\Delta:F_n\to\bbR^{2n+1}$ of~\S\ref{sssec: mom map} 
embeds each analytically connected equivalence class as a (not necessarily 
closed) convex polytope of some dimension~$m\le n$.  In particular, each of
these equivalence classes is contractible.
\end{remark}

\begin{corollary}\label{cor: complete fibration}
If~$(M,g,\Omega)$ is a complete, connected Bochner-K\"ahler structure, 
then the reduced momentum mapping~$h':M\to C(p_D,\mu)$ is surjective;
the submersion~$h':M^\circ\to C(p_D,\mu)^\circ$ is 
a fibration; and the fibers of~$h'$ in~$M$ are connected.
\end{corollary}

\begin{proof}
Fix~$x\in M$.  By Theorem~\ref{thm: analytically connected classes}, 
to prove that~$h'$ is surjective it suffices
to show that if~$(\tilde M, \tilde g,\tilde \Omega)$ is any Bochner-K\"ahler
structure containing an~$\tilde x$ with~$f(x)=\tilde f(\tilde x)$, then
$\tilde h'(\tilde M)$ is a subset of~$h'(M)$.  

Consider any~$\tilde y\in \tilde M$ and choose a smooth 
path~$\tilde c:[0,1]\to\tilde M$ with~$\tilde c(0) = \tilde x$ 
and~$\tilde c(1) = \tilde y$.  Choose a~$\tilde u_0\in \tilde P_{\tilde x}$
and let~$\tilde u:[0,1]\to \tilde P$ be the parallel transport 
of~$\tilde u_0$ along~$\tilde c$.  Thus~$(\tilde u)^*(\tilde \phi) = 0$
while~$(\tilde u)^*(\tilde \omega) = v(s)\,ds$ for some~$v:[0,1]\to\C{n}$.

Since~$f(x)=\tilde f(\tilde x)$ by hypothesis, there exists a~$u_0\in P_x$
so that~$H(u_0) = \tilde H(\tilde u_0)$, $T(u_0) = \tilde T(\tilde u_0)$, 
and $V(u_0) = \tilde V(\tilde u_0)$.  Since the metric on~$M$ is complete,
there will exist a unique curve~$u:[0,1]\to P$ satisfying the
initial condition~$u(0)=u_0$ and the ordinary differential equations
\begin{equation*}
u^*(\omega) = v(s)\,ds,\qquad\qquad u^*(\phi) = 0.
\end{equation*}
I.e., $u$ is the parallel transport of~$u_0$ along the 
curve~$c=\pi{\circ} u$ in~$M$.  
The structure equations~\eqref{eq: structure equations ii}
and the Chain Rule now imply that the two curves defined on~$[0,1]$
\begin{equation*}
\bigl(H{\circ}u, T{\circ}u, V{\circ}u \bigr)
\qquad\text{and}\qquad
\bigl(\tilde H{\circ}\tilde u, 
\tilde T{\circ}\tilde u, \tilde V{\circ}\tilde u \bigr)
\end{equation*}
in~$i\euu(n){\oplus}\C{n}{\oplus}\bbR$ satisfy the same initial conditions
and system of ordinary differential equations, so they are equal on~$[0,1]$. 
Now, setting~$s=1$ yields~$\tilde f(\tilde y) = f\bigl(c(1)\bigr)\in h'(M)$.  
Since~$\tilde y\in \tilde M$ was arbitrary,~$\tilde h'(\tilde M)$
is a subset of~$h'(M)$.

To prove the second part, suppose first that~$M$ is simply-connected. 
Then any local isometry~$\psi:U\to M$ defined 
on an open subset~$U\subset M$ extends uniquely to a global 
isometry of~$M$.%
\footnote{Simply choose~$x\in U$ and extend~$\psi$ so that it commutes
with the exponential map at~$x$, i.e., so that $\psi\bigl(\exp_x(v)\bigr)
=\exp_{\psi(x)}\bigl(\psi'(x)(v)\bigr)$. The  completeness and analyticity 
of the metric and the simple connectivity of~$M$ 
imply that such an extension of~$\psi$ to all of~$M$ exists 
and is an isometry, as desired.} 
Thus~$I(M)$, the global holomorphic isometry group of~$M$, 
acts transitively on the fibers of~$h'(M)$. 

Now~$(h')^{-1}\bigl(C(p_D,\mu)^\circ\bigr) = M^\circ$ and, by the first part 
of the proof,~$h':M^\circ\to C(p_D,\mu)^\circ$ is surjective.
Since~$h':M^\circ\to C(p_D,\mu)^\circ$ is also a submersion whose
fibers are $I(M)$-orbits, it follows
that~$h':M^\circ\to C(p_D,\mu)^\circ$ is a fibration.

To treat the case where~$M$ is not simply-connected, pass to the
universal cover~$\tilde M$ and note that~$\tilde h'$ is invariant under 
the deck transformations~$\Delta$ of the cover~$\tilde M\to M$ 
(which form a discrete subgroup of~$I(\tilde M)$).  Since~$\tilde h':
\tilde M^\circ\to C(p_D,\mu)^\circ$ is a fibration, 
dividing by the (free) action of~$\Delta$ yields 
that $h':M^\circ\to C(p_D,\mu)^\circ$ is also a fibration.

To prove the connectedness of the fibers of~$h'$ it suffices to treat
the case where~$M$ is simply-connected, so assume this.
Again,~$I(M)$ acts transitively on the fibers of~$h'$.   
Since~$M^\circ$ is the complement of a codimension~1 complex subvariety
in~$M$, it is connected.  The exact sequence of the
fibration~$h':M^\circ\to C(p_D,\mu)^\circ$ and the contractibility of
$C(p_D,\mu)^\circ$ then imply that the fibers 
of~$h':M^\circ\to C(p_D,\mu)^\circ$ are connected.  
By Corollary~\ref{cor: sym stab alg},
the $I(M)$-stabilizer of any point in~$M$ is a product of unitary groups
and hence is connected.  Since the $I(M)$-orbit of any~$x\in M^\circ$
has been shown to be connected, it follows that~$I(M)$ must be
connected as well.  Consequently, all of its orbits in~$M$ are connected
and these are the fibers of~$h'$.
\end{proof}

\subsection{A Riemannian submersion}
\label{ssec: Riem submer}

The mapping~$h':M\to C(p_D,\mu)\subset\bbR^m$ 
can be used to give more detailed information
about the Bochner-K\"ahler metric.  

\subsubsection{The cell metric}\label{sssec: cell metric}

The proof of Corollary~\ref{cor: complete fibration} 
shows that, at least when~$M$ is complete, 
there is a metric on~$C(p_D,\mu)^\circ$ for 
which~$h':M^\circ\to C(p_D,\mu)^\circ$ is a Riemannian submersion.  
It turns out that this metric exists even when~$M$ is not complete and can
be identified explicitly.

\begin{theorem}\label{thm: cell metric from pD} 
Given a monic polynomial $p_D(t)$ of degree~$m{+}2$
that falls into one of the cases of 
Proposition~{\upshape\ref{prop: pD root cases}}, there exist 
rational functions~$R^{ij}_D=R^{ji}_D$ on~$\bbR^m$ with the property 
that the quadratic form
\begin{equation}\label{eq: cell metric}
R_D = R^{ij}_D(u)\,du_i\,du_j
\end{equation}
restricts to be positive definite on the interior of each possible 
momentum cell for~$p_D(t)$ and moreover, so that for any Bochner-K\"ahler 
manifold with reduced momentum polynomial~$p_D(t)$, the reduced momentum 
mapping~$h':M\to\bbR^m$ is a Riemannian submersion when restricted 
to~$M^\circ$.
\end{theorem}

\begin{proof}
Recall the bundle~$P_2\subset P_1$ that was introduced in the proof of
Proposition~\ref{prop: ph'' roots are pD roots}
and the notation introduced there.  On~$P_2$, the 
matrix~$H$ is diagonal and~$H_{i\bar\imath} = \lambda_i$ for $1\le i\le m$.
The structure equation for~$dH_{i\bar\imath}$ then becomes
\begin{equation}\label{eq: lambda-T relation}
d\lambda_i = T_i(\omega_i+\overline{\omega_i}).
\end{equation}
By equation~\eqref{eq: pC to ph'' ratio},
\begin{equation*}
{T_i}^2 = \frac{p_D(\lambda_i)}{\prod_{j\not=i}(\lambda_i-\lambda_j)},
\end{equation*}
so \eqref{eq: lambda-T relation} can be written in the form
\begin{equation*}
\Re(\omega_i) 
= \frac12\sqrt{\frac{\prod_{j\not=i}(\lambda_i-\lambda_j)}{p_D(\lambda_i)}}\,
d\lambda_i\,.
\end{equation*}
In other words,
\begin{equation*}
\sum_{i=1}^m \Re(\omega_i) ^2
= \frac14\,\sum_{i=1}^m
\frac{\prod_{j\not=i}(\lambda_i-\lambda_j)}{p_D(\lambda_i)}\,
{d\lambda_i}^2\,.
\end{equation*}
Now suppose that~$h'(M)$ lies in~$C(p_D,\mu)$ and that
\begin{equation*}
\lambda\bigl(C(p_D,\mu)\bigr) 
= I_1\times I_2\times\cdots\times I_m\subset \bbR^m_{\ge}
\end{equation*}
as in \S\ref{sssec: prod rep}.  
Since~$(-1)^{i-1}p_D(y_i)>0$ for $y_i\in I^\circ_i$, 
the quadratic form
\begin{equation}\label{eq: S from y}
S =  \frac14\,\sum_{i=1}^m
\frac{\prod_{j\not=i}(y_i-y_j)}{p_D(y_i)}\,
{dy_i}^2
\end{equation}
is positive definite 
on~$I^\circ_1\times I^\circ_2\times\cdots\times I^\circ_m\subset\bbR^m_>$.
Since~$S$ has rational coefficients, is well-defined 
on~$\bbR^m$ minus the union of the hyperplanes~$p_D(y_i)=0$, 
and is invariant under permutations 
of the coordinates~$y_i$, it follows that there are unique 
rational functions~$R^{ij}_D=R^{ji}_D$ on~$\bbR^m$ so that
\begin{equation*}
S = R^{ij}_D\bigl(\sigma(y)\bigr)
\,d\bigl(\sigma_i(y)\bigr)\,d\bigl(\sigma_j(y)\bigr).
\end{equation*}
I.e., setting~$R_D = R^{ij}_D(u)\,du_i\,du_j$, the positive definite 
quadratic form~$S$ defined 
on~$I^\circ_1\times I^\circ_2\times\cdots\times I^\circ_m$
is of the form~$\sigma^*(R_D)$.  
Since~$\sigma:I^\circ_1\times I^\circ_2\times\cdots\times I^\circ_m
\to C(p_D,\mu)^\circ$ is a diffeomorphism,~$R_D$ is 
positive definite on~$C(p_D,\mu)^\circ$.  
Moreover, the above formula on~$M^\circ$ can now be 
written as
\begin{equation*}
\sum_{i=1}^m \Re(\omega_i) ^2  = (h')^*\bigl(R_D\bigr),
\end{equation*}
showing that~$h':M^\circ\to C(p_D,\mu)^\circ$ is a Riemannian submersion
when the target is given the Riemannian metric~$R_D$. \end{proof}

\subsubsection{Explicit formulae}\label{sssec: explicit forms}

When~$m=2$, with~$u_1 = y_1+y_2$ and~$u_2 = y_1y_2$, 
the relation between~$S$ and~$R$ is expressible 
in the intermediate form
\begin{equation*}
\begin{split}
4S = \frac{y_1-y_2}{p_D(y_1)}\,{dy_1}^2 + \frac{y_2-y_1}{p_D(y_2)}\,{dy_2}^2
&= \frac{\bigl({y_1}^2p_D(y_2){-}{y_2}^2p_D(y_1)\bigr)}
        {(y_1{-}y_2)p_D(y_1)p_D(y_2)}\,{du_1}^2\\
&\quad - 2\frac{\bigl(y_1p_D(y_2){-}y_2p_D(y_1)\bigr)}
               {(y_1{-}y_2)p_D(y_1)p_D(y_2)}\,du_1\,du_2\\
&\quad  + \frac{p_D(y_2){-}p_D(y_1)}{(y_1{-}y_2)p_D(y_1)p_D(y_2)}\,{du_2}^2.\\
\end{split}
\end{equation*}
Each of the coefficients on the right is visibly a symmetric rational
function of~$y_1$ and~$y_2$ and so can be written as a rational 
function of~$u_1$ and~$u_2$.  

Notice that the expression~$p_D(y_1)p_D(y_2)$
is a common denominator of all these rational expressions.  In the 
case where~$p_D(t)$ has all four of its roots real, this can be written
in the form
\begin{equation*}
p_D(y_1)p_D(y_2) = \prod_{\alpha=0}^3(y_1-r_\alpha)(y_2-r_\alpha) 
                 = \prod_{\alpha=0}^3(u_2-r_\alpha\,u_1+{r_\alpha}^2),
\end{equation*}
so that the denominator is actually a product of linear functions
on~$\bbR^2$, indeed, the very linear functions whose vanishing
defines the faces of possible momentum cells.

 For any~$m$, in Case~4, in which all of the roots of~$p_D$ are real and 
distinct, this generalizes, leading to an expression for the metric~$R_D$ that
will turn out to be very useful.  As usual, let
\begin{equation*}
r_0>r_1>\cdots>r_m>r_{m+1}
\end{equation*} 
be the real roots of~$p_D$. Since the roots are real and distinct,
 $(-1)^\alpha p'_D(r_\alpha)>0$ for~$0\le\alpha\le m{+}1$.
Note the identity
\begin{equation*}
p'_D(r_\alpha) = \prod_{\beta\not=\alpha}(r_\alpha{-}r_\beta),
\end{equation*}
as well as the classical identities
\begin{equation*}
\sum_{\alpha=0}^{m+1} \frac{{r_\alpha}^k}{p'_D(r_\alpha)}
= \begin{cases} {\displaystyle\frac{(-1)^{m+1}}{r_0r_1\cdots r_{m+1}}} 
                                      & \text{when $k = -1$}; \\
         \hfil 0 & \text{when $0\le k\le m$}; \\ 
         \hfil 1 & \text{when $k = m{+}1$;} \\ 
         r_0 +\cdots + r_{m+1} & \text{when $k = m{+}2$.} \\ 
  \end{cases}
\end{equation*}
(The cases $k=-1$ and $k=m{+}2$ will be used in a later section.)
Using coordinates~$(u_1,\ldots,u_m)$
on~$\bbR^m$ as above, define linear functions%
\footnote{The reader will note that the formula for~$l_\alpha$ makes
sense in general as long as $r_\alpha$ is a simple root of~$p_D(t)$.
Accordingly, $l_\alpha$ will taken to be defined 
by~\eqref{eq: define l-alpha} in
this more general case.}
\begin{equation}\label{eq: define l-alpha}
l_\alpha 
= -\frac{\,(r_\alpha^m{-}r_\alpha^{m-1}\,u_1+\cdots+(-1)^mu_m)\,}
        {p'_D(r_\alpha)}.
\end{equation}
for~$0\le \alpha\le m{+}1$, so that the equations~$l_\alpha=0$
define the hyperplanes that are the faces of the various possible
momentum cells for~$p_D$.  Note that the above classical identities
in the range $0\le k\le m{+}1$ are equivalent to the equations
\begin{equation}\label{eq: l relations}
\sum_{\alpha=0}^{m+1} l_\alpha = 0
\qquad\text{and}\qquad
\sum_{\alpha=0}^{m+1} r_\alpha\, l_\alpha = -1,
\end{equation}
which are the only linear relations among the~$l_\alpha$.  Note also that
\begin{equation*}
\sigma^*(l_\alpha) 
          = \frac{\prod_{i=1}^m(y_i-r_\alpha)}
        {\,\prod_{\beta\not=\alpha}(r_\beta{-}r_\alpha)\,}.
\end{equation*}

The metric~$R_D$ then has the simple expression
\begin{equation}\label{eq: R in l-terms}
R_D =  \sum_{\alpha=0}^{m+1} \frac{{dl_\alpha}^2}{4\,l_\alpha}.
\end{equation}
Indeed,
\begin{equation*}
\begin{split}
\sigma^*\left(\sum_{\alpha=0}^{m+1}
   l_\alpha\,\left(\frac{dl_\alpha}{l_\alpha}\right)^2\right)
&= \sum_{\alpha=0}^{m+1}
   \frac{\prod_{i=1}^m(y_i-r_\alpha)}
        {\,\prod_{\beta\not=\alpha}(r_\beta{-}r_\alpha)\,}\,
       \left(\sum_{j=1}^m\frac{dy_j}{(y_j{-}r_\alpha)}\right)^2\\
&= \sum_{j,k=1}^m \left(\sum_{\alpha=0}^{m+1}
   \frac{\prod_{i=1}^m(y_i{-}r_\alpha)}
        {\,(y_j{-}r_\alpha)(y_k{-}r_\alpha)\,
            \prod_{\beta\not=\alpha}(r_\beta{-}r_\alpha)\,}\right)\,
             dy_j\,dy_k\\
&=\sum_{i=1}^m
\frac{\prod_{j\not=i}(y_i{-}y_j)}{p_D(y_i)}\,
{dy_i}^2\,,\\
\end{split}
\end{equation*}
where the last equality follows from the classical identities
above and the Lagrange interpolation identity.  

The expression~\eqref{eq: R in l-terms}
 can also be written in Hessian form as
\begin{equation*}
R_D = R^{ij}_D\,du_i\,du_j 
= \frac{\partial^2G}{\partial u_i\,\partial u_j}\,du_i\,du_j 
\end{equation*}
where the potential function~$G$ has the form
\begin{equation*}
G = \frac{1}4\, \sum_{\alpha=0}^{m+1}
            l_\alpha\,\bigl(\log|l_\alpha|-1\bigr),
\end{equation*}
a fact that will be useful below.

The formula for~$R_D$ is evidently singular along the 
hyperplanes~$l_\alpha=0$, but this singularity is mild and can be 
`resolved' with little difficulty.  For simplicity, and since this case
will be useful in the analysis below, I will illustrate this for the
`lowest' cell, i.e., Subcase 4-0.  The spectral intervals in this subcase
are~$I_i = [r_i,r_{i+1}]$ and the functions $l_1,\ldots l_{m+1}$ are all
nonnegative on this cell~$C(p_D,\mu)$, which is an $m$-simplex.   
In fact, \eqref{eq: l relations}
 shows that~$l_0 = -(l_1+\cdots+l_{m+1})$ and that
\begin{equation}\label{eq: l-relations i}
1 = \sum_{\alpha=1}^{m+1} (r_0{-}r_\alpha)\,l_\alpha \,,
\end{equation}
so that the functions~$(r_0-r_\alpha)\,l_\alpha$ for~$1\le\alpha\le m{+}1$
can be regarded as homogeneous affine coordinates on this simplex.   Now
let~$E\subset\bbR^{m+1}$ be the $m$-dimensional ellipsoid defined by
\begin{equation}\label{eq: define ellipsoid i}
1 = \sum_{\alpha=1}^{m+1} (r_0{-}r_\alpha)\,{p_\alpha}^2 \,.
\end{equation}
There is then a unique smooth map~$s:E\to C(p_D,\mu)$ defined by
$s^*(l_\alpha) = p_\alpha^2$ for~$1\le \alpha\le m{+}1$.
Since~$s^*(l_0) = -({p_1}^2 + \cdots + {p_{m+1}}^2)$,
the $s$-pullback metric is
\begin{equation}\label{eq: s-pullback metric}
s^*(R_D) = s^*\left(\sum_{\alpha=0}^{m+1} \frac{{dl_\alpha}^2}
                                              {4\,l_\alpha}\right)
= \sum_{\alpha=1}^{m+1}{dp_\alpha}^2
    \ \ -\ \frac{(p_1\,dp_1 + \cdots + p_{m+1}\,dp_{m+1})^2}
       {({p_1}^2 + \cdots + {p_{m+1}}^2)}.
\end{equation}
The quadratic form on the right hand side is well-defined on~$\bbR^{m+1}$
minus the origin and is positive semidefinite there, 
with the null space of the quadratic form being spanned
by the radial vector at each point.%
\footnote{In fact, this metric is just the tangential part~$r^2\,d\sigma^2_m$
of the expression for the standard metric in polar coordinates~$dr^2 + r^2\,
d\sigma^2_m$, where~$d\sigma^2_m$ is the standard metric on the $m$-sphere.
A curious consequence of this fact is that the metric~$R_D$ is 
conformally flat.}
  Thus, this quadratic form is 
positive definite and smooth on~$E$, thereby providing the desired
`resolution' of the singularities of~$R_D$ on~$C(p_D,\mu)$.  Note that the
rank of the mapping~$s$ at~$p = (p_\alpha)\in E$ is equal to one
less than the number of nonzero entries~$p_\alpha$.  This will be useful 
below.

The analysis of~$R_D$ in Cases 1, 2, and 3 can be
derived from the Case 4 analysis by either regarding two of the roots 
as complex conjugates and combining the corresponding terms in the above sums
to obtain real expressions (Case 1) or collecting two or three of the
terms and taking the limit as the corresponding roots come together
(Cases 2 and 3).  This will only be needed in Case 3-1$b$ below, 
so I will do this case and leave the others to the interested reader.

Case 3-1 can be regarded as the limit of Case 4 as the root~$r_0$
approaches~$r_1$ while the roots~$r_1$ through~$r_{m+1}$ remain fixed.
Thus, the relations~\eqref{eq: l relations}
 can be solved for~$l_0$ and $l_1$ in the
form
\begin{equation*}
(r_1{-}r_0)\,l_0 = 1 - \sum_{\alpha=2}^{m+1} (r_1{-}r_\alpha)\,l_\alpha
\qquad\qquad
(r_0{-}r_1)\,l_1 = 1 - \sum_{\alpha=2}^{m+1} (r_0{-}r_\alpha)\,l_\alpha
\end{equation*}
Using these formulae, one computes the limit
\begin{equation*}
\lim_{r_0\to r_1}\left(\frac{{dl_0}^2}{l_0} + \frac{{dl_1}^2}{l_1}\right)
= \frac{t\,{da}^2}{a^2} - \frac{2\,da\,dt}{a},
\end{equation*}
where
\begin{equation}\label{eq: define a and t}
a = 1 - \sum_{\alpha=2}^{m+1} (r_1{-}r_\alpha)\,l_\alpha
\qquad\text{and}\qquad
t = l_2 + l_3 + \ldots + l_{m+1},
\end{equation}
and where the formulae~\eqref{eq: define l-alpha}
 for~$l_\alpha$ for~$2\le\alpha\le m{+}1$
remain valid.  Thus, the formula for~$R_D$ in Case~3-1 is 
\begin{equation}\label{eq: R in case 3-1}
R_D = \frac{t\,{da}^2}{4a^2} - \frac{da\,dt}{2a} 
       + \sum_{\alpha=2}^{m+1}\frac{{dl_\alpha}^2}{4l_\alpha}\,.
\end{equation}

In Case 3-1$b$, all of the quantities~$a,l_2,\ldots,l_{m+1}$ are
non-negative.  In fact, the quantity~$a>0$ together with the quantities
$(r_1{-}r_\alpha)\,l_\alpha$ for~$2\le\alpha\le m{+}1$ can
be regarded as affine homogeneous coordinates on the momentum
cell~$C(p_D,\mu)$. 

Set~$\rho_\alpha = r_1-r_\alpha>0$ and~$\rho=(\rho_2,\ldots,\rho_{m+1})$.
Let~$E_\rho\subset \bbR^m$ be the ellipsoidal domain defined by 
\begin{equation*}
\sum_{\alpha=2}^{m+1} \rho_\alpha\,{p_\alpha}^2 < 1.
\end{equation*}
Define a surjective map~$s:E_\rho\to C(p_D,\mu)$ by~$s^*l_\alpha 
= p_\alpha^2$ for~$2\le\alpha\le m{+}1$.  This map~$s$ satisfies
\begin{equation*}
\bar a = s^*(a) =  1 - \sum_{\alpha=2}^{m+1} \rho_\alpha\,{p_\alpha}^2
\qquad\text{and}\qquad
\bar t = s^*(t)  = {p_2}^2 + \ldots + {p_{m+1}}^2.
\end{equation*}
Thus~$R_\rho = s^*(R_D)$ has the form
\begin{equation*}
R_\rho = \frac{\bar t\,{d\bar a}^2}{4\bar a^2} 
          - \frac{d\bar a\,d\bar t}{2\bar a} 
       + \sum_{\alpha=2}^{m+1}{dp_\alpha}^2\,.
\end{equation*}
Since~$\bar t$ vanishes to second order at $p=0$ (the center of~$E_\rho$),
the quadratic form~$R_\rho$ is visibly positive definite and smooth
on a neighborhood of~$p=0$.  Let~$\delta>0$ be less than any
$1/\sqrt{\rho_\alpha}$ and consider the annular region~
$A_\delta\subset E_\rho$ 
defined by~$\bar t \ge \delta^2$.  On this region, $\bar a$ and $\bar t$ 
are both positive and~$R_\rho$ can be written in the form
\begin{equation*}
R_\rho 
= \frac{\bar t}4
  \left(\frac{d{\bar a}}{\bar a}-\frac{d{\bar t}}{\bar t} \right)^2
          - \frac{{d\bar t\,}^2}{4{\bar t\,}^2} 
       + \sum_{\alpha=2}^{m+1}{dp_\alpha}^2\,
= \frac{\bar t}4
  \left(\frac{d{\bar a}}{\bar a}-\frac{d{\bar t}}{\bar t} \right)^2
 + R^*_\rho\,,
\end{equation*}
where~$R^*_\rho$ is defined by this last equality.  Since~$\bar t = |p|^2$,
it follows without difficulty that $R^*_\rho$ is positive semidefinite
on $\bbR^m$ minus the origin (where it is singular) and that its null space
at each point is one dimensional and is spanned by the radial vector.  Since
the 1-form~$\rho = d\bar a/\bar a - d\bar t/\bar t$ is evidently nonvanishing
on the radial vector field, it follows that~$R_\rho$ is positive definite
(and smooth) everywhere on~$E_\rho$.   Thus,~$R_\rho$ on~$E_\rho$ provides
the desired resolution of the boundary singularities of~$R_D$ on the
momentum cell~$C(p_D,\mu)$. 

Moreover, on~$A_\delta$, whose outer boundary is defined by~$\bar a = 0$,
the inequality
\begin{equation*}
R_\rho \ge \frac{\delta^2}4
  \left(\frac{d{\bar a}}{\bar a}-\frac{d{\bar t}}{\bar t} \right)^2
 = \left(\frac\delta2\,
      d\left(\log\left(\frac{\bar t}{\bar a}\right)\,\right)\,\right)^2
\end{equation*}
holds.  Since~$\log(\bar t/\bar a)$ is proper on~$A_\delta$,
it follows that~$R_\rho$ is complete on~$E_\rho$.  

This result will be needed in~\S\ref{sssec: leaf metric}, 
when completeness is being discussed.  For use in that
section, I will point out that the above formulae define a
convex domain~$E_\rho \subset\bbR^m$ and a complete metric~$R_\rho$ 
on~$E_\rho$ for any~$\rho=(\rho_2,\ldots,\rho_{m+1})$ 
satisfying~$\rho_\alpha\ge0$ for all~$\alpha$.  (Recall that the metrics
that arise as resolutions of singular metrics on~$C(p_D,\mu)$ satisfy
$0<\rho_2<\cdots<\rho_{m+1}$.)

The metric~$R_\rho$ is flat only when~$\rho_2=\cdots=\rho_{m+1}=0$,
in which case~$E_0=\bbR^m$ and~$R_0$ is the standard flat metric.
Moreover, the above formulae show that $R_\rho$ is always conformally flat,
with~$(E_\rho,R_\rho)$ being globally conformal to~$(E_0,R_0)$.

\subsubsection{Necessary conditions for completeness}
\label{sssec: ness conds for complete}

It turns out that most of the possible momentum cells cannot be the
reduced momentum image of a complete Bochner-K\"ahler manifold.

\begin{proposition}\label{prop: metric complete implies cell bounded}
If there is a complete Bochner-K\"ahler~$(M,g,\Omega)$ 
whose reduced momentum mapping has image in~$C(p_D,\mu)$,
then~$C(p_D,\mu)$ is bounded.
\end{proposition}

\begin{proof}
Suppose that~$(M,g,\Omega)$ is connected and complete, 
with characteristic polynomials~$p_C$ and~$p_D$ but that its reduced momentum
mapping takes values in an unbounded momentum cell~$C(p_D,\mu)$.
Let~$I_1{\times}\cdots{\times}I_m$ be the corresponding 
spectral product.  The unboundedness of~$C(p_D,\mu)$ implies that~$I_1$ 
is either~$[r,\infty)$ or~$(r,\infty)$ where~$r$ is the largest real
root of~$p_D$.  

Again, let~$P_2\subset P_1$ be the bundle over~$M^\circ$ constructed 
in the course of the proof of Proposition~\ref{prop: ph'' roots are pD roots}.  
Because the structure group 
of~$P_2$ is~$\I_m\times G_\Lambda$, it follows that the 1-forms~$\omega_i$ 
for~$1\le i\le m$ are actually well-defined on~$M^\circ$.  Let~$E_1$
be the vector field on~$M^\circ$ that is $g$-dual to~$\Re(\omega_1)$.
Then, by the relation~$d\lambda_i = T_i(\omega_i+\overline{\omega_i})$,
it follows that $d\lambda_i(E_1)=0$ for~$1<i\le m$ and that
\begin{equation*}
d\lambda_1(E_1) = 2T_1 
= 2\sqrt{\frac{p_D(\lambda_1)}{\prod_{j\not=1}(\lambda_1-\lambda_j)}} >0.
\end{equation*}
In particular, along an integral curve of~$E_1$
the functions~$\lambda_j$ for~$1<j\le m$ are constant while~$\lambda_1$
is strictly increasing.  

Fix~$x\in M^\circ$ and let~$a:[0,T)\to M$ be the maximal forward integral
curve of~$E_1$ with~$a(0)=x$.  I claim that~$T$ cannot be finite.
If it were, the fact that~$E_1$ is a unit speed vector field and that~$M$
is complete would imply that~$a(t)$ approaches a limit~$y\in M$ as~$t$
approaches~$T$ (after all,~$d\bigl(a(t),a(s)\bigr)\le |t-s|$).  
The limit point~$y$ could not lie in~$M^\circ$ since then~$[0,T)$ 
would not be maximal. By continuity, $\lambda(y)$ must not lie 
in~$I^\circ_1{\times}\cdots{\times}I^\circ_m$.  However,~$\lambda_i(y)
=\lambda_i(x)$ for~$1<i\le m$ while~$\lambda_1(y)>\lambda_1(x)$.
Since~$I^\circ_1 = (r,\infty)$ this forces~$\lambda(y)$ to lie 
in~$I^\circ_1{\times}\cdots{\times}I^\circ_m$, a contradiction
since~$h'(M^\circ)$ lies in~$C(p_D,\mu)^\circ$.  Thus~$T=\infty$,
as claimed.  In particular, the forward flow of~$E_1$ exists for
all time on~$M^\circ$.

However, this leads to a contradiction:
Along~$a$, the element of arc is given by 
\begin{equation*}
ds = \frac12
      \sqrt{\frac{\prod_{j=2}^m
        \bigl(\lambda_1-\lambda_j(x)\bigr)}{p_D(\lambda_1)}}
     \,d\lambda_1\,.
\end{equation*}
Let~$\lim_{t\to\infty}\lambda_1\bigl(a(t)\bigr) = \lambda_\infty\le \infty$. 
Since~$a:[0,\infty)\to M$ has unit speed, the integral
\begin{equation*}
\int_{\lambda_1(x)}^{\lambda_\infty} \sqrt{\frac{\prod_{j=2}^m
           \bigl(\xi-\lambda_j(x)\bigr)}{p_D(\xi)}}
     \,d\xi
\end{equation*}
must be infinite.  However, this integral is bounded by
\begin{equation*}
\int_{\lambda_1(x)}^\infty \sqrt{\frac{\prod_{j=2}^m
           \bigl(\xi-\lambda_j(x)\bigr)}{p_D(\xi)}}
     \,d\xi\,
\end{equation*}
which converges, since~$p_D(t)$ has degree~$m{+}2$.

This contradiction implies that $(M,g)$ could not have been complete. 
\end{proof}

\begin{remark}[Bounded momentum cells]
The discussion in~\S\ref{sssec: pos mo cells} 
shows that there are only two cases in which the momentum cell is bounded:

The first case is SubCase~3-1$b$, i.e.,~$p_D$ has $r_1$ as a double root
and~$\mu_1=0$.  The spectral bands are~$I_1 = 
[r_2,r_1)$ and~$I_j = [r_{j+1},r_j]$ for~$1<j\le m$. This cell
is bounded but not compact.

The second case is SubCase 4-0, i.e., $p_D(t)$ has $m{+}2$ simple
roots~$r_0>\cdots>r_{m+1}$ 
and the spectral bands are~$I_j = [r_{j+1},r_j]$ for~$1\le j\le m$.  
This cell is compact.  However, as the next proposition shows, 
Subcase 4-0 never contains a complete example when~$m>0$.
\end{remark}

\begin{proposition}\label{prop: distinct roots is not complete}
When~$m>0$, there is no complete Bochner-K\"ahler manifold
whose reduced characteristic polynomial~$p_D$ has $m{+}2$ 
distinct roots.
\end{proposition}

\begin{proof}
In view of Proposition~\ref{prop: metric complete implies cell bounded} 
and the remark above, what has to be shown
is that SubCase 4-0 cannot occur for a complete Bochner-K\"ahler
manifold when~$m>0$.  This will involve an interesting examination
of the fixed points of the flow of the canonical torus action.

Thus, suppose, to the contrary, that~$(M,g,\Omega)$ is a complete
Bochner-K\"ahler structure with~$m>0$ and and that
\begin{equation*}
p_D(t) = (t-r_0)(t-r_1)\cdots(t-r_{m+1})
\end{equation*}
where~$r_0>\cdots> r_{m+1}$.
By Proposition~\ref{prop: metric complete implies cell bounded}, 
the momentum cell~$C(p_D,\mu)$ must be bounded, 
which implies that~SubCase 4-0 obtains, 
namely~$(-1)^{i-1}p_{h'}(r_i)\ge0$ for $1\le i\le m{+}1$.

For~$1\le\alpha\le m{+}1$, let~$F_\alpha\subset C(p_D,\mu)$ 
be the $\alpha$-th face of this $m$-simplex, i.e., 
the intersection of~$C(p_D,\mu)$ with the hyperplane~$l_\alpha = 0$ 
(where the functions~$l_\alpha$ 
are as defined in~\eqref{eq: define l-alpha}).  Let~$N_\alpha
= (h')^{-1}(F_\alpha)$ be the preimage of~$F_\alpha$.  Evidently, 
each~$N_\alpha$ is a closed, analytic subset of~$M$ 
and the union of the~$N_\alpha$ is the complex 
subvariety~$N\subset M$.  Thus,~$N_\alpha$ is a
(non-empty) complex subvariety of~$M$.

For~$0\le \alpha\le m{+}1$, define functions~$w_\alpha=(h')^*(l_\alpha)\ge0$
on~$M$ and then define vector fields~$W_\alpha\in\euz$ 
by~$W_\alpha\lhk\Omega = -dw_\alpha$.   
By~\eqref{eq: l relations}, the~$W_\alpha$ satisfy
\begin{equation}\label{eq: W-relations}
\sum_{\alpha=0}^{m+1} W_\alpha 
= \sum_{\alpha=0}^{m+1} r_\alpha\,W_\alpha = 0.
\end{equation} 
Moreover, any~$m$ of these vector fields are linearly independent
on~$M^\circ$.  Note that since~$w_\alpha$ reaches its minimum of~$0$
along~$N_\alpha$, the vector field~$W_\alpha$ vanishes
along~$N_\alpha$.  Since~$M$ is complete, the flows of the vector 
fields~$W_\alpha$ are complete.  

I am going to show that the flow of each vector field~$W_\alpha$ 
is periodic of period~$\pi$ by examining the rotation~$\nabla W_\alpha$
along the fixed hypersurface~$N_\alpha$.

Now, equation~\eqref{eq: differential of ph} can be written as
\begin{equation}\label{eq: char pol diff}
t^n\,d\bigl(p_h(t^{-1})\bigr) = -t^n\,p_h(t^{-1})
\bigl(t\,\alpha_1 + t^2\,\alpha_2 + \cdots\bigr),
\end{equation}
where I have replaced~$t$ by~$-t$ and am regarding~$t$ as a parameter, taken
to be sufficiently small so that the series converges in a neighborhood
of any given compact domain in~$M$.
Using~\eqref{eq: define alpha}, this can be written in the form
\begin{equation}\label{eq: char pol diff ii}
d\bigl(p_h(t^{-1})\bigr) = -t\,p_h(t^{-1})
\sum_{k=0}^{\infty} \left(T^*\,(tH)^k\,\omega + \omega^*\,(tH)^k\,T\right)
\end{equation}
and the series can then be summed, yielding the equation
\begin{equation*}
d\bigl(p_h(t^{-1})\bigr) = -t\,p_h(t^{-1})
\left( T^*\,(\I_n{-}tH)^{-1}\,\omega + \omega^*\,(\I_n{-}tH)^{-1}\,T\right).
\end{equation*}
Replacing~$t$ by~$t^{-1}$, this becomes
\begin{equation}\label{eq: char pol diff iii}
\begin{split}
d\bigl(p_h(t)\bigr) 
&= -p_h(t)\left( T^*\,(t\,\I_n{-}H)^{-1}\,\omega
                    +\omega^*\,(t\,\I_n{-}H)^{-1}\,T\right)\\
&= - \left(T^*\,\Cof(t\,\I_n{-}H)\,\omega
                    +\omega^*\,\Cof(t\,\I_n{-}H)\,T\right).
\end{split}
\end{equation}
The final expression is valid for all~$t$, while the middle
expression is valid away from the locus~$p_h(t)=0$ in~$P\times\bbR$.

Since~$p_h(t) = p_{h''}(t)p_{h'}(t)$, and since~$p_{h''}(t)$ has
constant coefficients, \eqref{eq: char pol diff iii} implies
\begin{equation}\label{eq: char pol diff iv}
d\bigl(p_{h'}(t)\bigr) 
= -p_{h'}(t)\left( T^*\,(t\,\I_n{-}H)^{-1}\,\omega
                    +\omega^*\,(t\,\I_n{-}H)^{-1}\,T\right)
\end{equation}
away from the locus~$p_{h}(t)\not=0$.  

Define a vector field~$W(t)$ on~$M$ by $W(t)\lhk\Omega 
= -d\bigl(p_{h'}(t)\bigr)$.  This vector field depends polynomially
on~$t$ and lies in~$\euz$ for all~$t$.  
In fact, comparison with~\eqref{eq: define l-alpha},
the definition of~$l_\alpha$,
shows that~$p_{h'}(r_\alpha) = -p'_D(r_\alpha)\,w_\alpha$,
so it follows that
\begin{equation}\label{eq: W at r-alpha}
W(r_\alpha) = -p'_D(r_\alpha)\,W_\alpha\,,
\qquad\qquad 0\le \alpha\le m{+}1.
\end{equation}
By~\eqref{eq: char pol diff iv}, the vector field~$W(t)$ has representative
function~$w(t):P\to\C{n}$ given by
\begin{equation}\label{eq: define w(t)}
w(t) = -2i\,p_{h'}(t)\,(t\,\I_n{-}H)^{-1}T.
\end{equation}
The expression on the left is polynomial in~$t$, so the expression
on the right must be also.  Since the flow of~$W(t)$ is a holomorphic 
isometry, it follows that
\begin{equation*}
d\bigl(w(t)\bigr) + \phi\,w(t) = w'(t)\,\omega,
\end{equation*}
where~$w'(t)$ takes values in~$\euu(n)$.  
In fact, by~\eqref{eq: structure equations ii},
\begin{equation}\label{eq: w' formula}
\begin{split}
w'(t) &= 2i\,p_{h'}(t)\,(t\,\I_n{-}H)^{-1}\,
        \Bigl[TT^*(t\,\I_n{-}H)^{-1} - T^*(t\,\I_n{-}H)^{-1}T\,\I_n\\
 &\qquad\qquad\qquad\qquad\qquad\qquad\qquad - H^2 - h_1\,H-V\,I_n\Bigr]\\
\end{split}
\end{equation}
and the matrix on the right is visibly skew-Hermitian.
When~$T(u)=0$, formula~\eqref{eq: w' formula}
 simplifies to the form in which it will
be the most useful:
\begin{equation}\label{eq: eigen w}
w'(t)(u) = -2i\,p_{h'(u)}(t)\,\bigl(t\,\I_n{-}H(u)\bigr)^{-1}\,
        \Bigl[ H(u)^2 + h_1(u)\,H(u)+V(u)\,I_n\Bigr].
\end{equation}

Now fix~$\beta$ in the range~$1\le\beta\le m{+}1$ and 
let~$k_\beta\in C(p_D,\mu)$ be the vertex that lies on the intersection
of the faces~$F_\alpha$ for~$\alpha\not=0,\beta$, i.e., $k_\beta$ is
the vertex that lies \emph{opposite} the face~$F_\beta$.  
Applying Corollary~\ref{cor: complete fibration}, 
choose~$x_\beta\in M$ to satisfy~$h'(x_\beta)=k_\beta$ 
and then let~$u_\beta\in P$ satisfy~$\pi(u_\beta) = x_\beta$.
Then $T(u_\beta)=0$ since the differential of~$h'$ vanishes
at~$x_\beta$.  In particular, $x_\beta$ is a zero of~$W_\alpha$ 
for all~$\alpha$.

Now, $r_\alpha$ is a root of~$p_{h'(u_\beta)}(t)$
for all~$\alpha\not=0,\beta$ since~$h'(u_\beta)$ lies on each~$F_\alpha$
with~$\alpha\not=\beta$.  Thus, the 
set~$\{\lambda_1(x_\beta),\ldots,\lambda_m(x_\beta)\}$ consists 
of the~$r_\alpha$ where~$\alpha\not=0,\beta$.  
Consequently, since~\eqref{eq: pC to ph'' ratio}
now simplifies to
\begin{equation*}
\begin{split}
\prod_{\alpha=0}^{m+1}(t{-}r_\alpha) = p_D(t) 
&= p_{h'(u_\beta)}(t)\,\bigl(t^2 + h_1(u_\beta)\,t + V(u_\beta)\bigr)\\
&= \left(\prod_{\alpha\not=0,\beta}^{m+1}(t{-}r_\alpha)\right)\,
      \bigl(t^2 + h_1(u_\beta)\,t + V(u_\beta)\bigr),\\
\end{split}
\end{equation*}
it follows that~$\bigl(t^2 + h_1(u_\beta)\,t + V(u_\beta)\bigr) 
= (t - r_0)(t-r_\beta)$.  In particular, \eqref{eq: eigen w} becomes
\begin{equation*}
w'(t)(u_\beta) 
 = -2i\,\,\left(\prod_{\alpha\not=0,\beta}(t{-}r_\alpha)\right)\,
   \bigl[ H(u_\beta){-}r_0\,\I_n\bigr]
   \bigl[H(u_\beta){-}r_\beta\,I_n\bigr]
   \,\bigl[t\,\I_n{-}H(u_\beta)\bigr]^{-1}\,.
\end{equation*}

Now, any eigenvalue of~$H(u_\beta)$ is a root 
of~$p_{h(u_\beta)}(t) = p_{h''}(t) p_{h'(u_\beta)}(t)$ and so, 
by Proposition~\ref{prop: ph'' roots are pD roots}, 
must be of the form~$r_\gamma$ for 
some~$\gamma=0,\ldots,m{+}1$.  Let~$V_{\beta,\gamma}\subset\C{n}$
denote the eigenspace of~$H(u_\beta)$ belonging to the eigenvalue~$r_\gamma$.
Then the above formula implies that~$w'(t)(u_\beta)$ 
annihilates~$V_{\beta,\beta}$ and~$V_{\beta,0}$
and that, for~$v\in V_{\beta,\gamma}$ with~$\gamma\not=0,\beta$,
\begin{equation}\label{eq: eigen w ii}
w'(t)(u_\beta)\,v = -2i\,(r_\gamma{-}r_0)(r_\gamma-r_\beta)
        \left(\prod_{\alpha\not=0,\beta,\gamma}(t{-}r_\alpha)\right)\,\,v.
\end{equation}
Since the right hand side of~\eqref{eq: eigen w ii} is a polynomial in~$t$, 
it now makes sense to substitute~$t = r_\alpha$ for any~$\alpha$. 
When~$\alpha\not=0,\beta,\gamma$, this gives~$w'(r_\alpha)(u_\beta)\,v = 0$ 
for~$v\in V_{\beta,\gamma}$, while, if~$\gamma\not=0,\beta$, this gives
\begin{equation*}
w'(r_\gamma)(u_\beta)\,v 
= -2i\,(r_\gamma{-}r_0)(r_\gamma-r_\beta)
   \left(\prod_{\alpha\not=0,\beta,\gamma}(r_\gamma{-}r_\alpha)\right)\,\,v
= -2i\,p'_D(r_\gamma) v.
\end{equation*}

In other words, for~$\beta\not=\gamma$ in the 
range~$1\le\beta,\gamma\le m{+}1$,
\begin{equation}\label{eq: eigen w iii}
w'(r_\gamma)(u_\beta) = -2i\,p'_D(r_\gamma)\,E_{\beta,\gamma}
\end{equation}
where~$E_{\beta,\gamma}:\C{n}\to V_{\beta,\gamma}$ is the orthogonal 
projection onto this eigenspace.  Thus, the flow of~$W(r_\gamma)$ 
is periodic of period~$\pi/p'_D(r_\gamma)$ and so, 
by~\eqref{eq: W at r-alpha}, 
the flow of~$W_\gamma$ is periodic of period~$\pi$, as claimed.

Now, further information can be got by evaluating~$w'(r_\beta)$ 
at~$u_\beta$ itself.  Indeed, in the above formula, if~$v$ lies 
in~$V_{\beta,\gamma}$ with~$\gamma\not=0,\beta$, then putting~$t = r_\beta$
gives
\begin{equation*}
\begin{split}
w'(r_\beta)(u_\beta)(v) 
&= -2i\,(r_\gamma{-}r_0)(r_\gamma-r_\beta)
    \left(\prod_{\alpha\not=0,\beta,\gamma}(r_\beta{-}r_\alpha)\right)\,\,v\\  
&= 2i\,\frac{(r_\gamma-r_0)}{(r_\beta-r_0)}\,
     \left(\prod_{\alpha\not=\beta}(r_\beta{-}r_\alpha)\right)\,\,v 
  = 2i\,p'_D(r_\beta)\,\frac{(r_\gamma-r_0)}{(r_\beta-r_0)}\, v,\\
\end{split}
\end{equation*}
In other words, using the projection notation already introduced,
\begin{equation*}
w'(r_\beta)(u_\beta) = 2i\,p'_D(r_\beta)
    \sum_{\gamma\not=0,\beta}\frac{(r_\gamma-r_0)}{(r_\beta-r_0)}
                                \,E_{\beta,\gamma}.
\end{equation*}
Since the flow of~$W_\beta$ has period~$\pi$, 
each of the ratios~$(r_\gamma{-}r_0)/(r_\beta{-}r_0)$ must be an
integer for~$1\le \beta\not=\gamma\le m{+}1$.
Since~$r_\beta\not=r_\gamma$ when~$\beta\not=\gamma$, these ratios
cannot be~$+1$.  Thus, as the inverse of each such ratio is another such
ratio, these integers must all be~$-1$.  However, this is equivalent
to~$\frac12(r_\beta+r_\gamma)=r_0$, which is impossible, since~$r_0$ 
is greater than either~$r_\beta$ or~$r_\gamma$.  This contradiction
establishes the proposition.  
\end{proof}

\begin{corollary}\label{cor: compact BK is loc sym}
The only connected compact Bochner-K\"ahler manifolds are the 
compact quotients of the symmetric Bochner-K\"ahler 
manifolds~$M^p_c\times M^{n-p}_{-c}$.
\end{corollary}

\begin{proof}
A compact Bochner-K\"ahler manifold is necessarily complete and its
reduced momentum image is necessarily compact.  
Proposition~\ref{prop: distinct roots is not complete}, 
Corollary~\ref{cor: complete fibration},
and the fact that only SubCase 4-0 has a compact momentum cell
imply that $m>0$ is impossible.  When $m=0$, the momentum mapping is
constant and the metric is therefore locally homogeneous, so that,
by Proposition~\ref{prop: BK loc symm}, its simply-connected cover 
(which is complete) must be isometric to~$M^p_c\times M^{n-p}_{-c}$,
as claimed.  
\end{proof}

\begin{remark}[Orbifolds]
While Proposition~\ref{prop: distinct roots is not complete} rules out the 
existence of a Bochner-K\"ahler manifold
in Subcase 4-0, it does not rule out the existence of orbifolds.  In fact,
the argument in Proposition~\ref{prop: distinct roots is not complete}
 implies that, if there is a complete orbifold
with reduced characteristic polynomial~$p_D(t)$ as in the proof, then
the ratios~$(r_\gamma{-}r_0)/(r_\beta{-}r_0)$ must all be rational
for~$1\le \beta,\gamma\le m{+}1$.  A little algebra then 
leads to the formulae
\begin{equation*}
\begin{split}
p_D(t) &= (t-r_0)\,(t-r_1)\,\cdots\,(t-r_{m+1})\\
p_C(t) &= (t-r_0)^{\nu_0+1}\,(t-r_1)^{\nu_1+1}\,
                 \cdots\,(t-r_{m+1})^{\nu_{m{+}1}+1}\\
\end{split}
\end{equation*}
with
\begin{equation*}
r_\beta = r\,\sum_{\alpha=0}^{m+1} (\nu_\alpha+1)(p_\alpha-p_\beta),
\qquad 0\le\beta\le m{+}1
\end{equation*}
where~$r>0$ is real, $0=p_0<p_1<p_2<\cdots<p_{m+1}$ is
a strictly increasing sequence of integers with no common divisor, and
$\nu_0,\ldots,\nu_{m+1}$ are nonnegative integers satisfying
\begin{equation*}
n = m + \nu_0+\cdots+\nu_{m+1}.
\end{equation*}

While I have not done all of the necessary calculations, 
it appears that, for each choice of~$r$, 
$p=(p_1,\ldots,p_{m+1})$, and~$\nu = (\nu_0,\ldots,\nu_{m+1})$
satisfying the above conditions, there exists a complete orbifold
with characteristic polynomials~$p_C$ and $p_D$ as above that fits
into SubCase 4-0.  The case~$n=1$ has already been verified 
in~\S\ref{sssec: explicit dim 1},
and the cases with~$n=m$ will be verified in~\S\ref{sssec: wtd proj spaces}.

The parameter~$r$ can be normalized to~$1$ by scaling the metric. 
Thus, up to scaling, there exists a countable family of complete 
Bochner-K\"ahler orbifolds in each dimension whose momentum cells 
are compact. It appears that these orbifolds are weighted projective 
space in most cases.  In fact, it will be seen that every weighted 
projective space carries a Bochner-K\"ahler metric.
\end{remark}

By the same methods as employed in the proof of 
Proposition~\ref{prop: distinct roots is not complete}, one
can prove the following more general periodicity result.  
Details will be left to the reader.

\begin{proposition}\label{prop: periods of W-alpha}
Let~$r_\alpha$ be a simple root of~$p_D$, 
let~$w_\alpha = (h')^*(l_\alpha)$, and let~$W_\alpha$ be the vector field
in~$\euz$ defined by~$W_\alpha\lhk\Omega = -dw_\alpha$.  Then on 
a neighborhood of any zero of~$W_\alpha$, the flow of~$W_\alpha$ is
periodic, with period~$\pi$. 
\end{proposition}

\begin{remark}[Locality]
One must restrict to a neighborhood of a fixed point for the conclusion
of Proposition~\ref{prop: periods of W-alpha}.  
In the first place, without some completeness
assumptions, there is no reason to believe that the flow of~$W_\alpha$
is even defined for all time except near a fixed point.  In the second
place, even if the flow is defined for all time, by removing the zero
locus of~$W_\alpha$ and passing to a covering space, one could conceivably
arrange that~$W_\alpha$ have no closed orbits at all.
\end{remark}

\subsection{A geodesic foliation}
\label{ssec: geod fol}

By Theorem~\ref{thm: generators of center}, 
the vector fields~$Z_2,\ldots, Z_{m+1}$ 
are linearly independent (over~$\bbC$) on~$M^\circ$ 
and satisfy~$[Z_i,Z_j]=0$ for~$2\le i,j\le m{+}1$. Moreover, since these
vector fields are the real parts of holomorphic vector fields, they satisfy
\begin{equation*}
[Z_i,Z_j] = [Z_i,JZ_j] = [JZ_i,JZ_j] = 0.
\end{equation*} 
Since the $2m$ vector fields~$Z_2,\ldots, Z_{m+1},JZ_2,\ldots, JZ_{m+1}$
Lie-commute and are linearly independent on~$M^\circ$,
they are tangent to a foliation~$\cF$ of~$M^\circ$ whose leaves
are complex submanifolds of~$M$ (of complex dimension~$m$).  Moreover,
the vector fields~$JZ_2,\ldots, JZ_{m+1}$ are tangent to a foliation~$\cE$
of~$M^\circ$ whose tangent spaces are the orthogonal complements
to fibers of~$h':M^\circ\to C(p_D,\mu)^\circ$.

\subsubsection{Geometry of the leaves}
\label{sssec: geom leaves}

It turns out that the $\cF$-leaves are themselves rather interesting
objects.

\begin{proposition}\label{prop: leaf data}
The leaves of the foliation~$\cF$ are totally geodesic in~$M^\circ$ 
and the induced K\"ahler structure on each $\cF$-leaf 
is Bochner-K\"ahler of cohomogeneity~$m$.  
The characteristic and momentum polynomials of any~$\cF$-leaf are
\begin{equation*}
p^L_C(t) = p^L_D(t) = p_D(t-\lambda),
\qquad\text{and}\qquad
p^L_{h}(t) = p^L_{h'}(t) = p_{h'}(t-\lambda)
\end{equation*}
where the constant~$\lambda$ is defined so that~
$p_{h''}(t) = t^{n-m}-(m{+}2)\lambda\,t^{n-m-1} + \cdots$.
\end{proposition}

\begin{proof}
Return to the structure equations on~$P_2$
that were introduced in the proof of 
Proposition~\ref{prop: ph'' roots are pD roots}.  
Since~$(-1)^{i-1}p_D(\lambda_i)>0$ for~$1\le i\le m$ 
while~$p_D(\lambda_a) = 0$ for~$a>m$, 
it follows that~$\lambda_i-\lambda_a$ is nonvanishing on~$M^\circ$
and hence on~$P_2$.  Equation~\eqref{eq: off block vanishing} 
can thus be written as
\begin{equation*}
\phi_{a\bar\imath} = \frac{T_i}{\lambda_i-\lambda_a}\,\omega_a\,.
\end{equation*}
Let~$L^\circ\subset M^\circ$ be an $\cF$-leaf, and let~$P^L_2\subset P_2$
be the bundle~$\pi^{-1}(L^\circ)\cap P_2$, which is a~$G_\Lambda$-bundle 
over~$L^\circ$.  By the definition of
the bundle~$P_2$ and the foliation~$\cF$, the forms~$\omega_a$ vanish
when pulled back to~$P^L_2$, so, by the above equations, so do the
forms~$\phi_{a\bar\imath}$. 
Consequently, the complex $m$-manifold $L^\circ$ is totally geodesic 
in~$M^\circ$, as claimed.

Denoting pullback to~$P^L_2$ by a superscript~$L$, the formulae
\begin{equation*}
\omega^L = \begin{pmatrix}\tilde\omega\\0\end{pmatrix},
\qquad
\phi^L = \begin{pmatrix}\tilde\phi&0\\ 0&\phi''\end{pmatrix},
\qquad
H^L = \begin{pmatrix} H'&0\\0&\Lambda\end{pmatrix},
\qquad
T^L = \begin{pmatrix} T'\\0\end{pmatrix}
\end{equation*}
hold, where~$\tilde\omega$ takes values in~$\C{m}$ and~$\tilde\phi$ takes
values in~$\euu(m)$.  The notations~$H'$, $T'$, and~$\Lambda$ are as
previously established.  The K\"ahler form~$\Upsilon$ induced on~$L^\circ$
by pullback from~$\Omega$ then satisfies~$\pi^*(\Upsilon) 
= -{\frac{i}2}\,\tilde\omega^*\w\tilde\omega$.

The pullbacks of the structure equations to~$P^L_2$ then imply 
$d\tilde\omega = -\tilde\phi\w\tilde\omega$, so that~$\tilde\phi$ is
the connection matrix of the torsion-free K\"ahler connection of
the induced K\"ahler structure.  The pullbacks further imply
\begin{equation*}
\begin{split}
d\tilde\phi + \tilde\phi\w\tilde\phi 
&= H'\,\tilde\omega^*\w\tilde\omega 
 - H'\,\tilde\omega\w\tilde\omega^* 
 - \tilde\omega\w\tilde\omega^*\,H' 
   + \tilde\omega^*\,H'\,\tilde\omega\,\I_m\\
&\qquad\qquad 
   + (\tr H' + \tr\Lambda)\,\bigl(\tilde\omega^*\w\tilde\omega\,\I_m
       - \tilde\omega\w\tilde\omega^*\bigr)\\
&= \tilde H\,\tilde\omega^*\w\tilde\omega 
 - \tilde H\,\tilde\omega\w\tilde\omega^* 
 - \tilde\omega\w\tilde\omega^*\,\tilde H 
 + \tilde\omega^*\,\tilde H\,\tilde\omega\,\I_m\\
&\qquad\qquad 
      + (\tr \tilde H)\,\bigl(\tilde\omega^*\w\tilde\omega\,\I_m
                      - \tilde\omega\w\tilde\omega^*\bigr)\\
\end{split}
\end{equation*} 
where~$\tilde H = H'+\lambda\,\I_m$ and~$(m{+}2)\,\lambda=\tr\Lambda$.
Since~$p_{h''}(t) = \det\bigl(t I_{n-m}{-}\Lambda\bigr)$, this defines
$\lambda$ as in the statement of the proposition.

Thus, by definition, the induced metric on~$L^\circ$ is Bochner-K\"ahler
and has momentum polynomial
\begin{equation*}
p^L_h(t) = \det\bigl(t\,I_m - \tilde H\bigr) 
         = \det\bigl((t-\lambda)\,I_m - H'\bigr) 
         = p_{h'}(t-\lambda),
\end{equation*}
as claimed.  Moreover, the pullback of the $dH$ equation implies
\begin{equation*}
d\tilde H = -\tilde\phi\,\tilde H + \tilde H\,\tilde\phi
              + T'\,\tilde\omega^* + \tilde\omega\,(T')^*\,
\end{equation*}
so that~$\tilde T = T'$ is the $\C{m}$-valued function
defined by the structure equations for~$\Upsilon$.

Since the entries of~$\tilde T = T'$ are all nonzero
and the eigenvalues of~$\tilde H$ are distinct,
the rank of the momentum mapping for~$L^\circ$ is~$m$, 
implying that~$p^L_{h'}(t) = p^L_h(t)$.

Finally, the pullback of the identity for~$dT$ becomes
\begin{equation*}
d\tilde T = -\tilde\phi\,\tilde T 
+ \bigl((\tilde H)^2 + \tr(\tilde H)\,\tilde H
  + (V {-}\lambda\,\tr(H') {-} (m{+}1)\lambda^2 ) \bigr)\tilde\omega,
\end{equation*}
so that, setting~$\tilde V = V{-}\lambda\,\tr(H'){-}(m{+}1)\lambda^2$,
the structure function for the metric on~$L^\circ$ 
takes the form~$(\tilde H,\tilde T,\tilde V)$.

The formula for~$p^L_C(t)$ then becomes
\begin{equation*}
\begin{split}
p^L_C(t) &= \det(t\,I_m - \tilde H)
             \bigl(t^2+\tr(\tilde H)\,t + \tilde V\bigr)
            + (\tilde T)^*\Cof(tI_m-\tilde H)\,\tilde T\\
 &=  \det((t{-}\lambda)\,I_m - H')
             \bigl((t{-}\lambda)^2+\tr(H)\,(t{-}\lambda) + V\bigr)\\
 &\qquad\qquad\qquad
            + (\tilde T)^*\Cof\bigl((t{-}\lambda)I_m-H'\bigr)\,\tilde T\\
&= \frac{p_C(t{-}\lambda)}{p_{h''}(t{-}\lambda)} = p_D(t{-}\lambda)\\
\end{split}
\end{equation*}
where the last line uses the definition of~$p_C(t)$, the
identity~$p_h(t) = p_{h'}(t)p_{h''}(t)$, and the identity
\begin{equation*}\Cof(tI_n-H) 
= \begin{pmatrix} p_{h''}(t)\Cof(tI_m{-}H') & 0\\
           0 & p_{h'}(t)\Cof(tI_{n-m}{-}\Lambda)\\
  \end{pmatrix}.
\end{equation*}

These formulae establish the proposition. 
\end{proof}

\begin{corollary}\label{cor: reduction to leaf}
The momentum mapping~$h^L:L^\circ\to\bbR^m$ is equal to the
restriction of~$h':M^\circ\to\bbR^m$ to~$L^\circ$ followed 
by an invertible linear map~$\Phi_\lambda:\bbR^m\to\bbR^m$. 
The corresponding momentum cells satisfy~$C(p^L_D,\mu^L) 
= \Phi_\lambda\bigl(C(p_D,\mu)\bigr)$. 
\end{corollary}

\subsubsection{Completion and real slices}
\label{sssec: comp and real slices}

It will now be shown that the~$\cF$- and~$\cE$-leaves can
be extended through the locus~$N$ where the~$Z_k$ become dependent.

\begin{proposition}\label{prop: complete metric implies leaf complete}
If the metric on~$M$ is complete, then the closure 
of any $\cF$-leaf~$L^\circ$ is a complete, 
totally geodesic complex $m$-manifold~$L\subset M$.  
Moreover, the geodesic completion 
of any~$\cE$-leaf~$R^\circ\subset L^\circ$ 
is a totally geodesic real $m$-manifold~$R$ and
the mapping~$h':R\to C(p_D,\mu)$ is surjective.
\end{proposition}

\begin{proof}
Before beginning the proof, it will be useful to establish the
following fact.  If~$g$ is any real-analytic, complete metric on 
a manifold~$M$ and~$S\subset M$ is a connected, totally geodesic submanifold
of some dimension~$s$, then $S$ can be `completed':
There exists an $s$-manifold~$\bar S$ and a totally geodesic
immersion~$\iota:\bar S\to M$ whose image contains~$S$ and for which
the induced metric~$\bar g = \iota^*g$ on~$\bar S$ is complete.
This completion~$(\bar S,\iota)$ is unique up to diffeomorphism.

Here is a sketch of the proof:  Fix~$x\in S$ and consider, for 
every~$v\in T_xS$, the constant speed 
geodesic~$\gamma_v:\bbR\to M$ defined by~$\gamma_v(t) = \exp_x(tv)$. 
Let~$E(v)\subset T_{\gamma_v(1)}M$ be the parallel translation 
of~$E(0)=T_xS$ along $\gamma_v$ from~$t=0$ to~$t=1$, and 
let~$\bar S\subset \Gr_s(TM)$ be the set of all such~$E(v)$.
Since~$S$ is totally geodesic, when~$|v|$ is sufficiently
small~$E(v)$ is equal to~$T_{\gamma_v(1)}S$.  It follows, 
by the real-analyticity of~$g$, 
that~$\exp_{\gamma_v(1)}:E(v)\to M$ embeds~$B_\delta(0)\subset E(v)$
into~$M$ as a totally geodesic submanifold of~$M$ as long as~$\delta>0$
is less than the injectivity radius at~$\gamma_v(1)$.  From this, 
it is not hard to prove that~$\bar S$ is an embedded
submanifold of~$\Gr_s(TM)$.  Moreover, the basepoint 
projection~$\iota:\bar S \to M$ defined by~$\iota\bigl(E(v)\bigr) 
= \gamma_v(1)$ is a totally geodesic immersion.  
Since each of the geodesics~$\gamma_v$ lifts to a complete 
geodesic~$t\mapsto E(tv)$ in~$\bar S$ 
for the induced metric~$\bar g=\iota^*g$, 
all of the $\bar g$-geodesics in~$\bar S$ passing through~$E(0)$ are
complete.  Thus,~$\bar g$ is complete.  The completeness of the metric
then ensures that $\iota(\bar S)$ contains~$S$, as desired.

Now apply this result to the leaf~$L^\circ\subset M^\circ$ and 
consider~$\overline{L^\circ}\subset\Gr_{2m}(TM)$.  
Since the induced metric on $\overline{L^\circ}$ is real-analytic
and is Bochner-K\"ahler on an open set, it is Bochner-K\"ahler
everywhere.  Moreover, by Proposition~\ref{prop: leaf data}, 
it has cohomogeneity~$m$, equal to its complex dimension.
Let~$h^L:\overline{L^\circ}\to\bbR^m$ denote its momentum mapping.
By Corollary~\ref{cor: reduction to leaf} and real-analyticity,
$h^L = \Phi_\lambda\circ h'\circ\iota$, since this holds on~$L^\circ$.
Since the rank of~$\euz$ is~$m$, the proof of 
Theorem~\ref{thm: generators of center} coupled with
the remarks of~\S\ref{sssec: min sym} show that~$\euz$ accounts for all of the
infinitesimal symmetries of~$\overline{L^\circ}$, i.e., that the
full symmetry group of~$\overline{L^\circ}$ is generated by the
canonical torus action, even if one were to pass to its simply-connected
cover.  In particular, by Corollary~\ref{cor: complete fibration}, 
the fibers of~$h^L$ are connected and are the orbits of the 
canonical torus action.

Now,~$h^L$ is a submersion outside some closed complex 
submanifold~$K\subset \overline{L^\circ}$.  Since~$h'$ is a submersion
only when it is restricted to~$M^\circ$, it follows that
that $\iota\bigl(\overline{L^\circ}\setminus K\bigr)$ must 
lie in~$M^\circ$.  Since $\overline{L^\circ}\setminus K$
is connected, and since it contains the $\cF$-leaf~$L^\circ$,
it follows from analyticity that $\iota(\overline{L^\circ}\setminus K)$
must be equal to~$L^\circ$.   Thus, $\overline{L^\circ}\setminus K$
consists of the tangent planes to~$L^\circ$.   It follows that
$\iota$ is one-to-one on~$\overline{L^\circ}\setminus K$.  

If~$\iota$ were not one-to-one on~$K$, this would violate the connectedness
of the fibers of~$h^L$, so~$\iota$ is one-to-one everywhere.
In other words, $L = \iota(\overline{L^\circ})$ is 
a submanifold of~$M$, as claimed.  Obviously, $L$ is the closure
of~$L^\circ$ in~$M$.

Now, turning to the geometry of the leaves of the foliation~$\cE$, 
note that these leaves are defined by the equations~$\omega_a = 
\Im(\omega_i) = 0$ (since~$H$ and~$T$ are real on~$P_2$).  To prove
that these leaves are totally geodesic, it would suffice to show that
the imaginary part of~$\phi'$ vanishes when one restricts to such a
leaf.  Thus, write~$\tilde\omega = \xi + i\,\eta$ 
and~$\phi' = \theta + i\,\psi$, where~$\xi$, $\eta$, $\theta$, 
and~$\psi$ take values in~$\bbR^m$, $\bbR^m$, $\euso(m)$, and the
space of real symmetric matrices, respectively.  Since~$H$ and~$T$
are real-valued, the imaginary part of the equation for~$dH'$ becomes
\begin{equation*}
 0 = -H'\,\psi - \psi\,H' - T'\,{}^t\eta + \eta\,{}^tT'.
\end{equation*}
It then follows by linear algebra that on the open set~$U\subset M^\circ$ 
where $H'$ has no two eigenvalues that sum to zero, the components of~$\psi$
are linear combinations of the components of~$\eta$.  
By the structure equation for~$dH'$, the eigenvalues of
$H'$ are independent on each leaf of~$\cE$ 
since~$d\lambda_i = 2T_i,\xi$.  Thus, the open set~$U$ intersects
each~$\cE$-leaf in a dense open set.  Consequently, the components of~$\psi$ 
vanish on each~$\cE$-leaf, implying that each~$\cE$-leaf is totally
geodesic, as desired.  

Now, let~$R^\circ$ be an $\cE$-leaf and let~$\bar R\subset \Gr_m(TM)$
be its geodesic completion.  By construction, $R^\circ$ meets each
isometry orbit in~$M^\circ$ orthogonally.  
Thus, by Theorem~\ref{thm: cell metric from pD},
the map $h':R^\circ\to C(p_D,\mu)^\circ$ is an isometry when~$R^\circ$
is given the induced metric.  It remains to show that ~$h'\circ\iota:\bar R
\to C(p_D,\mu)$ is surjective.  The completeness and real-analyticity
of the induced metric~$\bar g$ on~$\bar R$ coupled with the analysis of the 
resolution of the singular cell metric in SubCase 3-1$b$ done
at the end of~\S\ref{sssec: explicit forms}, shows that $(\bar R,\bar g)$ 
must be an isometric quotient of~$(E_\rho,R_\rho)$ for some~$\rho$.  
The completeness of
this mapping and the surjectivity of this resolution imply the desired
surjectivity. 
\end{proof}

\subsubsection{The leaf metric}
\label{sssec: leaf metric}

The Bochner-K\"ahler metric induced on a leaf~$L^\circ$ can now be 
described rather explicitly in terms of the geometry of the momentum cell
associated to~$M$.

\begin{theorem}\label{thm: leaf metric}
Let~$\bigl(R_{ij}^D(u)\bigr)$ be the inverse matrix to the 
coefficient matrix~$\bigl(R^{ij}_D(u)\bigr)$ of the cell metric~$R_D$
on~$C(p_D,\mu)^\circ$.  Then, on the universal cover of~${L^\circ}$,
there exist functions~$\theta^1,\ldots, \theta^m$ for which the 
induced K\"ahler form and metric are
\begin{equation*}
\Upsilon = dh'_k\w d\theta^k
\qquad\text{and}\qquad
ds^2 = R^{jk}_D(h')\,dh'_j{\circ}dh'_k 
      + R_{jk}^D(h')\,d\theta^j{\circ}d\theta^k.
\end{equation*}
\end{theorem}

\begin{proof}
First, it will be useful to take a different
basis for~$\euz$.  Recall that the functions~$(h'_1,\ldots,h'_m)$
are constant linear combinations of the functions~$(h_1,\ldots, h_m)$
and vice versa.  By equation~\eqref{eq: Z-k is Ham},
$Z_{k+1}\lhk\Omega = -dh_k$ 
for~$1\le k\le m$, so the vector fields~$Z'_2,\ldots, Z'_{m+1}$ defined by
$Z'_{k+1}\lhk\Omega = -dh'_k$ for~$1\le k\le m$ are also a basis of~$\euz$.

Let~$L^\circ$ be an $\cF$-leaf. 
The vector fields~$\{Z'_k{-}iJZ'_k\mid2\le k\le m{+}1\}$ 
are a basis for the holomorphic vector fields on~$L^\circ$, 
so there are unique holomorphic 1-forms~$\zeta^1,\ldots,\zeta^m$ 
on~$L^\circ$ that satisfy
\begin{equation*}
\zeta^j\bigl(Z'_{k+1}-iJZ'_{k+1}\bigr)=\imath\,\delta^j_k
\end{equation*}
for~$1\le j,k\le m$.  (The introduction of the factor of~$\imath$ 
simplifies formulae to appear below.)
Because the vector fields~$Z'_k{-}iJZ'_k$ are Lie-commuting,
the 1-forms~$\zeta^j$ are closed.

Write~$\zeta^j=\xi^j+\imath\,\eta^j$ where~$\xi^j$ and~$\eta^j$ are
real $1$-forms.  Then the defining equation above is equivalent to
\begin{equation*}
\eta^j(Z'_{k+1})=-\xi^j(JZ'_{k+1}) = {\ts\frac12}\,\delta^j_k
\quad\text{and}\quad
\xi^j(Z'_{k+1})=\eta^j(JZ'_{k+1}) = 0.
\end{equation*}

Since the~$\zeta^j$ are a basis for the holomorphic 1-forms on~$L^\circ$,
the metric on~$L^\circ$ can be written in the form
\begin{equation*}
ds^2 = g_{jk}\,\zeta^j{\circ}\overline{\zeta^k},
\end{equation*}
where~$g_{jk} = \overline{g_{kj}}$ 
and where the pullback of~$\Omega$ to~$L^\circ$ is
\begin{equation*}
\Upsilon = {\ts\frac{i}2}\,g_{jk}\,\zeta^j\w \overline{\zeta^k}.
\end{equation*}
Now,
\begin{equation*}
dh'_k = -Z'_{k+1}\lhk\Upsilon 
= {\ts\frac12}\bigl(g_{kj}\,\overline{\zeta^j}
                    + g_{jk}\,\zeta^j\bigr),
\end{equation*}
or, equivalently,
\begin{equation*}
dh'_k = {\ts\frac12}(g_{kj}+g_{jk})\,\xi^j
      + {\ts\frac{i}2}(g_{kj}-g_{jk})\,\eta^j.
\end{equation*}
Since~$Z'_{j+1}$ for~$1\le j\le m$
is tangent to the fibers of~$h'$, the coefficient of~$\eta^j$ in the
above equation must vanish, i.e., $g_{kj} = g_{jk}=\overline{g_{kj}}$. Thus,
\begin{equation*}
dh'_k = g_{kj}\, \xi^j.
\end{equation*}

Define~$g^{ij}=g^{ji}$ so that~$g^{ij}g_{jk} = \delta^i_k$. 
Note, in particular, that~$\xi^j = g^{jk}\,dh'_k$.
The metric on~$L^\circ$ can now be written in the form
\begin{equation*}
\begin{split}
ds^2 &= g_{jk}\,\,\zeta^j{\circ}\overline{\zeta^k}
      = g_{jk}\,(\xi^j+i\,\eta^j){\circ}(\xi^k-i\,\eta^k)\\
     &= g_{jk}\bigl(\xi^j{\circ}\xi^k + \eta^j{\circ}\eta^k\bigr)\\
     &= g^{jk}\,dh'_j{\circ}dh'_k + g_{jk}\,\eta^j{\circ}\eta^k.\\
\end{split}
\end{equation*}
Since~$L^\circ$ is totally geodesic in~$M^\circ$, it follows from 
Theorem~\ref{thm: cell metric from pD} that
\begin{equation*}
g^{ij}\,dh'_i\,dh'_j = (h')^*(R_D) = R^{ij}_D(h')\,dh'_i\,dh'_j\,.
\end{equation*}
Thus,~$g^{ij} = R^{ij}_D(h')$ and so
\begin{equation*}
\xi^j = g^{jk}\,dh'_k = R^{jk}_D(h')\,dh'_k 
= (h')^*\bigl(R^{jk}_D(u)\,du_k \bigr) 
\end{equation*}

Now lift everything to the universal cover of~$L^\circ$.
Since the~$\eta^k$ are closed,  
there exist functions~$\theta^1,\ldots,\theta^m$ 
on this universal cover so that~$\eta^k = d\theta^k$.
Then
\begin{equation*}
ds^2 = R^{jk}_D(h')\,dh'_j{\circ}dh'_k 
    + R_{jk}^D(h')\,d\theta^j{\circ}d\theta^k,
\end{equation*}
where~$(R_{jk}^D)$ is the inverse matrix to~$(R^{jk}_D)$ and
\begin{equation*}
\Upsilon = R_{jk}^D(h')\, \xi^j\w d\theta^k = dh'_k\w d\theta^k.
\end{equation*}
These are the desired formulae.  
\end{proof}

\begin{remark}[K\"ahler metrics of Hessian type]
The reader may find the metric in the above form to be very familiar.
In fact, K\"ahler metrics of this form are well-known in the literature
as being of \emph{Hessian type}.  Their general form is as follows:
Let~$D\subset\bbR^m$ be an open domain in~$\bbR^m$, assumed to be
simply-connected for simplicity.  Let~$x_1,\ldots,x_m$ be any linear
coordinates on~$\bbR^m$ and suppose that~$g$ is a Riemannian metric 
on~$D$, written in the form
\begin{equation*}
g = g^{jk}(x)\,dx_j\circ dx_k\,.
\end{equation*}
Using the flat affine structure on~$\bbR^m$ restricted to~$D$, 
one gets a canonical metric on~$T^*D$ as follows:  Let~$y^1,\ldots,y^m$
be the coordinates that are linear on the fibers of~$T^*D\to D$ and
dual to the coordinates~$x_1,\ldots,x_m$ in the sense that the
tautological 1-form on~$T^*D$ is $y^j\,dx_j$.  Let~$g_{jk}=g_{kj}$ be
the functions on~$D$ so that~$g^{jk}g_{kl} = \delta^j_l$ and
define the metric
\begin{equation*}
\hat g = g^{jk}(x)\,dx_j\circ dx_k + g_{jk}(x)\,dy^j\circ dy^k.
\end{equation*}
This metric on~$T^*D$ does not depend on the choice of coordinates~$x_i$,
but only on~$g$ and the flat affine structure that~$D$ inherits from~$\bbR^m$.
Moreover, this metric is compatible with the symplectic
form on~$T^*D$ given by~$\Upsilon = -d\bigl(y^j\,dx_j\bigr) = dx_j\w dy^j$.

Thus, the metric~$\hat g$ and 2-form~$\Upsilon$ define an almost 
complex structure on~$T^*D$ for which the 1-forms
\begin{equation*}
\zeta^j = g^{jk}(x)\,dx_k + \imath\,dy^j
\end{equation*}
give a basis for the $(1,0)$-forms.  This
almost complex structure will be integrable if and only if 
the forms~$\zeta^j$ are closed.  In other words, the pair~$(\hat g,\Upsilon)$
defines a K\"ahler metric on~$T^*D$ if and only 
if~$d\bigl( g^{jk}(x)\,dx_k\bigr)=0$ for all~$1\le j\le m$.  Since~$D$ is
simply-connected, this closure condition holds if and only if there
exists a convex `potential' function~$G:D\to\bbR$ for which
\begin{equation*}
g^{jk} = \frac{\partial^2 G}{\partial x_j\,\partial x_k},
\end{equation*}
i.e., if and only if the metric~$g$ is of Hessian type.
For this reason, metrics of the form~$\hat g$ as above are often called
K\"ahler metrics of Hessian type.  Note that, for such a metric, 
translation in the $y$-variables defines a Hamiltonian torus action
that is holomorphic and whose momentum mapping is the 
basepoint projection~$T^*D\to D$.

For further investigation of these metrics, the reader might 
consult~\cite{Gu} and~\cite{Ab}.

In the case of Bochner-K\"ahler metrics, the formula for the potential
function~$G:C(p_D,\mu)^\circ\to\bbR$ has been indicated 
in~\S\ref{sssec: explicit forms}.  
The main problem with this representation is that it only
describes the leaf metric away from the singular locus~$N$.  More
work must now be done to analyze the metric near this locus.
\end{remark}

\subsubsection{A partial completion}
\label{sssec: part complt}

By Theorem~\ref{thm: leaf metric}, 
the $2$-form and Riemannian metric on~$C(p_D,\mu)^\circ\times\bbR^m$ 
defined by
\begin{equation*}
\Upsilon = du_k\w d\theta^k
\qquad\text{and}\qquad
ds^2 = R^{jk}_D(u)\,du_j{\circ}du_k 
    + R_{jk}^D(u)\,d\theta^j{\circ}d\theta^k
\end{equation*}
define a Bochner-K\"ahler structure on~$C(p_D,\mu)^\circ\times\bbR^m$.
The simply-connected cover
of any $\cF$-leaf~$L^\circ$ has an immersion into~ 
$C(p_D,\mu)^\circ\times\bbR^m$, canonical up to a
translation in the $\theta$-coordinates, that pulls this Bochner-K\"ahler
structure back to the induced one on~$L^\circ$.  In this sense, this
Bochner-K\"ahler structure on~$C(p_D,\mu)^\circ\times\bbR^m$ is
universal for Bochner-K\"ahler metrics associated to this reduced
momentum cell.  Under this immersion, which is a local
diffeomorphism, the vector field~$Z'_{k+1}$ is carried 
into~$\partial/\partial\theta^k$.  

Suppose now that~$r_\alpha$ is a simple root of~$p_D(t)$ such that
$l_\alpha=0$ defines a face of~$C(p_D,\mu)$. Then the
vector field~$W_\alpha$ is defined 
(Proposition~\ref{prop: periods of W-alpha}) and has the
expansion
\begin{equation*}
W_\alpha = \frac{1}{p'_D(r_\alpha)}
  \bigl({r_\alpha}^{m-1}Z'_2 - {r_\alpha}^{m-2}Z'_3 
         + \cdots + (-1)^{m-1}Z'_{m+1}\bigr).
\end{equation*}
It follows that, under the canonical immersion, $W_\alpha$ is carried
over to the vector field
\begin{equation*}
\Theta_\alpha = \frac{1}{p'_D(r_\alpha)}
  \left(  {r_\alpha}^{m-1}\frac{\partial\hfill}{\partial\theta^1} 
        - {r_\alpha}^{m-2}\frac{\partial\hfill}{\partial\theta^2} 
         + \cdots 
        + (-1)^{m-1}\frac{\partial\hfill}{\partial\theta^m} \right).
\end{equation*}

Suppose now that $M$ contains points that satisfy~$w_\alpha=0$, i.e.,
the image of the reduced momentum mapping contains points that lie on
the face~$l_\alpha=0$.  
Then by Proposition~\ref{prop: periods of W-alpha}, near such points 
the flow of the vector field~$W_\alpha$ is periodic with period~$\pi$.
This suggests considering the vector
\begin{equation*}
\tau_\alpha = \frac{\pi}{p'_D(r_\alpha)}
  \bigl( (-r_\alpha)^{m-1},  (-r_\alpha)^{m-2}, \cdots , 1 \bigr)\in\bbR^m.
\end{equation*}

The above Bochner-K\"ahler structure is well-defined 
on~$C(p_D,\mu)^\circ\times\bigl(\bbR^m/(\bbZ\,\tau_\alpha)\bigr)$.
It is not hard to see that there exists a simply-connected complex 
$m$-manifold~$X_\alpha$ endowed with a Bochner-K\"ahler structure
and a totally geodesic hypersurface~$Y_\alpha\subset X_\alpha$ so that,
first,  the Bochner-K\"ahler structure on~$X_\alpha\setminus Y_\alpha$ 
is isomorphic to the above Bochner-K\"ahler structure 
on~$C(p_D,\mu)^\circ\times\bigl(\bbR^m/(\bbZ\,\tau_\alpha)\bigr)$
and, second, the image of the momentum mapping on~$X_\alpha$ is
equal to~$C(p_D,\mu)^\circ$ union the interior of the face~$l_\alpha=0$.

The argument for this `face-wise' extension is based on 
Theorem~\ref{thm: analytically connected classes}, which
shows that, for any point~$v\in C(p_D,\mu)$,
there must exist \emph{some} Bochner-K\"ahler metric
in the given analytically connected equivalence class whose reduced
momentum mapping assumes the value~$v$.  Taking~$v$ to lie in the
interior of the face~$l_\alpha=0$ and applying uniqueness, one sees
that the extension must exist locally.  A simple patching argument
then allows one to produce the extension~$X_\alpha$.  Details are
left to the reader, but see the next section, where the extension 
is computed explicitly in a couple of cases of interest.

\begin{remark}[Guillemin's completion]
The reader should also compare Guillemin's description~\cite{Gu}
of a K\"ahler metric constructed from the data of a polytope, since
the issue of completion across the facets is much the same.  However,
one big difference in the present case is that the polytopes involved
here are not necessarily closed.  Another is that they do not generally
satisfy the rationality requirements for the global existence theorems
that Guillemin is able to cite.  Instead, in the present case
the singular loci corresponding
to the faces are `filled in' with `patches' whose existence stem from
Theorem~\ref{thm: analytically connected classes}.  
\end{remark}

This construction generalizes:  
If~$A = \{\alpha_1,\ldots,\alpha_k\}$ is a set of~$k\le m$ distinct 
simple roots of~$p_D(t)$ for which each hyperplane~$l_{\alpha_j}=0$ defines
a face of~$C(p_D,\mu)$, then the vectors~$v_{\alpha_j}$ as defined above
generate a discrete subgroup~$\Lambda_A\subset\bbR^m$ and the
Bochner-K\"ahler structure 
descends to~$C(p_D,\mu)^\circ\times\bigl(\bbR^m/\Lambda_A\bigr)$.  
Moreover, there is a simply-connected complex $m$-manifold~$X_A$
endowed with a Bochner-K\"ahler structure and a (reducible)
hypersurface~$Y_A\subset X_A$ (whose irreducible components are 
totally geodesic) so that, first,~$X_A\setminus Y_A$ is isomorphic as 
a Bochner-K\"ahler manifold to 
$C(p_D,\mu)^\circ\times\bigl(\bbR^m/\Lambda_A\bigr)$, and, second,
the image of the momentum mapping on~$X_A$ is
equal to~$C(p_D,\mu)^\circ$ union the faces~$l_{\alpha_j}=0$
($1\le j\le k$) and minus any faces omitted from this list.

In cases where $C(p_D,\mu)$ has at most~$m$ simple faces (which includes
all the cases except~SubCase 4-$i$ for~$i<m$), one can take~$A$ to be
the set of all the~$\alpha$ for which~$l_\alpha=0$ is a simple face of~
$C(p_D,\mu)$, and the result will be~$X_A$, whose momentum image is
the entire momentum cell.  This is, in some sense, `the maximally 
complete' Bochner-K\"ahler structure of dimension~$m$ with the given 
momentum cell as momentum image.  By 
Propositions~\ref{prop: metric complete implies cell bounded}~
and~\ref{prop: leaf data}, however, 
$X_A$ cannot be metrically complete unless the cell is bounded and has 
exactly~$m$ simple faces.  As has been seen, when~$m>0$ this can only 
happen in SubCase 3-1$b$.  As will be seen in the next section, 
the metric on~$X_A$ does turn out to be complete in this SubCase.

In SubCase~4-$i$ for~$i<m$, the hyperplane~$l_\alpha=0$ is a simple
face of~$C(p_D,\mu)$ for all~$\alpha\not=i$.  The relations
\begin{equation*}
\sum_\alpha W_\alpha = \sum_\alpha r_\alpha\,W_\alpha =0,
\end{equation*}
then imply the relation
\begin{equation*}
\sum_{\alpha\not=i} (r_i{-}r_\alpha)\, W_\alpha = 0
\end{equation*}
among the vector fields that would have to be periodic if one were
going to be able to complete the metric across all $m{+}1$ of the faces 
simultaneously.  This, in turn implies that
\begin{equation*}
\sum_{\alpha\not=i} (r_i{-}r_\alpha)\, \tau_\alpha = 0
\end{equation*}
and this is the unique linear relation among 
the~$\{\tau_\alpha\,|\,\alpha\not=i\}$.  These
$m{+}1$ vectors generate a discrete lattice~$\Lambda_i$ in~$\bbR^m$ 
if and only if the ratios~$(r_i{-}r_\alpha)/(r_i{-}r_\beta)$ are rational 
for all~$\alpha,\beta\not=i$.  

However, this rationality condition is not sufficient
for the Bochner-K\"ahler structure on~
$C(p_D,\mu)^\circ\times\bigl(\bbR^m/\Lambda_i\bigr)$ to complete to 
a smooth manifold.  In fact, the necessary condition for this is
that these ratios all be integers, which an elementary argument
shows not to be possible.  Instead, the rationality is sufficient to
ensure that the metric extends to a smooth orbifold whose 
momentum mapping is onto~$C(p_D,\mu)$.  
By Propositions~\ref{prop: metric complete implies cell bounded}~
and~\ref{prop: leaf data},
this metric is complete only in SubCase~4-0, and this returns to the 
orbifold discussion at the end of~\S\ref{sssec: ness conds for complete}.

\subsubsection{Complete examples}
\label{sssec: complete examples}

I am now going to describe a formula that defines an $n$-parameter family 
of complete Bochner-K\"ahler metrics on~$\C{n}$.  I will then state
a theorem about these metrics and follow this with a discussion that
motivates the derivation of this (rather unlikely looking) formula.

First, fix~$\rho = (\rho_1,\ldots,\rho_n)\in\bbR^n$ where each~$\rho_i$ 
is a non-negative real number. I claim that there is a real-analytic 
function~$s:\C{n}\to[0,\infty)$ so that
\begin{equation*}
s(z) - \sum_{i=1}^n e^{-\rho_i s(z)}\,|z_i|^2 = 0.
\end{equation*}
for all~$z\in\C{n}$.  (This claim will be justified below.)  Of course,
the function~$s$ depends on~$\rho$, but I will not notate this.  
When~$\rho=0$, one has~$s(z) = |z|^2$, but otherwise this is not an
elementary function.  By construction, the function~$s$ is invariant 
under the standard~$n$-torus action on~$\C{n}$ defined by
\begin{equation*}
(e^{i\theta_1},\ldots,e^{i\theta_n})\cdot (z_1,\ldots,z_n)
= \bigl(e^{i\theta_1}z_1,\ldots,e^{i\theta_n}z_n\bigr).
\end{equation*}

Now set
\begin{equation*}
S(z) = 1 + \sum_{i=1}^n \rho_i\,e^{-\rho_i s(z)}\,|z_i|^2 \ge 1.
\end{equation*}
Define an Hermitian symmetric positive definite matrix~$G(z)
=\bigl(G^{ij}(z)\bigr)$ by
\begin{equation*}
G^{ij}(z) 
=  S(z) \left(\delta^{ij}e^{\rho_is(z) }
  + \bigl(\rho_i+\rho_j+\rho_i\rho_j\,s(z)\bigr)\overline{z_i}\,z_j\right).
\end{equation*}
Write~$G(z)^{-1} = \bigl(G_{ij}(z)\bigr)>0$ and define the Hermitian metric
\begin{equation*}
g_\rho = G_{ij}(z)\,dz_i\circ d\overline{z_j}.
\end{equation*}
This metric is evidently 
invariant under the standard $n$-torus action defined above.  Of course,
$g_0$ is the standard flat metric on~$\C{n}$.

\begin{theorem}\label{thm: complete metrics}
For every choice of~$\rho_i\ge0$, the metric~$g_\rho$ is Bochner-K\"ahler
and complete on~$\C{n}$.  Conversely, every simply-connected, complete
Bochner-K\"ahler manifold in dimension~$n$ is either
homogeneous or is isometric to~$(\C{n},g_\rho)$
for some~$\rho$ with~$\rho_i\ge0$.   When the $\rho_i$ are distinct and 
positive, the only symmetries of this metric belong to the standard 
$n$-torus action on~$\C{n}$. 
\end{theorem}

\begin{proof}
The structure of the proof will be as follows:  I will first
assume that I have a complete Bochner-K\"ahler metric that is not
locally homogeneous and consider the induced metric on a 
completed~$\cF$-leaf.  Knowing by earlier discussions that the only 
possibility for this is in SubCase 3-1$b$, I will use knowledge of
the form of~$p_D$ and the momentum cell to choose a particularly
good basis for~$\euz$, one for which each of the vector fields of the
basis has a periodic flow of period~$\pi$.  I will then attempt to find
global holomorphic coordinates on the leaf that will carry these
vector fields into the vector fields that generate the standard $m$-torus
action defined above.  Using these calculations as a guide and then
comparing with the discussion at the end of~\S\ref{sssec: explicit forms} 
of the `resolution'
of boundary singularities of the cell metric in SubCase~3-1$b$, I will
finally arrive at a candidate for the metric in these good coordinates
and finish by showing how completeness and real-analyticity give the
conclusions of the theorem.

Thus, suppose that~$M^n$ is simply connected and  has a complete 
Bochner-K\"ahler metric of cohomogeneity~$m>0$.  
As has already been remarked, the momentum cell must fall 
into SubCase 3-1$b$, so that
\begin{equation*}
p_D(t) = (t-r_1)^2(t-r_2)\cdots(t-r_{m+1}), 
 \qquad\qquad (r_1>\cdots>r_{m+1}).
\end{equation*}
For notational simplicity, I will use the index 
range~$2\le\alpha,\beta,\gamma\le m{+}1$ and the 
abbreviation~$\rho_\alpha$ ($>0$) for~$r_1{-}r_\alpha$ in what follows.
 Use~\eqref{eq: define l-alpha} as the definition of the linear 
function~$l_\alpha:\bbR^m\to\bbR$ for~$\alpha\ge 2$, as before,
and set
\begin{equation*}
a = 1 - \sum_\alpha \rho_\alpha\,l_\alpha
\quad\qquad\text{and}\qquad\quad
t = l_2 + \cdots + l_{m+1}\,,
\end{equation*}
as was done in~\S\ref{sssec: explicit forms} during the analysis of 
the metric~$R_D$
on this momentum cell~$C(p_D,\mu)\subset\bbR^m$, which is defined
by the inequalities~$l_\alpha\ge0$ and~$a>0$.  Note that this
momentum cell contains only one vertex, namely the point~$k_1$ 
where all of the $l_\alpha$ vanish.  Also as before, let~$F_\alpha
\subset C(p_D,\mu)$ be the face defined by~$l_\alpha=0$ 
for~$2\le\alpha\le m{+}1$.  

As was done in the proof of 
Proposition~\ref{prop: distinct roots is not complete}, 
set~$w_\alpha = (h')^*(l_\alpha)$ and consider the~$m$ vector 
fields~$W_\alpha\in\euz$ defined by~$W_\alpha\lhk\Omega = -dw_\alpha$. 
These vector fields are a basis of~$\euz$ and are linearly independent 
on~$M^\circ$.  

Fix~$q\in M^\circ$, let~$L\subset M$ be the leaf of the foliation~$\cF$
passing through~$q$, and let~$E\subset L$ be the leaf of the 
foliation~$\cE$ passing through~$q$.  On~$E^\circ$, the 
map~$h':E\to C(p_D,\mu)$ is a local isometry when~$C(p_D,\mu)$ is
given the metric~$R_D$.  

Recall the discussion and notation at the 
end of~\S\ref{sssec: explicit forms} about the metric~$R_\rho$ 
on the ellipsoidal 
domain~$E_\rho\subset\bbR^m$.  The map~$s:E_\rho\to C(p_D,\mu)$
is surjective and is isometric and submersive away from the 
hyperplanes~$p_\alpha=0$.  Letting~$p\in E_\rho$ be the point with
coordinates~$p_\alpha = \sqrt{w_\alpha(q)}>0$, it follows that 
there is a real-analytic map~$\psi$ from a neighborhood of~$p\in E_\rho$ 
to a neighborhood of~$q\in E$ satisfying~$\psi(p)=q$ and~$h'{\circ}\psi
= s$.  This map is an isometry when~$E_\rho$ is endowed with the 
metric~$R_\rho$.  Since~$E_\rho$ is simply-connected and the 
metric~$R_\rho$ is both real-analytic and complete, it follows 
that~$\psi$ can be extended uniquely as a global 
isometry~$\psi:E_\rho\to E$ and that it satisfies~$h'{\circ}\psi = s$.  
Since the rank of the differential of~$s:E_\rho\to C(p_D,\mu)$
at any~$p = (p_\alpha)$ is equal to the number of nonzero
entries~$p_\alpha$, it follows that the rank of the differential
of~$h'$ at any~$x\in E$ is equal to $m$ minus the number of
faces~$F_\alpha$ on which~$h'(x)$ lies.  Since the rank of the 
differential of~$h'$ at~$x$ is equal to the dimension of the span
of~$\{W_\alpha(x)|2\le\alpha\le m{+}1\}$, it follows from this
discussion that for any~$x$, the nonzero elements in the 
list~$\bigl(W_2(x),\ldots,W_{m+1}(x)\bigr)$ are linearly independent.
This observation will be useful below.

Each~$W_\alpha$ vanishes 
on~$N_\alpha = (h')^{-1}(F_\alpha)$, which, since the flow of
$W_\alpha$ is isometric, is a totally geodesic submanifold 
of~$M$ and, moreover, is a complex submanifold of~$M$ as well (since
$W_\alpha$ is the real part of a holomorphic vector field on~$M$).
It also follows from the discussion in the 
previous paragraph that~$W_\alpha$ is nonzero off of~$N_\alpha$.  
Let $L_\alpha = N_\alpha\cap L$.  Then $L_\alpha$ is a totally
geodesic complex hypersurface in~$L$.  

One of the goals of this argument is to show that there are
holomorphic coordinates~$z_2,\ldots z_{m+1}$ on~$L$ for which
\begin{equation*}
W_\alpha-iJW_\alpha 
= 4\,i\,z_\alpha\,\frac{\partial\hfill}{\partial z_\alpha}\,,
\qquad\qquad 2\le \alpha\le m{+}1
\end{equation*}
and to find the expression for the induced K\"ahler metric on~$L$ 
in these coordinates.  (The choice of the coefficient~$4$ is dictated
by the fact that the flow of the vector field~$W_\alpha$ has 
period~$\pi$.  The proof of this periodicity follows the same lines as 
the corresponding proof in the SubCase 4-0 situation analyzed in 
Proposition~\ref{prop: distinct roots is not complete}.  
Since it only differs in details from that proof, 
the argument will be left to the reader.)

Accordingly, let~$\zeta^2,\ldots,\zeta^{m+1}$ be the holomorphic
1-forms on~$L^\circ$ that satisfy
\begin{equation*}
\zeta^\alpha(W_\beta-iJW_\beta) = 4i\,\delta^\alpha_\beta
\end{equation*}
These forms extend meromorphically to~$L$, with simple poles along
the hypersurfaces~$L_\alpha$.  Since the vector fields~$W_\alpha$
Lie-commute, it follows that each~$\zeta^\alpha$ is closed. 
Note that, if the coordinates~$z_\alpha$ are to exist as claimed,
it will have to be true that~$\zeta^\alpha = dz_\alpha/z_\alpha$.

Writing~$\zeta^\alpha = \xi^\alpha + i\,\eta^\alpha$, the above
equations are equivalent to
\begin{equation*}
  \xi^\alpha(W_\beta) = \eta^\alpha(JW_\beta) = 0,
\qquad\qquad
-\xi^\alpha(JW_\beta) = \eta^\alpha(W_\beta) = 2\,\delta^\alpha_\beta,
\end{equation*}
Again, if the coordinates~$z_\alpha$ exist as claimed, it will follow
that~$2\xi^\alpha = d\bigl(\log|z_\alpha|^2\bigr)$.  

Since the~$\zeta^\alpha$ are a basis for the holomorphic 1-forms 
on~$L^\circ$, the metric on~$L^\circ$ can be written in the form
\begin{equation*}
ds^2 = g_{\alpha\beta}\,\zeta^\alpha{\circ}\overline{\zeta^\beta},
\end{equation*}
where~$g_{\alpha\beta} = \overline{g_{\beta\alpha}}$ 
and where the pullback of~$\Omega$ to~$L^\circ$ is
\begin{equation*}
\Upsilon 
= {\ts\frac{i}2}\,g_{\alpha\beta}\,\zeta^\alpha\w \overline{\zeta^\beta}.
\end{equation*}
Now, the identity~$\zeta^\alpha(W_\beta) = 2i\,\delta^\alpha_\beta$ implies
\begin{equation*}
dw_\alpha = -W_\alpha\lhk\Upsilon 
= g_{\alpha\beta}\,\overline{\zeta^\beta}
                    + g_{\beta\alpha}\,\zeta^\beta,
\end{equation*}
or, equivalently,
\begin{equation*}
dw_\alpha = (g_{\alpha\beta}+g_{\beta\alpha})\,\xi^\beta
      - i\,(g_{\alpha\beta}-g_{\beta\alpha})\,\eta^\beta.
\end{equation*}
Since~$W_\alpha$ is tangent to the fibers of~$h'$, and since~$w_\alpha = 
(h')^*(l_\alpha)$ is constant on those fibers, it follows that
the coefficient of~$\eta^\beta$ in the
above equation must vanish, i.e., $g_{\alpha\beta} 
= g_{\beta\alpha}=\overline{g_{\alpha\beta}}$. Thus,
\begin{equation*}
dw_\alpha = 2g_{\alpha\beta}\, \xi^\beta.
\end{equation*}

Define~$g^{\alpha\beta}=g^{\beta\alpha}$ so 
that~$g^{\alpha\beta}g_{\beta\gamma} = \delta^\alpha_\gamma$. 
Note, in particular, that~$\xi^\alpha = \frac12 g^{\alpha\beta}\,dw_\beta$.
The metric on~$L^\circ$ can now be written in the form
\begin{equation*}
\begin{split}
ds^2 &= g_{\alpha\beta}\,\,\zeta^\alpha{\circ}\overline{\zeta^\beta}
      = g_{\alpha\beta}\,
       (\xi^\alpha+i\,\eta^\alpha){\circ}(\xi^\beta-i\,\eta^\beta)\\
     &= g_{\alpha\beta}\bigl(\xi^\alpha{\circ}\xi^\beta 
                        + \eta^\alpha{\circ}\eta^\beta\bigr)\\
     &= {\ts\frac14}g^{\alpha\beta}\,dw_\alpha{\circ}dw_\beta 
          + g_{\alpha\beta}\,\eta^\alpha{\circ}\eta^\beta.\\
\end{split}
\end{equation*}
Since~$L^\circ$ is totally geodesic in~$M^\circ$, it follows from 
Theorem~\ref{thm: cell metric from pD} that
\begin{equation*}
{\ts\frac14}g^{\alpha\beta}\,dw_\alpha\,dw_\beta 
= (h')^*(R_D) 
=(h')^*\bigl(T^{\alpha\beta}\bigr)\,dw_\alpha\,dw_\beta\,,
\end{equation*}
where~$T^{\alpha\beta} = T^{\beta\alpha}$ for~$2\le\alpha,\beta\le m{+}1$ 
is defined on~$C(p_D,\mu)^\circ$ so that the formula
\begin{equation*}
\sum_{\alpha,\beta=2}^{m+1} T^{\alpha\beta}\,dl_\alpha\,dl_\beta 
= R_D = \frac{t\,{da}^2}{4a^2} - \frac{da\,dt}{2a} 
       + \sum_{\alpha=2}^{m+1}\frac{{dl_\alpha}^2}{4l_\alpha}
\end{equation*}
holds.  Thus,~$g^{\alpha\beta}=4(h')^*\bigl(T^{\alpha\beta}\bigr)$.

Using the definitions of~$a$ and~$t$, it follows from the formula 
for~$R_D$ that
\begin{equation*}
4T^{\alpha\beta}
=   \frac{\delta_{\alpha\beta}}{l_\alpha} 
  + \frac{(\rho_\alpha{+}\rho_\beta)}{a}
  + \frac{\rho_\alpha\rho_\beta\,t}{a^2}\,,
\end{equation*}
so that, in particular,
\begin{equation*}
\begin{split}
4T^{\alpha\beta}\,dl_\beta  
&= \frac{dl_\alpha}{l_\alpha} - \frac{da}{a}
   + \rho_\alpha\left(\frac{dt}{a}-\frac{t\,da}{a^2}\right)\\
&= d\left(\log{\frac{l_\alpha}a} +\rho_\alpha\,\frac{t}{a}\right).
\end{split}
\end{equation*}
Meanwhile, if the coordinates~$z_\alpha$ exist as claimed, this will
imply that
\begin{equation*}
\frac{d|z_\alpha|^2}{|z_\alpha|^2} = 2\xi^\alpha 
= g^{\alpha\beta}\,dw_\beta 
= (h')^*\bigl(4T^{\alpha\beta}\,dl_\beta \bigr)
= d\left((h')^*\left(\log{\frac{l_\alpha}a} 
   +\rho_\alpha\frac{t}{a}\right)\right),
\end{equation*}
i.e., there will exist constants~$c_\alpha>0$ so that
\begin{equation*}
|z_\alpha|^2 
= c_\alpha (h')^*\left(\frac{l_\alpha}a\,e^{\rho_\alpha\,t/a}\right).
\end{equation*}
Since~$z_\alpha$ would only be determined up to a multiplicative
constant anyway by the above normalizations, one might as well
take~$c_\alpha=1$, which will normalize the~$z_\alpha$ up to a
phase.

These calculations suggest the following construction of a
candidate for the leaf metric:  Consider the system of equations
\begin{equation*}
p_\alpha = \frac{y_\alpha}b\,e^{\rho_\alpha\,s/b},
\qquad\text{where}\qquad
b = 1-\sum_{\beta=2}^{m+1} \rho_\beta\,y_\beta
\quad\text{and}\quad
s = \sum_{\beta=2}^{m+1} y_\beta,
\end{equation*}
relating the $m$ variables~$y_2,\ldots,y_{m+1}$ to the 
variables~$p_2,\ldots,p_{m+1}$.  These formulae define
a real-analytic mapping~$\mathbf{p}$ from the open 
halfspace~$H_y\subset\bbR^m$
defined by~$b>0$ in~$y$-space into $p$-space.  

I claim that the mapping~$\mathbf{p}$ is a diffeomorphism from~$H_y$ onto its 
image~$D_p\subset\bbR^m$ and that this open image contains the primary
orthant~$O_p\subset\bbR^m$, i.e., the closed domain defined by~$p_\alpha\ge0$.  
Consequently, $\mathbf{p}$ has an inverse~$\mathbf{y}:D_p\to H_y$, i.e., 
the above equations can be solved real-analytically in the form
\begin{equation*}
y_\alpha = {\mathbf{y}}_\alpha(p_2,\ldots, p_{m+1}).
\end{equation*}
Moreover, this inverse~$\mathbf{y}$ takes~$O_p$ diffeomorphically 
onto the partially open simplex~$\Sigma\subset H_y$ 
defined by the inequalities~$y_\alpha\ge0$
and~$\sum_\beta \rho_\beta\,y_\beta <1$.

To prove this claim, it is helpful to introduce the intermediate quantities
\begin{equation*}
u_\alpha = \frac{y_\alpha}{(1-\sum_\beta \rho_\beta\,y_\beta)}\,.
\end{equation*}
These equations can be inverted in the form
\begin{equation*}
y_\alpha = \frac{u_\alpha}{(1+\sum_\beta \rho_\beta\,u_\beta)}\,,
\end{equation*}
thus showing that they define a diffeomorphism from~$H_y$ to the
half-space~$H_u\subset\bbR^m$ defined 
by~$1+\sum_\beta \rho_\beta\,u_\beta>0$.  Then the claim above amounts
to showing that the mapping defined by
\begin{equation*}
p_\alpha = u_\alpha e^{\rho_\alpha(u_2+\cdots+u_{m+1})}
\end{equation*}
is invertible on the domain~$H_u$ and that its image has the desired
properties.  

Consider the function~$f$ on~$\bbR\times\bbR^m$ defined by
\begin{equation*}
f(r,p) = r - \sum_\alpha e^{-\rho_\alpha r}p_\alpha\,.
\end{equation*}
Now, $\partial f/\partial r 
= 1 + \sum_\alpha \rho_\alpha e^{-\rho_\alpha r}p_\alpha$ is
positive on~$\bbR\times O_p$, so that~$r\mapsto f(r,\bar p)$
is a strictly increasing function on~$\bbR$ for every~$\bar p\in O_p$.
Note that~$f(0,\bar p)<0$ for~$\bar p\in O_p$ and
that~$\lim_{r\to\infty}f(r,\bar p)=\infty$ for~$\bar p\in O_p$ 
(since each of the $\rho_\alpha$ is positive).  
It then follows by the intermediate value
theorem and the implicit function theorem that 
the equation~$f(r,p)=0$ can be solved uniquely and real-analytically 
for~$r\ge0$ on an open set~$O^*_p\subset\bbR^m$ containing the domain~$O_p$. 

Thus, let~$\mathbf{r}:O^*_p\to\bbR$ satisfy~$f\bigl(\mathbf{r}(p),p\bigr)=0$
and  set
\begin{equation*}
u_\alpha = p_\alpha e^{-\rho_\alpha \mathbf{r}(p)}.
\end{equation*} 
Then, by construction,
\begin{equation*}
\sum_\alpha u_\alpha = \sum_\alpha  p_\alpha e^{-\rho_\alpha \mathbf{r}(p)}
= \mathbf{r}(p),
\end{equation*}
so that~$p_\alpha = u_\alpha e^{\rho_\alpha(u_2+\cdots+u_{m+1})}$.  
Moreover, when~$p$ lies in~$O_p$,
\begin{equation*}
1+\sum_\beta \rho_\alpha\,u_\alpha 
= 1 + \sum_\alpha \rho_\alpha e^{-\rho_\alpha \mathbf{r}(p)}p_\alpha
= \frac{\partial f}{\partial r}\bigl(\mathbf{r}(p),p) >0,
\end{equation*}
so the image point lies in~$H_u$.  The
inversion of the original system is therefore
\begin{equation*}
y_\alpha = \frac{p_\alpha\,e^{-\rho_\alpha \mathbf{r}(p)}}
 {1 + \sum_\beta \rho_\beta\,p_\beta\,e^{-\rho_\beta \mathbf{r}(p)}}
= \mathbf{y}_\alpha\bigl(p_2,\ldots,p_{m+1}\bigr),
\end{equation*}
as was desired.

Now define a metric on~$\C{m}$ as follows:  
First, define functions on~$\C{m}$ by
\begin{equation*}
G^{\alpha\beta}(\mathbf{z})
=   \overline{z_\alpha}\,z_\beta\,
    \left( \frac{\delta_{\alpha\beta}}
      { \mathbf{y}_\alpha\,\bigl(|z_2|^2,\ldots,|z_{m+1}|^2\bigr) } 
  + \frac{(\rho_\alpha{+}\rho_\beta)}{\mathbf{a}}
  + \frac{\,\rho_\alpha\rho_\beta\,\mathbf{t}\,}{\mathbf{a}^2}\right)\,,
\end{equation*}
where
\begin{equation*}
\mathbf{a} 
=1-\sum_\beta\rho_\beta\,
\mathbf{y}_\beta\bigl(|z_2|^2,\ldots,|z_{m+1}|^2\bigr)
\qquad\text{and}\qquad
\mathbf{t} 
=\sum_\beta\mathbf{y}_\beta\bigl(|z_2|^2,\ldots,|z_{m+1}|^2\bigr).
\end{equation*}
Note that~$\mathbf{a}$ is strictly positive on~$\C{m}$.  Moreover, 
since~$\mathbf{y}_\alpha = p_\alpha \mathbf{y}^*_\alpha$ 
where~$\mathbf{y}^*_\alpha$ is strictly positive on~$\C{m}$, the formula
for~$G^{\alpha\beta} = \overline{G^{\beta\alpha}}$ defines a smooth 
function on~$\C{m}$ for all~$\alpha$ and~$\beta$.  Moreover, the
inequalities satisfied by the~$\mathbf{y}_\alpha$ show that the Hermitian
matrix~$G(\mathbf{z}) = \bigl(G^{\alpha\beta}(\mathbf{z})\bigr)$ is 
positive definite for all~$\mathbf{z}\in\C{m}$.  
Let~$G_{\alpha\beta}(\mathbf{z})$
denote the components of the inverse matrix and define
\begin{equation*}
d\mathbf{s}^2 = G_{\alpha\beta}(\mathbf{z})\,dz_\alpha\,d\overline{z_\beta}.
\end{equation*}
This is an Hermitian metric on~$\C{m}$.  It is
visibly invariant under the torus action generated by the
real parts of the holomorphic vector fields
\begin{equation*}
Z_\alpha = 4i\,z_\alpha\,\frac{\partial\hfill}{\partial z_\alpha}.
\end{equation*}
Setting~$\zeta^\alpha = dz_\alpha/z_\alpha$ yields
$\zeta^\alpha(Z_\beta) = 4i\,\delta^\alpha_\beta$.
Tracing through the construction above, one sees that, 
away from the complex hyperplanes~$z_\alpha = 0$, the metric can 
be written in the form
\begin{equation*}
d\mathbf{s}^2 
= f_{\alpha\beta}(\mathbf{z})\,\zeta^\alpha\circ\overline{\zeta^\beta},
\end{equation*}
where the inverse matrix~$f^{\alpha\beta}$ has the form
\begin{equation*}
f^{\alpha\beta}(\mathbf{z})
=   \frac{\delta_{\alpha\beta}}
      { \mathbf{y}_\alpha\bigl(|z_2|^2,\ldots,|z_{m+1}|^2\bigr)} 
  + \frac{(\rho_\alpha{+}\rho_\beta)}{\mathbf{a}}
  + \frac{\rho_\alpha\rho_\beta\,\mathbf{t}}{\mathbf{a}^2}\,,
\end{equation*}
with
\begin{equation*}
\mathbf{a} 
=1-\sum_\beta\rho_\beta\,\mathbf{y}_\beta\bigl(|z_2|^2,\ldots,|z_{m+1}|^2\bigr)
\qquad\text{and}\qquad
\mathbf{t} 
=\sum_\beta\mathbf{y}_\beta\bigl(|z_2|^2,\ldots,|z_{m+1}|^2\bigr).
\end{equation*}
Thus, the map from~$\C{m}$ to~$C(p_D,\mu)$ defined by~$l_\alpha 
= \mathbf{y}_\alpha\bigl(|z_2|^2,\ldots,|z_{m+1}|^2\bigr)$ is
a Riemannian submersion from the complement of the 
hyperplanes~$z_\alpha=0$ onto~$C(p_D,\mu)$. 

It not difficult now to trace through the construction and see that the 
restriction of the metric~$d\mathbf{s}^2$ to~$\bbR^m\subset\C{m}$ is 
isometric to the metric~$R_\rho$ on~$E_\rho$ and is hence complete.  
It now follows without difficulty that~$d\mathbf{s}^2$ is complete on~$\C{m}$.
Note that this completeness is a consequence of the completeness of
the metric~$R_\rho$ on~$E_\rho$ and so, by the discussion 
in~\S\ref{sssec: explicit forms},
is valid for any~$\rho$ all of whose entries are non-negative, and
not just for those whose entries are positive and strictly increasing.    

Moreover, looking back at the formula for the metric on~$L$ 
and comparing terms, one sees that~$\bigl(\C{m},d\mathbf{s}^2\bigr)$ 
is locally and (therefore, by completeness) globally holomorphically 
isometric to~$L$ with its induced metric and that, 
under this isomorphism, the K\"ahler form corresponding to~$d\mathbf{s}^2$ 
is simply
\begin{equation*}
\Upsilon 
= {\ts\frac{i}2}\,f_{\alpha\beta}(\mathbf{z})\,
          \zeta^\alpha\w\overline{\zeta^\beta}
= {\ts\frac{i}2}\,G_{\alpha\beta}(\mathbf{z})\,
          dz_\alpha\w\overline{dz_\beta}\,.
\end{equation*}
This provides the desired `explicit' formula for the metric on the 
leaf~$L$.  (Bear in mind, though, that the function~$\mathbf{r}$, which
is the crucial ingredient in the recipe for the metric, was found 
by abstractly solving an implicit equation.) 

As the reader can verify, the formula given above simplifies (after 
some trivial changes in notation) to the formula for~$g_\rho$ given
before the statement of Theorem~\ref{thm: complete metrics}.

The argument to this point shows that the metric~$g_\rho$ defined 
before the statement of Theorem~\ref{thm: complete metrics} 
is Bochner-K\"ahler for any~$\rho$ 
all of whose entries are positive and distinct.  However,
the property of being Bochner-K\"ahler is preserved in the limit
as any of the entries of~$\rho$ vanish or become equal.  (The 
curvature tensor is evidently analytic in~$\rho$.)  Consequently,
the metric~$g_\rho$ is Bochner-K\"ahler and complete for any~$\rho$ with 
non-negative entries.

Finally, returning to the notation used before the statement of
Theorem~\ref{thm: complete metrics}, 
suppose $\rho = (\rho_1,\ldots,\rho^n)$ with
\begin{equation*}
0\le\rho_1\le\rho_2\le\cdots\le\rho_n\,.
\end{equation*}
Suppose first that these inequalities are strict and set
\begin{equation*}
r_1 = \frac1{(n{+}2)}\bigl(\rho_1+\cdots+\rho_n\bigr)
\end{equation*}
and then~$r_\alpha = r_1-\rho_{\alpha-1}$ for~$2\le\alpha\le n{+}1$,
so that
\begin{equation*}
2r_1 + r_2 + \cdots + r_{n+1} = 0
\end{equation*} 
and~$r_1>r_2>\cdots>r_{n+1}$.  From the construction in the first part 
of the proof, it follows that the metric~$g_\rho$ satisfies
\begin{equation*}
p_C(t) = p_D(t) = (t-r_1)^2(t-r_2)\cdots(t-r_{n+1})
\end{equation*}
and has cohomogeneity~$n$.  Since the metric is complete, 
by Propositions~\ref{prop: metric complete implies cell bounded} 
and~\ref{prop: distinct roots is not complete}, 
the momentum cell must fall into the SubCase~3-1$b$.
Moreover, since any strictly decreasing sequence~$(r_1,\ldots,r_{n+1})$ 
satisfying~$2r_1 + r_2 + \cdots + r_{n+1}=0$ can be written in the above form
for a unique~$\rho$ with~$0<\rho_1<\cdots<\rho_n$, it follows that such
parameters account for all of the $n$-dimensional reduced momentum cells
in SubCase 3-1$b$.  Thus, by Theorem~\ref{thm: analytically connected classes}, 
this formula gives all of the 
complete, simply-connected cohomogeneity~$n$ Bochner-K\"ahler metrics
of dimension~$n$. Note that the origin is the unique fixed point of
all of the vector fields in~$\euz$, and it follows 
from~\eqref{eq: define pC} that
\begin{equation*}
p_{h(0)}(t) = (t-r_2)\cdots(t-r_{n+1})
\qquad\text{and}\qquad
\bigl(t^2 + h_1(0)\,t + V(0)\bigr) = (t-r_1)^2.
\end{equation*}
 From this, it follows from Proposition~\ref{prop: dim sym alg} 
that the Lie algebra of
the symmetry group is~$\euz$.  Since the group of symmetries is
necessarily connected, it follows that the flows in~$\euz$ generate
the entire symmetry group.

Now consider what happens as the~$\rho_i$ vary.  The metric~$g_\rho$
depends analytically on~$\rho$, so the formulae for~$p_{h(0)}(t)$ and~$V(0)$ 
must remain true for all values of~$\rho$.  The vector fields
in~$\euz$ all vanish at~$0$, so it follows that~$B_3 = |T|^2$ must
vanish at~$0$.  Now, applying Theorem~\ref{thm: point data to polys}, 
one sees that, as~$\rho$ varies 
through~$\bbR^n$ satisfying~$0\le \rho_1\le\cdots\le\rho_n$, 
the values of the moduli pass through all of the values 
that can give rise to momentum cells in SubCase 3-1$b$, with the one
exception of~$\rho=0$, since, in this case, there is no such cell.
Consequently, as $\rho$ varies in the primary orthant, 
the~$g_\rho$ account for all of the possible analytically connected 
equivalence classes that can contain a complete metric.  Since these
metrics are all complete, it follows from 
Theorem~\ref{thm: analytically connected classes}, and 
Propositions~\ref{prop: metric complete implies cell bounded}
and~\ref{prop: distinct roots is not complete}, 
that these contain all of the inhomogeneous complete Bochner-K\"ahler 
metrics on simply connected manifolds.  
\end{proof}

\begin{remark}[Existence]
Interestingly, the argument above justifies the original assumption
that there exists a complete Bochner-K\"ahler metric that is not locally
homogeneous by producing such examples at the end.
\end{remark}

\subsubsection{Weighted projective spaces}
\label{sssec: wtd proj spaces}

A construction similar to that in the SubCase 3-1$b$ can be used to 
express the leaf metric in complex coordinates in SubCase 4-0.  
Since the details are similar to those in the 
proof of Theorem~\ref{thm: complete metrics}, I will be brief.

Consider a momentum cell~$C(p_D,\mu)$ in SubCase 4-0, with
\begin{equation*}
p_D(t) = (t-r_0)(t-r_1)\cdots(t-r_{m+1}), 
 \qquad\qquad (r_0>\cdots>r_{m+1}).
\end{equation*}
The cell~$C(p_D,\mu)$ is defined by the inequalities~
$l_\alpha\ge0$ for~$1\le\alpha\le m{+}1$.

The first task is to produce holomorphic coordinates on the completion~$X_A$
when~$A = \{2,\ldots,m{+}1\}$.   In fact, the argument to follow will show 
that~$X_A$ is biholomorphic to~$\C{m}$.  By Proposition~\ref{prop: leaf data}, 
it suffices consider the case~$n=m$, for one can always reduce to this 
case by simultaneously translating all of the~$r_\alpha$ 
until~$r_0+\cdots+r_{m+1}=0$.  So assume that this has been done.

For simplicity, use the 
abbreviation~$\rho_\alpha$ ($>0$) for~$r_0{-}r_\alpha$ when~$\alpha\ge1$.  
Use~\eqref{eq: define l-alpha} as the definition of the linear 
function~$l_\alpha:\bbR^m\to\bbR$ as before,
and note that the relations~\eqref{eq: l relations} can be written as
\begin{equation*}
 \rho_1\,l_1 = 1 - \sum_{\alpha>1} \rho_\alpha\,l_\alpha
\quad\qquad\text{and}\qquad\quad
-\rho_1\,l_0 = 1 - \sum_{\alpha>1} (\rho_\alpha{-}\rho_1)\,l_\alpha\,.
\end{equation*}
The functions~$l_2,\ldots,l_{m+1}$ are nonnegative coordinates on the cell,
the function~$l_1$ is nonnegative, and the function~$l_0$ is strictly
negative.  Of course, the function~$l_1$ is strictly positive on the 
cell minus the face~$l_1=0$, and this will be important below.
In what follows, whenever repeated indices invoke the summation convention,
the range will be assumed to be~$2\le\alpha,\beta\le m{+}1$ unless stated
otherwise.

As was done in the proof of 
Proposition~\ref{prop: distinct roots is not complete}, 
set~$w_\alpha = (h')^*(l_\alpha)$ and consider the vector 
fields~$W_\alpha\in\euz$ defined by~$W_\alpha\lhk\Omega = -dw_\alpha$. 
The vector fields~$W_2,\ldots,W_{m+1}$ are a basis of~$\euz$ and are 
linearly independent on~$X_A^\circ$.  The map~$h':X_A\to C(p_D,\mu)$ is a 
Riemannian submersion on~$X_A^\circ$ when~$C(p_D,\mu)$ 
is given the metric~$R_D$.  The image~$h'(X_A)$ is equal to~$C(p_D,\mu)$
minus the face~$l_1=0$.  

Each~$W_\alpha$ vanishes 
on~$N_\alpha = (h')^{-1}(F_\alpha)$, which, since the flow of
$W_\alpha$ is isometric, is a totally geodesic complex hypersurface 
in~$X_A$.  Moreover,~$W_\alpha$ is nonzero off of~$N_\alpha$ 
for~$\alpha\ge2$.

As before, I will show that there are
holomorphic coordinates~$z_2,\ldots z_{m+1}$ on~$X_A$ for which
\begin{equation*}
W_\alpha-iJW_\alpha 
= 4\,i\,z_\alpha\,\frac{\partial\hfill}{\partial z_\alpha}\,,
\qquad\qquad 2\le \alpha\le m{+}1,
\end{equation*}
and find the expression for the Bochner-K\"ahler metric on~$X_A$ 
in these coordinates. 

Accordingly, let~$\zeta^2,\ldots,\zeta^{m+1}$ be the holomorphic
1-forms on~$X_A^\circ$ that satisfy
\begin{equation*}
\zeta^\alpha(W_\beta-iJW_\beta) = 4i\,\delta^\alpha_\beta
\end{equation*}
These forms extend meromorphically to~$X_A$, with simple poles along
the hypersurfaces~$N_\alpha$.  Since the vector fields~$W_\alpha$
Lie-commute, it follows that each~$\zeta^\alpha$ is closed. 
As before, if the coordinates~$z_\alpha$ are to exist as claimed,
it will have to be true that~$\zeta^\alpha = dz_\alpha/z_\alpha$.

Writing~$\zeta^\alpha = \xi^\alpha + i\,\eta^\alpha$, the above
equations are equivalent to
\begin{equation*}
  \xi^\alpha(W_\beta) = \eta^\alpha(JW_\beta) = 0,
\qquad\qquad
-\xi^\alpha(JW_\beta) = \eta^\alpha(W_\beta) = 2\,\delta^\alpha_\beta,
\end{equation*}
Again, if the coordinates~$z_\alpha$ exist as claimed, it will follow
that~$2\xi^\alpha = d\bigl(\log|z_\alpha|^2\bigr)$.  

Since the~$\zeta^\alpha$ are a basis for the holomorphic 1-forms 
on~$X_A^\circ$, the metric on~$X_A^\circ$ can be written in the form
\begin{equation*}
ds^2 = g_{\alpha\beta}\,\zeta^\alpha{\circ}\overline{\zeta^\beta},
\end{equation*}
where~$g_{\alpha\beta} = \overline{g_{\beta\alpha}}$ 
and where the pullback of~$\Omega$ to~$X_A^\circ$ is
\begin{equation*}
\Upsilon 
= {\ts\frac{i}2}\,g_{\alpha\beta}\,\zeta^\alpha\w \overline{\zeta^\beta}.
\end{equation*}
Now, the identity~$\zeta^\alpha(W_\beta) = 2i\,\delta^\alpha_\beta$ implies
\begin{equation*}
dw_\alpha = -W_\alpha\lhk\Upsilon 
= g_{\alpha\beta}\,\overline{\zeta^\beta}
                    + g_{\beta\alpha}\,\zeta^\beta,
\end{equation*}
or, equivalently,
\begin{equation*}
dw_\alpha = (g_{\alpha\beta}+g_{\beta\alpha})\,\xi^\beta
      - i\,(g_{\alpha\beta}-g_{\beta\alpha})\,\eta^\beta.
\end{equation*}
Since~$W_\alpha$ is tangent to the fibers of~$h'$, and since~$w_\alpha = 
(h')^*(l_\alpha)$ is constant on those fibers, it follows that
the coefficient of~$\eta^\beta$ in the
above equation must vanish, i.e., $g_{\alpha\beta} 
= g_{\beta\alpha}=\overline{g_{\alpha\beta}}$. Thus,
\begin{equation*}
dw_\alpha = 2g_{\alpha\beta}\, \xi^\beta.
\end{equation*}

Define~$g^{\alpha\beta}=g^{\beta\alpha}$ so 
that~$g^{\alpha\beta}g_{\beta\gamma} = \delta^\alpha_\gamma$. 
Note, in particular, that~$\xi^\alpha = \frac12 g^{\alpha\beta}\,dw_\beta$.
The metric on~$X_A^\circ$ can now be written in the form
\begin{equation*}
\begin{split}
ds^2 &= g_{\alpha\beta}\,\,\zeta^\alpha{\circ}\overline{\zeta^\beta}
      = g_{\alpha\beta}\,
       (\xi^\alpha+i\,\eta^\alpha){\circ}(\xi^\beta-i\,\eta^\beta)\\
     &= g_{\alpha\beta}\bigl(\xi^\alpha{\circ}\xi^\beta 
                        + \eta^\alpha{\circ}\eta^\beta\bigr)\\
     &= {\ts\frac14}g^{\alpha\beta}\,dw_\alpha{\circ}dw_\beta 
          + g_{\alpha\beta}\,\eta^\alpha{\circ}\eta^\beta.\\
\end{split}
\end{equation*}
Since~$h'$ is a Riemannian submersion on~$X_A^\circ$, it follows that
\begin{equation*}
{\ts\frac14}g^{\alpha\beta}\,dw_\alpha\,dw_\beta 
= (h')^*(R_D) 
=(h')^*\bigl(T^{\alpha\beta}\bigr)\,dw_\alpha\,dw_\beta\,,
\end{equation*}
where~$T^{\alpha\beta} = T^{\beta\alpha}$ for~$2\le\alpha,\beta\le m{+}1$ 
is defined on~$C(p_D,\mu)^\circ$ so that the formula
\begin{equation*}
\sum_{\alpha,\beta=2}^{m+1} T^{\alpha\beta}\,dl_\alpha\,dl_\beta 
= R_D = \sum_{\alpha=0}^{m+1}\frac{{dl_\alpha}^2}{4l_\alpha}
\end{equation*}
holds.  Thus,~$g^{\alpha\beta}=4(h')^*\bigl(T^{\alpha\beta}\bigr)$.

Using the relations above that express~$l_0$ and $l_1$
in terms of~$l_2,\ldots,l_{m+1}$, 
it follows from the formula for~$R_D$ that
\begin{equation*}
4T^{\alpha\beta}
=   \frac{\delta_{\alpha\beta}}{l_\alpha} 
  + \frac{(\rho_\alpha{-}\rho_1)(\rho_\beta{-}\rho_1)}{{\rho_1}^2\,l_0}
  + \frac{\rho_\alpha\rho_\beta}{{\rho_1}^2\,l_1}\,,
\end{equation*}
so that, in particular,
\begin{equation*}
\begin{split}
4T^{\alpha\beta}\,dl_\beta  
&= \frac{dl_\alpha}{l_\alpha} 
   + \frac{\rho_\alpha{-}\rho_1}{\rho_1}\frac{dl_0}{l_0}
   - \frac{\rho_\alpha}{\rho_1}\frac{dl_1}{l_1}\\
&= d\left(\log\left(l_\alpha\,(-l_0)^{(\rho_\alpha{-}\rho_1)/\rho_1}
              (l_1)^{-\rho_\alpha/\rho_1}  \right) \right).
\end{split}
\end{equation*}
Meanwhile, if the coordinates~$z_\alpha$ exist as claimed, this will
imply that
\begin{equation*}
\frac{d|z_\alpha|^2}{|z_\alpha|^2} = 2\xi^\alpha 
= g^{\alpha\beta}\,dw_\beta 
= (h')^*\bigl(4T^{\alpha\beta}\,dl_\beta \bigr)
= d\left((h')^*\left(\log
        \frac{l_\alpha\,(-l_0)^{(\rho_\alpha{-}\rho_1)/\rho_1}}
             {(l_1)^{\rho_\alpha/\rho_1}}
     \right)\right),
\end{equation*}
i.e., there will exist constants~$c_\alpha>0$ so that
\begin{equation*}
|z_\alpha|^2 
= c_\alpha (h')^*\left(\frac{l_\alpha\,(-l_0)^{(\rho_\alpha{-}\rho_1)/\rho_1}}
             {(l_1)^{\rho_\alpha/\rho_1}}\right).
\end{equation*}
Since~$z_\alpha$ would only be determined up to a multiplicative
constant anyway by the above normalizations, one might as well
take~$c_\alpha=1$, which will normalize the~$z_\alpha$ up to a
phase.  Writing~$x_\alpha = -l_\alpha/l_0\ge0$ for $\alpha>0$, this
equation can be written more simply as
\begin{equation*}
|z_\alpha|^2 
= (h')^*\left(\frac{x_\alpha}
             {(x_1)^{\rho_\alpha/\rho_1}}\right),
\qquad\qquad \alpha = 2,3,\ldots,m{+}1,
\end{equation*}
where the~$x_\alpha\ge0$ satisfy the relation~$x_1+\cdots+x_{m+1}=1$.

Consider the equation
\begin{equation*}
s + \sum_{\alpha=2}^{m+1} |z_\alpha|^2 s^{\rho_\alpha/\rho_1} = 1
\end{equation*}
on~$\bbR\times\C{m}$.  An analysis similar to the one performed in the 
proof of Theorem~\ref{thm: complete metrics} 
shows that when~$\rho_\alpha/\rho_1\ge0$ for~$\alpha\ge2$
there is a unique real-analytic function~$\mathbf{s}: \C{m}\to (0,\infty)$
that satisfies
\begin{equation*}
\mathbf{s}(z) 
+ \sum_{\alpha=2}^{m+1} 
|z_\alpha|^2 \bigl(\mathbf{s}(z)\bigr)^{\rho_\alpha/\rho_1} = 1
\end{equation*}
for all~$z\in\C{m}$.  Note that the function~$\mathbf{s}$ is invariant
under the standard $m$-torus action on~$\C{m}$ and is algebraic if
and only if all of the ratios~$\rho_\alpha/\rho_1$ are rational.
Using the function~$\mathbf{s}$, the equations above can be
solved in the form~$(h')^*(x_1) = \mathbf{s}(z)$ and
\begin{equation*}
(h')^*(x_\alpha) = |z_\alpha|^2 \bigl(\mathbf{s}(z)\bigr)^{\rho_\alpha/\rho_1},
\qquad (\alpha>1),
\end{equation*}
whence, by algebra, follows the formula
\begin{equation*}
w_\alpha = (h')^*(l_\alpha) 
= \frac{|z_\alpha|^2 \bigl(\mathbf{s}(z)\bigr)^{\rho_\alpha/\rho_1}}
  {\rho_1 + \sum_{\beta=2}^{m+1}
    (\rho_\beta{-}\rho_1)\,|z_\beta|^2 
        \bigl(\mathbf{s}(z)\bigr)^{\rho_\beta/\rho_1} },
\qquad (2\le\alpha\le m{+}1).
\end{equation*}

This motivates defining a metric on~$\C{m}$ as follows:  Set
\begin{equation*}
\mathbf{w}_\alpha(z)
= \frac{|z_\alpha|^2 \bigl(\mathbf{s}(z)\bigr)^{\rho_\alpha/\rho_1}}
  {\rho_1 + \sum_{\beta=2}^{m+1}
    (\rho_\beta{-}\rho_1)\,|z_\beta|^2 
        \bigl(\mathbf{s}(z)\bigr)^{\rho_\beta/\rho_1} }
\end{equation*}
for $2\le\alpha\le m{+}1$ 
and define functions~$\mathbf{w}_1$ and~$\mathbf{w}_0$ on~$\C{m}$ by
\begin{equation*}
 \rho_1\,\mathbf{w}_1(z) 
= 1 - \sum_{\alpha>1} \rho_\alpha\,\mathbf{w}_\alpha(z)
\qquad\text{and}\qquad
-\rho_1\,\mathbf{w}_0(z) 
= 1 - \sum_{\alpha>1} (\rho_\alpha{-}\rho_1)\,\mathbf{w}_\alpha(z)\,.
\end{equation*}
Then~$\mathbf{w}_1$ and $-\mathbf{w}_0$ are strictly positive on~$\C{m}$.
For $2\le\alpha,\beta\le m{+}1$, define functions on~$\C{m}$ by
\begin{equation*}
G^{\alpha\beta}(z)
=   \overline{z_\alpha}\,z_\beta\,\left(
    \frac{\delta_{\alpha\beta}}{\mathbf{w}_\alpha(z)} 
  + \frac{(\rho_\alpha{-}\rho_1)(\rho_\beta{-}\rho_1)}
           {{\rho_1}^2\,\mathbf{w}_0(z)}
  + \frac{\rho_\alpha\rho_\beta}{{\rho_1}^2\,\mathbf{w}_1(z)}\,\right)\,.
\end{equation*}
It is not difficult to show that the Hermitian
matrix~$G(z) = \bigl(G^{\alpha\beta}(\mathbf{z})\bigr)$ is 
positive definite for all~$z\in\C{m}$.  
Let~$G_{\alpha\beta}(z)$
denote the components of the inverse matrix and define
\begin{equation*}
d\mathbf{s}^2 = G_{\alpha\beta}(z)\,dz_\alpha\,d\overline{z_\beta}.
\end{equation*}
This is an Hermitian metric on~$\C{m}$.  It is
visibly invariant under the torus action generated by the
real parts of the holomorphic vector fields
\begin{equation*}
Z_\alpha = 4i\,z_\alpha\,\frac{\partial\hfill}{\partial z_\alpha}.
\end{equation*}
Setting~$\zeta^\alpha = dz_\alpha/z_\alpha$ yields
$\zeta^\alpha(Z_\beta) = 4i\,\delta^\alpha_\beta$.
Tracing through the construction above, one sees that, 
away from the complex hyperplanes~$z_\alpha = 0$, the metric can 
be written in the form
\begin{equation*}
d\mathbf{s}^2 
= f_{\alpha\beta}(z)\,\zeta^\alpha\circ\overline{\zeta^\beta},
\end{equation*}
where the inverse matrix~$f^{\alpha\beta}$ has the form
\begin{equation*}
f^{\alpha\beta}(z)
=  \frac{\delta_{\alpha\beta}}{\mathbf{w}_\alpha(z)} 
  + \frac{(\rho_\alpha{-}\rho_1)(\rho_\beta{-}\rho_1)}
           {{\rho_1}^2\,\mathbf{w}_0(z)}
  + \frac{\rho_\alpha\rho_\beta}{{\rho_1}^2\,\mathbf{w}_1(z)}\,.
\end{equation*}
In particular, the map from~$\C{m}$ to~$C(p_D,\mu)$ defined by~$l_\alpha 
= \mathbf{w}_\alpha\bigl(|z_2|^2,\ldots,|z_{m+1}|^2\bigr)$ is
a Riemannian submersion from the complement of the 
hyperplanes~$z_\alpha=0$ onto~$C(p_D,\mu)^\circ$ endowed with the 
metric~$R_D$. 

It not difficult now to trace through the construction and see that the 
restriction of the metric~$d\mathbf{s}^2$ to~$\bbR^m\subset\C{m}$ is 
isometric to the metric~$R_D$ on~$E$ as defined 
in~\S\ref{sssec: explicit forms}.   

Moreover, looking back at the formulae for the metric on~$X_A^\circ$ 
and comparing terms, one sees that~$\bigl(\C{m},d\mathbf{s}^2\bigr)$ 
must be globally holomorphically 
isometric to~$X_A$ with its Bochner-K\"ahler metric and that, 
under this isomorphism, the K\"ahler form corresponding to~$d\mathbf{s}^2$ 
is simply
\begin{equation*}
\Upsilon 
= {\ts\frac{i}2}\,f_{\alpha\beta}(\mathbf{z})\,
          \zeta^\alpha\w\overline{\zeta^\beta}
= {\ts\frac{i}2}\,G_{\alpha\beta}(\mathbf{z})\,
          dz_\alpha\w\overline{dz_\beta}\,.
\end{equation*}
This provides the desired explicit formula for the metric on the 
leaf~$X_A$.

Although the derivation provided the inequalities~
$0<\rho_1<\cdots<\rho_{m+1}$, the recipe given for the metric
only needs the assumption~$\rho_\alpha>0$ for~$1\le\alpha\le m{+}1$.
This means, for example, that the metric makes sense when all of
the~$\rho_\alpha$ are equal.  In this case, the reader can verify
that the metric~$d\mathbf{s}^2$ on~$\C{m}$ is simply the Fubini-Study 
metric on~$\bbC\bbP^m$ restricted to the complement of a hyperplane. 

Suppose now that all of the ratios~$\rho_\alpha/\rho_1$ are rational
and let~$r>0$ be such that~$\rho_\alpha = (m{+}2)r\, p_\alpha$ where the
numbers~$0<p_1<\ldots<p_{m+1})$ are integers with no common divisor.
This uniquely defines~$r$ and the integers~$p_\alpha$.  Moreover,
the equations~$r_0-r_\alpha = \rho_\alpha = rp_\alpha$ 
and~$r_0+r_1+\cdots+r_{m+1}=0$ imply
\begin{equation*}
r_\alpha 
= r\left(\sum_{\beta=0}^{m+1} p_\beta - (m{+}2)p_\alpha\right)
\end{equation*}
where, for notational symmetry, I have set~$p_0 = 0$.

Recalling 
that~$\rho_1\,W_1 + \rho_2\,W_2 + \cdots + \rho_{m+1}W_{m+1}=0$, it
follows that
\begin{equation*}
p_1\bigl(W_1-iJW_1\bigr) = -4i\left(
 p_2\,z_2\,\frac{\partial\hfill}{\partial z_2}
+\cdots+
 p_{m+1}\,z_{m+1}\,\frac{\partial\hfill}{\partial z_{m+1}}
\right).
\end{equation*}

Now, set~$[p]=[p_1,\ldots,p_{m+1}]$ and consider the weighted projective 
space~$\bbC\bbP^{[p]}$ one gets by taking the quotient of~$\C{m+1}$
minus the origin by the~$\bbC^*$-action
\begin{equation*}
\lambda\cdot(Z_1,Z_2,\ldots,Z_{m+1})
= \bigl(\lambda^{p_1}Z_1,\lambda^{p_2}Z_2,
          \ldots,\lambda^{p_{m+1}}Z_{m+1}\bigr).
\end{equation*}
This is an orbifold and not a manifold except when~$p_1 = \cdots = p_{m+1}$
(in which case, this is $\bbC\bbP^m$).   
Let~$[Z_1,\ldots,Z_{m+1}]\in \bbC\bbP^{[p]}$ denote 
the orbit of~$(Z_1,\ldots,Z_{m+1})\in\C{m}$.
Consider the holomorphic mapping~$\Phi_1:\C{m}\to \bbC\bbP^{[p]}$
defined by
\begin{equation*}
\Phi_1(Z_2,\ldots,Z_{m+1}) = [1,Z_2,\ldots,Z_{m+1}].
\end{equation*}
This mapping is a $p_1$-fold branched covering of its image and 
the above considerations show that the metric~$d\mathbf{s}^2$ extends
to be a smooth orbifold metric on~$\bbC\bbP^{[p]}$.  The end result is
the following:

\begin{theorem}\label{thm: wtd proj space}
Every weighted projective space~$\bbC\bbP^{[p]}$ supports a
Bochner-K\"ahler metric.  
\end{theorem}

\begin{remark}[Uniqueness]
Of course, in the classical case of projective space, the Bochner-K\"ahler
metric so constructed is a constant multiple of the Fubini-Study metric.
By Corollary~\ref{cor: compact BK is loc sym}, 
this is the unique Bochner-K\"ahler metric on~$\bbC\bbP^m$,
up to a constant multiple.  I suspect, though I have not checked
all of the details, that this uniqueness holds for all of the 
weighted projective spaces.  
\end{remark}

\begin{remark}[Reduction]
The reader will recall that one way of constructing the Fubini-Study metric 
is to apply reduction to the flat K\"ahler metric under the
diagonal $S^1$-action.  Given this, one might suspect that the 
Bochner-K\"ahler metric on $\bbC\bbP^{[p]}$ is got from the flat K\"ahler
metric by applying reduction to the weighted $S^1$-action described
above.  However, calculation shows that, except in the classical case,
the reduction metric is \emph{not} Bochner-K\"ahler.
\end{remark}

\subsection{Reduction and the full metric}
\label{ssec: red and full metric}

Theorem~\ref{thm: leaf metric} 
provides a formula for the induced metric on the~$\cF$-leaves
of a Bochner-K\"ahler metric.  In the case of maximal cohomogeneity, i.e.,
$m=n$, the regular set~$M^\circ$ constitutes a single~$\cF$-leaf, so
this formula determines the metric completely.  In this section, 
I will indicate how one can reconstruct the full metric from the knowledge 
of the metric on the~$\cF$-leaves.  Thus, for the rest of this section,
I will assume that~$M^n$ is endowed with a Bochner-K\"ahler
metric of cohomogeneity~$m$ satisfying~$0<m<n$, 
since otherwise, there is nothing to do.  

Let~$p_C(t)$ and~$p_D(t)$
be the characteristic polynomials of the Bochner-K\"ahler structure.
Write~
\begin{equation*}
p_{h''}(t) = (t-\lambda_{m+1})\cdots (t-\lambda_{n})
\end{equation*}
where, by Proposition~\ref{prop: ph'' roots are pD roots}, 
the roots~$\lambda_{m+1}\ge\cdots\ge\lambda_{n}$ 
are also roots of~$p_D(t)$.  Let~$\pi:P_2\to M^\circ$ be the 
$G_\Lambda$-bundle as described in the proof of 
Proposition~\ref{prop: ph'' roots are pD roots} and return
to that notation, particularly the index ranges.  Recall that the
$\lambda_i$ are distinct and not equal to any of the~$\lambda_a$, and
that the~$T_i$ are positive and real and satisfy
\begin{equation*}
{T_i}^2 = \frac{p_D(\lambda_i)}{\prod_{j\not=i}(\lambda_i{-}\lambda_j)}.
\end{equation*}
Also, recall the relations
\begin{equation*}
\phi_{a\bar\imath} = \frac{T_i\,\omega_a}{\lambda_i{-}\lambda_a}\,,
\end{equation*} 
which followed from the structure equations applied to~$h_{a\bar\imath}=0$.
The structure equations applied to~$h_{a\bar b} = \delta_{ab}\lambda_a$ 
yield
\begin{equation*}
(\lambda_a - \lambda_b)\phi_{a\bar b} = 0,
\end{equation*}
so that~$\phi_{a\bar b} = 0$ when~$\lambda_a\not=\lambda_b$.  The
structure equations applied to~$h_{i\bar\jmath}=0$ for~$i\not=j$ yield
the relations
\begin{equation*}
\phi_{i\bar\jmath} 
= -\frac{T_i\,\overline{\omega_j}+T_j\,\omega_i}{\lambda_i-\lambda_j},
\qquad i\not=j,
\end{equation*}
while the structure equations applied to~$h_{i\bar\imath}=\lambda_i$ yield
\begin{equation*}
d\lambda_i = T_i(\omega_i + \overline{\omega_i}).
\end{equation*}
Meanwhile, the structure equations applied to~$T_i$ yield
\begin{equation*}
\begin{split}
dT_i &= -\phi_{i\bar\jmath}\,T_j 
          + ({\lambda_i}^2{+}h_1\,\lambda_i{+}V)\,\omega_i\\
&= -\phi_{i\bar\imath}\,T_i 
   +\sum_{j\not=i}\left(\frac{T_i\,\overline{\omega_j}+T_j\,\omega_i}
                       {\lambda_i-\lambda_j}\right)T_j
   + ({\lambda_i}^2{+}h_1\,\lambda_i{+}V)\,\omega_i\\
\end{split}
\end{equation*}
which can be rearranged to give
\begin{equation*}
\frac{dT_i}{T_i} 
= -\phi_{i\bar\imath}
 +\sum_{j\not=i}\frac{T_j\,\overline{\omega_j}}{\lambda_i-\lambda_j}
 +\left(({\lambda_i}^2{+}h_1\,\lambda_i{+}V)
    +\sum_{j\not=i}\frac{{T_j}^2}{\lambda_i-\lambda_j} \right)
       \,\frac{\omega_i}{T_i}\,.
\end{equation*}

Now, the structure equations for~$\omega_a$ are (summation over $i$ and $b$)
\begin{equation*}
\begin{split}
d\omega_a &= -\phi_{a\bar\imath}\w\omega_i - \phi_{a\bar b}\w\omega_b
           = \frac{T_i\,\omega_a}{\lambda_a{-}\lambda_i}\,\omega_i 
                      -\phi_{a\bar b}\w\omega_b\\
&= -\left(\phi_{a\bar b} 
      +\delta_{a\bar b} \,\frac{T_i\,\omega_i}{\lambda_a{-}\lambda_i}\right)
     \w\omega_b
  = -\left(\varphi_{a\bar b} 
     +\frac12\delta_{a\bar b}\,\frac{d\lambda_i}{\lambda_a{-}\lambda_i}\right)
      \w\omega_b\,,\\
\end{split}
\end{equation*}
where
\begin{equation*}
\varphi_{a\bar b} = -\overline{\varphi_{b\bar a}}
 = \phi_{a\bar b} 
      +\frac12\delta_{a\bar b}\,
   \sum_i\frac{T_i\,(\omega_i-\overline{\omega_i})}
                                              {\lambda_a{-}\lambda_i}\,.
\end{equation*}
Since~$p_{h'}(\lambda_a) = \prod_i(\lambda_a{-}\lambda_i)$
and since~$\phi_{a\bar b}=\varphi_{a\bar b} =0$ 
when~$\lambda_a\not=\lambda_b$, setting
\begin{equation*}
\eta_a = |p_{h'}(\lambda_a)|^{-1/2}\omega_a
\end{equation*}
yields~$d\eta_a = -\varphi_{a\bar b}\w\eta_b$.  This implies that, for
each root~$r$ of~$p_{h''}(t)$, the quadratic form
and 2-form
\begin{equation*}
g_r = \sum_{\{a:\lambda_a = r\}} \eta_a\circ\overline{\eta_a}
\qquad\text{and}\qquad
\Omega_r = {\frac\imath2}\,\sum_{\{a:\lambda_a = r\}} 
            \,\eta_a\w\overline{\eta_a}
\end{equation*}
define a K\"ahler structure on the space of leaves of the 
system~$\{\eta_a = 0\ \vrule\ \lambda_a = r\,\}$ in any open set
in~$M^\circ$ where this leaf space is Hausdorff.  If $r$ has
multiplicity~$\nu>0$, this leaf space has complex dimension~$\nu$.

To compute the curvature of this leaf space, one needs to compute the
2-forms
\begin{equation*}\Phi_{a\bar b} = d\varphi_{a\bar b} 
+ \varphi_{a\bar c}\w\varphi_{c\bar b},
\end{equation*}
so I now turn to this task.  Since~$\varphi_{a\bar b} = 
\phi_{a\bar b}$ when~$a\not=b$, the structure equations for~$a\not=b$
yield (summation on $i$ and $c$)
\begin{equation*}
\begin{split}
\Phi_{a\bar b} &= d\phi_{a\bar b} + \phi_{a\bar c}\w\phi_{c\bar b}
  = -\phi_{a\bar i}\w\phi_{i\bar b} 
       -(\lambda_a + \lambda_b + h_1)\,\omega_a\w\overline{\omega_b}\\
 & = \left[ \frac{{T_i}^2}{(\lambda_a{-}\lambda_i)(\lambda_b{-}\lambda_i)}
             -(\lambda_a + \lambda_b + h_1)\right]
                   \,\omega_a\w\overline{\omega_b}\\
 & = \left[ \frac{p_D(\lambda_i)}
        {(\lambda_a{-}\lambda_i)(\lambda_b{-}\lambda_i)
                    \prod_{j\not=i}(\lambda_i{-}\lambda_j)}
             -(\lambda_a + \lambda_b + h_1)\right]
                   \,\omega_a\w\overline{\omega_b}
\end{split}
\end{equation*}
Rather miraculously, when~$\lambda_a\not=\lambda_b$, the classical identities
of~\S\ref{sssec: explicit forms} show that this expression is zero, as should 
have been expected.
On the other hand, if~$\lambda_a = \lambda_b = r_i$ (but still~$a\not=b$), 
the same classical identities show that this expression simplifies to
\begin{equation*}
\Phi_{a\bar b} = \frac{p'_D(r_i)}{p_{h'}(r_i)}\,\omega_a\w\overline{\omega_b} 
    = \frac{p'_D(r_i)|p_{h'}(r_i)|}{p_{h'}(r_i)}\,\eta_a\w\overline{\eta_b}
    = (-1)^{\mu_i}p'_D(r_i)\,\eta_a\w\overline{\eta_b}\,.
\end{equation*}
(Recall that~$(-1)^{\mu_i}p_{h'}(r_i)>0$ on~$M^\circ$.)  It remains to
compute the quantities~$\Phi_{a\bar a}$.  This computation is greatly
simplified by first observing that $\Phi_{a\bar a}$ must be a sum of terms
of the form~$\omega_b\w\overline{\omega_c}$ where~$\lambda_a = \lambda_b
=\lambda_c$.  Thus, in carrying out the expansion from the definitions,
all other terms can be ignored.  For simplicity, I will use $\equiv$
to denote equality modulo the ideal generated by the 1-forms~$\omega_i$
and~$\overline{\omega_i}$ for~$1\le i\le m$ and the 1-forms~$\phi_{a\bar b}$ 
for~$m<a,b,\le n$.  
Then, first of all (summation over $j$ and~$b$),
\begin{equation*}
d\omega_i = -\phi_{i\bar j}\w\omega_j -\phi_{i\bar b}\w\omega_b
   \equiv  \frac{T_i\,\omega_b\w\overline{\omega_b}}
                      {\lambda_b{-}\lambda_i}\,.
\end{equation*}
Using this and the identities quoted above, 
the calculation of~$\Phi_{a\bar a}$ follows
from the structure equations goes as (summation over~$j$ and~$b$)
\begin{equation*}
\begin{split}
\Phi_{a\bar a} &= d\varphi_{a\bar a} + \varphi_{a\bar b}\w\varphi_{b\bar a}
\equiv d\left(\phi_{a\bar a} 
      +\frac12\,
   \frac{T_j\,(\omega_j-\overline{\omega_j})}
                                              {\lambda_a{-}\lambda_j}\right)\\
&\equiv -\phi_{a\bar\jmath}\w\phi_{j\bar a} 
       - (2\lambda_a + h_1)\,\omega_a\w\overline{\omega_a}
      - (\lambda_a{+}\lambda_b+h_1)\,\omega_b\w\overline{\omega_b} 
      + \frac{{T_j}^2\,\omega_b\w\overline{\omega_b}}
             {(\lambda_a{-}\lambda_j)(\lambda_b{-}\lambda_j)}\\
&=  \frac{p'_D(r_i)}{p_{h'}(r_i)}\,
     \left(\omega_a\w\overline{\omega_a}  \quad
          + \sum_{\{b:\lambda_a = r_i\}}\omega_b\w\overline{\omega_b}\right)\\
& = (-1)^{\mu_i}\,p'_D(r_i)\,
     \left(\eta_a\w\overline{\eta_a}  \quad
          + \sum_{\{b:\lambda_a = r_i\}}\eta_b\w\overline{\eta_b}\right).   
\end{split}
\end{equation*}
These formulae imply that the K\"ahler structure defined by~$g_{r_i}$
and~$\Omega_{r_i}$ actually has constant holomorphic sectional curvature
equal to~$(-1)^{\mu_i}\,4\,p'_D(r_i)$.  

\subsubsection{Reduction}
\label{sssec: red}

Since~$h':M\to C(p_D,\mu)$ is the momentum mapping of the infinitesimal
torus action defined by the basis~$Z'_1,\ldots,Z'_m$ of~$\euz$, it is 
natural to consider the effect of applying symplectic reduction.  Since
the torus action is not assumed to be globally defined (because no
completeness assumptions have been made about the metric), 
this can only be done locally.  

For simplicity, I will only consider reduction at a regular value of the 
reduced momentum mapping~$h':M\to C(p_D,\mu)$.  Recall that~$h':M^\circ\to
C(p_D,\mu)^\circ$ is a submersion, let~$x\in M^\circ$ 
be fixed and let~$\kappa = h'(x)\in C(p_D,\mu)^\circ$.  The method of
symplectic reduction then consists of the following:  Consider the 
codimension~$m$ submanifold~$(h')^{-1}(\kappa)\subset M^\circ$.  
This submanifold is foliated by $m$-dimensional leaves whose tangent spaces 
are spanned by the vector fields~$Z_2,\ldots,Z_{m+1}$.  Suppose that this 
foliation is simple, i.e., its leaf space~$M_\kappa$ is Hausdorff.  
(This can always be arranged by restricting to a suitable open neighborhood 
of~$x$.)  Then the pullback of~$\Omega$ to~$(h')^{-1}(\kappa)\subset M^\circ$
is a closed 2-form that is the pullback to~$(h')^{-1}(\kappa)$ of a 
symplectic form~$\Omega_\kappa$ on~$M_\kappa$.  
The symplectic manifold~$(M_\kappa,\Omega_\kappa)$ is then called the
\emph{symplectic reduction} of~$(M,\Omega)$ at~$\kappa$.  

\begin{proposition}\label{prop: reduction metric} 
Fix~$x\in M^\circ$ and let
$\kappa = h'(x)\in C(p_D,\mu)^\circ$.
There is a unique metric~$g_\kappa$
on~$M_\kappa$ for which the leaf projection~$(h')^{-1}(\kappa)\to M_\kappa$
is a Riemannian submersion.  

The data~$(M_\kappa,g_\kappa,\Omega_\kappa)$
defines a K\"ahler structure on~$M_\kappa$ that is locally isomorphic
to a product of complex space forms.  Specifically, for each root~$r$
of~$p_{h''}(t)$ of multiplicity~$\nu$, the local product contains
a $\nu$-dimensional complex space form of constant holomorphic
sectional curvature
\begin{equation*}
c(r,\kappa) = \frac{4\,p'_D(r)}{p_{h(x)}(r)}
\end{equation*} 
and these are all of the factors.
\end{proposition}

\begin{proof}
Let~$P_2(\kappa)\subset P_2$ be the part of~$P_2$ that lies 
over~$(h')^{-1}(\kappa)\subset M^\circ$.  The structure equations 
on~$P_2(\kappa)$ are the same as those on~$P_2$ with the difference that,
after restriction to~$P_2(\kappa)$ the functions~$\lambda_i$ and~$T_i$
become constant and the 1-forms~$\omega_i$ become purely imaginary.
Note that the equations~$\omega_a = 0$ define the foliation by the
torus leaves.

Now going back to the structure equations, just derived above, one sees
that, on~$P_2(\kappa)$, the equations
\begin{equation*}
d\omega_a = -\varphi_{a\bar b}\w\omega_b
\end{equation*}
hold, where~$\varphi = (\varphi_{a\bar b}) = -\varphi^*$ is blocked
according to the multiplicities in the descending string of eigenvalues
\begin{equation*}
\lambda_{m+1}\ge \lambda_{m+2} \ge \cdots \ge \lambda_{n}\,.
\end{equation*}
It follows, of course, that quadratic form~$g_\kappa 
= \omega_a\circ\overline{\omega_a}$ is well-defined on the leaf 
space~$M_\kappa$ and that this metric and the symplectically reduced
$2$-form~$\Omega_\kappa = \frac\imath 2\,\omega_a\w\overline{\omega_a}$
define a K\"ahler structure on~$M_\kappa$.  

Finally, the computation of the curvature forms above shows that 
this K\"ahler structure is a product of the type described in the 
proposition. 
\end{proof}

\subsubsection{The general metric}
\label{sssec: gen metric}

As the preceding formulae and Proposition~\ref{prop: reduction metric} 
 now make clear, a recipe for
any Bochner-K\"ahler metric on its regular locus~$M^\circ$ can be 
constructed as a generalized warped product over a momentum cell,
where the fibers are products of so-called Sasakian space forms, i.e.,
the canonical circle bundles over complex space forms of constant
holomorphic sectional curvature.  In other words, once the leaf
metric has been found, as in Theorem~\ref{thm: leaf metric}, 
the full metric can be
constructed by group theoretic means.  This is to be expected,
since, after all, the pseudo-group of local symmetries of a connected
Bochner-K\"ahler metric acts transitively on the fibers of the
momentum mapping.  

The explicit formula does not appear to be of great interest. 
For brevity, I will not go into details.

\section{Final remarks}
\label{sec: final remarks}

In this last section, I will make some remarks about related geometries.

\subsection{Pseudo-K\"ahler geometry}
\label{ssec: psdoK geom}

When a complex $n$-manifold~$M$ is endowed with a pseudo-K\"ahler 
structure, i.e., a pseudo-Riemannian metric~$g$ that is invariant under
the complex structure and whose associated 2-form~$\Omega$ is closed, 
the structure group of the geometry is~$\Un(p,q)$ for some~$p,q>0$ with
$p{+}q=n$.  Since this group is simply a different real form of the
group~$\Un(n)$, one would expect a similar decomposition of the 
curvature tensor.  Indeed, this is what happens, the curvature tensor
again breaking into the sum of three irreducible tensors.  For simplicity
of terminology, I will still refer to these as the scalar curvature, 
the traceless Ricci curvature, and the Bochner curvature and will
refer to pseudo-K\"ahler structures for which the Bochner curvature
vanishes as Bochner-K\"ahler.  

The differential analysis of~\S\ref{ssec: diff analysis} extends essentially 
without change to the pseudo-K\"ahler case; it is just a matter of changing 
a few signs. Theorems~\ref{thm: existence} 
through~\ref{thm: generators of center} 
generalize with essentially no change as 
well.  However, past this point, the analysis becomes somewhat more 
complicated because the orbit structure of the action of~$\Un(p,q)$ 
on~$\euu(p,q){\oplus}\C{n}{\oplus}\bbR$ is considerably more complicated
than before.  One must deal with non-diagonalizable elements, nilpotent
orbits, and a host of other problems.  It seems unlikely that the simple
description of the analytically connected equivalence classes found in the
K\"ahler case can be carried through in the pseudo-K\"ahler case.

\subsection{A split-form analog}
\label{ssec: split form}

There is another `real form' of K\"ahler geometry that has an 
analog of the Bochner-K\"ahler condition.  

A K\"ahler structure
can be thought of as a symplectic manifold~$(M^{2n},\Omega)$ endowed 
with an $\Omega$-skew endomorphism~$J:TM\to TM$ that 
satisfies~$J^2 = -\I$ and a torsion-free connection~$\nabla$
with respect to which both~$\Omega$ and~$J$ are parallel.  

A different geometry results if one starts with a symplectic 
manifold as above and considers an $\Omega$-skew endomorphism~$K:TM\to TM$
that satisfies~$K^2 = +\I$ together with a torsion-free 
connection~$\nabla$ with respect to which both~$\Omega$ and~$K$
are parallel.  Some authors~\cite{Pu1,Pu2} call the 
data~$(M,\Omega,K,\nabla)$ a \emph{hyperbolic K\"ahler structure}, 
though this terminology seems likely to invite confusion. 

Since the null plane fields of~$K\pm\I$ are  necessarily $\Omega$-Lagrangian 
plane fields and since the hypothesis that there be a torsion-free connection 
with respect to which they are parallel implies that these two plane fields 
are integrable, such a structure endows the symplectic manifold with a 
pair of transverse, $\Omega$-Lagrangian foliations~$\cF_\pm$.  

Conversely, any symplectic manifold~$(M^{2n},\Omega)$ endowed with a pair 
of transverse, $\Omega$-Lagrangian foliations~$\cF_\pm$ has an
$\Omega$-skew endomorphism~$K:TM\to TM$ so that the tangent plane fields
to the two foliations are the kernels of~$K\pm\I$ and a unique 
torsion-free connection~$\nabla$ with respect to which both~$\Omega$
and~$K$ are parallel.  Thus, it makes sense to call such a structure
a \emph{bi-Lagrangian} structure, which I will do for the rest of this
subsection.

Let~$\bbR_n$ denote the space of \emph{row} vectors of length~$n$
whose entries are real numbers, so that the natural matrix multiplication
$\bbR_n\times \bbR^n\to\bbR$ is a non-degenerate pairing and~$\bbR_n$
is thus identified as the dual space of~$\bbR^n$. Endow~$\bbR_n\oplus\bbR^n$ 
with its natural induced symplectic structure.  Let~$\GL(n,\bbR)$
act on~$\bbR_n\oplus\bbR^n$ on the left by
\begin{equation*}
A\cdot(\xi, x) = \bigl(\xi\,A^{-1},\ A\,x\bigr).
\end{equation*}
This action preserves the symplectic structure on~$\bbR_n\oplus\bbR^n$
and its bi-Lagrangian
splitting into~$L_-=\bbR_n\oplus0$ and~$L_+ = 0\oplus\bbR^n$.  In fact,
$\GL(n,\bbR)$ is the largest subgroup of~$\Aut(\bbR_n\oplus\bbR^n)$
that preserves these structures.

Now let~$(M^{2n},\Omega,\cF_\pm)$ be a bi-Lagrangian manifold. 
Say that a coframe~$u:T_xM\to\bbR_n\oplus\bbR^n$
is \emph{adapted} if $u$ is a symplectic isomorphism,
identifies~$T_x\cF_-$ with~$\bbR_n\oplus 0$, and identifies~$T_x\cF_+$
with~$0\oplus\bbR^n$.  The bundle~$\pi:P\to M$ of adapted coframes is
then naturally a right~$\GL(n,\bbR)$-bundle with the action defined by
\begin{equation*}
(u\cdot A)(v) = A^{-1}\cdot u(v).
\end{equation*}
The tautological 1-form of this~$\GL(n,\bbR)$-structure can be
written in the form~$(\eta,\omega)$ where~$\eta$ takes values in~$\bbR_n$
and~$\omega$ takes values in~$\bbR^n$.  
One then has the formula~$\pi^*(\Omega) = \eta\w\omega$.

The existence of a torsion-free connection with respect to which~$\Omega$
and~$K$ are parallel is equivalent to the
existence of a~$\eugl(n,\bbR)$-valued 1-form~$\phi$ on~$P$ satisfying
the equations
\begin{equation}\label{eq: biL first str eq}
d\eta = - \eta\w\phi,\qquad\qquad d\omega =-\phi\w\omega.
\end{equation}
This is the \emph{first structure equation} of Cartan.  The 2-form
$\Phi = d\phi + \phi\w\phi$ then satisfies the \emph{first Bianchi
identities}
\begin{equation}
\eta\w\Phi = \Phi\w\omega = 0.
\end{equation}
These identities imply that there is a 
function~$R:P\to\Hom\bigl(\bbR_n\otimes\bbR^n,\eugl(n,\bbR)\bigr)$ so
that the \emph{second structure equation} holds:
\begin{equation*}
\Phi = d\phi + \phi\w\phi = R(\eta\otimes\omega)
\end{equation*}
and, moreover, that~$R$ can be interpreted
as taking values in a `curvature space'~$\cK$ that is isomorphic
as a~$\GL(n,\bbR)$-module to~$S^2(\bbR_n)\otimes S^2(\bbR^n)$.  
Applying the trace (or `contraction') maps
\begin{equation*}
S^2(\bbR_n)\otimes S^2(\bbR^n)\longrightarrow
\bbR_n\otimes\bbR^n\longrightarrow\bbR,
\end{equation*}
then yields, as in the K\"ahler case, a decomposition of~$\cK$ into three 
irreducible, inequivalent $\GL(n,\bbR)$-modules and a corresponding 
decomposition of the curvature tensor of any bi-Lagrangian structure 
into three parts.  For simplicity, I will refer to these three parts as the 
scalar curvature, the traceless Ricci curvature and the Bochner curvature.
(In~\cite{Pu1,Pu2}, the latter curvature is called the ``$HB$-tensor''.)

When the Bochner curvature vanishes, the bi-Lagrangian structure will be
said to be \emph{Bochner-bi-Lagrangian}.  This vanishing condition
is equivalent to the existence of a 
function~$S:P\to\eugl(n,\bbR)\simeq\bbR^n\otimes\bbR_n$
that satisfies
\begin{equation}\label{eq: biL second str eq}
d\phi = -\phi\w\phi + S\,\eta\w\omega - S\,\omega\w\eta - \omega\w\eta\,S
           + \eta\,S\w\omega\,\,\I_n\,.
\end{equation}
The reader will note the analogy with the second structure equation
for Bochner-K\"ahler structures.

The same sort of analysis as in~\S\ref{ssec: diff analysis} 
shows that there exist
functions~$F:P\to\bbR^n$ and~$G:P\to\bbR_n$ so that
\begin{equation}\label{eq: biL third str eq}
dS = -\phi\,S +S\,\phi +F\,\eta+\omega\ G
            +{\ts\frac12}(G\,\omega +\eta\,F)\,I_n\,;
\end{equation}
that there exists a function~$Q:P\to\bbR$ so that
\begin{equation}\label{eq: biL fourth str eq}
dF = -\phi\,F + \bigl(Q\,\I_n + S^2\bigr)\,\omega\,,
\qquad\qquad
dG =  G\,\phi + \eta\,\bigl(Q\,\I_n + S^2\bigr)\,;
\end{equation}
and that
\begin{equation}\label{eq: biL fifth str eq}
dQ = GS\,\omega + \eta\,SF.
\end{equation}
Moreover, the exterior derivatives of 
equations~(\ref{eq: biL first str eq}--\ref{eq: biL fifth str eq}) 
are identities.

Thus, the system of structure equations~
(\ref{eq: biL first str eq}--\ref{eq: biL fifth str eq}) 
satisfies the conditions for Cartan's Theorem~\ref{thm: Cartan existence}
to apply (see the appendix).
In particular, the analog of Theorem~\ref{thm: existence} 
 holds for Bochner-bi-Lagrangian
structures and there is a finite-dimensional moduli space of germs
of such structures.  The analog of Theorem~\ref{thm:  BK constants} 
will hold as well, in
that there will be $n{+}1$ polynomials in the functions~$S$, $F$, $G$, 
and~$Q$ that are constant on each connected Bochner-bi-Lagrangian
structure bundle and the rank of the mapping
$(S,F,G,Q):P\to\eugl(n,\bbR)\oplus\bbR^n\oplus\bbR_n\oplus\bbR$ is
never more than~$n$, implying that the `group' of local isometries
of the structure always acts with local cohomogeneity at most~$n$.

In principle, one could describe the analytically connected equivalence 
classes for this type of structure and examine completeness questions, 
and so on.  This project is made much more complicated than its
K\"ahler analog by the fact
that the $\GL(n,\bbR)$-invariant polynomials on~
$\eugl(n,\bbR)\oplus\bbR^n\oplus\bbR_n\oplus\bbR$ do not separate the
$\GL(n,\bbR)$-orbits.  This is potentially interesting, since it means
that one could possibly have continuous families of Bochner-bi-Lagrangian 
structures all with the same coarse moduli.  Whether this really does
happen is an interesting question.

\subsection{Self-dual K\"ahler metrics}
\label{ssec: self-dual K metrics}

\emph{This section was added after P. Gauduchon sent me the
preprint~\cite{AG}.  I thank the authors for bringing it to my attention.} 

The reader will recall that, when~$n=2$, the Bochner tensor is
the same as the anti-self-dual part of the Weyl tensor.  I.e., 
when $n=2$, Bochner-K\"ahler metrics are the same
as self-dual K\"ahler metrics.  The self-dual part of the Weyl
curvature in this case is essentially the scalar curvature~$s$.  
In particular, the squared norm of the Weyl curvature is the same 
as~$s^2$, up to a universal constant factor.

 From this point of view, some of the results in this article
in the case of dimension~2 had already been obtained.  
For example, in~\cite[Theorem~1]{De}
(which also follows from earlier work by B.-Y. Chen~\cite{Ch}) asserts
that there are no compact self-dual K\"ahler manifolds other than 
the locally symmetric ones.  Of course, this is the ~$n=2$ case of
Corollary~\ref{cor: compact BK is loc sym} of the present article.  
Their proofs use non-trivial
global results about complex surfaces, while the proof in the present
article is essentially self-contained.  It is also interesting to
note that, in view of Theorem~\ref{thm: wtd proj space}, 
their proofs must make essential 
use of the hypothesis that the domain of definition of the metric 
is a compact manifold, rather than just a compact orbifold.

After the initial version of this article was posted to the arXiv, 
I was contacted by P. Gauduchon, who explained that he and V. Apostolov 
had recently obtained a local classification of self-dual
Hermitian-Einstein metrics and that this implied a
local classification of self-dual K\"ahler metrics.  In particular,
they had also proved that such metrics always have local cohomogeneity
at most 2.  For more information about their version of the local 
classification, the reader should consult their preprint~\cite{AG}.   
In particular, their work provides an independent alternative to the 
classification derived in this article when~$n=2$. 

In fact, a remarkable relation between self-dual K\"ahler metrics 
and Einstein metrics follows from the work of Derdzinski~\cite{De}
and Apostolov and Gauduchon~\cite{AG}.  The interested reader should
consult~\cite{AG} for details, but I will summarize some of their results
here as preparation for the remarks I want to make at the end of 
this subsection.

If~$g$ is a self-dual K\"ahler metric on a complex 2-manifold~$M$ with 
scalar curvature~$s$ not identically zero, then~$g$ is not
conformally flat.  Apostolov and Gauduchon show that on the open set~$M^*$ 
where~$s$ is nonzero, the Hermitian metric~$g^* = s^{-2}g$ 
is Einstein (as well as being self-dual).
Of course, unless~$s$ is constant (which only happens when~$g$
is locally symmetric), $g^*$ will not be K\"ahler.
 
Conversely, Apostolov and Gauduchon show that any self-dual Hermitian 
Einstein metric that is not conformally flat is of the form~$g^*$ for 
a unique self-dual K\"ahler metric~$g$ with non-zero scalar curvature.

However, from the point of view in~\cite{AG}, completeness issues 
for either self-dual Hermitian Einstein metrics or self-dual K\"ahler
metrics appear not to be easily resolvable.  
For example, they did not know%
\footnote{P. Gauduchon, private communication.}
whether or not there were any 
complete examples of cohomogeneity~2.  
Using the description in this article, however,
it is easy to see that there are many complete examples of 
self-dual Hermitian Einstein metrics with cohomogeneity~2.  

Before discussing these examples, here are three general observations
that will be useful:  Let~$M$ be a connected complex surface endowed with 
a Bochner-K\"ahler metric~$g$ and characteristic 
polynomial~$p_C(t) = t^4 {+} C_2\,t^2 {+} C_3\,t {+} C_4$ 
and momentum mapping~$h = (h_1,h_2):M\to\bbR^2$.
First, the scalar curvature of~$g$ is~$s = -24h_1$.
Second, the Einstein constant of~$g^*$ is~$-6912C_3$.  
Third, the squared norm of the self-dual part of the Weyl curvature 
of~$g^*$ is~$c\,s^6>0$ for some universal constant~$c>0$. 

Now, consider the complete cohomogeneity~2 metrics on~$\C{2}$ provided 
by Theorem~\ref{thm: complete metrics}, 
where the parameters~$\rho_1$ and $\rho_2$
satisfy~$0<\rho_1<\rho_2$.  The characteristic polynomials are
\begin{equation*}
p_C(t) = p_D(t) = (t-r_1)^2(t-r_2)(t-r_3)
\end{equation*}
where
\begin{equation*}
r_1 = {\ts\frac14}(\rho_1{+}\rho_2),\qquad
r_2 = {\ts\frac14}(\rho_2{-}3\rho_1),\qquad
r_3 = {\ts\frac14}(\rho_1{-}3\rho_2).
\end{equation*}
The momentum cell~$C(p_D,\mu)$ is the bounded cell 
of SubCase~3-1$b$ (see Figure~\ref{fig: three roots}). 
Since the momentum mapping~$h=h':\C{2}\to C(p_D,\mu)$
is surjective and since the eigenvalues of~$H$ satisfy~$r_2\le\lambda_1<r_1$
and~$r_3\le\lambda_2\le r_2$, it follows that~$h_1 
= \tr H = \lambda_1{+}\lambda_2$ varies between an infimum 
of~$r_2{+}r_3$ (achieved only at~$0\in\C{2}$) and a supremum of~$r_1{+}r_2$ 
(not achieved). Thus, since~$s=-24h_1$, the scalar curvature 
satisfies the bounds
\begin{equation*}
-12(\rho_2{-}\rho_1) < s \le 12(\rho_2{+}\rho_1) = s(0).
\end{equation*}
Moreover, since~$C(p_D,\mu)$ has only one vertex and neither of
its two faces is vertical, it follows that $dh_1$ vanishes only at~$0$.
Consequently, the equation~$s=0$ 
defines a smooth hypersurface~$S\subset\C{2}$.  This hypersurface is 
unbounded because the $u_2$-axis (i.e., $u_1=0$) cuts through the 
omitted face of~$C(p_D,\mu)$. 

Since~$s$ is bounded, it follows that~$g^*_\rho = s^{-2}g_\rho$ 
is complete on each of the two domains~$D_+$ (where~$s>0$) and~$D_-$
(where~$s<0$).  The domain~$D_+$ is contractible, while~$D_-$
has the homotopy type of a circle.  Thus, this one example of
a self-dual K\"ahler metric gives rise to two \emph{distinct} complete, 
self-dual Hermitian Einstein manifolds.  

As another interesting example, consider the Bochner-K\"ahler metric of 
SubCase~4-1, where~$r_0>r_1>r_2>r_3$ are chosen so that~$r_0{+}r_3<0$.
Choosing the `completion'~$X_2\subset\C{2}$ obtained by omitting
the face~$l_2=0$, one sees that
the domain~$D_+\subset\C{2}$ defined by~$s>0$ is bounded, 
with boundary a smooth compact hypersurface~$S\subset\C{2}$.  
Again, the corresponding self-dual Hermitian-Einstein metric~$g^*$ 
is complete on~$D_+$.  This case is interesting because~$(D_+,g^*)$ exists
even when the ratios of the roots~$r_i$ are not rational, so that
any attempt to `complete' the corresponding self-dual 
K\"ahler metric~$(D_+,g)$ to a maximal domain leads inevitably to worse 
than orbifold singularities, i.e., to a non-Hausdorff complex space.

By considering Case 1, one can construct an example of a
self-dual Hermitian-Einstein manifold~$(M,g^*)$ that is maximally
extended and the corresponding self-dual K\"ahler metrics~$(M,g)$
is maximally extended, but such that neither~$g$ nor~$g^*$ is complete.
Neither can be extended because the scalar curvature
of~$g$ is proper on~$M$ and tends to $-\infty$ while the squared norm
of the Weyl curvature of~$g^*$ is proper on~$M$ and tends to ~$+\infty$.

What is perhaps more interesting are the Case 4-0 examples, which 
include the weighted projective planes~$\bbC\bbP^{[p_1,p_2,p_3]}$
where~$0<p_1<p_2<p_3$ are integers with greatest common divisor equal
to~1.  For the Bochner-K\"ahler metric~$g$ on this orbifold, the scalar
curvature is everywhere positive as long as $p_3 < p_1{+}p_2$ and the
corresponding Hermitian Einstein metric has positive Einstein constant.
When~$p_3 = p_1{+}p_2$, the scalar curvature is positive except at
one point (a singular orbifold point) and the corresponding Hermitian
Einstein metric has vanishing Einstein constant and is complete on the
(orbifold) complement of this point.  Finally, when~$p_3 > p_1{+}p_2$, 
the scalar curvature vanishes along 
a hypersurface~$S\subset\bbC\bbP^{[p_1,p_2,p_3]}$.
The complement consists of two open sets~$\bbC\bbP^{[p_1,p_2,p_3]}_\pm$
(labeled according to the sign of~$s$), each endowed with
a complete Hermitian Einstein metric with negative Einstein constant.
One of these two pieces, $\bbC\bbP^{[p_1,p_2,p_3]}_-$, 
can be `unfolded' to become a smooth, complete, Hermitian Einstein manifold 
that is biholomorphic to a bounded domain in~$\C{2}$, 
while the other, $\bbC\bbP^{[p_1,p_2,p_3]}_+$, 
has unremovable orbifold singularities.

\appendix

\section{Cartan's Generalization of Lie's Third Theorem}
\label{sec: appendix}

This appendix is an exposition of the passage~
\cite[Chapter II, \S\S17--29]{Ca} from Cartan's work on
a generalization of Lie's Third Fundamental Theorem to the 
`intransitive case' together with a few comments of an elementary 
nature designed to extend the applicability of Cartan's
results to the smooth category and to a `semi-global' setting.  
(In~\cite{Ca}, Cartan worked almost entirely in
what would now probably be called the category 
of real-analytic germs.)  These results have, in modern times, 
been incorporated into the theory of local Lie algebras, Lie
algebroids, and Lie groupoids.  For references and surveys of
this modern work the reader might consult~\cite{SW} and~\cite{Mack}.
The point of view that I take in this appendix is decidedly not
modern; instead I follow Cartan's exposition and development.
I do this since Cartan's version of the result
is more suited for the application in this article.

\subsection{Cartan's Problem}
One is given the following data:
\begin{enumerate}
\item a nonempty open set~$X\subset\bbR^s$ 
               (with coordinates~$x = (x^a)$ on~$\bbR^s$),
\item an integer~$n\ge1$, and
\item functions~$F^a_i$ and~$C^i_{jk}=-C^i_{kj}$ on~$X$, 
               for $1\le i,j,k\le n$ and $1\le a\le s$.
\end{enumerate}
The goal is to describe the solutions to the following 
`realization problem':  Find
\begin{enumerate}
\item a manifold~$N^n$,
\item a coframing~$\eta =(\eta^i)$ of~$N$, and 
\item a mapping~$h = (h^a):N\to X\subset\bbR^s$
\end{enumerate}
satisfying
\begin{equation}\label{eq: Cartan str eqs}
d\eta^i = {\ts\frac12} C^i_{jk}(h)\,\eta^j\w\eta^k,\qquad
dh^a = F^a_i(h)\,\eta^i\,.
\end{equation}

\begin{example}[Lie's Third Fundamental Theorem]
Consider the simple case  where the~$F^a_i$ are all zero. Then
the mapping~$h:N\to X$ of any realization must be constant, say~$h=\bar h$.
A necessary condition on the constants~$\bar C^i_{jk} = C^i_{jk}(\bar h)$ 
can then be found by computing the exterior derivatives of the equations
\begin{equation*}
d\eta^i = {\ts\frac12} \bar C^i_{jk}\,\eta^j\w\eta^k.  
\end{equation*}
These give~$0 = \bar C^i_{jl}\,d\eta^j\w\eta^l$, which, in view of the
above relations, can be rewritten (after an index substitution 
and skewsymmetrization) in the form
\begin{equation*}
 0 = {\ts\frac12}\,\bar C^i_{pl}\bar C^p_{jk}\,\eta^j\w\eta^k\w\eta^l
  ={\ts\frac16}\bigl(\bar C^i_{pj}\bar C^p_{kl}
    +\bar C^i_{pk}\bar C^p_{lj}+\bar C^i_{pl}\bar C^p_{jk}\bigr)
         \,\eta^j\w\eta^k\w\eta^l.
\end{equation*}
Using the linear independence of the~$\eta^i$, one derives
the \emph{Jacobi conditions}
\begin{equation*}
\bar C^i_{pj}\bar C^p_{kl}
    +\bar C^i_{pk}\bar C^p_{lj} + \bar C^i_{pl}\bar C^p_{jk} = 0
\end{equation*}
as necessary conditions for the existence of a solution to the problem.
In other words, any realization~$(N,\eta,h)$ must have~$h$ be constant
and take values in the locus~$X'\subset X$ defined by the equations
\begin{equation*}
C^i_{pj}C^p_{kl}+C^i_{pk}C^p_{lj}+C^i_{pl}C^p_{jk} = 0.
\end{equation*}

Conversely, Lie's Third Fundamental Theorem asserts 
that the Jacobi conditions suffice to ensure the existence of a solution
to the realization problem. I.e., if~$\bar h$ lies in~$X'$,
then there exists a realization~$(N,\eta, h)$ with~$h\equiv\bar h$.
Moreover, any two realizations assuming the same value~$\bar h$
are locally equivalent in the obvious sense. 
\end{example}

\subsection{Differential conditions in the general case}
Even when the~$F^a_i$ are not assumed to be zero, exterior differentiation
of the equations~\eqref{eq: Cartan str eqs}
 of a realization~$(N,\eta,h)$ yields a  set of 
necessary conditions on the map~$h:N\to X$. Namely, it must satisfy
\begin{equation*}
F^b_i(h){{\partial F^a_j}\over{\partial x^b}}(h)
-F^b_j(h){{\partial F^a_i}\over{\partial x^b}}(h)
= - C^l_{ij}(h)\,F^a_l(h)
\end{equation*}
(which is equivalent to~$d(dh^a)=0$) and
\begin{multline*}
 F^a_j(h){{\partial C^i_{kl}}\over{\partial x^a}}(h)
+F^a_k(h){{\partial C^i_{lj}}\over{\partial x^a}}(h)
+F^a_l(h){{\partial C^i_{jk}}\over{\partial x^a}}(h)\\
=
-\bigl(C^i_{mj}(h)C^m_{kl}(h)
         +C^i_{mk}(h)C^m_{lj}(h)+C^i_{ml}(h)C^m_{jk}(h)\bigr)
\end{multline*}
(which is equivalent to $d(d\eta^i)=0$).
Unless these equations are identities,
they place restrictions on the range of~$h$.

\subsection{Cartan's existence theorem}
On the other
hand, if the above equations \emph{are} identities on the functions~$F^a_i$
and~$C^i_{jk}$, then one might hope to find realizations 
of~\eqref{eq: Cartan str eqs} 
without placing any further restrictions on the range of~$h$.

In~\cite{Ca}, Cartan proved%
\footnote{It would be more accurate to say that Cartan only outlined the 
proof of this result. However, the reader knowledgeable
about Cartan-K\"ahler theory will have no trouble supplying the
details.  Also, while Cartan does not always explicitly state the assumption 
of real-analyticity, it is clear from context that he intended this 
assumption to be in force.}
just such a result in the real-analytic category.

\begin{apptheorem}[Cartan]\label{thm: Cartan existence}
Suppose that~$X\subset\bbR^s$ is an open set and suppose that~$F^a_i$
and $C^i_{jk}=-C^i_{kj}$ for~$1\le a\le s$ and~$1\le i,j,k\le n$ 
are real-analytic functions on~$X$ that satisfy
\begin{equation}\label{eq: Cartan cond one}
F^b_i{{\partial F^a_j}\over{\partial x^b}}
-F^b_j{{\partial F^a_i}\over{\partial x^b}}
= - F^a_lC^l_{ij}\,.
\end{equation}
and
\begin{equation}\label{eq: Cartan cond two}
 F^a_j{{\partial C^i_{kl}}\over{\partial x^a}}
+F^a_k{{\partial C^i_{lj}}\over{\partial x^a}}
+F^a_l{{\partial C^i_{jk}}\over{\partial x^a}}
= -\bigl(C^i_{mj}C^m_{kl}+C^i_{mk}C^m_{lj}+C^i_{ml}C^m_{jk}\bigr)
\end{equation}
Then for every~$h_0\in X$, there exists a real-analytic 
realization~$(N,\eta,h)$ satisfying the structure equations
\begin{equation*}
d\eta^i = {\ts\frac12} C^i_{jk}(h)\,\eta^j\w\eta^k,\qquad
dh^a = F^a_i(h)\,\eta^i\,.
\end{equation*}
and a~$p_0\in N$ for which~$h(p_0)=h_0$.  

Moreover, this realization is locally unique in the following sense:  
Given any other real-analytic realization~$(\tilde N,\tilde\eta,\tilde h)$ 
satisfying the corresponding structure equations
that contains a point~$\tilde p_0\in\tilde N$ satisfying
$\tilde h(\tilde p_0)=h_0$, there exists a $p_0$-neighbor\-hood $U\subset N$,
a $\tilde p_0$-neighbor\-hood $\tilde U\subset\tilde  N$, 
and a real-analytic diffeomorphism~$\phi:\tilde U\to U$ so 
that
\begin{equation*}
\phi(\tilde p_0) = p_0,\qquad \phi^*(\eta)=\tilde\eta, 
\qquad \text{and}\ \ \phi^*(h)=\tilde h.
\end{equation*}
\end{apptheorem}

\begin{remark}[A Paraphrase]
Informally, one can state Cartan's result in the following way:  There
is a `solution' of the structure equations~\eqref{eq: Cartan str eqs}
 provided that the
exterior derivatives of these equations are identities, i.e.,
$d^2=0$ is a formal consequence of~\eqref{eq: Cartan str eqs}.  A solution
is uniquely specified by choosing the values of the `invariants'
$h=(h^a)$ at one point in the domain of the solution.
\end{remark}

\subsubsection{Real-analyticity}
The full theorem that Cartan proves in the cited passage 
is more general than Theorem~\ref{thm: Cartan existence} 
and has to do with existence of so-called 
`infinite groups' (nowadays called pseudo-groups) satisfying a given set 
of structure equations.  
However, Theorem~\ref{thm: Cartan existence} is all that is needed 
in this article.  Cartan's proof is via the Cartan-K\"ahler Theorem, which is 
only valid in the real-analytic category.  While the general theorem that 
Cartan proves really does need real-analyticity, the special case being
discussed here as Theorem~\ref{thm: Cartan existence} 
can be proved without recourse to the 
Cartan-K\"ahler Theorem. Indeed, it can be proved using only the Frobenius 
Theorem, the Poincar\'e Lemma, and Lie's Third Fundamental Theorem 
(the classical one).   See the work of Pradines~\cite{Pra1,Pra2}
for this development.

Thus, the above theorem (both existence and uniqueness) is actually 
valid in the smooth category.  However, note that, in the case
where~$F$ and $C$ actually are real-analytic and satisfy~
\eqref{eq: Cartan cond one} and~\eqref{eq: Cartan cond two},
it follows from Cartan's uniqueness result that any sufficiently
differentiable realization of~\eqref{eq: Cartan str eqs} 
is real-analytic in suitable coordinates.

\begin{example}[Application]
In this article, Theorem~\ref{thm: Cartan existence} 
will be applied to the equations~\eqref{eq: structure equations ii}. 
In that case, the functions~$F$ and~$C$
are polynomial (in fact, either linear or quadratic) in the 
linear coordinates on~$X = \bbR^s = i\euu(n){\oplus}\C{n}{\oplus}\bbR$ 
(here~$s = n^2{+}2n{+}1$).  Thus, the realizations are all real-analytic
in this case.
\end{example}

\subsection{A coordinate-free reformulation}
Cartan's conditions can be recast into a somewhat more geometric form
as follows:  Suppose there are given functions~$F^a_i$ 
and~$C^i_{jk}=-C^i_{kj}$ on a domain~$X\subset\bbR^s$.  
Define $n$ vector fields on~$X$ by
\begin{equation*}
F_i = F^b_i\,{{\partial \hfill}\over{\partial x^b}}
\end{equation*}
for~$1\le i\le n$.  
Then~\eqref{eq: Cartan cond one} can be written in terms
of the Lie bracket as
\begin{equation}\label{eq: alt Cartan cond one}
[F_i,F_j] = -C^k_{ij}\,F_k\,.
\end{equation}
Also, \eqref{eq: Cartan cond two} can be written as
\begin{equation}\label{eq: alt Cartan cond two}
F_jC^i_{kl}+F_kC^i_{lj}+F_lC^i_{jk} 
= -\bigl(C^i_{mj}C^m_{kl}+C^i_{mk}C^m_{lj}+C^i_{ml}C^m_{jk}\bigr).
\end{equation}
When the vector fields~$F_i$ are linearly independent, 
\eqref{eq: alt Cartan cond two}
follows directly from~\eqref{eq: alt Cartan cond one}; 
it is simply the Jacobi identity for the
Lie bracket.  However, when the~$F_i$ are everywhere linearly dependent
(as is the case in the application in this article) the equations~
\eqref{eq: alt Cartan cond two}
are not consequences of~\eqref{eq: alt Cartan cond one}.

In~\eqref{eq: alt Cartan cond one} and~
\eqref{eq: alt Cartan cond two}, no explicit reference
is made to the coordinates~$x^a$ on~$X$.  Thus, it makes sense to 
speak of systems~$(X,F,C)$ satisfying \eqref{eq: alt Cartan cond one}
 and~\eqref{eq: alt Cartan cond two} where~$X$
is any smooth manifold, the~$F_i$ are smooth vector fields on~$X$
and the~$C^k_{ij}$ are smooth functions on~$X$.  
Such a system~$(X,F,C)$ is an example of what has since become known 
as a \emph{local Lie algebra}~\cite{Sh} or a \emph{Lie algebroid}~
\cite{SW,Mack}.  The notion of 
a realization~$(N,\eta,h)$ generalizes as well, with the formula
for~$d\eta^i$ remaining the same but the formula~$dh = F_i\,\eta^i$ now
being interpreted as a formula for~$dh :TN\to TX$ in the obvious sense.  
This `coordinate free' formulation of Cartan's problem will not be needed 
in this article, so I will not discuss it any further here.

\subsection{The leaves of~$F$}
For each~$x\in X$, let~$r(x)$ be the dimension of the span of
the vectors~$\{F_i(x)\}_{1\le i\le n}$.  When~$r$ is a constant function 
on~$X$ and~\eqref{eq: alt Cartan cond one} holds, the Frobenius
theorem asserts that the vector fields~$F_i$ are tangent to a foliation
of~$X$ whose leaves have dimension~$r$.

In most applications, however, the function~$r$ is not constant on~$X$.
(Indeed, it is not constant for the 
system~\eqref{eq: structure equations ii}.)
Nevertheless, there is a simple generalization of the Frobenius theorem
that does hold whenever~\eqref{eq: alt Cartan cond one} holds.
 
Say that a smooth curve~$\xi:[a,b]\to X$ is an~\emph{$F$-curve}%
\footnote{When~$r$ is not constant, this condition is \emph{a priori}
stronger than the mere condition that~$\xi'(t)$ lie
in the span of~$\{F_i\bigl(\xi(t)\bigr)\}_{1\le i\le n}$ 
for all~$t\in[a,b]$.}
if there exist smooth functions~$v^i$ on~$[a,b]$ for which
$\xi'(t) = v^i(t)\,F_i\bigl(\xi(t)\bigr)$ and say that~$x_1$
and~$x_2$ in~$X$ are \emph{$F$-equivalent}
if they can be joined by a smooth~$F$-curve.  

The generalized Frobenius theorem says that, if
the vector field system~$F$ satisfies~\eqref{eq: alt Cartan cond one}, 
then the $F$-equivalence
class~$[x]_F$ of~$x\in X$ is a smooth, connected submanifold of~$X$ 
of dimension~$r(x)$.  It is called the~\emph{$F$-leaf} through~$x$.
This generalized notion of a foliation is sometimes known as a \emph{Stefan
foliation} in the literature.  For further discussion of this 
singular leaf structure, which is virtually the same
as the sort of singular leaf structure that one encounters
in the theory of Poisson manifolds, see~\cite{SW}.

\subsection{The rank of a realization}
Suppose now that~$(X,F,C)$ satisfies~
\eqref{eq: Cartan cond one} and \eqref{eq: Cartan cond two} 
(or, equivalently, \eqref{eq: alt Cartan cond one}
 and \eqref{eq: alt Cartan cond two}).  

Then, for any realization~$(N,\eta,h)$ of the structure
equations~\eqref{eq: Cartan str eqs}
with $N$ connected, the map~$h:N\to X$ will have its 
image lie in a single~$F$-leaf~$L\subset X$, whose dimension will 
be~$r(N)=r\bigl(h(p)\bigr)$ for some (and hence any)~$p\in N$.  
Moreover, the structure equations~\eqref{eq: Cartan str eqs} imply that 
the map~$h:N\to L$ has constant rank~$r(N)$ and hence is a 
submersion onto its (open) image in~$L$.  

In~\cite{Ca}, Cartan
assumes these results without proof or remark.  It is not clear 
whether he knew these facts (which, even in the real-analytic
case, require argument it seems to me) 
or merely assumed that he was in some `generic' case where they held.
In any case, he does not make an issue of it.

The integer~$r(N)$ 
will be referred to as the \emph{rank} of the realization~$(N,\eta,h)$.

\begin{example}[Application]
For the system~\eqref{eq: structure equations ii}, 
the dimension of an $F$-leaf can be as low as~$0$ or as high as~$n(n{+}1)$.  
\end{example}

\subsection{The symmetry algebra of a leaf}
Let~$(X,F,C)$ satisfy~
\eqref{eq: Cartan cond one} and \eqref{eq: Cartan cond two}. 
Let~$L\subset X$
be an~$F$-leaf of rank~$r$, and let~$(N,\eta,h)$ be a realization
of the structure equations~\eqref{eq: Cartan str eqs}
 whose image~$h(N)$ is an open subset
of~$L$.  Then by Theorem~\ref{thm: Cartan existence}, 
given any~$\bar h\in h(N)$ and any
two points~$p_1$ and $p_2$ in the fiber~$h^{-1}(\bar h)\subset N$,
there is a locally defined `symmetry' of the realization that carries~$p_1$
to~$p_2$.  This locally defined symmetry is unique in a neighborhood 
of~$p_1$.  

Cartan might have expressed this fact by saying something like `the group of 
symmetries of the system~$(\eta,h)$ acts simply transitively on the fibers 
of~$h$'.  In the modern literature, this sort of vagueness about the domain 
of the `group' of `local symmetries' of such data is usually avoided by 
giving a more precise statement using the language of (finite-dimensional) 
pseudo-groups.  Rather than introduce this sort of terminology, 
I will give the corresponding infinitesimal formulation, which is simpler.

\begin{apptheorem}\label{thm: symmetry group}
If $N$ is connected and simply-connected 
and~$(N,\eta,h)$ is a realization of~\eqref{eq: Cartan str eqs} 
of rank~$r$, then the subset~$\euh\subset\euX(N)$ consisting 
of the vector fields on~$N$ whose {\upshape(}local\/{\upshape)} flows on~$N$ 
preserve~$\eta$ and~$h$ is a Lie algebra of dimension~$n{-}r$.
Moreover, for any~$x\in N$, the evaluation map~$e_x:\euh\to T_xN$
is a vector space isomorphism onto the kernel of~$h'(x):T_xN\to \bbR^s$.
\end{apptheorem}

Up to isomorphism, the Lie algebra~$\euh$ depends only on the leaf~$L$
that contains~$h(N)$.  
It will be referred to as the \emph{symmetry algebra} of~$L$.

It is useful to note that the symmetry algebra of a leaf~$L$ can be 
computed without actually having to find a realization~$(N,\eta,h)$ 
with~$h(N)\subset L$.  In fact, for any~$\bar h\in L$, define a
skewsymmetric bilinear pairing~$[,]_{\bar h}:\bbR^n\times\bbR^n\to\bbR^n$ by
\begin{equation*}
[E_i,E_j]_{\bar h} = C^k_{ij}(\bar h)\,E_k\,,
\end{equation*}
where~$E_i$ is the standard basis of~$\bbR^n$.  
Let~$\euh_{\bar h}\subset\bbR^n$ be the subspace 
that is the kernel of the (surjective) linear 
map~$\lambda_{\bar h}:\bbR^n\to T_{\bar h}L$ 
that satisfies~$\lambda_{\bar h}(E_i) = F_i(\bar h)$.  
Then the restriction of~$[,]_{\bar h}$ to~$\euh_{\bar h}$
defines a Lie algebra structure on~$\euh_{\bar h}$.  One can
verify that, up to isomorphism, this Lie algebra does not
depend on the choice of~$\bar h\in L$ and that this is indeed
the symmetry algebra of~$L$.

\subsection{A semi-global realization}
With these concepts,  a `semi-global' version of Cartan's existence 
and uniqueness result can be stated.
For lack of space, I will not discuss the (relatively straightforward)
proof, which, in any case, can be found in the above cited references.

\begin{apptheorem}\label{thm: realization}
Let~$(X,F,C)$ satisfy~\eqref{eq: Cartan cond one} 
and~\eqref{eq: Cartan cond two}, let~$L\subset X$ be
an $F$-leaf with symmetry algebra~$\euh$. Let~$H$ be a
Lie group whose Lie algebra is~$\euh$.  

Then over any contractible open subset~$U\subset L$  
there exists a principal left $H$-bundle~$(h^a) = h:N\to U$
together with an $H$-invariant coframing~$\eta = (\eta^i)$ on~$N$
so that~$(N,\eta,h)$ satisfies~\eqref{eq: Cartan str eqs}.
This realization is unique up to isomorphism.
\end{apptheorem}

Simple examples show that existence and/or uniqueness can fail
when~$U$ has nontrivial homotopy groups.  In fact, this is the source
of the orbifold singularities encountered 
in~\S\ref{sssec: ness conds for complete}.

When~$H$ is abelian, the obstruction to global existence
on a leaf~$L$ can be formulated as the vanishing of an element of
an appropriate cohomology group on~$L$.  When~$H$ is nonabelian,
there is still a cohomological condition, but it takes values
in a certain nonabelian cohomology set.  Since this refinement will
not be needed in this article, it will not be discussed.

\bibliographystyle{amsplain}

\begin{thebibliography}{10}

\bibitem{Ab} 
 M. Abreu, 
\textit{K\"ahler geometry of toric varieties and extremal metrics},
Internat. J. Math. \textbf{9} (1998), 641--651.
\href{http://www.ams.org/mathscinet-getitem?mr=99j:58047}{MR 99j:58047}

\bibitem{AG}
 V. Apostolov and P. Gauduchon, 
\textit{Self-dual Einstein Hermitian four-manifolds}
preprint 2000, {\tt arXiv:math.DG/0003162}.

\bibitem{Be} 
A. Besse, \emph{Einstein Manifolds}, Springer-Verlag, New York, 1987.
\href{http://www.ams.org/mathscinet-getitem?mr=88f:53087}{MR 88f:53087}

\bibitem{Bo} 
 S. Bochner, \textit{Curvature and Betti numbers, II},
Ann. Math. \textbf{50} (1949), 77--93.
\href{http://www.ams.org/mathscinet-getitem?mr=10:571f}{MR 10:571f}

\bibitem{BY} 
 S. Bochner and K. Yano,
\emph{Curvature and Betti Numbers},
Annals of Math. Studies, No.~32,
Princeton University Press Princeton, 1953.
\href{http://www.ams.org/mathscinet-getitem?mr=15:989f}{MR 15:989f}

\bibitem{Ca}
 \'E. Cartan, 
\textit{Sur la structure des groupes inifinis de transformations},
Ann. \'Ec. Norm. \textbf{3} (1904), 153--206.
(Reprinted in Cartan's Collected Works, Part II.)

\bibitem{Ch}
 B. Y. Chen, \textit{Some topological obstructions to Bochner-K\"ahler metrics 
       and their applications},
J. Diff. Geom. \textbf{13} (1978), 547--558.
\href{http://www.ams.org/mathscinet-getitem?mr=81f:32037}{MR 81f:32037}

\bibitem{DSV}
 J. Deprez, et al, 
\textit{Classifications of K\"ahler manifolds satisfying 
        some curvature conditions},
Sci. Rep. Niigata Univ. Ser. A \textbf{24} (1988), 1--12.
\href{http://www.ams.org/mathscinet-getitem?mr=89d:53045}{MR 89d:53045} 

\bibitem{De}
  A. Derdzi\'nski,
\textit{Self-dual K{\"a}hler manifolds 
         and Einstein manifolds of dimension four},
Compositio Math. \textbf{49 } (1983), 405--433.
\href{http://www.ams.org/mathscinet-getitem?mr=84h:53060}{MR 84h:53060} 

\bibitem{Ej}
 N. Ejiri, \textit{Bochner-K\"ahler metrics},
Bull. Sci. Math. \textbf{108} (1984), 423--436.
\href{http://www.ams.org/mathscinet-getitem?mr=86g:53026}{MR 86g:53026} 

\bibitem{FH}
 W. Fulton and J. Harris,
\textit{Representation Theory},
Graduate Texts in Mathematics, no. 129, Springer-Verlag, New York, 1991.
\href{http://www.ams.org/mathscinet-getitem?mr=93a:20069}{MR 93a:20069}

\bibitem{Gu}
 V. Guillemin, \textit{Kaehler metrics on toric varieties},
J. Diff. Geom. \textbf{40} (1994), 285--309.
\href{http://www.ams.org/mathscinet-getitem?mr=95h:32029}{MR 95h:32029}

\bibitem{He}
 S. Helgason,
\emph{Differential Geometry, Lie Groups, and Symmetric Spaces},
Academic Press, Princeton, 1978.
\href{http://www.ams.org/mathscinet-getitem?mr=80k:53081}{MR 80k:53081}

\bibitem{KK}
 U.-H. Ki and B. H. Kim, 
\textit{Manifolds with Kaehler-Bochner metric},
Kyungpook Math. J. \textbf{32} (1992), 285--290.
\href{http://www.ams.org/mathscinet-getitem?mr=93m:53075}{MR 93m:53075}

\bibitem{LPV}
 J. Leysen, et al 
\textit{Some curvature conditions in Bochner-Kaehler manifolds},
Atti Accad. Peloritana Pericolanti Cl. Sci. Fis. Mat. Natur.
 \textbf{65} (1987), 85--94.
\href{http://www.ams.org/mathscinet-getitem?mr=90b:53027}{MR 90b:53027}

\bibitem{KN}
 S. Kobayashi and K. Nomizu,
\emph{Foundations of Differential Geometry, vol.~II}, 
John Wiley \&~Sons, New York, 1963.
\href{http://www.ams.org/mathscinet-getitem?mr=38:6501}{MR 38:6501}

\bibitem{Mack}
K. Mackenzie,
Lie algebroids and Lie pseudoalgebras,
{\em Bull. London Math. Soc. } {\bf 27} (1995), 97--147.  
\href{http://www.ams.org/mathscinet-getitem?mr=MR96i:58183}{MR 96i:58183}

\bibitem{Ma}
 M. Matsumoto, \textit{On K\"ahlerian spaces with parallel or vanishing 
          Bochner curvature tensor},
Tensor (N.S.) \textbf{20} (1969), 25--28.
\href{http://www.ams.org/mathscinet-getitem?mr=39:3433}{MR 39:3433}

\bibitem{MT}
 M. Matsumoto and S. Tanno, 
\textit{K\"ahlerian spaces with parallel 
       or vanishing Bochner curvature tensor},
Tensor (N.S.) \textbf{27} (1973), 291--294.
\href{http://www.ams.org/mathscinet-getitem?mr=49:7943}{MR 49:7943} 

\bibitem{Pra1}
J. Pradines, \textit{Th\'eorie de {L}ie pour 
les groupo\"\i des diff\'erentiables. {C}alcul diff\'erenetiel 
dans la cat\'egorie des groupo\"\i des
infinit\'esimaux,}  C. R. Acad. Sci. Paris S\'er. A-B \textbf{264} (1967), 
A245--A248. 
\href{http://www.ams.org/mathscinet-getitem?mr=MR35:7242}{MR 35:7242} 

\bibitem{Pra2} 
J. Pradines, \textit{Troisi\`eme th\'eor\`eme de Lie 
sur les groupo\"{\i}des diff\'e\-ren\-tiables,} 
 C. R. Acad. Sci. Paris S\'er. A-B \textbf{267} (1968), 
A21--A23. 
\href{http://www.ams.org/mathscinet-getitem?mr=MR37:6969}{MR 37:6969} 

\bibitem{Pr}
 C. Procesi, \textit{The invariant theory of $n\times n$ matrices},
Advances in Math. \textbf{19} (1976), 306--381.
\href{http://www.ams.org/mathscinet-getitem?mr=54:7512}{MR 54:7512}

\bibitem{Pu1}
 N. Pu\v{s}i\'c, 
\textit{On an invariant tensor of a conformal transformation of a hyperbolic
Kaehlerian space},
Zb. Rad. \textbf{4} (1990), 55--64.
\href{http://www.ams.org/mathscinet-getitem?mr=92j:53014 }{MR 92j:53014 }

\bibitem{Pu2}
 N. Pu\v{s}i\'c, \textit{On $HB$-flat hyperbolic Kaehlerian spaces},
Mat. Vesnik \textbf{49} (1997), 35--44.
\href{http://www.ams.org/mathscinet-getitem?mr=98i:53044}{MR 98i:53044}

\bibitem{Sh}
 K. Shiga, \textit{Cohomology of Lie algebras over a manifold. I \& II},
J. Math. Soc. Japan, \textbf{26} (1974), 324--361, 587--607.
\href{http://www.ams.org/mathscinet-getitem?mr=51:4267}{MR 51:4267}
\href{http://www.ams.org/mathscinet-getitem?mr=51:4268}{MR 51:4268}

\bibitem{SW}
A. Cannas da Silva and A. Weinstein,
\textit{Geometric Models for Noncommutative Algebras},
University of Californiat at Berkeley Lecture Notes,
American Mathematical Society, 1999.

\bibitem{TL}
 S. Tachibana and R. Liu, 
\textit{Notes on Kaehlerian metrics with vanishing Bochner curvature tensor},
Kodai Math. Sem. Rep. \textbf{22} (1970), 313--321.
\href{http://www.ams.org/mathscinet-getitem?mr=42:1030}{MR 42:1030}

\bibitem{TW}
 H. Takagi and Y. Watanabe, 
\textit{K\"ahlerian manifolds with vanishing Bochner curvature tensor 
           satisfying $R(X,\,Y)áR_1=0$},
Hokkaido Math. J. \textbf{3} (1974), 129--132.
\href{http://www.ams.org/mathscinet-getitem?mr=49:3736}{MR 49:3736}

\bibitem{VV}
 D. Van Lindt and L. Verstraelen, 
\textit{A survey on axioms of submanifolds 
          in Riemannian and K\"ahlerian geometry},
Colloq. Math. \textbf{54} (1987), 193--213.
\href{http://www.ams.org/mathscinet-getitem?mr=89h:53115}{MR 89h:53115}

\end{thebibliography}

\end{document}